\documentclass[12pt]{article}

\usepackage[utf8]{inputenc}

\usepackage{sky}
\usepackage{minjae}
\usepackage{wrapfig}
\usepackage{graphicx}
\usepackage{pgf}
\usepackage{color}
\usepackage{subfig}
\usepackage{bm}
\usepackage{caption}

\newcommand{\revision}[1]{{#1}}

\newcommand{\pointstobrauerfn}{F}
\newcommand{\trace}{\mathrm{tr}}

\newcommand{\symgrp}{\mrm{S}}
\newcommand{\bpm}{\mrm{DBM}}
\newcommand{\Wg}{\mrm{Wg}}
\newcommand{\epe}{\mrm{EPE}}
\newcommand{\area}{\mrm{area}}
\newcommand{\splitting}{\mbb{S}}
\newcommand{\merge}{\mbb{M}}
\newcommand{\deform}{\mbb{D}}

\newcommand{\numcomp}{\#\mrm{comp}}
\newcommand{\cycles}{\#\mrm{cycles}}
\newcommand{\orthogonal}{\mrm{O}}
\newcommand{\proj}{\mrm{P}}
\newcommand{\fsfm}{\mrm{FEPE}}
\newcommand{\numswaps}{\#\mrm{swaps}}
\newcommand{\numintloops}{\#\mrm{loops}}

\def\@rst #1 #2other{#1}
\newcommand\MR[1]{\relax\ifhmode\unskip\spacefactor3000 \space\fi
\MRhref{\expandafter\@rst #1 other}{#1}}
\newcommand{\MRhref}[2]{\href{http://www.ams.org/mathscinet-getitem?mr=#1}{MR#2}}

\title{Random surfaces and lattice Yang-Mills}

\author{ 
    \begin{tabular}{c}{Sky Cao}\\[-4pt]\small MIT\end{tabular}
    \begin{tabular}{c}{Minjae Park}\\[-4pt]\small University of Chicago \end{tabular}
    \begin{tabular}{c}{Scott Sheffield}\\[-4pt]\small MIT\end{tabular}
}

\begin{document}

\maketitle

\begin{abstract}
We study Wilson loop expectations in lattice Yang-Mills models with a compact Lie group $G$.
% and similar Lie groups such as $\SON$, $\SUN$ and $\SphN$. 
Using tools recently introduced in a companion paper~\cite{park2023wilson}, we provide alternate derivations, interpretations, and generalizations of several recent theorems about Brownian motion limits (Dahlqvist), lattice string trajectories (Chatterjee and Jafarov) and surface sums (Magee and Puder). We show further that one can express Wilson loop expectations as sums over embedded planar maps in a manner that applies to any matrix dimension $N \geq 1$, any inverse temperature $\beta>0$, and any lattice dimension $d \geq 2$. The surface expansions we use for this purpose have not (as far as we know) appeared in the literature before.

When $G=\UN$, the embedded maps we consider are pairs $(\mathcal M, \phi)$ where $\mathcal M$ is a planar (or higher genus) map and $\phi$ is a graph homomorphism from $\mathcal M$ to a lattice such as $\mathbb Z^d$. The faces of $\mathcal M$ come in two partite classes: {\em edge-faces} (each mapped by $\phi$ onto a single edge) and {\em plaquette-faces} (each mapped by $\phi$ onto a single plaquette). The weight of a lattice edge $e$ is the Weingarten function applied to the partition whose parts are given by half the boundary lengths of the faces in $\phi^{-1}(e)$. (The Weingarten function becomes quite simple in the $N\to \infty$ limit.) The overall weight of an embedded map is proportional to $N^\chi$ (where $\chi$ is the Euler characteristic) times the product of the edge weights. We establish analogous results for $\SUN$, $\ON$, $\SON$, and $\SphN$, where the embedded surfaces and weights take a different form. There are several variants of these constructions. In this context, we present a list of relevant open problems spanning several disciplines: random matrix theory, representation theory, statistical physics, and the theory of random surfaces, including random planar maps and Liouville quantum gravity.

\end{abstract}

\maketitle

\tableofcontents

\section{Introduction}\label{section:intro}

\subsection{Overview} \label{subec::introoverview}
On a heuristic level, Euclidean Yang-Mills theory is a ``probability measure'' defined by
\begin{equation*}\label{eqn-the-YM-measure}
 d\mu_{\text{YM}}(\omega) = \frac{1}{Z}e^{-\frac{1}{2g^2}S_{\text{YM}}(\omega)}d\omega 
\end{equation*}
where $\omega$ ranges over a space $\mathcal A$ of Lie-algebra-valued connection forms on some Riemannian manifold, the Yang-Mills action $S_{\text{YM}}$ is the $L^2$-norm of the curvature of $\omega$, $g$ is a coupling constant, and $d\omega$ is a ``Lebesgue measure'' on $\mathcal A$. Making precise sense of the heuristic definition above is a famous open problem that we will not solve here \cite{Jaffe2006a}.

Instead, we will study {\em lattice Yang-Mills theory} (a.k.a.\ {\em lattice gauge theory}), an approximation to the continuum theory introduced in 1974 by Wilson~\cite{Wilson1974} who also credits Polyakov and Smit for similar ideas \cite{wilson2004origins}. An online search for scholarly work on ``lattice gauge theory'' turns up tens of thousands of articles in physics and mathematics, and we cannot cover all of the variants and applications here. Wilson's memoir and Chatterjee's recent survey for probabilists are good places to start \cite{wilson2004origins, chatterjee2016a}. See also Yang's account of his early work with Mills in 1954 \cite{wilson2001mathematical}.

Lattice Yang-Mills assigns a random $N$-by-$N$ matrix from some compact Lie group $G$ --- usually $\UN$, $\ON$, $\SUN$, $\SON$, or $\SphN$ --- to each directed edge of a graph $\Lambda$, which is usually $\mathbb Z^d$ or a finite induced subgraph of $\mathbb Z^d$. We require this assignment to have an edge-reversal symmetry: if $Q_e$ is the matrix assigned to a directed edge $e=(v,w)$, then $Q_{(w,v)} = Q_{(v,w)}^{-1}$. If $p= (e_1, e_2, \ldots, e_k)$ is a directed path, then we write $Q_p = Q_{e_1} Q_{e_2} \ldots Q_{e_k}$. A {\em loop} is a directed cycle $\ell$ defined modulo cyclical reordering (which amounts to repositioning the starting point of the loop). We define a set $\mathcal P$ of directed loops in $\Lambda$ that we call {\em plaquettes}. Usually $\mathcal P$ is the set of directed unit squares in $\Lambda$ (i.e., directed cycles with four distinct vertices), but in principle $\mathcal P$ can be any collection of loops that is closed under reversal (i.e.\ $p \in \mathcal P$ implies that the orientation reversal of $p$ is in $\mathcal P$).

Let $M$ be one of the aforementioned classical Lie groups. Define the \textbf{normalized trace} by $\trace(M) := \frac{1}{N} \Tr(M) = \frac{1}{N} \Bigl( \sum_{j=1}^N M_{j,j}\Bigr)$ and write $\mrm{Re}(z)$ for the real part of $z$.  Note that if $M$ is the identity, then $\mrm{Re}\bigl(\trace(M)\bigr) = 1$ and $\mrm{Re} \bigl(\trace(-M)\bigr) = -1$. In some sense $\mrm{Re} \bigl( \trace(M) \bigr)\in [-1,1]$ is a measure of how close $M$ is to the identity matrix. It is large (close to 1) if $M$ is near the identity. If $\ell$ is a loop then $\trace(Q_\ell)$ is well-defined because the conjugacy class of $Q_{e_1}Q_{e_2} \ldots Q_{e_k}$ (and hence the trace) does not change if we cyclically reorder the $e_i$. If $\ell^{-1}$ is the orientation reversal of $\ell$ then $\trace(Q_\ell) = \overline{\trace(Q_{\ell^{-1}})}$. This is because inverting a matrix inverts its eigenvalues, and (for matrices in compact Lie groups) each eigenvalue $z$ satisfies $|z|^2=z \overline z=1$ so that $1/z = \overline z$. This also implies that for the matrices $M$ in our compact Lie groups, we can write $\frac12 \bigl( \trace(M) + \trace(M^{-1}) \bigr) = \mrm{Re} \bigl(\trace(M)\bigr)$.

For finite graphs $\Lambda$, the lattice Yang-Mills measure is the probability measure
\begin{equation} \label{eqn::latticeym} Z^{-1} \prod_{p\in \mathcal P} \exp \bigl( N\beta \Tr(Q_p) \bigr) \prod_{e \in E_\Lambda^+} d Q_e\end{equation}
where $\beta>0$ is an inverse temperature, $Z$ is a normalizing constant, each $dQ_e$ is Haar measure on the compact Lie group $G$, and (to avoid counting an undirected edge twice) $E_\Lambda^+$ is the set of oriented edges of $\Lambda$ for which the endpoint is lexicographically after the starting point. This is a positive measure because $\mathcal P$ is closed under direction-reversal --- this direction-reversal property implies that we can define a set $\mathcal P^+$ of ``positively oriented plaquettes'' containing exactly one element of $\{\ell, \ell^{-1} \}$ for each $\ell \in \mathcal P$, and then rewrite \eqref{eqn::latticeym} as
\begin{equation}\label{eqn::latticeymreal} Z^{-1} \prod_{p\in \mathcal P^+} \exp \bigl( 2N\beta\, \mrm{Re}(\Tr(Q_p)) \bigr) \prod_{e \in E_\Lambda^+} d Q_e.\end{equation}

\begin{remark}\label{remark:beta-vs-two-beta}
We note here that the above action differs from some previous work~\cite{Chatterjee2019a, chatterjee2016, shen2022new} by a factor of $2$: where we have $2\beta$ the previous works have just $\beta$. This slightly simplifies many of our formulas later on, where $\beta$ appears instead of $\frac{\beta}{2}$.
\end{remark}

Informally, the Yang-Mills measure on $(Q_e)$ configurations corresponds to i.i.d.\ {\em Haar measure} (one instance of Haar measure for each positively directed edge of $\Lambda$) modified by a {\em weighting} that favors configurations for which $Q_p$ is close to the identity whenever $p \in \mathcal P$. A {\em Wilson loop observable} is a quantity of the form
\begin{equs}
W_{s}(Q) := \prod_{\ell \in s} \trace(Q_\ell), 
\end{equs}
where $s$ is some finite collection of loops in $\Lambda$. A
{\em Wilson loop expectation} is a quantity of the form $$
\langle W_{s} \rangle_{\Lambda, \beta} :=
\mathbb E \big[ W_{s}(Q) \big],$$ 
where in the right\revision{-}hand side above we take expectation with respect to the lattice Yang-Mills measure defined in \eqref{eqn::latticeym}.

\begin{remark}
In contrast to some previous works \cite{Chatterjee2019a, chatterjee2016, shen2022new}, our Wilson loops are defined with the normalized trace rather than the trace. Thus, our Wilson loop expectations are $N^{-|s|}$ (here $|s|$ denotes the number of loops in $s$) times the Wilson loop expectations that appear in the works mentioned above. This is a cosmetic distinction; the scaling we use is natural when taking large $N$ limits.
\end{remark}

The fundamental goal of lattice Yang-Mills theory is to understand these quantities. That is, one seeks to compute \begin{equation}
\label{eqn::latticeym3}  Z^{-1} \int \prod_{\ell \in s} \trace(Q_\ell)  \prod_{p\in \mathcal P} \exp \bigl( N\beta \Tr(Q_p) \bigr) \prod_{e \in E_\Lambda^+} d Q_e,\end{equation} which we can Taylor expand and write as
\begin{equation}
\label{eqn::latticeym4} Z^{-1}\int \prod_{\ell \in s} \trace(Q_\ell)\ \prod_{p\in \mathcal P} \Bigl( \sum_{k=0}^\infty \frac{(N\beta)^k}{k!} \Tr(Q_p)^k \Bigr) \prod_{e \in E_\Lambda^+} d Q_e.\end{equation}
Given $K: \mathcal P \to \mathbb N$, write $\displaystyle K!=\prod_{p \in \mathcal P}K(p)!$ and $\displaystyle \beta^{K} = \prod_{p \in \mathcal P} \beta^{K(\rho)}$. Using this notation, write \eqref{eqn::latticeym4} as
\begin{equation}
\label{eqn::latticeymexpanded} Z^{-1} \sum_{K: \mathcal P \to \N} \frac{(N\beta)^K}{K!} \int \prod_{\ell\in s} \trace(Q_\ell)\, \prod_{p\in \mathcal P} \big(  \Tr(Q_p) \big)^{K(p)} \prod_{e \in E_\Lambda^+} d Q_e.\end{equation}

\begin{remark}\label{remark:K-gives-plaquette-copies}
We will later interpret $K$ as specifying the number of copies of each (oriented) plaquette $p \in \mc{P}$. Thus we will often refer to $K$ as the ``plaquette count". 
% It may help the reader to keep this in mind.
\end{remark}

This leads to a classical problem in random matrix theory, which is somehow at the heart of this subject. How can we best compute and understand the individual summands in \eqref{eqn::latticeymexpanded}, which can be described in words as ``expected products of traces of products of matrices---each of which comes from a set of i.i.d.\ Haar-distributed matrices and their inverses''? This question is expressed more carefully in Section~\ref{subsec::randommatrixintro}.

Variants of this question have a long history, beginning with the foundational work of 't Hooft and Br\'ezin et al and Itzykson, Parisi, Zuber from the 1970's \cite{Hooft1973, brezin1978planar} and expanding greatly over subsequent decades, encompassing various types of random matrices, including Gaussian ensembles (such as GUE or GOE) as well as Haar measure on compact Lie groups \cite{MR2858440, itzykson1980planar, mehta1981method, gross1991two, di19952d, zvonkin1997matrix,bessis1980quantum, okounkov2000random,bouttier2002census, 
zinn2003some, guionnet2005combinatorial, guionnet2005combinatorial, maurel2006high,
collins2007second,
eynard2008algebraic,
collins2009asymptotics, eynard2011formal, guionnet2015asymptotics, eynard2016counting}. These papers make connections to the random planar map theory developed by Tutte (and many others) \cite{tutte1968enumeration} and the continuum random surface theory developed by Polyakov (and many others) \cite{polyakov1981quantum}. The third author's recent random surface survey contains many additional references on both sides~\cite{she2022}.

Despite these decades of work, fundamental advances continue to be made. For example, the precise question above was recently addressed in the groundbreaking work of Magee and Puder \cite{Magee2019, Magee2019a}, as explained further in Section~\ref{subsec::randommatrixintro} below. The analog of $\Lambda$ in their setting is a ``blossom graph'' which contains a single vertex and an edge set $E_\Lambda$ that consists of finitely many (distinct and labeled) self-loops. This is somehow the most general setting, because in this scenario {\em any} element of the free group generated by $(Q_e)$ can be written as $Q_\ell$ for a loop $\ell$ in $\Lambda$. Their analysis treats this as a fundamental random matrix question (not necessarily motivated by Yang-Mills) and builds on the classical work of Collins and \'{S}niady \cite{collins2006integration} on the so-called {\em Weingarten calculus} which in turn builds on earlier work by Weingarten himself \cite{weingarten1978asymptotic}.
% See also the representation-theoretic ideas due to Dahlqvist and others \cite{Collins2018, Dahlqvist2016}.
Further analysis on this theme appears in recent work by Buc-d'Alch{\'e} which in particular describes the $N\to \infty$ asymptotic behavior of Wilson loop expectations in terms of so-called {\em unitary maps} \cite{bucdalche2023topological} while also considering generalizations to mixtures of deterministic and unitary matrices.

A series of groundbreaking papers by Chatterjee and/or Jafarov has provided a different approach to identities involving \eqref{eqn::latticeymexpanded} including the {\em Makeenko-Migdal/Master Loop/Schwinger-Dyson equations}. This approach enables them to describe the $N \to \infty$ behavior of \eqref{eqn::latticeymexpanded} in terms of so-called {\em lattice string trajectories}~\cite{chatterjee2016,chatterjee2016a, jafarov2016, Chatterjee2019a}, see also other recent derivations by \cite{shen2022new, AN2023} and several generalizations due to Diez and Miaskiwskyi \cite{diez2022expectation}. These works build on a vast literature in this area, including early works of Makeenko and Migdal \cite{makeenko1979} (see also the recent physics paper \cite{kazakov2023bootstrap} which combines these equations with the bootstrap method in order to numerically compute Wilson loop expectations). Although they work in the setting where $\Lambda$ is an induced subgraph of $\mathbb Z^d$, one may also recall a standard ``gauge fixing'' argument that allows one to reduce to the case that $Q_e$ is fixed to be the identity for all $e$ within some spanning tree of $\Lambda$ (see e.g. \cite[Section9]{chatterjee2016leading}). This is equivalent to identifying that entire tree with a single vertex, which reduces $\Lambda$ to a blossom graph that (in the case $\beta = 0$) agrees exactly with the setting discussed by Magee and Puder. 

We will provide an alternate derivation of some of the blossom graph results of Magee and Puder \cite{Magee2019, Magee2019a} as well as the master field and string trajectory results of Chatterjee and Jafarov \cite{chatterjee2016,chatterjee2016a, jafarov2016, Chatterjee2019a}. Our statements along these lines will be in several ways more general than those in previous works. \revision{We also mention that there is an old paper \cite{BGDF1986} which considers the special cases $G = \mathbb{Z}_2, \unitary(1), \mrm{SU}(2)$, although the methods used therein are quite different from the methods of the present paper.~The resulting surface sum looks quite different as well.}

\begin{enumerate}
\item {\bf General $N$:}
We allow for any matrix dimension $N \geq 1$ (the results in \cite{Magee2019, Magee2019a} are stated for $N$ sufficiently large; one has to use a slightly different definition of the Weingarten function for smaller $N$).
\item {\bf General graphs:} We consider general $\Lambda$ and $\mathcal P$ in our derivation of the Makeenko-Migdal/Master loop/Schwinger-Dyson relations (in \cite{chatterjee2016,chatterjee2016a, jafarov2016, Chatterjee2019a} the plaquettes $\mathcal P$ are taken to be squares, though this is not fundamental to the argument).
\item {\bf More general recurrence formula:} We also derive a more general form of the above-mentioned relations. To roughly explain the distinction, recall that the Makeenko-Migdal/Master loop/Schwinger-Dyson relation in \cite[Theorem 3.6]{chatterjee2016} expresses the Wilson loop expectation of a string $s$ in terms of strings $s'$ obtained by applying local moves to $s$. A stronger result \cite[Theorem 8.1]{chatterjee2016} uses only the $s'$ obtained from local moves involving a {\em single fixed edge} $e \in \Lambda$. Our slightly stronger result uses only the $s'$ obtained from local moves involving a {\em single fixed edge} of $s$. The distinction is that there may be many edges in $s$ that correspond to the same $e \in \Lambda$.  We refer to this as the {\it single-location} Makeenko-Migdal/Master loop/Schwinger-Dyson equation.
\item {\bf General matrix families:} We also include analogs of our result for the most fundamental Lie group families (namely $\UN$, $\ON$, $\SUN$, $\SON$, $\SphN$) while some of the earlier papers focused on one or two such groups. While~\cite{chatterjee2016,chatterjee2016a, jafarov2016, Chatterjee2019a} first frame their results in terms of $\SON$ and $\SUN$ we will frame our results and discussion in terms of $\UN$, which from our point of view is the simplest case.  We then extend the theory to $\ON$, $\SUN$, $\SON$, and $\SphN$ in Section~\ref{sec::othergroups}. This is the longest and most technically challenging part of the paper, as each group family comes with its own interesting set of challenges.
\end{enumerate}
Another straightforward generalization of our result would be to include some deterministic matrices in the words; this type of generalization is considered e.g.\ in~\cite{bucdalche2023topological}. We expect that this should be possible in our setting as well, but we will not discuss this here. %We later explain how our results change if one instead considers $\SUN$, $\SON$, or $\SphN$.
% See Remark \ref{remark:other-groups} for more discussion.

In all of the settings described above, we will explain how to express Wilson loop expectations in terms of random lattice-embedded planar maps, which give rise to convergent sums for any $\Lambda$, any $N \geq 1$, and any $\beta$. These are closely related to both the topological surface sums in \cite{Magee2019, Magee2019a} and the string trajectories in \cite{chatterjee2016,chatterjee2016a, jafarov2016, Chatterjee2019a}, but our derivation and planar map interpretation will be rather different. The main point we want to stress in this paper is that there are powerful ways to express Wilson loop expectations as sums over embedded planar maps. Some settings are more challenging than others (certain symmetries that apply in one setting may not apply in all settings) but we will nonetheless develop a framework that is very general, and that we hope will lead to progress on some of the open problems listed in Section~\ref{sec::open}.

\begin{remark} One of the long-term goals of this theory is to construct and understand a continuum scaling limit of quantities like \eqref{eqn::latticeymexpanded} as $\beta \to \infty$ and the lattice mesh size simultaneously goes to zero at an appropriate rate.  Thus, ideally one desires an understanding of the terms of \eqref{eqn::latticeymexpanded} that is sufficiently robust that it allows one to make predictions about these limits.
\end{remark}

\begin{remark} When $\beta$ is large, the function $x \to \exp (2 \beta x)$, defined for $x \in [-1,1]$, is largest for $x$ near $1$ and {\em much smaller} in the rest of $[-1,1]$. In principle, one could replace the $\exp$ in \eqref{eqn::latticeym3} by a different function with this property: say $x \to \frac12(x^b + x^{b+1})$ for some large $b$. If we took this approach, then the analog of \eqref{eqn::latticeymexpanded} would have only finitely many summands, but we would still expect it to have a similar scaling limit behavior as $b \to \infty$ and the lattice mesh size simultaneously goes to zero. Somehow $b$ is playing the role of $\beta$ here: instead of taking the number of plaquettes of a given type to be {\em a priori} Poisson with parameter $\beta$ we can take the number to be either $b$ or $b+1$ (each with probability $1/2$). Alternatively, one can replace \eqref{eqn::latticeym} with
\begin{equation} \label{eqn::latticefixedfacenum} Z^{-1} \Bigl[ \Bigl(|\mathcal P|^{-1}  \sum_{p\in \mathcal P} \trace(Q_p) \Bigr)^b +  \Bigl( |\mathcal P|^{-1} \sum_{p\in \mathcal P} \trace(Q_p) \Bigr)^{b+1} \Bigl] \prod_{e \in E_\Lambda^+} d Q_e\end{equation}
which somehow fixes the {\em total} number of plaquettes to be $b$ or $b+1$. This approach might also have a similar scaling limit if $b \to \infty$ at the right rate. If one is working toward the goal of ``constructing a candidate continuum theory'' one is allowed to use whatever approach turns out to be most computationally tractable.
\end{remark}

{\em Acknowledgements.} We thank Bjoern Bringmann, Jacopo Borga, Sourav Chatterjee, Doron Puder, \revision{Yahui Qu}, Hao Shen, Jasper Aaron Shogren-Knaak, and Tom Spencer for helpful comments and conversations. \revision{We also thank the anonymous referees whose suggestions have greatly improved this paper.} We thank the Institute for Advanced Study for hosting us while this work was completed. The first author was supported by the Minerva Research Foundation while at IAS, as well as by NSF Award: DMS 2303165. The third author is supported by NSF Award: DMS 2153742.

%\begin{enumerate}
%\item How can we understand the expected trace of a word in $U_1, \ldots, U_n$ where the $U_i$ are i.i.d.\ elements of a Lie group each chosen from Haar measure? For example, what is the expected trace of $U_1 U_2 U_3^2 U_1^{-1} U_3^{-1} U_2^{-1} U_3^{-1}$?
%\item How does the answer change if instead of sampling from Haar measure, we obtain each $U_i$ by running Brownian motion on a Lie group for $t_i$ units of time, starting at the identity?
%\end{enumerate}
%Since Brownian motion is mixing, the answer to the section question is the $t_i \to \infty$ of the answer to the first question.
\subsection{Random matrices} \label{subsec::randommatrixintro}

At the heart of our analysis are two classical questions about the traces of random matrices. The first is the one we discussed in Section~\ref{subec::introoverview} and the second is a close variant.

\begin{enumerate}
\item Suppose $M_1, M_2, \ldots M_k$ are i.i.d.\ samples from Haar measure on $\UN$ (or a similar Lie group) and that $W_1, \ldots, W_m$ are words in the $M_i$ and $M_i^{-1}$. Can we compute the expectation $$\mathbb E \Bigl[\prod_{i=1}^k \trace(W_i)\Bigr]$$ in a ``nice'' way? For example, can we express 
\[ \mathbb E \Bigl[ \trace \bigl( M_1 M_2 M_1^{-1} M_2^{-3}\bigr) \trace(M_1^3 M_2 M_1^{-1} M_2^{-1}) \Bigr] \]
as a simple function of $N$?
\item How does the answer to the previous question change if instead of sampling from Haar measure, we obtain each $M_i$ by running a Brownian motion on the Lie group for $t_i$ units of time, starting with the identity?
\end{enumerate}

The second question can be understood as an ``external field'' version of the first question. This is because when the $t_i$ are small, the $M_i$ are more likely to be close to the identity, and this ``bias toward the identity'' is similar in spirit to the ``bias toward positive spin'' imposed in e.g.\ an Ising model with an external field.  The second question also arises naturally in {\em two-dimensional} Yang-Mills theory and has been heavily studied in that context \cite{Gross1989, Driver1989, fine1991, Witten91, migdal1996recursion, sengupta1997, levy2003, levy2010a, Levy2011a, Driver2017, shen2018, che19, driver2019functional, CCHS2D, dahlqvist2020, dahlqvist2022large}.

In two dimensions, the fine-mesh scaling limit of Yang-Mills theory is well understood, but if one attempts to compute the Wilson loop expectation for a complicated collection of loops (perhaps with many intersections and self-intersections) one obtains precisely an instance of the second problem above---indeed, the problems are equivalent since one can obtain {\em any} instance of the second problem for {\em some} two-dimensional loop.

In a recent companion paper~\cite{park2023wilson} (including some of the authors of this paper), it was shown that the answer to the second question can be expressed as an expectation with respect to a certain Poisson point process. In this paper, we will explain how the answer to the first question can be derived directly from the analysis in~\cite{park2023wilson} by taking $t_i \to \infty$. This is our first main result, which we state informally as follows. For a precise version, see Theorem \ref{thm:weingarten-recovery}. \revision{Later, we will show analogs of this theorem for orthogonal and symplectic groups, see Theorem \ref{thm:os-weingarten-recovery} for the precise statement.}

\begin{theorem}[Recovery of Weingarten calculus via Brownian motion]\label{thm:informal-weignarten-recovery}
The expectations of traces of words of Unitary Brownian motion converge as the time parameter goes to infinity to an explicit limit given in terms of the Weingarten function. Similar results hold for the other classical Lie groups.
\end{theorem}

\begin{remark}\label{remark:weingarten-recovery-remark}
We note that Theorem \ref{thm:weingarten-recovery} has previously appeared in \cite{Dahlqvist2017}, albeit stated in slightly different (but equivalent) terms -- see Sections 4 and 5 of the paper. Dahlqvist's proof relies heavily on representation theory. On the other hand, we believe that our proof may be easier to read for those who have a probability background but perhaps are not as familiar with representation theory. Additionally, our proof technique differs from Dahlqvist's in an essential way, which allows us to obtain the more general version of the Makeenko-Migdal/Master loop/Schwinger-Dyson equation for lattice Yang-Mills that we previously alluded to (Theorem \ref{thm:informal-master-loop}). See Remark \ref{remark:dahlqvist-comparision} for more discussion on the differences between the two arguments.
% Dahlqvist's proof relies heavily on representation theory. On the other hand, our proof is essentially entirely probabilistic (involving the analysis of a certain stochastic process), and we will explain all the representation theory prerequisites that are needed. For these reasons, we believe that our proof may be easier to read for those who have a solid probability background but perhaps are not that familiar with representation theory.
\end{remark}

% \scott{What does the following sentence refer to now? Should Remark 2.8 be brought sooner? Is this also rederived in Magee and Puder?} This is precisely the ``alternate route'' by which we re-derive some of the key results of Magee and Puder, and of Chatterjee and Jafarov.
Our approach to this result is in some sense very straightforward. The analysis in~\cite{park2023wilson} notes that when all of the $t_i$ are less than infinity, the noise generating the Lie group Brownian motion is a Gaussian white noise on a Lie algebra; because all randomness is Gaussian, all of the relevant quantities can be easily deduced from Wick's formula and planar maps (see the overview of these techniques \cite{zvonkin1997matrix}) which leads to a Poisson point process formulation of the theory. The analysis in this paper begins with the Poisson point process formulation obtained in~\cite{park2023wilson} and shows that geometric cancellations simplify in the $t_i \to \infty$ limit, so that the {\em Weingarten function} (as originally introduced in \cite{weingarten1978asymptotic}) appears naturally without any difficult computation. This approach also provides other insights -- for instance, certain single-edge analogs of the string-exploration steps in \cite{chatterjee2016} can be interpreted in terms of the so-called Jucys-Murphy elements \cite{jucys1974symmetric, murphy1981new, zinn2009jucys, novak2010jucys, Applebaum2011, matsumoto2013jucys}.

\subsection{Continuum Yang-Mills}
The famous continuum {\em Yang-Mills problem} \cite{Jaffe2006a} is (roughly speaking) to construct and understand the basic properties of a continuum analog of the lattice models described above, which should somehow make rigorous sense of the measure in~\ref{eqn-the-YM-measure}. This problem remains open for $d \geq 3$ and its solution for $d=4$ would in some sense also yield a solution to the quantum Minkowski version of Yang-Mills that forms the basis of the standard model in physics, see the Millennium Prize description \cite{Jaffe2006a}.
%\begin{center}
%\includegraphics[width=.36\textwidth]{matrixgrid.pdf}
%\end{center}

This paper is focused on understanding a {\em lattice} version of Yang-Mills theory in terms of sums over surfaces, with the aim of gaining insight into a possible continuum theory. It is not clear what kind of fine-mesh scaling limit one should expect the lattice models to have, but our hope is that the lattice analysis presented here will provide some clues, and we present several open problems along these lines in Section~\ref{sec::open}.

We remark that a number of purely continuum approaches to this problem are also being actively pursued.  For example, there is an {\em SPDE-based approach} which aims to construct a dynamical version of continuum Yang-Mills (on a torus, say) and show that it converges to a stationary law in the large time limit. One can take as the initial value a ``Lie-algebra-valued Gaussian free field connection'' that one expects to approximate the correct continuum theory at small scales and try to argue that the behavior at large scales converges to a limit over time. See e.g. \cite{CCHS2D, CCHS3D,
ChevyrevSurvey, chevyrev2023invariant, cao2021state,
cao2021yang,
BC2023, BC2024}. There has been some significant recent progress in this area, especially in two and three dimensions.

Alternatively, one can also work directly in the continuum {\em without} attempting to understand a dynamical process. One might regularize the continuum model in some other way---perhaps starting with a continuum Gaussian. Some form of this was implemented by Magnen, Rivasseau, and S\'{e}n\'{e}or \cite{MRS1993}. Some approaches along these lines might also be amenable to the type of random surface analysis discussed in this paper; see Section~\ref{sec::open}.

\subsection{Lattice models and planar maps}\label{section:edge-plaquette-embedding}
%There are three natural quantities to consider:

%\begin{enumerate}
%\item Expected Wilson loop trace in a lattice Yang-Mills model.
%\item Expected cycle-number power in a Poisson swap model.
%\item Surface weight sum in a random planar map model
%\end{enumerate}

%The above-mentioned companion paper proved the equivalence of the first and second items~\cite{park2023wilson}.  This paper will prove the equivalence of the second and third items. But we should also stress that the previous paper dealt with the finite $t$ case and we will focus here on the limiting $t \to \infty$ case. The {\bf random surfaces} treated here will be planar or higher genus maps that are ``embedded in'' (or mapped to) $d$-dimensional space.  Special embedded maps called {\bf edge-plaquette embeddings} will play a key role.

%In this paper, we will show what happens to the formulation in~\cite{park2023wilson} in the $t_i \to \infty$ limit, where the random matrix laws converge to Haar measure.  Our techniques rely primarily on Poisson point processes, planar maps and Wick's formula, whereas many other works have relied more heavily on representation theory, Weingarten calculus, and Haar measure integration by parts. We do obtain the {\em Weingarten function} in our end result, but our route to obtaining that function is somewhat different from the more representation-theoretic approaches.

\revision{Recall that a \textit{map} is a collection of polygons with a complete pairing of the set of edges, which should be thought of as 
a collection of specified gluings (see e.g. \cite{curien2016planar}). One should think of a map as a discretized surface. In particular, the ``surface sums" of the present paper will be sums over maps.}

\begin{figure}[ht!]
    \centering
\includegraphics[width=.8\textwidth]{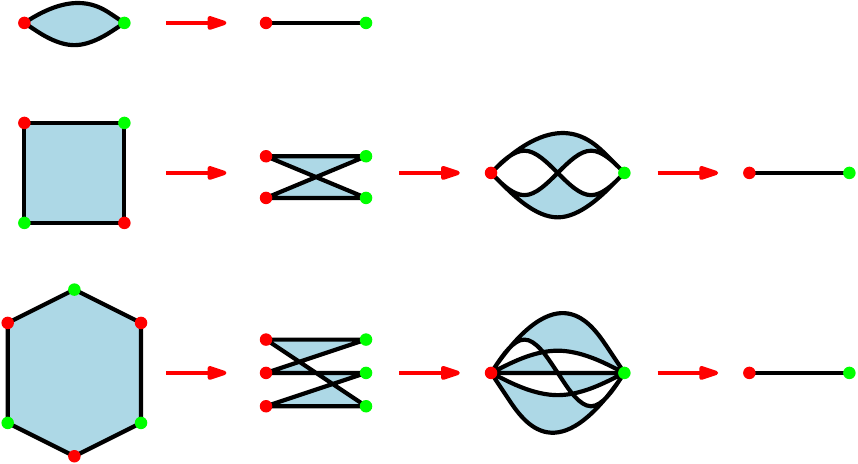}
    \caption{In an edge-plaquette embedding, we can imagine that each blue face is ``twisted and collapsed'' onto a single edge, see Figure~\ref{fig::folding}. In the sequence above, we first twist, then collapse matching vertices, then collapse edges. \revision{This is meant to illustrate condition 4 in the definition of edge-plaquette-embedding. Also, due to condition 2, blue faces are even degree, and so we can also color its vertices in alternating colors (red and green in the present example).}}
    \label{fig::flatten}
\end{figure}

\begin{figure}[ht!]
    \centering
\includegraphics[width=.6\textwidth]{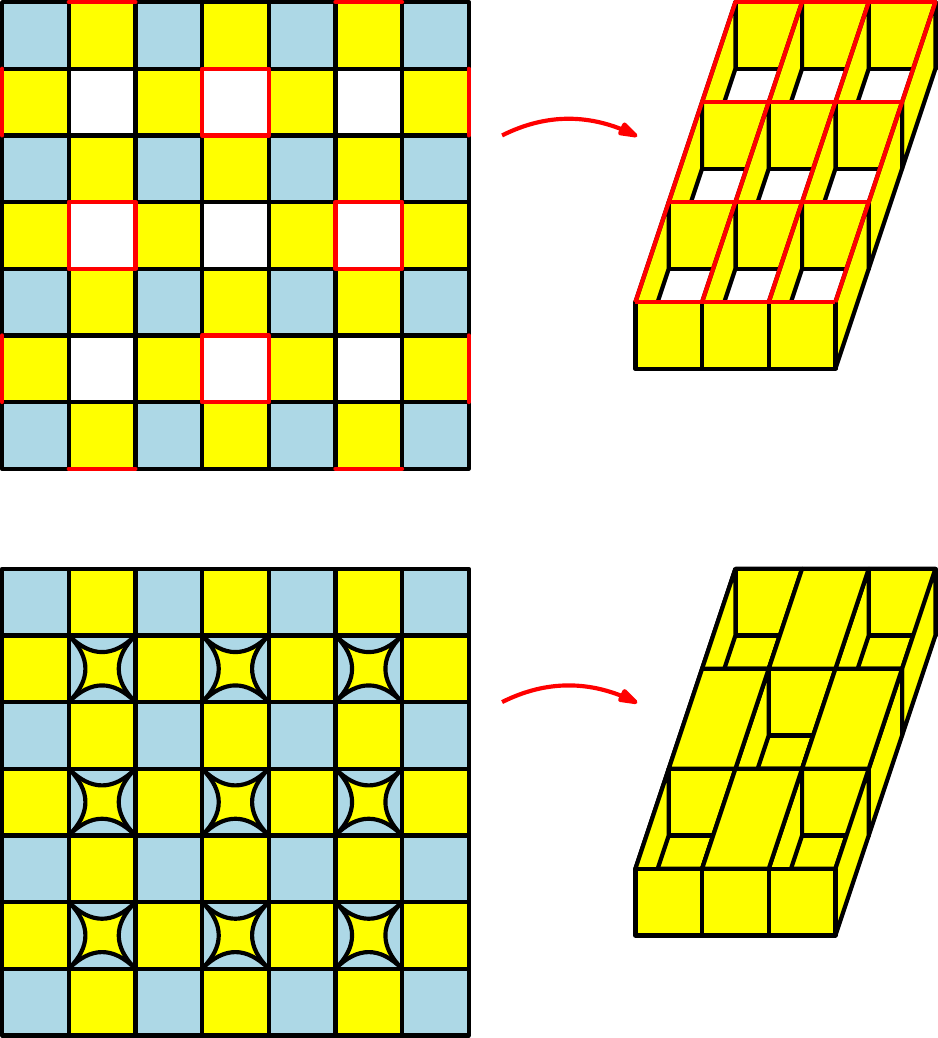}
    \caption{{\bf Edge-plaquette embedding example:} Each of the 16 blue faces on the upper left gets mapped to a single vertical edge in the upper right, while each yellow face on the upper left gets mapped to a vertical yellow face on the upper right---the edge colored red is the one mapped to the top. On the lower left, additional yellow faces are added; their images on the right alternate between upper and lower layers in checkerboard fashion. Going from left to right requires ``folding up'' the blue squares and collapsing the blue 2-gons. \revision{We emphasize that the pictures on the right are drawn in the lattice $\Lambda$.}}
    \label{fig::folding}
\end{figure}

\begin{figure}[ht!]
    \centering
\includegraphics[width=.4\textwidth]{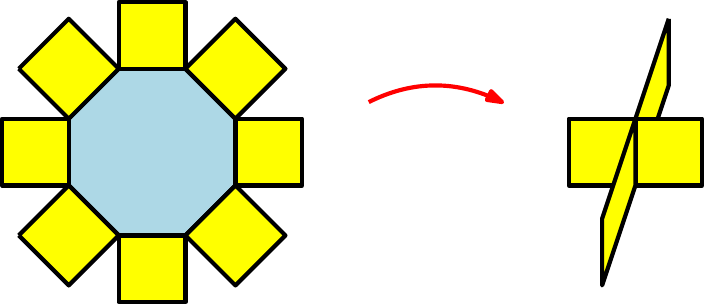}
    \caption{{\bf Edge-plaquette embedding example:} If the blue face is an octagon, then there will be 8 yellow plaquettes meeting at the corresponding edge. In this example shown, the pre-image of each yellow face on the right may consist of two yellow faces on the left. In other words, there are two ``copies'' of each of the four plaquettes shown on the right. \revision{We emphasize that the pictures on the right are drawn in the lattice $\Lambda$.}}
    \label{fig::octagon}
\end{figure}

\begin{figure}[ht!]
    \centering
\includegraphics[width=.8\textwidth]{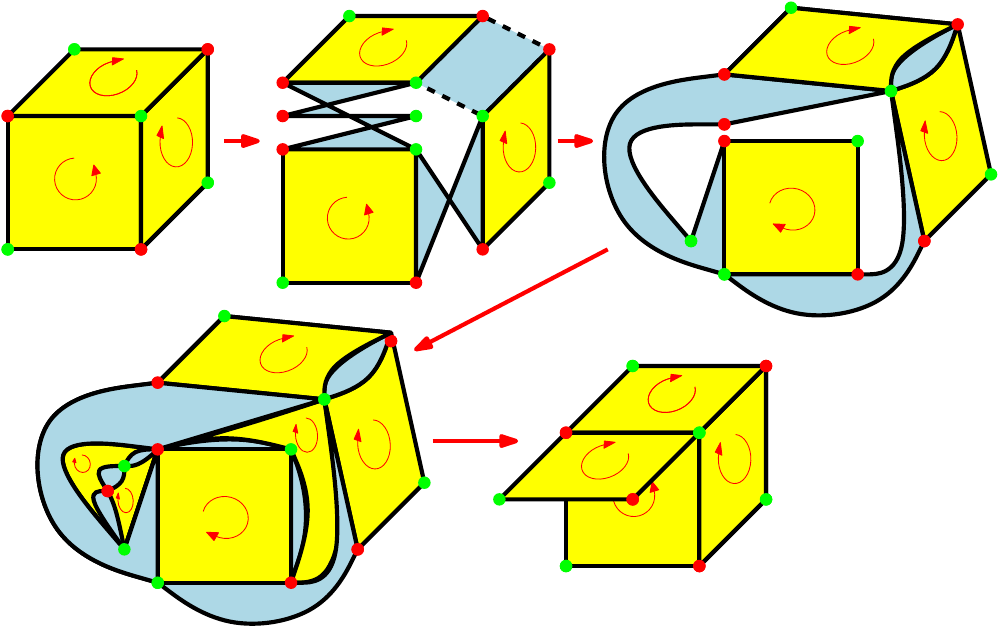}
    \caption{{\bf Edge-plaquette embedding example showing orientations:} (1) Three oriented plaquette images in $\phi(\mc M)$. (2) The blue faces connecting them have different types. (3) ``Untwist'' by flipping the lower-left plaquette across its red-red diagonal so that the three red and three blue faces are orientably embedded in the plane. (4) Add some new faces (three yellow squares and five blue $2$-gons) to fill in the hole. Interpret the resulting colored map as a portion of $\mathcal M$ orientably embedded in the plane. (5) Map this portion back into the lattice. Not all six yellow plaquettes are visible on the right because some overlap each other.}
    \label{fig::twistedfaces}
\end{figure}

Consider a pair $(\mathcal M, \phi)$ where $\mathcal M$ is a planar (or higher genus) map and $\phi:\mathcal M \to \Lambda$ is a graph homomorphism.\footnote{In other words, if two vertices $v,w\in V(\mathcal M)$ are adjacent in the graph $\mathcal M$, then $\phi(v),\phi(w)$ are adjacent in the graph $\phi(\mathcal M)\subset\Lambda$.} We call this pair an {\bf edge-plaquette embedding} if the following hold:
\begin{enumerate}
\item The dual graph of $\mathcal M$ is bipartite. The faces of $\mathcal M$ in one partite class are designated as ``edge-faces'' (shown blue in figures) and those in the other class are called ``plaquette-faces'' (shown yellow in figures).
\item \revision{All edge-faces have even degree.}
\item $\phi$ maps each plaquette-face of $\mathcal M$ isometrically {\em onto} a plaquette in $\mathcal P$.
\item $\phi$ maps each edge-face of $\mathcal M$ onto a single edge of $\Lambda$.
\end{enumerate}
See Figures~\ref{fig::flatten}-\ref{fig::twistedfaces} for examples and intuition.

\revision{Given an edge plaquette embedding $(\mathcal{M}, \phi)$, we denote by $\phi^{-1} : \mathcal{P} \ra \N$ the associated ``plaquette  count", i.e. the function which assigns to each oriented plaquette $p \in \mathcal{P}$ the number of copies of $p$ in $(\mathcal{M}, \phi)$. In other words, $\phi^{-1}(p)$ is the number of plaquette-faces of $\mathcal{M}$ which are mapped isometrically onto the plaquette $p$. We also write $(\phi^{-1})! := \prod_{p \in \mathcal{P}} (\phi^{-1}(p))!$.}

In order to construct a model of random edge-plaquette-embedding that is useful in Yang-Mills theory\revision{, w}e will need to assign a ``weight'' to every plaquette $p$ of $\Lambda$ (depending on the number of copies of $p$, i.e. the value of $K(p)$ in \eqref{eqn::latticeymexpanded}) and every edge (depending on the number and type of blue faces there). This weight is closely related to the so-called Weingarten function, which we discuss next.

\subsubsection{Weingarten function}\label{section:weingarten-discussion}

% In this section, we recall the definition of Weingarten functions that defines the weight of each edge-plaquette embedding. To avoid technical details in representation theory, we adapt the definition of Weingarten functions from~\cite{zinn2009jucys} in terms of the Gram matrix. For alternative definitions in terms of character expansion, e.g.\ see~\cite{collins2006integration, collins2022weingarten, kostenberger2021weingarten}.

Throughout this paper, let $\symgrp_n$ be the symmetric group on $n$ elements. A complex-valued function on $\symgrp_n$ can be identified as an element in the group algebra $\mathbb C[\symgrp_n]$, that is $\sigma\mapsto f(\sigma)$ is identified as $\sum_{\sigma\in \symgrp_n} f(\sigma)\sigma$. 
% Let $\mathbb Q[N]\subset \mathbb C[N]$ be the field of rational functions with rational coefficients in the variable $N$.
When $N \geq n$, the \textbf{Weingarten function} $\Wg_N$\footnote{We omit the dependence on $n$ for brevity.~} can be defined as the inverse in \revision{the group algebra $\C[\symgrp_n]$ of the function $\sigma \mapsto N^{\#\textrm{cycles} (\sigma)}$.}
% the group ring of $\mathbb Q(N)[\symgrp_n]$ of the function $\sigma \to N^{\#\textrm{cycles} (\sigma)}$. 
(There is a slightly different definition for $N < n$, see Section~\ref{sec::mainresults}.)
Note that $\Wg_N(\sigma)$ depends only on the conjugacy class of $\sigma$, i.e.\ on the cycle structure of $\sigma$. We can order cycles from biggest to smallest, represent this by a Young diagram, and interpret $\Wg_N$ as a function on Young diagrams. It is not the simplest function, and one explicit formula which works for general (\revision{i.e., even if $N < n$}) values of $N$ (see \cite[Equation (9)]{collins2006integration} and \cite[Section 2.1]{Magee2019})
% when $N \geq n$ (see the overview and additional references in \cite{Magee2019}) 
is as follows:
\begin{equs}\label{eq:weingarten-character-sum}
\Wg_N(\sigma) = \frac{1}{n!} \sum_{\substack{\lambda \vdash n \\ \ell(\lambda) \leq N }} \Bigl[ \chi_\lambda(\id) \chi_\lambda(\sigma) \prod_{(i,j) \in \lambda}(N+j-i)^{-1} \Bigr]
\end{equs}
where $\id$ is the identity permutation, $\lambda \vdash n$ denotes that $\lambda$ is a partition of $n$, $\ell(\lambda)$ is the number of rows of the Young diagram $\lambda$, $\chi_\lambda(\sigma)$ is the character (trace of $\sigma$ in the irreducible representation indexed by $\lambda$), \revision{and $(i, j) \in \lambda$ denotes a box of the Young diagram $\lambda$, where $i$ is the row index and $j$ is the column index. As an example, for the Young diagram consisting of a single row of length $n$, the boxes are indexed by (in order from left to right) $(1, 1), (1, 2), \ldots, (1, n)$.} Observe that $\lambda \vdash n$ implies $\ell(\lambda) \leq n$, and so if $N \geq n$, the condition $\ell(\lambda) \leq N$ is automatically satisfied.

Alternatively, letting $f \in \C[\symgrp_n]$ be the function $\sigma \mapsto N^{\cycles(\sigma)}$, observe that we may write $f = N^n(\id + g)$, where $g = N^{-n}(f - \id)$. Thus the inverse of $f$ may be formally expanded as
\begin{equation} \label{eqn::expandedinverse}
\Wg_N = f^{-1} = N^{-n} \big(\id + (-g) + (-g)^2 + \cdots \big). \end{equation} 
This observation plays a key role in the derivation of \cite[Theorem 2.9]{Magee2019}.

As will be explained in Section \ref{section:wilson-loop-expectation-epe}, we interpret $\sigma$ as a collection of blue faces (one blue face of length $2k$ for each cycle of $\sigma$ of length $k$). Then $\Wg_N(\sigma)$ is essentially the {\em weight} associated to given collection of blue faces at an edge. Actually, as we detail in Section \ref{section:wilson-loop-expectation-epe}, the edge weights are given by the {\em normalized Weingarten function}, which we define as
\begin{equs}
\ovl{\Wg}_N(\sigma) := N^{2n - \cycles(\sigma)} \Wg_N(\sigma).
\end{equs}
This is the normalization which leads to a nontrivial $N \toinf$ limit (see Remark \ref{remark:weingarten-large-N-limit}).

\revision{By \eqref{eq:weingarten-character-sum}, we see that $\Wg_N$ is a linear combination of characters, and thus $\Wg_N$ is a class function. Recalling that the conjugacy classes of $\symgrp_n$ are indexed by partitions $\mu$ of $[n]$, we may think of $\Wg_N$ as a function on the set of partitions of $[n]$, whose value $\Wg_N(\mu)$ at a given partition $\mu$ is defined to be $\Wg_N(\sigma)$, where $\sigma \in \symgrp_n$ is any permutation whose cycle structure is given by $\mu$ (since two permutations are in the same conjugacy class if and only if they have the same cycle structure).}

Given an edge-plaquette embedding $(\mc{M}, \phi)$ and an edge $e$ of our lattice $\Lambda$, we will write $\ovl{\Wg}_N(e)$ as shorthand for $\ovl{\Wg}_N(\mu_e(\mc{M}, \phi))$, where $\mu_e(\mc{M}, \phi)$ is the partition given by half the degrees of the edge-faces mapped to $e$. \revision{To expand on this, recall from the definition of edge-plaquette embedding at the beginning of Section \ref{section:edge-plaquette-embedding} that all edge-faces have even degree. Thus, if we take the collection of numbers given by half the degrees of the edge-faces mapped into a given lattice edge $e$, we obtain a partition $\mu_e(\phi)$ of some positive integer, say $n_e$. Then by the preceding paragraph, we may consider $\ovl{\Wg}_N(\mu_e(\phi))$, i.e. the Weingarten function applied to a partition of $[n_e]$.}

\subsection{Main results}

We already informally stated the first of our main results -- recall Theorem \ref{thm:informal-weignarten-recovery}. In this subsection, we proceed to state the remaining main results of this paper. The next result is the main conceptual contribution of the paper. We state it informally. See Theorem~\ref{thm:wilson-loop-expectation-sum-over-epe} for the corresponding precise statement. \revision{In the following, recall the notation $(\phi^{-1})!$ from Section \ref{section:edge-plaquette-embedding}.}

\begin{theorem}[Surface-sum representation of Wilson loop expectations]\label{thm:informal-wilson-loop-expectation-sum-over-epe}
\revision{Let $\Lambda$ be a finite lattice.} When the gauge group is $\UN$, the Wilson loop expectation $\langle W_{s} \rangle_{\Lambda, \beta}$ is proportional to 
\[ \sum_{\ptl (\mc{M}, \phi) = s} \frac{\beta^{\mrm{area(\mc{M}, \phi)}}}{(\phi^{-1})!} \Big(\prod_e \ovl{\Wg}_N(e)\Big) \cdot N^{\chi(\mc{M}) - k},\]
where the sum is over edge-plaquette embeddings $(\mc{M}, \phi)$ with boundary $s$, $\mrm{area}(\mc{M}, \phi)$ is the total number of plaquettes in the edge-plaquette embedding, $(\phi^{-1})!$ is a combinatorial factor depending only on the plaquette counts, $\chi(\mc{M})$ is the Euler characteristic of $\mc{M}$, and $k = |s|$ is the number of loops. \revision{The sum converges absolutely in the following sense: the set of $(\mc{M}, \phi)$ with boundary $s$ may be partitioned based on the value of the plaquette count $\phi^{-1} = K$, and we have that
\[ \sum_{K : \mc{P} \ra \N}  \Bigg| \sum_{\substack{\ptl (\mc{M}, \phi) = s \\ \phi^{-1} = K}} \frac{\beta^{\mrm{area(\mc{M}, \phi)}}}{K!} \Big(\prod_e \ovl{\Wg}_N(e)\Big) \cdot N^{\chi(\mc{M}) - k} \Bigg| < \infty.\]}
\end{theorem}

% First, when computing Wilson loop expectations, we imagine the simplest setting in which we fix the number of yellow faces of each type (i.e.\ assign weight 1 to that number and 0 to all others). This corresponds to focusing on a single summand in \eqref{eqn::latticeymexpanded}, or equivalently to taking $\beta = 0$ in \eqref{eqn::latticeymexpanded}. In this case we have the following:

% \begin{theorem}[Surface-sum representation for word expectations]\label{thm:informal-word-expectation-as-sum-over-bipartite-maps}
% When the gauge group is $\UN$, the expected trace product is proportional to $\sum \big(\prod_e \ovl{\Wg}_N(e)\big) \cdot N^{\chi - 2k}$ where the sum is over spanning edge-plaquette embeddings with given plaquette numbers, $\chi$ is the Euler characteristic and $k = |\mc{L}|$ is the number of loops. 
% % (Variants apply to other gauge groups.)
% \end{theorem}

% For a precise statement of this theorem, see Proposition \ref{prop:word-expectation-as-sum-over-bipartite-maps}. 

\begin{remark}\label{remark:main-point}
Theorem \ref{thm:informal-wilson-loop-expectation-sum-over-epe} is the main conceptual contribution of this paper, because it introduces the new concept of an edge-plaquette embedding, and gives a fundamentally new description of Wilson loop expectations in terms of random planar maps\footnote{This is rather loose terminology, as our surface sums are signed, and in general higher genus surfaces may appear.}, thereby connecting two very different areas of research. Ultimately, we hope to prove new results about lattice gauge theories via analysis of these random planar maps, in particular building on the many advances in their understanding -- see Section \ref{sec::open} for some open problems. See also Remark \ref{remark:taggi}.
\end{remark}

\begin{remark}
The results of Magee and Puder \cite{Magee2019, Magee2019a} could also be applied here to give a surface sum representation of the terms in \eqref{eqn::latticeymexpanded}. However, the relation to random planar maps is not as clear in their formulation. As mentioned in Remark \ref{remark:main-point}, this is the main point of our result. For more comparison with Magee and Puder, see Section \ref{section:discussion-of-magee-puder}.
\end{remark}

\begin{remark}
We will explain in Section~\ref{sec::othergroups} the variants of this theorem that apply when $\UN$ is replaced with another compact Lie groups. 
Even in the $\UN$ case there are several variants to this result. The various ``string trajectory moves'' in \cite{chatterjee2016} can be interpreted in terms of the exploration of a surface built out of blue $2$-gons and $4$-gons and yellow squares. One can also interpret the individual Jucys-Murphy elements in these terms.
\end{remark}

\begin{remark}
\revision{Due to the gauge symmetry present in Yang--Mills, Wilson string expectations do not change even if in the integral \eqref{eqn::latticeymexpanded} we enforce that all the edges in a spanning tree of the lattice are assigned the identity matrix. See e.g. \cite[Section 9]{chatterjee2016leading} for a proof of this for a particular choice of spanning tree (but the proof directly extends for general spanning trees). In terms of the surface sum in Theorem \ref{thm:informal-wilson-loop-expectation-sum-over-epe}, this gauge symmetry implies that for any choice of spanning tree, the surface sum is equal to a modified surface sum where all edges in the spanning tree are collapsed to a single vertex. To directly see this from the surface sum, one would have to use certain properties of the Weingarten weights $\ovl{\Wg}_N$ to argue that collapsing all edges in a spanning tree leaves the total surface sum unchanged.}
\end{remark}

% By applying Theorem \ref{thm:informal-word-expectation-as-sum-over-bipartite-maps} to every term in the series appearing in equation \eqref{eqn::latticeymexpanded}, we obtain that Wilson loop expectations may be expressed as a weighted sum over edge-plaquette embeddings. We state this informally as the following corollary.

% We note that $\mrm{area}$ and $\chi$ only depend on $\mathcal{M}$, which can be considered a covering space of the embedded map. On the other hand, $K!$ and the product weight also depend on the embedding $\phi$.

\begin{remark}\label{remark:taggi}
Recently, Taggi and coauthors~\cite{LeesTaggi2020, LeesTaggi2021, QuitmannTaggi2023} have succeeded in proving various results about spin $O(n)$ and related models by analyzing a certain related random path (or random loop) model. Starting from the spin $O(n)$ model, they arrive at their random path model in exactly an analogous manner as how we arrive at Theorem~\ref{thm:informal-wilson-loop-expectation-sum-over-epe}. Namely, starting from the action for the spin $O(n)$ model, which at a single edge is of the form $\exp(\beta \sigma_x \cdot \sigma_y)$, where $\sigma_x, \sigma_y \in S^n$ (here $S^n$ is the unit $n$-sphere), they expand 
\[ \exp(\beta (\sigma_x \cdot \sigma_y)) = \sum_{k=0}^\infty \frac{\beta^k}{k!} (\sigma_x \cdot \sigma_y)^k\] for each edge $(x, y)$, and then compute the resulting $S^n$-integrals. The $S^n$-integrals may be easily computed, with the resulting expressions only involving very explicit quantities such as factorials and Gamma functions (see~\cite[equation (2.12)]{LeesTaggi2021}). This is one simplification compared to our setting, where the $\unitary(N)$-integrals lead to the appearance of the Weingarten function, which is much more complicated to understand. Another key difference is that while the $S^n$-integrals are always positive, the $\unitary(N)$-integrals may be both positive and negative. Thus the random path model of Taggi et al.\ may be interpreted as a genuine probability measure, while our surface sums may only be interpreted as signed measures. 
% Nevertheless, we introduce our surface sum model for Wilson loop expectations with the ultimate goal of using them to understand lattice Yang-Mills. 
\end{remark}

Next, we give an informal statement of the Makeenko-Migdal/Master loop/Schwinger-Dyson equations satisfied by Wilson loop expectations. The corresponding precise statement is Theorem~\ref{thm:master-loop}.
% We note that our recursion is slightly different from and more general than the existing literature -- see Remark \ref{remark:schwinger-dyson}.

\begin{theorem}[Single-location Makeenko-Migdal/Master loop/Schwinger-Dyson equation]\label{thm:informal-master-loop}
Wilson loop expectations satisfy the following recursion:
\begin{equs}
\langle W_{s} \rangle_{\Lambda, \beta} = \mrm{splitting} + \mrm{merger} + \mrm{deformation}
\end{equs}
\end{theorem}
Here, the splitting, merger, and deformation terms correspond to certain types of operations we may apply to a given collection of loops $s$ to obtain a new collection of loops. They will be precisely defined in Section \ref{section:master-loop}.

\begin{remark}
As previously mentioned, versions of this recursion for various Lie groups have previously appeared \cite{chatterjee2016,chatterjee2016a, jafarov2016, Chatterjee2019a, shen2022new}. We note that the precise form of our recursion is slightly different from (and more general than) the existing literature -- see Remarks \ref{remark:schwinger-dyson} and \ref{remark:schwiner-dyson-2}. Ultimately, the reason for this difference is due to our proof method. Whereas previous approaches are based on integration-by-parts\footnote{The argument in \cite{shen2022new} essentially reduces to integration by parts, as explained in \cite[Appendix A.2]{AN2023}. To sketch the argument, the proof uses the fact that if the lattice Langevin dynamics is started at stationarity,  then the expectation of any observable must be constant in time. Then, applying It\^{o}'s formula to Wilson loop observables, one obtains an identity saying that the drift term must have expectation zero.
This identity is precisely integration by parts.}, our approach is essentially equivalent to applying a certain recursion that is satisfied by the Weingarten function (see e.g.\ \cite[Proposition 2.2]{Collins2017}), which comes out of our particular approach towards proving Theorem \ref{thm:informal-weignarten-recovery} (recall Remark \ref{remark:weingarten-recovery-remark}, and see also Remark \ref{remark:dahlqvist-comparision}).
% , although we don't phrase our argument in this way -- we prefer to proceed more probabilistically via our aforementioned Poisson point process formulation.

\revision{After the present article first appeared, the single-location master loop equation has since found application in the work of \cite{SSZ2024}, where it was crucially used to derive the continuum master loop equation as a scaling limit of lattice master loop equations.}
\end{remark}

% Having completed the statements of our main results, we next remark how we expect our results to change when we consider other Lie groups such as $\SUN$, $\SON$, or $\SphN$. In a forthcoming revision, we will include results for these other cases.

% \begin{remark}[Other Lie groups]\label{remark:other-groups}
% First, the results of \cite{park2023wilson} also include the cases of $\SUN$, $\SON$, or $\SphN$. That is, expectations of products of words of independent Brownian motions are expressed as certain Poisson process expectations. In particular, \cite{park2023wilson} explains the key differences in terms of how one associates weights (depending on the Lie group) to a given realization of a Poisson process. Once these different weights are accounted for, the proof of Theorem \ref{thm:informal-master-loop} would proceed in exactly the same manner. On the other hand, with regards to Theorem \ref{thm:informal-word-expectation-as-sum-over-bipartite-maps} and Corollary \ref{thm:informal-wilson-loop-expectation-sum-over-epe}, the analogs of the Weingarten function one obtains (and the corresponding planar map interpretation) are different -- see e.g. \cite{collins2006integration, Dahlqvist2016, Magee2019a}.
% \end{remark}

\subsubsection{Discussion of Magee and Puder}\label{section:discussion-of-magee-puder}
The vocabulary in~\cite{Magee2015, Magee2019} is somewhat different from ours, but the results can be expressed in similar terms. We won't give a detailed account of those results, but let us briefly outline a couple of key ideas to assist readers trying to compare their approach to ours. The approach in~\cite{Magee2015, Magee2019} makes heavy use of commutator words. Suppose a loop $\ell$ in $s$ corresponds to a commutator word $ABA^{-1}B^{-1}$ (where $A$ and $B$ could in principle describe paths of length longer than one).  Imagine then that we have a surface $S$ with a single boundary loop, whose boundary is mapped to $\ell$. We can turn this surface into a closed surface in two ways. First, we can identify the boundary of an ordinary disk (with circular boundary) with the boundary of $S$, thereby gluing a circular disk onto $S$. Second, we can glue the boundary of $S$ to itself by first gluing the pre-images of the $A$ and $A^{-1}$ segments to each other and then gluing the pre-images of the $B$ and $B^{-1}$ segments to each other---which somehow turns the disk bounded by $\ell$ into a torus. It is not hard to see that the second approach produces a surface whose genus is $1$ higher than the surface produced by the first approach: it effectively ``adds a handle'' to the surface. If we write a long loop $\ell$ as a product of $n$ commutator words, then those words provide us a recipe for turning a disk bounded by $\ell$ into an $n$-holed torus (by performing gluings of the type mentioned above for each commutator).

An intermediate result that we use to prove Theorem~\ref{thm:informal-wilson-loop-expectation-sum-over-epe} (see Proposition \ref{prop:word-expectation-as-sum-over-bipartite-maps}) is closely related to \cite[Theorem 2.8]{Magee2019}. We remark that one could also interpret \cite[Theorem 2.9]{Magee2019} in terms of embedded maps (somehow involving multiple layers of blue faces). We note that \cite[Theorem 2.9]{Magee2019} is in some ways simpler than \cite[Theorem 2.8]{Magee2019} (it does not involve the Weingarten function) and in other ways more complicated (it leads to \cite[Theorem 1.7]{Magee2019}, which involves another quantity called the $L^2$ Euler characteristic, which is in general not so trivial). We note that \cite[Theorem 2.9]{Magee2019} is derived from \cite[Theorem 2.8]{Magee2019}. We will not give an alternate derivation of this step, aside from remarking that the expansion in \eqref{eqn::expandedinverse} plays a role.

\subsection{Summary of paper and reading guide}

We close this section off with a summary of the rest of the paper. In Section \ref{sec::mainresults}, we introduce the notation and background material that will be needed in the rest of the paper. In Section \ref{section:wilson-loop-expectation-epe}, we derive our surface-sum representation of Wilson loop expectations. In Section \ref{section:poisson-exploration}, we show how to recover the Weingarten calculus by taking limits of Unitary Brownian motion, using a certain strand-by-strand exploration that we introduce in the section. In Section \ref{section:master-loop}, we apply our strand-by-strand exploration to obtain the single-location Makeenko-Migdal/Master loop/Schwinger-Dyson equation for Wilson loop expectations. Finally, in Section \ref{sec::othergroups}, we adapt our results to the cases of $G = \orthogonal(N), \SphN, \SUN, \SON$.

To the reader who wants to understand our surface-sum representation of Wilson loop expectations as quickly as possible, we recommend the following expedited reading strategy. First, read enough of Section \ref{section:strand-diagrams-and-weingarten-caculus} to understand the statement of the Unitary Weingarten calculus (Theorem \ref{thm:weingarten}). Then, proceed directly to Section \ref{section:wilson-loop-expectation-epe} to see how this theorem is applied to obtain the surface-sum representation. 
% This is roughly ten pages of material.

\section{Notation and background} \label{sec::mainresults}

We introduce some basic notation and background that will be needed throughout this paper.

\begin{itemize}
 \item For $n \in \mathbb{N}$, let $[n] := \{1, \ldots, n\}$.
 \item For $a,b \in \Z$, $a < b$, we denote $(a:b] := \{a+1, \ldots, b\}$. So $[n] = (0:n]$.
 % \item For a set $A$, we let $\binom{A}{2}$ denote the unordered set of ordered pairs of elements of $A$.
\end{itemize} 

\subsection{Strand diagrams and Weingarten calculus}\label{section:strand-diagrams-and-weingarten-caculus}

In this section, we describe how to interpret the Unitary Weingarten calculus in terms of diagrammatic sums. This is an alternate interpretation of the usual way of stating the Weingarten calculus in terms of expectations of products of matrix entries (as is e.g. done in the survey \cite{collins2022weingarten}). We find this alternate interpretation particularly useful for obtaining surface-sum representations of Wilson loop expectations.

Let $\bm \Gamma = (\Gamma_1, \Gamma_2, \dots, \Gamma_k)$ be an (ordered) collection of words $\Gamma_i$ of lengths $m_i$ on letters $\{\lambda_1,\cdots,\lambda_L\}$, where
\begin{equs}\label{eq:general-word-form}
\Gamma_i = \lambda_{c_i(1)}^{\ep_i(1)}\cdots \lambda_{c_i(m_i)}^{\ep_i(m_i)}
\end{equs}
for some $c_i: [m_i]\to [L]$ and $\ep_i: [m_i]\to \{-1,1\}$. For each letter $\lambda$, let $n_+(\lambda)$ (resp. $n_-(\lambda)$ be the total number of times $\lambda$ (resp $\lambda^{-1}$) appears in $\bm \Gamma$.

The {\it strand diagram} corresponding to $\bm \Gamma$, denoted $\mrm{SD}(\bm \Gamma)$, is a directed graph with two types of edges, constructed as follows. To start, for each letter $\lambda$, we draw $n_+(\lambda)$ right-directed edges and $n_-(\lambda)$ left-directed edges. Such edges are drawn as solid black edges, and we refer to them as strands. Define $n_{\pm}(\lambda) := n_+(\lambda) + n_-(\lambda)$.

\begin{figure}[h]
    \centering
    \includegraphics[page = 1, width=10cm]{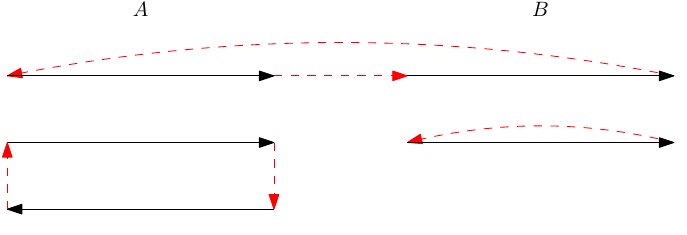}
    \caption{In this example, we have two letters $\{A, B\}$, and three words $(AB, AA^{-1}, B)$. The dashed red edges are drawn so that after ``shrinking" them away, what is left is exactly the three words.}\label{figure:strand_diagram_example}
\end{figure}

Next, we add edges of a different type to the graph, which we draw as dashed red edges, so that the connected components of the resulting graph exactly correspond to the words in $\bm \Gamma$. This is perhaps best explained with an example -- see Figure \ref{figure:strand_diagram_example}.

We refer to the dashed red edges as the {\it exterior connections} of the strand diagram, because they may connect vertices of strands associated to different letters.

\begin{remark}
In principle, there are many different ways one could place the dashed red edges according to $\bm \Gamma$. For instance, in Figure \ref{figure:strand_diagram_example}, we could permute the two strands corresponding to the letter $B$. Let us just say that for each $\bm \Gamma$, we fix one way of placing the dashed red edges, so that there is a unique strand diagram $\mrm{SD}(\bm \Gamma)$ associated to each $\bm \Gamma$.
\end{remark}

Next, we describe another collection of edges that we will add to the strand diagram $\mrm{SD}(\bm \Gamma)$. We refer to these edges as the {\it interior connections} of the strand diagram, because they only connect vertices of strands associated to the same letter. When drawn on top of the strand diagram, we will always draw them as blue edges, to distinguish from the (black) strands and (red) exterior connections. When drawing these edges by themselves, we will default to using the color black.

\begin{figure}[h]
    \centering
    \includegraphics[width=8cm]{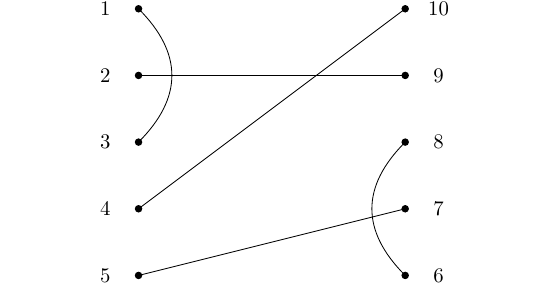}
        \caption{$\pi = \{\{1, 3\}, \{2, 9\}, \{4, 10\}, \{5, 7\}, \{6, 8\}\}$}\label{figure:matching-example}
\end{figure}

\begin{definition}\label{def:matchings}
Let $n \geq 1$ be an integer. Let $\mc{M}(n)$ be the set of matchings $\pi$ of $[2n]$, i.e. partitions of $[2n]$ into two-element sets. We view matchings pictorially as in Figure \ref{figure:matching-example}.
\end{definition}

\begin{figure}[h]
    \centering
    \includegraphics[page = 2, width=10cm]{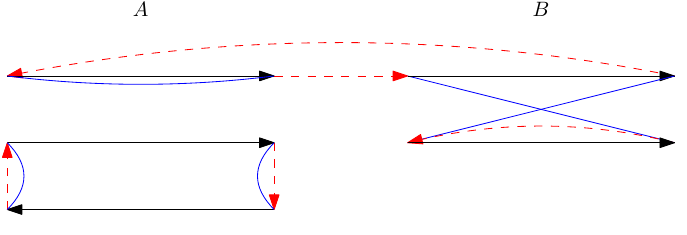}
    \caption{Some particular choice of matchings $\bm \pi = (\pi_A, \pi_B)$ may lead to this graph, starting from Figure \ref{figure:strand_diagram_example}.}\label{figure:strand_diagram_example_2}
\end{figure}

Suppose we have a collection of matchings $\bm \pi = (\pi_\ell, \ell \in [L])$, where for each $\ell \in [L]$, $\pi_\ell \in \mc{M}(n_{\pm}(\lambda))$ is a matching on the vertices of the  strand diagram associated to $\lambda_\ell$. For instance, we might obtain Figure \ref{figure:strand_diagram_example_2} from the strand diagram in Figure \ref{figure:strand_diagram_example}.

\begin{figure}[h]
    \centering
    \includegraphics[page = 3, width=10cm]{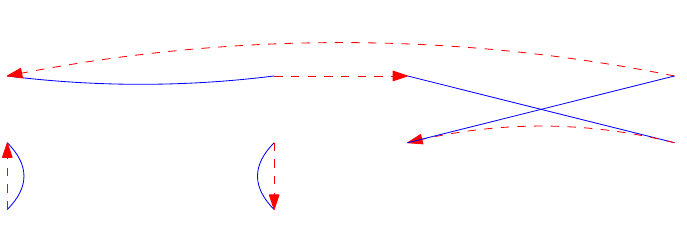}
    \caption{By deleting the strands from Figure \ref{figure:strand_diagram_example_2}, we obtain the above graph, which we see has two components. Thus in this case, $\numcomp(\bm \Gamma, \bm \pi) = 2$.}\label{figure:strand_diagram_example_3}
\end{figure}

Given $(\bm \Gamma, \bm \pi)$, let $\numcomp(\bm \Gamma, \bm \pi)$ be the number of components of the graph obtained by starting with the strand diagram $\mrm{SD}(\bm \Gamma)$, adding in the interior connections $\bm \pi$, and deleting the strands. For example, from Figure \ref{figure:strand_diagram_example_2}, we obtain Figure \ref{figure:strand_diagram_example_3}.

We are now almost ready to state the Unitary Weingarten calculus. We first need to make some assumptions on the collection of words $\bm \Gamma$ and interior connections $\bm \pi$.

\begin{definition}
We say that a collection of words $\bm \Gamma$ is {\it balanced} if $n_+(\lambda) = n_-(\lambda)$ for all letters $\lambda$.
\end{definition}

In terms of the strand diagram, a balanced collection of words $\bm \Gamma$ leads to a strand diagram where for each letter, the number of right-directed edges is equal to the number of left-directed edges. We now introduce a special collection of interior connections on such a strand diagram.

\begin{figure}[h]
    \centering
    \includegraphics[page=3]{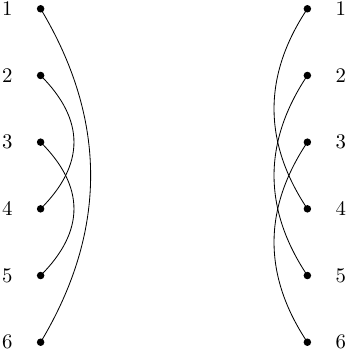}
    \caption{Example of an element of $\mc{M}(2n)$ defined by a pair of bijections $(\sigma, \tau)$. Here $n = 3$,  $\sigma$ maps $1 \mapsto 6$, $2 \mapsto 4$, $3 \mapsto 5$, and $\tau$ maps $1 \mapsto 4$, $2 \mapsto 5$, $3 \mapsto 6$. }\label{figure:sigma-tau-matching-example}
\end{figure}

\begin{definition}\label{def:left-right-bijection}
Let $n \geq 1$. Let $\sigma, \tau : [n] \ra (n : 2n]$ be bijections. The pair $(\sigma, \tau)$ defines a matching in $\mc{M}(2n)$ by viewing $\sigma$ as specifying a partition of the left vertices into two-element sets, and $\tau$ as specifying a partition of the right vertices into two-element sets. We denote the matching induced in this manner by $[\sigma ~ \tau] \in \mc{M}(2n)$. See Figure \ref{figure:sigma-tau-matching-example} for a pictorial example.
\end{definition}

Observe that if $\sigma, \tau : [n] \ra (n: 2n]$ is a pair of bijections, then $\sigma^{-1} \tau : [n] \ra [n]$ is a bijection, i.e. an element of $\symgrp_n$.

We can now finally state the Unitary Weingarten calculus. Let $(U(\lambda_\ell), \ell \in [L])$ be an i.i.d. collection of random Haar-distributed $\UN$ matrices. Define the shorthand (recalling \eqref{eq:general-word-form})
\begin{equs}
\Tr(U(\bm \Gamma)) &:=  \Tr(U(\Gamma_1)) \cdots \Tr(U(\Gamma_k)), \\
\Tr(U(\Gamma_i)) &:= (U(\lambda_{c_i(1)}))^{\varep_i(1)} \cdots (U(\lambda_{c_i(m_i)}))^{\varep_i(m_i)}, ~~ i \in [k].
\end{equs}
Due to invariance of Haar measure, $\E[\Tr(U(\bm \Gamma))] = 0$ if $\bm \Gamma$ is not balanced (multiply one of the $U(\lambda_\ell)$ by $e^{\icomplex \theta} I$ for some appropriately chosen $\theta \in [0, 2\pi]$). Thus, we will assume that $\bm \Gamma$ is balanced. Recall the Weingarten function $\Wg_N$, defined in equation \eqref{eq:weingarten-character-sum}.

\begin{theorem}[Unitary Weingarten calculus]\label{thm:weingarten}
Let $\bm \Gamma = (\Gamma_1, \ldots, \Gamma_M)$ be a balanced collection of words on $\{\lambda_1, \ldots, \lambda_L\}$. Then
\begin{equs}
\E[\Tr(U(\bm \Gamma))] = \sum_{ (\sigma_\ell, \tau_\ell), \ell \in [L]} \bigg(\prod_{\ell \in [L]} \Wg_N(\sigma_\ell^{-1} \tau_\ell) \bigg) N^{\numcomp(\bm \Gamma, \bm \pi)}.
\end{equs}
Here, the sum in the right\revision{-}hand side is over pairs of  bijections $\sigma_\ell, \tau_\ell : [n_+(\lambda_\ell)] \ra (n_+(\lambda_\ell) : 2n_+(\lambda_\ell)]$, $\ell \in [L]$, and we define $\bm \pi = ([\sigma_\ell ~ \tau_\ell], \ell \in [L])$.
\end{theorem}

\begin{remark}
We make the following remarks.

\begin{enumerate}
    \item Theorem \ref{thm:weingarten} says that expectations of products of traces of words of independent Haar Unitaries can be computed as a diagrammatic sum. We find this form of the Weingarten calculus to be particularly clean, and hope it is of benefit to future readers who want to learn about this topic.

    \item We could not find the exact statement of Theorem \ref{thm:weingarten} anywhere in the literature. For instance, \cite[Corollary 2.4]{collins2006integration} gives the matrix-entry version of the Weingarten calculus. Of course, Theorem \ref{thm:weingarten} follows by an application of \cite[Corollary 2.4]{collins2006integration}, since one may expand out $\Tr(U(\bm \Gamma))$ into a giant sum of products of matrix entries. Theorem \ref{thm:weingarten} also directly follows from our Theorem \ref{thm:weingarten-recovery}, which we will prove in Section \ref{section:poisson-exploration}.

    \item \revision{As a hint towards how we will eventually apply Theorem \ref{thm:weingarten} to lattice Yang--Mills, recall from \eqref{eqn::latticeymexpanded} that lattice Yang--Mills can be thought of as a model where one starts with an independent Haar unitary at each lattice edge, and recall Section \ref{subsec::randommatrixintro} that we claimed that the computation of Wilson string expectations will reduce to computing expectations of traces of products of independent Haar Unitaries. This claim will follow because we will take the set of letters $\{\lambda_1, \ldots, \lambda_L\}$ to be the set of edges of the lattice, and the set of words $\mbf{\Gamma}$ to be loops in the lattice.}
\end{enumerate}
\end{remark}

\begin{remark}
Certain representation theory concepts underlie all of this discussion. For instance, the matchings $\pi \in \mc{M}(n)$ should be interpreted as elements of $\mc{B}_n$, the Brauer algebra. We will introduce the Brauer algebra and other relevant representation theory concepts in Section \ref{section:rep-theory}. We chose not to do so yet because we want to give as quick and direct a statement of the Unitary Weingarten calculus as possible. This way, the reader who wants to understand our surface sum representation of Wilson loop expectations can do so as quickly as possible. The representation theory material is only needed for understanding the large-time limits of Lie group Brownian motion.
\end{remark}

We finish off this subsection with an instructive example which illustrates how one may go about computing $\numcomp(\bm \Gamma, \pi)$ for a given $\bm \Gamma, \pi$.

\begin{example}\label{ex:computing-connected-components}
% We give an example of computing the number of components in a collection of strand diagrams. For instance,
As in Figure \ref{figure:strand_diagram_example}, suppose our letters are $\{A, B\}$, and our words are $\bm \Gamma = (AB, AA^{-1}, B)$. Suppose the collection of matchings $\bm \pi = (\pi_A, \pi_B)$ is as in Figure \ref{figure:strand_diagram_example_2}. To compute $\numcomp(\bm \Gamma, \bm \pi)$, we may add vertex labels to Figure \ref{figure:strand_diagram_example_3}. This gives Figure \ref{figure:strand_diagram_example_4}.

\begin{figure}[h]
    \centering
    \includegraphics[page=4, width=10cm]{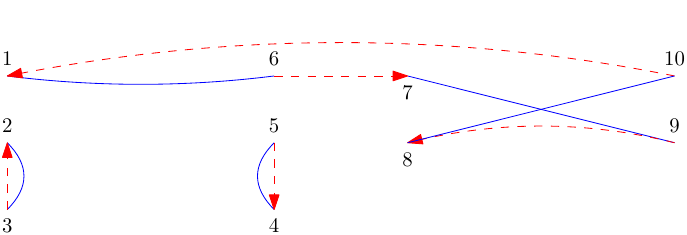}
    \caption{}\label{figure:strand_diagram_example_4}
\end{figure}

It really does not matter how we label the vertices, as long as each vertex gets a unique label. The point is that by using the labels, we can define a permutation whose cycles are in 1-1 correspondence with the the connected components of the graph obtained by combining the interior connections $\bm \pi$ and the exterior connections specified by $\bm \Gamma$. The cycles are obtained by starting at a given vertex and alternately following the (blue) interior connections and (red) exterior connections. In the example of Figure \ref{figure:strand_diagram_example_4}, we obtain the permutation $(1~ 6 ~7~ 9 ~8 ~10) (2 ~3) (4 ~5)$.
\end{example}

\subsection{Poisson point process on strand diagrams}\label{section:poisson-process-intro}

In this section, we review a result in the companion paper~\cite{park2023wilson} that is necessary for this paper. In particular, we express the expected product of traces of words of Unitary Brownian motions in terms of a certain Poisson point process on the strand diagram. We recall that the Unitary Brownian motion started from a unitary matrix $A$ is defined as the Markov process started from $A$ with generator one-half the Laplace-Beltrami operator associated to the Lie \revision{group} $\UN$, \revision{where the Laplace-Beltrami operator is defined using the metric $\langle X, Y \rangle_{\mathfrak{u}(N)} := N\Tr(X^* Y)$ on the Lie algebra $\mathfrak{u}(N)$ of $\UN$.} For more background, see \cite[Section 3.1]{park2023wilson} or \cite[Section 2]{Dahlqvist2017}.

As in Section \ref{section:strand-diagrams-and-weingarten-caculus}, let $\bm \Gamma = (\Gamma_1, \ldots, \Gamma_k)$ be an ordered collection of words on letters $\{\lambda_1, \ldots, \lambda_L\}$, where each $\Gamma_i$ has the form \eqref{eq:general-word-form}. Let $(B_T(\lambda_\ell), \ell \in [L])$ be an i.i.d. collection of Unitary Brownian motions run for time $T$ started at the identity. As in the Haar-distributed case, define
\begin{equs}
\Tr(B_T(\bm \Gamma)) &:= \Tr(B_T(\Gamma_1)) \cdots \Tr(B_T(\Gamma_k)), \\
\Tr(B_T(\Gamma_i)) &:= (B_T(\lambda_{c_i(1)}))^{\varep_i(1)} \cdots (B_T(\lambda_{c_i(m_i)}))^{\varep_i(m_i)}, ~~ i \in [k].
\end{equs}
We now begin to define the Poisson point process on the strand diagram. We denote this process by $\Sigma$. Formally, $\Sigma$ is just a rate-$1$ Poisson point process on the space
\begin{equs}
\mc{D} := \bigsqcup_{\ell \in [L]} \mc{D}(\lambda_\ell), ~~ \text{ where } ~~\mc{D}(\lambda_\ell) := \bigsqcup_{\substack{i, j \in [n_{\pm}(\lambda)] \\ i < j}} [0, \infty), ~~ \text{ for each $\ell \in [L]$}.
\end{equs}
We may write 
\begin{equs}
\Sigma = \bigsqcup_{\ell \in [L]} \Sigma(\lambda_\ell),
\end{equs}
where $\Sigma(\lambda_\ell)$ is a rate-1 Poisson point process on $\mc{D}(\lambda_\ell)$ for each $\ell \in [L]$. For $T \in [0, \infty)$, we also define the finite-interval counterparts
\begin{equs}
\mc{D}_T := \bigsqcup_{\ell \in [L]} \mc{D}_T(\lambda_\ell), ~~ \text{ where } ~~\mc{D}_T(\lambda_\ell) := \bigsqcup_{\substack{i, j \in [n_{\pm}(\lambda)] \\ i < j}} [0, T], ~~ \text{ for each $\ell \in [L]$}.
% \mc{D}_T := \bigsqcup_{\ell \in [L]} \mc{D}_T(\lambda_\ell) := \bigsqcup_{\ell \in [L]} \bigsqcup_{\substack{i, j \in [n_{\pm}(\lambda)] \\ i < j}} [0, T],
\end{equs}
and 
\begin{equs}
\Sigma_T := \bigsqcup_{\ell \in [L]} \Sigma_T(\lambda_\ell), ~~ \text{ where } ~~ \Sigma_T(\lambda_\ell) := \Sigma(\lambda_\ell) \cap \mc{D}_T(\lambda_\ell) \text{ for each $\ell \in [L]$}.
\end{equs}
In words, one should think of $\Sigma(\lambda_\ell)$ as specifying a Poisson collection of points for each pair of strands in the portion of the strand diagram of $\bm \Gamma$ corresponding to the letter $\lambda_\ell$. Next, we explain how to visually interpret these points. In what follows, fix $T \in [0, \infty)$ and a letter $\lambda$.

Let $P \sse \mc{D}_T(\lambda)$ be a finite collection of points. Let $i, j \in [n_{\pm}(\lambda)]$, $i < j$. We say that a point $x \in P$ {\it connects} strands $i$ and $j$ if $x$ comes from the interval $[0, T]$ which is indexed by the pair $\{i, j\}$. In this case, we also say that the strands $i$ and $j$ are {\it connected} by $x$. 

Now $P$ may be visualized as follows. The reader may wish to follow Figure \ref{figure:point-process-example} for a visual example as they read. For each point $x$ in $P$, we draw a vertical line between the two strands which are connected by $x$. If the two strands connected by $x$ have the same direction, then we use the color green to draw the line, and we call $x$ a {\it swap}. If the two strands connected by $x$ have opposite directions, then we use the color blue to draw the line, and we call $x$ a {\it turnaround}. The choices for these names will become clear later. The locations at which we draw these vertical lines are proportional to $x$ (which by assumption is in $[0, T]$). I.e., if we imagine the strands as having unit length, then we draw the line corresponding to $x$ at location $\frac{x}{T}$. Alternatively, if we imagine the strands as having length $T$, then we draw the line at location $x$. It will not matter what length we choose for the strands,
% what matters is that only vertical lines are drawn at locations in proportion to the values of the corresponding points. 
% This is
because of how we are going to interpret these vertical lines, which we begin to describe next.

\begin{figure}[ht!]
    \centering
\includegraphics[width=7cm]{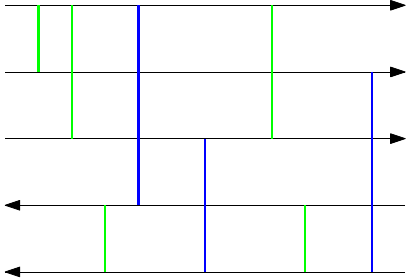}
    \caption{A realization of a finite subset of $\mc{D}_T(\lambda)$ in terms of lines connecting strands.}\label{figure:point-process-example}
\end{figure}

The visual representation of $P \sse \mc{D}_T(\lambda_\ell)$ in terms of the strands defines a matching $\pi \in \mc{M}(n_{\pm}(\lambda))$ in the following manner. Suppose we take some vertex of some strand and want to know which other vertex this is matched to by $\pi$. To do so, we begin ``exploring" the strand emanating from this vertex, and every time we come upon a point (i.e. a swap or turnaround), we ``jump" to the other strand. If the point was a swap, we continue exploring in the same direction as before, while if the point was a turnaround, we actually continue exploring in the opposite direction as before (hence the name ``turnaround"). We continue this process until we reach some other vertex of some strand. We define $\pi$ to pair the two vertices which form the beginning and end points of this exploration. See Figure \ref{figure:point-process-example-2} for a visual example of how this exploration works.

\begin{figure}[ht!]
    \centering
\includegraphics[page = 2, width=7cm]{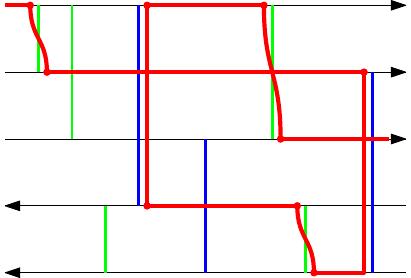}
    \caption{Continuation of the example in Figure \ref{figure:point-process-example}. To determine what the top-left vertex gets matched to, we explore in the indicated manner.}\label{figure:point-process-example-2}
\end{figure}

\begin{figure}[ht!]
    \centering
\includegraphics[page = 3, width=7cm]{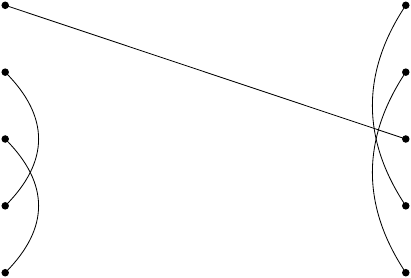}
    \caption{The matching corresponding to Figure \ref{figure:point-process-example}.}\label{figure:point-process-example-3}
\end{figure}

If we separately perform this exploration starting from each vertex of each strand, we obtain the matching $\pi$. As an example, see Figure \ref{figure:point-process-example-3} for the matching one obtains from Figure \ref{figure:point-process-example}.

\begin{notation}
For finite subsets $P \sse\mc{D}_T(\lambda_\ell)$, let $\pi(P)$ be the matching obtained from $P$ in the manner described.
\end{notation}

\begin{figure}[ht!]
    \centering
\includegraphics[page = 1, width=7cm]{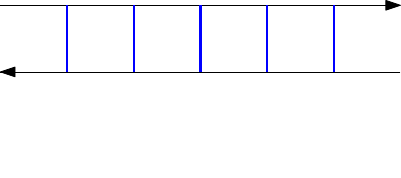}
\includegraphics[page = 2, width=7cm]{figures/point_process_interior_components_example.pdf}
    \vspace{-15mm}
    \caption{Left: a collection of points which results in interior loops. Right: the four interior loops in this example are indicated red. Technically, the vertical parts of the red loops should be drawn on top of the blue lines, but we avoid doing this for visual reasons.}\label{figure:point_process_interior_components_example}
\end{figure}

There is an additional quantity defined in terms of $P$ which will play a role in relating the Poisson point process to expectations of Unitary Brownian motion. The visual representation of $P$ in general determines not only a matching $\pi(P)$ but also a certain number of ``interior loops". These interior loops may arise when there are multiple turnarounds between two strands. See Figure \ref{figure:point_process_interior_components_example} for a simple example to follow along. Each interior loop is obtained by starting at some location on a given strand and then exploring, following swaps and turnarounds until one winds up back at the starting point. Note here we do not need to start at a vertex of a strand (in fact if we start a vertex, then it is impossible to form a loop), rather we allow ourselves to start at (say) the middle of a strand. These interior loops are viewed modulo starting point, so there are only a finite number for any given finite $P$.

\begin{figure}[ht!]
    \centering
\includegraphics[page = 3, width=7cm]{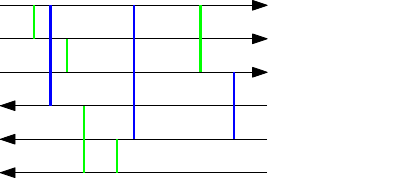}
\includegraphics[page = 4, width=7cm]{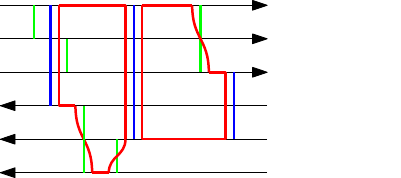}
    \caption{Left: the visual representation of some subset $P$. Right: the two interior loops in this example are indicated red. }\label{figure:point_process_interior_components_example_3}
\end{figure}

See Figure \ref{figure:point_process_interior_components_example_3} for a less trivial example of a subset $P$ with a nonzero number of interior loops.

\begin{definition}
Given a finite subset $P \sse \mc{D}_T(\lambda_\ell)$, let $\numintloops(P)$ be the number of interior loops of $P$. Similarly, for finite subsets $P \sse \mc{D}_T$, which we may express as $P = \bigsqcup_{\ell \in [L]} P(\lambda_\ell)$, where $P(\lambda_\ell)$ is a finite subset of $\mc{D}_T(\lambda_\ell)$, let 
\begin{equs}
\numintloops(P) := \sum_{\ell \in [L]} \numintloops(P(\lambda_\ell)).
\end{equs}
\end{definition}

% We are now almost ready to state how expectations of Unitary Brownian motion relate to the Poisson point process $\Sigma$. First, some more notation.

We define one more piece of notation.

\begin{definition}
For finite subsets $P \sse\mc{D}_T(\lambda_\ell)$, let $\numswaps(P)$ be the number of swaps that $P$ contains. Similarly, for finite subsets $P \sse \mc{D}_T$, which we may express as $P = \bigsqcup_{\ell \in [L]} P(\lambda_\ell)$, where $P(\lambda_\ell)$ is a finite subset of $\mc{D}_T(\lambda_\ell)$, let 
\begin{equs}
\numswaps(P) := \sum_{\ell \in [L]} \numswaps(P(\lambda_\ell)).
\end{equs}
\end{definition}

We now give a precise statement of how Unitary Brownian motion expectations reduce to certain diagrammatic sums. We quote the following result from \cite{park2023wilson}.

\begin{lemma}[Expected trace as Poisson sums~\cite{park2023wilson}]
\label{lem-poisson-sum}
    Let $\bm \Gamma$ be a collection of words on $\{\lambda_1, \ldots, \lambda_L\}$ and $T>0$. Let $\Sigma_T$ be the Poisson point process corresponding to $\bm \Gamma$. Let $\bm \pi (\Sigma_T) := (\pi(\Sigma_T(\lambda_\ell)), \ell \in [L])$, and let $|\Sigma_T|$ be the number of points of $\Sigma_T$.
    % Let $\Sigma$ be the Poisson point process on $\mathcal D_T$. Consider the strand diagram $S$ for $(\bm \Gamma, \Sigma)$.
    Then
    \begin{equs}
    \E \big[\Tr(B_T(\bm \Gamma))\big] = \exp\bigg(T \sum_{\ell \in [L]} &\bigg(\binom{n_{\pm}(\lambda_\ell)}{2} - \frac{n_{\pm}(\lambda_\ell)}{2}\bigg) \bigg) ~\times \\
    &\E\Bigl[(-1)^{\numswaps(\Sigma_T)} N^{-|\Sigma_T|} N^{\numintloops(\Sigma_T)} N^{\numcomp(\bm \Gamma, \bm \pi(\Sigma_T))}\Bigr].
    \end{equs}
\end{lemma} 

\begin{remark}
This lemma gives a probabilistic viewpoint (in terms of expectations over Poisson point processes) for expressing expectations of Unitary Brownian motion. When phrased in these terms, certain ideas become natural. Section \ref{section:poisson-exploration} gives such an example, where we analyze the large-$T$ limit of Unitary Brownian motion expectations by defining an exploration of the Poisson point process with unexpectedly nice properties. 

Of course, there are existing formulas for expectations of Unitary Brownian motion in terms of representation theory (see e.g. \cite{levy2008schur, Dahlqvist2017}), and in particular in terms of the Brauer algebra (which we will define in Section \ref{section:rep-theory}). We will see at the beginning of Section \ref{section:poisson-exploration} how these two viewpoints are related. The Brauer algebra provides a convenient language in which to phrase our proofs. On the other hand, the idea of ``exploring" a point process is not as natural from the representation-theoretic point of view, which is why we make the effort to introduce the Poisson point process formulation.
\end{remark}

We slightly restate Lemma \ref{lem-poisson-sum} to make it appear more similar to the Weingarten calculus (Theorem \ref{thm:weingarten}). For a collection of matchings $\bm \pi = (\pi_\ell, \ell \in [L])$, define
\begin{equs}
w_T&(\bm \pi) := \prod_{\ell \in [L]} w_T(\pi_\ell), 
\end{equs}
where
\begin{equs}\label{eq:w-T-def}
w_T(\pi_\ell) := \exp\bigg(&\binom{n_{\pm}(\lambda_\ell)}{2}T  - \frac{n_{\pm}(\lambda_\ell)}{2}T\bigg) ~\times \\
&\E\Big[ (-1)^{\numswaps(\Sigma_T(\lambda_\ell))} N^{-|\Sigma_T(\lambda_\ell)|} N^{\numintloops(\Sigma_T(\lambda_\ell))}  \ind(\bm \pi(\Sigma_T(\lambda_\ell))= \pi_\ell) \Big],
\end{equs}
which can be interpreted as the partition function of all Poisson point configurations which results in the collection of matchings $\bm \pi$. Lemma \ref{lem-poisson-sum} can then be restated in the following form.

\begin{lemma}\label{lemma:expectations-words-unitary-bm}
 Let $\bm \Gamma$ be a collection of words on $\{\lambda_1, \ldots, \lambda_L\}$. Let $\Sigma_T$ be the Poisson point process corresponding to $\bm \Gamma$. Then
\begin{equs}\label{eq:expectations-words-unitary-bm}
\E \big[ \Tr(B_T(\bm \Gamma))\big] = \sum_{\bm \pi = (\pi_1, \ldots, \pi_L)} w_T(\bm \pi) N^{\numcomp(\bm \Gamma, \bm \pi)} .
\end{equs}
\end{lemma}

\begin{remark}
Observe the similarity between this lemma and the Weingarten calculus (Theorem \ref{thm:weingarten}). In both cases, expectations of traces of products of words are expressed as diagrammatic sums. However, the weights of the sums differ. In the Haar case, the matchings which have nonzero weight are made out of combining a pair of left and right bijections, and the weight of such a matching is given by the Weingarten function of their ``difference". In the Brownian motion case, in principle any matching can have nonzero weight, and the weight can be interpreted as the partition function of some Poisson point process.
\end{remark}

% Lemma \ref{lemma:expectations-words-unitary-bm} says that in order to compute the expectations of traces of words of Unitary Brownian motion, we may perform a weighted sum over all partitions of the corresponding strand diagram, where the weights are given by $w_T(\pi)$, and the statistic we are averaging over is $N$ raised to the number of components of the diagram made from the exterior connections specified by the collection of words $\Gamma$ and the interior connections specified by the collection of partitions $\pi$. 

\begin{figure}[ht!]
    \centering
\includegraphics[width=.8\textwidth]{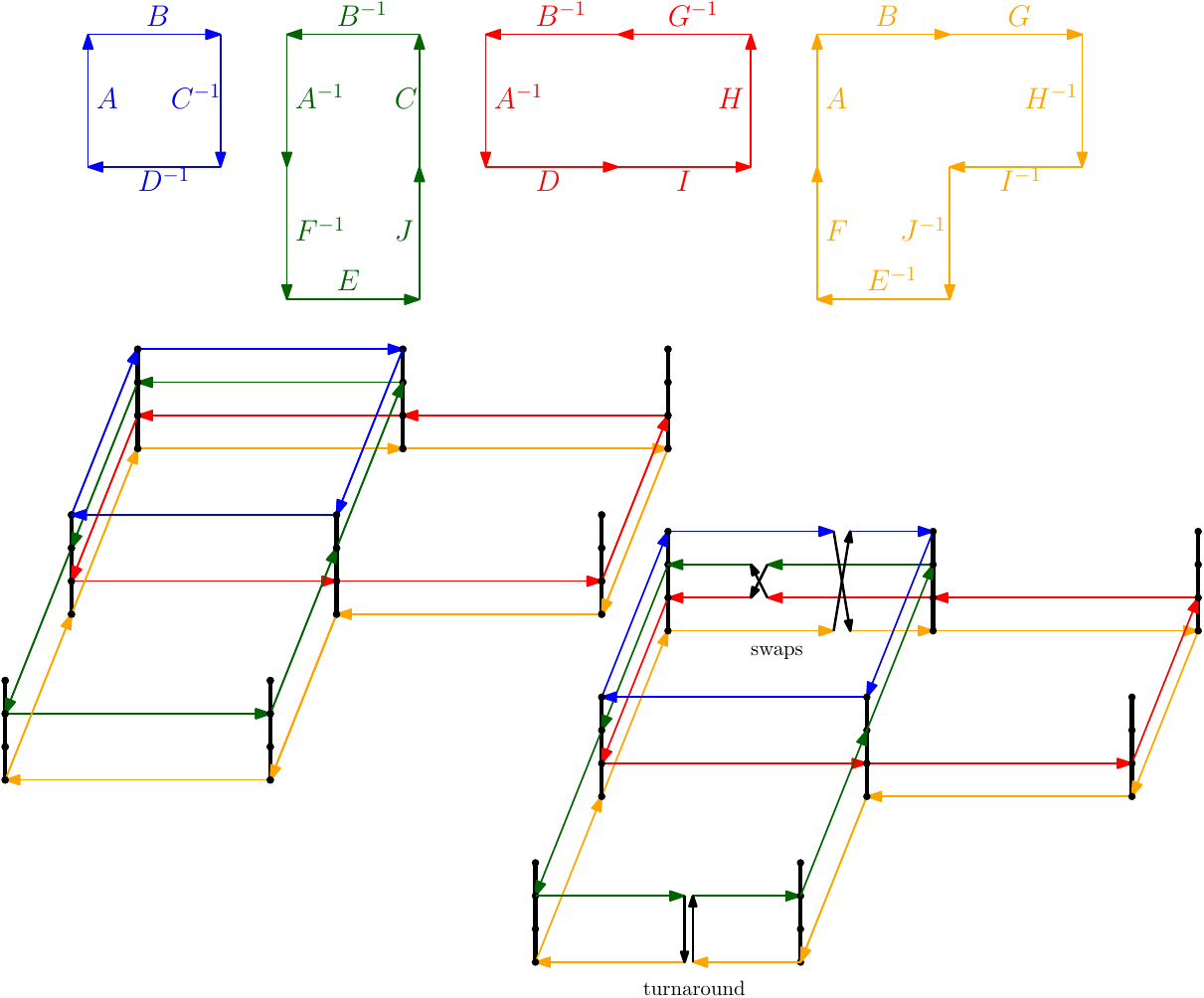}
    \captionsetup{singlelinecheck=off}
    \caption[.]{The quantity 
    \[ \E \Bigl[ \tr (ABC^{-1}D^{-1}) \tr(A^{-1} F^{-1} E J C B^{-1} ) \tr(A^{-1} D I H G^{-1} B^{-1}) \tr( B G H^{-1} I^{-1} E^{-1} F A) \Bigr] \]
    is the expected trace product of the four loops above, where the symbols are i.i.d.\ random elements of $\UN$. \cite{park2023wilson} explains one way to compute this quantity when each symbol has the law of the time-$T$ value of a Brownian motion on $\UN$ started at the identity. First one stacks the strands on top of each other (lower left) so that the arrows corresponding to the same symbol lie on top of each other. Then one chooses locations for points according to a certain Poisson point process (lower right). 
    % The desired expectation is a constant times $\mathbb E\Bigl[ (-1/N)^{\textrm{\#\,swaps}} N^{\textrm{\#\,components}}\Bigr] $.
    As $T \to \infty$ the expected number of points tends to infinity. In this paper we derive the $T \to \infty$ limit for the expectation by starting with the formulation above and applying simple geometric arguments and sign cancellations. 
    % Similar arguments apply to other compact Lie groups.
    }
    \label{fig::strands}
\end{figure}

We proceed to give a precise statement of Theorem \ref{thm:informal-weignarten-recovery}, which is we are able to obtain the $T \toinf$ limit of the right-hand side of~\eqref{eq:expectations-words-unitary-bm} and thus recover the Weingarten calculus (Theorem \ref{thm:weingarten}). The proof of the following theorem is the subject of Section \ref{section:poisson-exploration}.

% First, we make the following definition.

% \begin{definition}[Balanced collection of words]
% A collection of words $\bm \Gamma = (\Gamma_1, \ldots, \Gamma_k)$ on letters $\{\lambda_1, \ldots, \lambda_L\}$ is balanced if for each letter $\lambda_i$, the number of times that $\lambda_i$ appears in $\bm \Gamma$ is equal to the number of times $\lambda_i^{-1}$ appears in $\bm \Gamma$.
% \end{definition}

% Next, we describe a certain special set of partitions that plays a key role in the limiting formula. Suppose that $\bm\Gamma$ is balanced. Then $n_\ell$ is even for all $\ell \in [L]$. Given a pair of bijections $\sigma, \tau : [n_\ell/2] \ra (n_\ell /2:n_\ell]$, we may obtain a partition $[\sigma ~ \tau]$ of $[2n_\ell]$ as in Figure \ref{figure:walled-brauer-def-example-3}.
% \begin{figure}[h]
%     \centering
%     \includegraphics[page=3]{figures/walled-brauer-def-example.pdf}
%     \caption{Example when $n_\ell = 6$. Here $\sigma_\ell$ maps $1 \mapsto 6$, $2 \mapsto 4$, $3 \mapsto 5$, and $\tau_\ell$ maps $1 \mapsto 4$, $2 \mapsto 5$, $3 \mapsto 6$. }\label{figure:walled-brauer-def-example-3}
% \end{figure}

% Clearly, the set of all matchings that arise this way is a strict subset of the set of all matchings of $[2n_\ell]$. Observe that $\sigma \tau^{-1} : [n_\ell/2] \ra [n_\ell/2]$ is a bijection, and thus we may view $\sigma \tau^{-1} \in \symgrp_{n_\ell/2}$. It turns out that as $T \toinf$, these are the only matchings that have a non-vanishing weight. This is inherent in the following theorem. Its proof is the subject of Section \ref{section:poisson-exploration}.

\begin{theorem}\label{thm:weingarten-recovery}
Let $\bm \Gamma = (\Gamma_1, \ldots, \Gamma_M)$ be a balanced collection of words on $\{\lambda_1, \ldots, \lambda_L\}$. Then
\begin{equs}
\lim_{T \toinf} \E \big[ \Tr(B_T(\bm \Gamma))\big] =  \sum_{ (\sigma_\ell, \tau_\ell), \ell \in [L]} \bigg(\prod_{\ell \in [L]} \Wg_N(\sigma_\ell^{-1} \tau_\ell) \bigg) N^{\numcomp(\bm \Gamma, \bm \pi)}.
\end{equs}
Here, as in Theorem \ref{thm:weingarten}, the sum in the right\revision{-}hand side is over pairs of  bijections $\sigma_\ell, \tau_\ell : [n_+(\lambda_\ell)] \ra (n_+(\lambda_\ell) : 2n_+(\lambda_\ell)]$, $\ell \in [L]$, and we define $\bm \pi = ([\sigma_\ell ~ \tau_\ell], \ell \in [L])$.
\end{theorem}

See Figure 
\ref{fig::strands} for a visualization of Theorem \ref{thm:weingarten-recovery}.

\subsection{Representation theory  and other preliminaries}\label{section:rep-theory}

The concepts of Section \ref{section:poisson-process-intro} may naturally be phrased in terms of the Brauer algebra, which is a well-studied object in mathematics. We proceed to introduce the Brauer algebra because this will form a convenient language when phrasing our proofs. \revision{We emphasize here that all the concepts in this section are quite classical in representation theory.}

Recall from Definition \ref{def:matchings} that $\mc{M}(n)$ is the set of matchings $\pi$ of $[2n]$. 
% We view matchings pictorially as in Figure \ref{figure:brauer-def-example}.
%     \begin{figure}[ht!]
%         \centering
%          \includegraphics[width=10cm]{figures/brauer-def-example.pdf}
%         \caption{$\pi = \{\{1, 3\}, \{2, 9\}, \{4, 10\}, \{5, 7\}, \{6, 8\}\}$}\label{figure:brauer-def-example}
%     \end{figure}

% We refer to pairs that involve both a left and right vertex as ``left-right pairings", and pairs that involve two left elements or two right elements as ``same-side pairings". In the above picture, $\{1, 3\}, \{6, 8\}$ are same-side pairings, while $\{2, 9\}, \{4, 10\}, \{5, 7\}$ are left-right pairings.

\begin{definition}[Brauer algebra]\label{def:brauer-algebra}
Let $\mc{B}_n$ be the vector space of $\C$-valued functions on $\mc{M}(n)$. We will often view elements $f \in \mc{B}_n$ as formal sums $f = \sum_\pi f(\pi) \pi$, where $\pi$ ranges over $\mc{M}(n)$. Fix $\zeta \in \C$. We define a product $\pi_1 \pi_2 \in \mc{B}_n$ of matchings $\pi_1, \pi_2 \in \mc{M}(n)$ as in Figure \ref{figure:brauer-algebra-multiplication}.
    \begin{figure}[ht!]
        \centering
        \includegraphics[page=3]{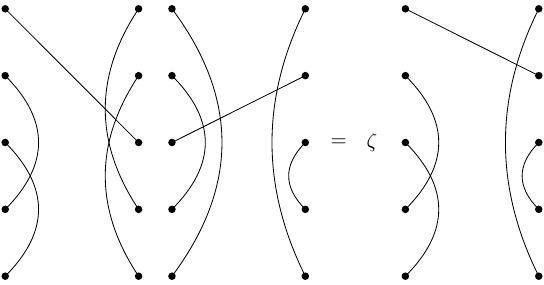}
        \caption{Example of multiplication in the Brauer algebra. In words, we put $\pi_1, \pi_2$ together side-by-side, and then follow the lines to obtain a new matching. Any closed loops incur a factor of $\zeta$. }\label{figure:brauer-algebra-multiplication}
    \end{figure}

Observe that this product induces a product on $\mc{B}_n$ which turns $\mc{B}_n$ into an algebra. Explicitly, if we represent $f, g \in \mc{B}_n$ by formal linear combinations $f = \sum_{\pi_1 \in \mc{M}(n)} f(\pi_1) \pi_1$, $g = \sum_{\pi_2 \in \mc{M}(n)} g(\pi_2) \pi_2$, then the product $fg$ is given by
    \begin{equs}
        f g = \sum_{\pi_1, \pi_2 \in \mc{M}(n)} f(\pi_1) g(\pi_2) \pi_1 \pi_2.
    \end{equs}
    We refer to $\mc{B}_n$ as the Brauer algebra.
\end{definition}

\begin{remark}\leavevmode \revision{We make the following remarks regarding the Brauer algebra.}
 
\begin{enumerate}
\item Clearly, we may view $\mc{M}(n) \sse \mc{B}_n$ by viewing a matching $\pi \in \mc{M}(n)$ as the delta function $\delta_\pi$ defined by $\delta_\pi(\pi') = \ind(\pi = \pi')$.
\item Typically, elements of $\mc{B}_n$ are drawn as top-bottom matchings, yet here we have chosen to draw them as left-right matchings.
\end{enumerate}
\end{remark}

In what follows, we always take $\zeta = N$. This is the choice of $\zeta$ which relates multiplication in the Brauer algebra with expectations of Unitary Brownian motion: note that the factor of $N$ that we incur when we form a loop exactly matches the factor of $N$ that we incur in the strand diagram when we add another interior loop. We will later give a precise statement of this in Lemma \ref{lemma:brauer-algebra-F-relation}.

We specify a norm on $\mc{B}_n$, which will enable us to later talk about convergence in $\mc{B}_n$.

\begin{definition}[Norm on $\mc{B}_n$]\label{def:brauer-norm}
For $f \in \mc{B}_n$, define $\|f\|$ to be the $L^1$ norm, i.e. $\|f\| := \sum_{\pi \in \mc{M}(n)} |f(\pi)|$.
\end{definition}

Next, we define a certain sub-algebra of the Brauer algebra, called the walled Brauer algebra. This arises naturally in computing expectations of Unitary Brownian motion, as it turns out that the strand diagrams of Section \ref{section:poisson-process-intro} are not only elements of the Brauer algebra, but even more they are elements of the walled Brauer algebra.

\begin{definition}[Walled Brauer algebra]\label{def:walled-brauer-algebra}
Let $n, m \geq 1$. Let $\mc{M}(n, m) \sse \mc{M}(n+m)$ be the subset of matchings of $[2(n+m)]$ such that every pairing of same-side vertices includes one vertex among the top $n$ vertices, and one vertex among the bottom $m$ vertices, and every pairing of opposite-side vertices is either between two top $n$ vertices, or between two bottom $m$ vertices. 

Pictorially, one imagines a dashed line separating the top $n$ vertices from the bottom $m$ vertices, and the only pairings which can cross this dashed line are pairings of same-side vertices. See Figure \ref{figure:walled-brauer-def-example} for an example when $n = m = 3$.
\begin{figure}[ht!]
        \centering
        \includegraphics{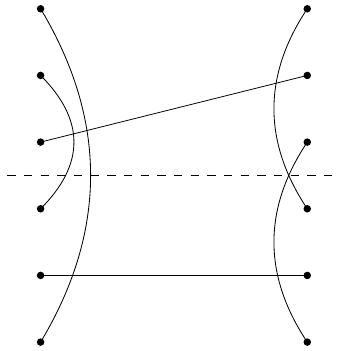}
        \caption{Element of $\mc{B}_{3, 3}$}\label{figure:walled-brauer-def-example}
    \end{figure}

The walled Brauer algebra $\mc{B}_{n, m}$ is the sub-algebra of $\mc{B}_{n+m}$ consisting of functions $f \in \mc{B}_{n+m}$ which are supported on the matchings $\mc{M}(n, m)$. One may check that given two matchings $\pi_1, \pi_2 \in \mc{M}(n, m)$, their product $\pi_1 \pi_2$ is proportional to a matching in $\mc{M}(n, m)$. This implies that the product on $\mc{B}_{n+m}$ descends to a product on $\mc{B}_{n, m}$.
\end{definition}

\begin{remark}
Recall from Definition \ref{def:left-right-bijection} that given two bijections $\sigma, \tau : [n] \ra (n:2n]$, we may obtain a matching $[\sigma ~ \tau] \in \mc{M}(2n)$. In fact, this matching is in $[\sigma ~ \tau] \in \mc{M}(n, n)$.
\end{remark}

\begin{figure}[ht!]
 \centering
    % \minipage{0.55\textwidth}
    \includegraphics[page=2,width=8cm]{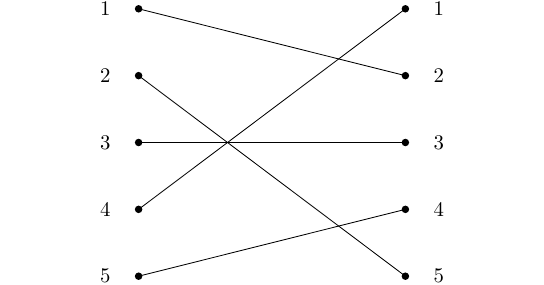}
        \caption{$\sigma = (1 ~ 2 ~ 5 ~ 4)$}\label{figure:brauer-def-example-2}
    % \endminipage
    % \hfill
    % \minipage{0.38\textwidth}
    %     \includegraphics[width=\linewidth]{figures/brauer_13.png}
    %         \caption{$(1 ~ 3)$}\label{figure:brauer-(1 3)}
    % \endminipage
    \end{figure}

The embedding $\mc{M}(n) \sse \mc{B}_n$ also restricts to an embedding $\mc{M}(n,m) \sse \mc{B}_{n,m}$. Moreover, observe that $\symgrp_n$ can be embedded in $\mc{M}(n) \sse \mc{B}_n$ as follows. Given $\sigma \in \symgrp_n$, we can view it as an element of $\mc{M}(n)$ as in Figure \ref{figure:brauer-def-example-2}.

We may also embed $\symgrp_n$ into $\mc{M}(n, n) \sse \mc{B}_{n, n}$ by connecting the top $n$ vertices on the left and right as we did to embed $\symgrp_n$ into $\mc{B}_n$, and then connecting the bottom $n$ vertices on the left and right by straight lines.

Next, we define the following notation for certain special elements of the walled Brauer algebra $\mc{B}_{n, m}$. These correspond to the swaps and turnarounds introduced in Section \ref{section:poisson-process-intro}.

\begin{definition}\label{def:brauer-algebra-generators}
Given $i, j \in [n]$ or $i, j \in (n:n+m]$, define $(i ~ j) \in \mc{M}(n, m) \sse \mc{B}_{n, m}$ to be the matching of $[2(n+m)]$ which swaps the $i, j$ vertices with their corresponding versions on the right, while keeping the other vertices fixed. This is best explained by the example in the left of Figure~\ref{figure:brauer-generator-examples}.
% \begin{figure}[h]
%   \centering
%         \includegraphics[width=5cm]{figures/brauer-(1 3).png}
%         \caption{$(1 ~ 3)$}
% \end{figure}

Given $i \in [n]$ and $j \in (n:n+m]$, let $\langle i ~ j \rangle$ be the matching which has a same-side pairing between $i, j$ on the left, as well as their corresponding versions on the right, and every other pairing is a straight line. See the right of Figure \ref{figure:brauer-generator-examples} for an example.
\begin{figure}[h]
  \centering
    \includegraphics[page=1, width=5cm]{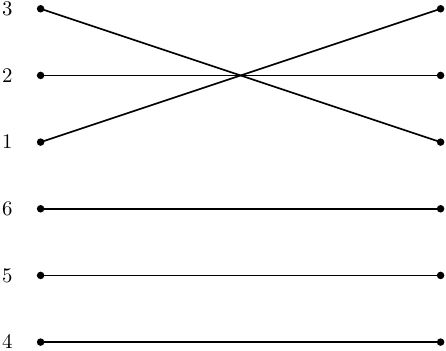}
    \hspace{10mm}
    \includegraphics[page = 2, width=5cm]{figures/brauer-element.pdf}
        \caption{Examples of the two types of generators when $n = m = 3$. Left: $(1 ~ 3)$. Right: $\langle 2 ~ 5 \rangle$}\label{figure:brauer-generator-examples}

%       \includegraphics[page=2,width=5cm]{figures/walled-brauer-def-example.pdf}
%     \caption{$\sigma = \{1 \mapsto 5, 2 \mapsto 6, 3 \mapsto 4\}, ~ \tau = \{1 \mapsto 4, 2 \mapsto 5, 3 \mapsto 6\}$}\label{figure:walled-brauer-def-example-2}
%    \endminipage\hfill
\end{figure}
\end{definition}

% Next, we define the following notation for another set of special elements of the walled Brauer algebra $\mc{B}_{n, n}$. These elements are matchings which have no left-right pairings.

% \begin{definition}
% Let $\sigma, \tau : [n] \ra (n:2n]$ be bijections. Define $[\sigma ~ \tau] \in \mc{B}_{n, n}$ to be the element of the walled Brauer algebra which is given by $\sigma$ on the left and $\tau$ on the right. See the Figure \ref{figure:walled-brauer-def-example-2} for an example when $n = 3$.
% % \begin{figure}
% %     \centering
% %     \includegraphics[page=2]{figures/walled-brauer-def-example.pdf}
% %     \caption{$\sigma = \{1 \mapsto 5, 2 \mapsto 6, 3 \mapsto 4\}, ~ \tau = \{1 \mapsto 4, 2 \mapsto 5, 3 \mapsto 6\}$}\label{figure:walled-brauer-def-example-2}
% % \end{figure}
% \end{definition}

% \begin{notation}
% We note here the notational difference between $(\sigma, \tau)$ and $[\sigma ~ \tau]$. Whenever we want to denote a pair of $\sigma$ and $\tau$, we use the notation $(\sigma, \tau)$. The notation $[\sigma ~ \tau]$ will only be used when we want to denote the element of $\mc{B}_{n, n}$ induced by the pair $(\sigma, \tau)$.
% \end{notation}

Note the particular way we have chosen to label the vertices in Figure \ref{figure:brauer-generator-examples}. From now on, this is how we will label vertices when working with the walled Brauer algebra $\mc{B}_{n, n}$. Ultimately the labeling will not matter, but we have chosen to label in this way to better relate to the Jucys-Murphy elements, which we next define.

\begin{definition}[Jucys-Murphy elements]
For $n \geq 2$, define the Jucys-Murphy element $J_n := (1 ~ n) + \cdots + (n-1 ~ n) \in \C[\symgrp_n]$. We also view $J_n \in \C[S_m]$ for any $n \leq m$. Define $J_1 := 0$. 
\end{definition}

\begin{remark}\label{remark:jucys-murphy-commute}
One may show that the Jucys-Murphy elements commute with each other.
\end{remark}

In the following, we will also view $J_1, \ldots, J_n$ as elements of $\mc{B}_{n, n}$, by using the previously mentioned embedding of $\symgrp_n \sse \mc{B}_{n, n}$. We will also need to refer to Jucys-Murphy elements which act on bottom elements rather than top elements. We define this next.

\begin{definition}\label{def:bottom-jucys-murphy}
Let $n, m \geq 1$. Define $J_1', \ldots, J_m' \in \mc{B}_{n, m}$ by 
\begin{equs}
J_k' := (n + 1 ~ n + k) + \cdots + (n+ k - 1 ~ n+k), ~~ k \in [m].
\end{equs}
\end{definition}

\begin{definition}[Norm on group algebra]
For $f \in \C[\symgrp_n]$, we define $\|f\|$ to be the $L^1$ norm, i.e. $\|f\| := \sum_{\pi \in \symgrp_n} |f(\pi)|$.
\end{definition}

Note that the norms on $\C[\symgrp_n]$ and $\mc{B}_n$ are compatible with our embeddings $\C[\symgrp_n] \sse \mc{B}_n$ and $\C[\symgrp_n] \sse \mc{B}_{n, n} \sse \mc{B}_{2n}$.

\begin{remark}\label{remark:inverse-in-group-algebra}
With this definition of the norm, we have that $\|fg\| \leq \|f\| \cdot \|g\|$ for $f, g \in \C[\symgrp_n]$. This implies that $\|e^f\| \leq e^{\|f\|}$, which further implies that if $\|f\| < 1$, then 
\begin{equs}
\int_0^\infty e^{-u(\id +  f)} du  \in \C[\symgrp_n]
\end{equs}
converges absolutely. Moreover, one has that
\begin{equs}
\int_0^\infty e^{-u(\id +  f)} du  = (\id + f)^{-1}.
\end{equs}
\end{remark}

Next, we discuss an alternate form of the Weingarten function $\Wg_N$ which arises naturally in the proof of Theorem \ref{thm:weingarten-recovery}. First, suppose $N \geq n$. The case of general $N$ will be addressed a bit later. Then $\Wg_N \in \C[\symgrp_n]$ is the following inverse:
\begin{equs}
\Wg_N := \bigg(\sum_{\sigma \in \symgrp_n} N^{\#\mrm{cycles}(\sigma)} \sigma\bigg)^{-1}.
\end{equs}
Jucys \cite{jucys1974symmetric} proved the following identity:
\begin{equs}\label{eq:jucys-identity}
\sum_{\sigma \in \symgrp_n} N^{\#\mrm{cycles}(\sigma)} \sigma = (N + J_n) \cdots (N + J_1).
\end{equs}
Note that when $N \geq n$, each $N + J_k$ for $k \in [n]$ is invertible because then $\|J_k\| = k -1 < N$ (recall Remark \ref{remark:inverse-in-group-algebra}), with inverse given by:
\begin{equs}
(N + J_k)^{-1} = N (\id + J_k / N)^{-1} = N \int_0^\infty e^{-u(\id + J_k/N)} du.
\end{equs}
% explicitly given by the convergent sum
% \begin{equs}
% (N + J_k)^{-1} = \frac{1}{N}\sum_{j=0}^\infty \bigg(\frac{-J_k}{N}\bigg)^j.
% \end{equs}
Since the $J_1, \ldots, J_n$ commute with each other, we have that (as observed by \cite{novak2010jucys})
\begin{equs}\label{eq:weingarten-jucy-murphy-large-N}
\Wg_N = (N+J_n)^{-1} \cdots (N + J_1)^{-1}, ~~ \text{ when $N \geq n$.}
\end{equs}
The reason why we introduce this formula for the Weingarten function is because the terms $(N + J_k)^{-1}$, $k \in [n]$ will appear naturally in our argument.

% Next, we discuss the definition of the Weingarten function in case of general $N$. We follow \cite{collins2006integration}.

% \begin{definition}[Weingarten function]\label{def:weingarten-general-N}
% Let $N, n \geq 1$. Define $\Wg_N \in \C[\symgrp_n]$ (as usual, we omit the dependence on $n$) by
% \begin{equs}\label{eq:weginarten-def-general}
% \Wg_N (\sigma) := \frac{1}{n!} \sum_{\substack{\lambda \vdash n \\ \ell(\lambda) \leq N} }\Bigl[ \chi_\lambda(\id) \chi_\lambda(\sigma) \prod_{(i,j) \in \lambda}(N+j-i)^{-1} \Bigr].
% \end{equs}
% Here, $\ell(\lambda)$ is the number of rows of $\lambda$, i.e.\ the number of parts in the partition of $n$ given by $\lambda$.
% \end{definition}

% Compared with the formula \eqref{eq:weingarten-character-sum} for $N \geq n$, the only difference is in the restriction $\ell(\lambda) \leq N$ when summing over Young diagrams $\lambda$. Note that when $N \geq n$, every Young diagram with $n$ boxes has at most $n \leq N$ rows, and thus the definition \eqref{eq:weginarten-def-general} reduces to \eqref{eq:weingarten-character-sum} if $N \geq n$.

\subsubsection{Additional technicalities for the small \texorpdfstring{$N$}{N} case}\label{section:general-N-preliminaries}

The following material is only needed to prove Theorem \ref{thm:weingarten-recovery} in the case $N < 2\max_{\ell \in [L]} n_\ell$. We encourage the reader on a first reading to skip this subsection and continue on to Section \ref{section:poisson-exploration} to first read over the proof in the case $N \geq 2\max_{\ell \in [L]} n_\ell$, which already contains the main probabilistic ideas. The reader may come back to this section once they are ready to read Section \ref{section:poisson-exploration-general-N-proof}, where the results introduced here will be needed.

Let $e_1, \ldots, e_N$ denote the standard basis of $\C^N$. The tensor space $(\C^N)^{\otimes n}$ has a basis given by $(e_i, i = (i_1, \ldots, i_n) \in [N]^n)$, where $e_i := e_{i_1} \otimes \cdots  \otimes e_{i_n}$. The space $(\C^N)^{\otimes n}$ has a natural inner product which when restricted to basis elements is given by
\begin{equs}
\langle e_i, e_j \rangle = \delta_{ij} = \delta_{i_1 j_1} \cdots \delta_{i_k j_k}.
\end{equs}
\revision{To be clear, our convention is that the inner product is linear in the first variable and conjugate-linear in the second variable, i.e. $\langle \lambda v, w \rangle = \lambda \langle v, w \rangle$ and $\langle v, \lambda w \rangle = \bar{\lambda} \langle v, w \rangle$ for $\lambda \in \C$.}
Let $M \in \End((\C^N)^{\otimes n})$. One may think of $M$ as an $N^n \times N^n$ matrix, whose matrix entries are given by:
\begin{equs}
M_{ij} = \langle e_i, M e_j \rangle, ~~ i, j \in [N]^n.
\end{equs}
In particular, if $M_1, \ldots, M_n \in \End(\C^N)$, then the matrix entries of the tensor product $M = M_1 \otimes \cdots \otimes M_n$ are given by
\begin{equs}
M_{ij} = \langle e_i, (M_1 \otimes \cdots \otimes M_n) e_j \rangle &= \langle e_{i_1} \otimes \cdots e_{i_n}, (M_1 e_{j_1}) \otimes \cdots \otimes (M_n e_{j_n}) \rangle \\ 
&= \langle e_{i_1}, M_1 e_{j_1} \rangle \cdots \langle e_{i_n}, M_n e_{j_n} \rangle \\
&= (M_1)_{i_1 j_1} \cdots (M_n)_{i_n j_n},
\end{equs}
i.e. the product of the corresponding matrix entries of $M_1, \ldots, M_n$.

\begin{definition}\label{def:rho-plus}
Let $N, n \geq 1$. We define a representation $\rho_+$ of $\mc{B}_n$ as follows. Given a pairing $\pi$ of $[2n]$, define $\rho_+(\pi)$ to be the linear map in $\End((\C^N)^{\otimes n})$ whose matrix entries are given by:
\begin{equs}
(\rho_+(\pi))_{(i_1, \ldots, i_n), (i_{2n}, \ldots, i_{n+1})} := \prod_{\{a, b\} \in \pi} \delta^{i_a i_b} .
\end{equs}
\end{definition}

\begin{figure}[h]
    \centering
    \includegraphics{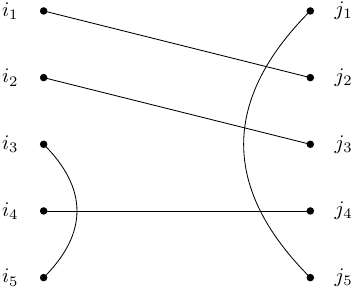}
    \caption{Visualization of matrix entries of $\rho_+(\pi)$.}\label{figure:brauer-representation-example}
\end{figure}

For notational brevity, we omit the dependence of $\rho_+$ on $N, n$. In the following, we mostly apply $\rho_+$ to elements of $\mc{B}_{n, n} \sse \mc{B}_{2n}$. The way one visualizes this definition is as follows. Suppose $n = 5$ and we are given the pairing displayed in Figure \ref{figure:brauer-representation-example}.

Then the matrix entry corresponding to indices $(i_1, \ldots, i_5), (j_1, \ldots, j_5)$ is simply $1$ if all constraints indicated by the pairing are satisfied (in this case, $i_1 = j_2$, $i_2 = j_3$, $i_3 = i_5$, $i_4 = j_4$, and $j_1 = j_5$), and $0$ otherwise.

\begin{remark}
There is an alternative definition of $\rho_+$ that one typically sees (e.g. \cite{Dahlqvist2016}). First, let $E_{ij} \in \End(\C^N)$ be the elementary matrix which has a $1$ in its $(i, j)$ entry and zeros everywhere else. We may write $E_{ij} = e_i e_j^T$. Then
\begin{equs}
\rho_+(\pi) = \prod_{\{a, b\} \in \pi} \delta^{i_a i_b} E_{i_1 i_{2n}} \otimes E_{i_2 i_{2n-1}} \otimes \cdots \otimes E_{i_n i_{n+1}}.
\end{equs}
Here and in the following, repeated indices are implicitly summed over. To see why this definition is equivalent, we may compute an arbitrary matrix entry:
\begin{equs}
\big(\rho_+(\pi)\big)_{(i_1, \ldots, i_n), (i_{2n}, \ldots, i_{n+1})} &= \prod_{\{a, b\} \in \pi} \delta^{j_a j_b}\langle e_{i_1} \otimes \cdots e_{i_n}, \big(E_{j_1 j_{2n}} \otimes \cdots \otimes E_{j_n j_{n+1}}\big) \big( e_{i_{2n}} \otimes \cdots e_{i_{n+1}} \big)\rangle \\
&= \prod_{\{a, b\} \in \pi} \delta^{j_a j_b} \langle e_{i_1}, e_{j_1} e_{j_{2n}}^T e_{i_{2n}} \rangle \cdots \langle e_{i_n}, e_{j_n} e_{j_{n+1}}^T e_{i_{n+1}} \rangle \\
&=  \prod_{\{a, b\} \in \pi} \delta^{j_a j_b}  \delta_{i_1 j_1} \delta_{j_{2n} i_{2n}} \cdots \delta_{i_n j_n} \delta_{j_{n+1} i_{n+1}} \\
&= \prod_{\{a, b\} \in \pi} \delta^{i_a i_b} .
\end{equs}
\end{remark}

\begin{example}
\revision{We record the following example calculation of $\rho_+$. Take $n = 2$, and consider the Brauer algebra element $\langle 1 ~ 2 \rangle \in \mc{B}_{1, 1} \sse \mc{B}_2$ (which recall was defined in Definition \ref{def:brauer-algebra-generators}). In this case, $\rho_+(\langle 1 ~ 2 \rangle) \in \End((\C^N)^{\otimes 2})$, and its action on pure tensors is given by
\begin{equs}
\rho_+(\langle 1 ~ 2 \rangle)(v \otimes w) &= \delta^{i_1 i_2} \delta^{i_3 i_4} (E_{i_1 i_4} v) \otimes (E_{i_2 i_3} w) \\
&= \delta^{i_1 i_2} \delta^{i_3 i_4} \langle v, e_{i_4} \rangle \langle w, e_{i_3} \rangle e_{i_2} \otimes e_{i_3} \\
&= \langle v, \bar{w} \rangle \delta^{ij} e_i \otimes e_j.
\end{equs}
In words, $\rho_+(\langle 1 ~ 2 \rangle)$ acts on pure tensors by contraction of the two tensor components, and its image is in the linear span of the vector $\delta^{ij} e_i \otimes e_j$. 
}

\revision{More generally, for $\langle i ~ j \rangle \in \mc{B}_{n, m} \sse \mc{B}_{n+m}$, the matrix $\rho_+(\langle i ~ j \rangle)$ acts on pure tensors by sending
\begin{equs}
v_1 \otimes \cdots v_i \otimes \cdots  \otimes v_n \otimes w_1 \otimes \cdots w_j \otimes w_m
\end{equs}
to the vector
\begin{equs}
\langle v_i, \bar{w}_j \rangle v_1 \otimes \cdots e_i \otimes \cdots \otimes v_n \otimes w_1 \otimes \cdots \otimes e_j \otimes \cdots \otimes w_m,
\end{equs}
where the $i$th component among the $v$s and the $j$th component among the $w$s are contracted. On all other components, the matrix acts as the identity.
}

\revision{
By an explicit computation, we have that
\begin{equs}
(U \otimes \bar{U}) \rho_+(\langle 1 ~ 2 \rangle) (v \otimes w) &= \langle v, \bar{w} \rangle \delta^{ij} U e_i \otimes \bar{U} e_j \\
&= \langle v, \bar{w} \rangle \delta^{ij} U_{\ell i} \bar{U}_{k j} e^\ell \otimes e^k \\
&= \langle v, \bar{w} \rangle \delta_{\ell k} e^\ell \otimes e^k \\
&= \rho_+(\langle 1 ~ 2 \rangle)(v \otimes w).
\end{equs}
Here, in the second-to-last line, we used that $U$ is a Unitary matrix, in particular that its rows are orthogonal, so that upon fixing $\ell, k$ and summing in $i, j$, we have that $\delta^{ij} U_{\ell i} \bar{U}_{kj} = \delta_{\ell k}$. Similarly, we may compute 
\begin{equs}
\rho_+(\langle 1 ~ 2 \rangle)(U \otimes \bar{U}) (v \otimes w) &= \langle v, e_i \rangle \langle w, e_j \rangle  U_{\ell i} \bar{U}_{kj}  \rho_+(\langle 1 ~ 2 \rangle)(e^i \otimes e^j) \\
&= \langle U v , \ovl{\bar{U} w} \rangle \delta^{ij} e_i \otimes e_j \\
&= \langle U v, U \bar{w} \rangle \delta^{ij} e_i \otimes e_j \\
&= \langle v, \bar{w} \rangle \delta^{ij} e_i \otimes e_j \\
&= \rho_+(\langle 1 ~ 2 \rangle)(v \otimes w),
\end{equs}
where in the fourth line, we used that $U$ is a Unitary matrix, so that it preserves inner products. These considerations directly generalize to give that
\begin{equs}\label{eq:unitary-rho-plus-identity}
(U^{\otimes n} \otimes \bar{U}^{\otimes m}) \rho_+(\langle i ~ j \rangle) = \rho_+(\langle i ~ j \rangle) = \rho_+(\langle i ~ j \rangle) (U^{\otimes n} \otimes \bar{U}^{\otimes m}). 
\end{equs}}
\revision{This identity will be needed for later remarks, but will strictly be needed for the later proofs.}
\end{example}

Recalling that we may view $\symgrp_n$ as embedded in $\mc{B}_{n}$, the restriction of the representation $\rho_+$ to $\symgrp_n$ defines a representation of $\symgrp_n$.

\begin{definition}\label{def:rho}
Let $N, n \geq 1$. Define the representation $\rho : \C[\symgrp_n] \ra \End((\C^N)^{\otimes n})$ to be the restriction of $\rho_+ : \mc{B}_n \ra \End((\C^N)^{\otimes n})$ to $\C[\symgrp_n] \sse \mc{B}_n$.
\end{definition}

Again, we omit the dependence of $\rho$ on $N, n$ for notational brevity. 

\begin{remark}
One may verify that $\rho$ has the following explicit form on pure tensors:
\begin{equs}
\rho(\sigma) (v_1 \otimes \cdots \otimes v_n) = v_{\sigma(1)} \otimes \cdots \otimes v_{\sigma(n)}, ~~ \sigma \in \symgrp_n, ~~ v_1, \ldots, v_n \in \C^N.    
\end{equs}
In words, $\rho(\sigma)$ acts by permutation of tensors.
\end{remark}

Next, we discuss how the formula \eqref{eq:weingarten-jucy-murphy-large-N} needs to be modified when $N$ is general. First,  recall that when $N \geq n$, the Weingarten function may also be defined as the inverse of $\sum_\sigma N^{\cycles(\sigma)} \sigma$ in $\C[\symgrp_n]$. For general $N$, this inverse may not exist. However, we quote the following result from \cite{collins2006integration}, which says that the Weingarten function can still be interpreted as an inverse, in a suitable sense.

\begin{lemma}[Section 2 of \cite{collins2006integration}]\label{lemma:general-N-weingarten-as-inverse}
Let $N, n \geq 1$. We have that $\rho\big(\sum_\sigma N^{\cycles(\sigma)} \sigma\big)$ is invertible, with inverse given by $\rho(\Wg_N)$.
\end{lemma}

\begin{remark}\label{remark:rho-motivation}
This is the whole point of introducing the representation $\rho$, in that the image $\rho(\sum_{\sigma} N^{\cycles(\sigma)} \sigma)$ is always invertible as a matrix, even though $\sum_\sigma N^{\cycles(\sigma)} \sigma$ may not be invertible in $\C[\symgrp_n]$. The simplest example of this difference is when $N = 1$ and $n = 2$, in which case $\sum_\sigma N^{\cycles(\sigma)} \sigma = \id + (1 ~ 2)$. Now clearly, $\id + (1 ~ 2)$ is not invertible in $\C[\symgrp_2]$, because the inverse would be given in general by $a \cdot \id + b \cdot (1 ~ 2)$, where $a, b \in \C$ solve the following system of equations:
\begin{equs}
a + b &= 1, \\
a + b &= 0.
\end{equs}
On the other hand, when $N = 1$, the space $(\C^N)^{\otimes n}$ is one-dimensional no matter the value of $n$. On this space, both $\rho(\id)$ and $\rho((1 ~ 2))$ are the identity operator. (Recall that $\rho((1 ~2)) (u \otimes v) = v \otimes u$. If $u, v  \in \C$, then $v \otimes u = u \otimes v$, so that $\rho((1 ~ 2))(u \otimes v)  = u \otimes v$.) Thus $\rho(\id + (1 ~ 2))$ acts as multiplication by $2$, and thus $\rho(\id + (1 ~ 2))^{-1}$ is multiplication by $1/2$.
\end{remark}

Similarly, we next show that the elements $N + J_k$, $k \in [n]$ are always invertible, if we apply the representation $\rho_+$. \revision{Following up on Remark \ref{remark:rho-motivation}, consider again the simple example $n = 2$. In this case, the Jucys-Murphy element $J_2 = \groupid + (1 ~ 2)$, which in the remark was shown not to be invertible in $\C[S_2]$, but on the other hand we saw that $\rho(\groupid + (1 ~ 2))$ is invertible as a matrix.}

\begin{lemma}\label{lemma:rho-plus-Jucy-Murphy-eigenvalues}
Let $N, n \geq 1$. Let $\rho_+ : \mc{B}_{n, n} \ra \End((\C^N)^{\otimes 2n})$ be the representation\footnote{Recall that $\mc{B}_{n, n} \sse \mc{B}_{2n}$. The representation $\rho_+$ is originally defined on $\mc{B}_{2n}$, here we restrict it to $\mc{B}_{n, n}$.} from Definition \ref{def:rho-plus}. For all $k \in [n]$, all eigenvalues of $\rho_+(J_k)$ are at least $-N + 1$.
\end{lemma}
\begin{proof}
Due to our embedding of $\symgrp_n$ into $\mc{B}_{n, n}$, $\rho_+(J_k)$ acts as the identity on the last $n$ coordinates of $(\C^N)^{\otimes 2n}$. On the first $n$ coordinates, $\rho_+(J_k)$ acts as $\rho(J_k)$ (as defined in Definition \ref{def:rho}) . Thus, it suffices to show that all eigenvalues of $\rho(J_k)$ are at least $-N + 1$. This follows from the combination of two classic results in the representation theory of the symmetric group:
\begin{enumerate}
    \item By Schur-Weyl duality (see e.g. \cite[Theorem 2.1]{collins2006integration}), we have that in the decomposition of \revision{the representation} $\rho$ into irreps, only those irreps corresponding to Young diagrams $\lambda$ with at most $N$ rows (i.e. $\ell(\lambda) \leq N$) appear.
    \item Let $\rho^\lambda$ be the irrep corresponding to $\lambda$. The eigenvalues of $\rho^\lambda(J_k)$ are explicitly known: for each Young tableaux with shape $\lambda$, let $(i, j)$ be the coordinates of the box which contains the integer $k$. Here, $i$ is the row index and $j$ the column index. Then $\rho^\lambda(J_k)$ has an eigenvalue equal to $j-i$. Moreover, all eigenvalues of $\rho^\lambda(J_k)$ arise this way. This result was proven by Jucys \cite{jucys1974symmetric} and independently later by Murphy \cite{murphy1981new}.
\end{enumerate}
The second fact implies that the every eigenvalue of $\rho^\lambda(J_k)$ is at least $-\ell(\lambda) + 1$, since the box with the most negative value of $j-i$ is $(\ell(\lambda), 1)$. Combining this with the first fact, the desired result now follows.
\end{proof}

This lemma shows that all for all $k \in [n]$, all eigenvalues of $\rho_+(N + J_k)$, are at least $1$, and thus $\rho_+(N +J_k)$ is invertible. Moreover, we have the following lemma, which generalizes \eqref{eq:weingarten-jucy-murphy-large-N} to the case of general $N$.

\begin{lemma}\label{lemma:rho-plus-weingarten-is-inverse-rho-plus-cycles}
Let $N, n \geq 1$. We have that
\begin{equs}
\rho_+(\Wg_N) = \rho_+(N + J_n)^{-1} \cdots \rho_+(N + J_1)^{-1}.
\end{equs}
\end{lemma}
\begin{proof}
Due to our embedding of $\symgrp_n$ into $\mc{B}_{n, n}$, for any element $f \in \C[\symgrp_n]$, the matrix $\rho_+(f) \in \End((\C^N)^{\otimes 2n})$ acts as the identity on the last $n$ coordinates of $(\C^N)^{\otimes 2n}$. On the first $n$ coordinates, $\rho_+(f)$ acts as $\rho(f) \in \End((\C^N)^{\otimes n})$. Thus it suffices to prove the claimed identity with $\rho_+$ replaced by $\rho$. Since $\rho$ is a representation, we have that (using that the Jucys-Murphy elements commute with each other and applying \eqref{eq:jucys-identity} in the final identity)
\begin{equs}
\rho(N + J_n)^{-1} \cdots \rho(N + J_1)^{-1} = \rho((N + J_1) \cdots (N + J_n))^{-1} = \rho\bigg(\sum_\sigma N^{\cycles(\sigma)} \sigma \bigg)^{-1}.
\end{equs}
The desired result now follows by Lemma \ref{lemma:general-N-weingarten-as-inverse}.
% Next, by definition, we have that $\Wg_N$ is the inverse of $\sum_\sigma N^{\cycles(\sigma)} \sigma$ in $\C_N[\symgrp_n]$, which implies that
% \begin{equs}
% \rho(\Wg_N) \rho\bigg(\sum_\sigma N^{\cycles(\sigma)} \sigma \bigg) = I, ~~ \text{ i.e. } \rho(\Wg_N) = \rho\bigg(\sum_\sigma N^{\cycles(\sigma)} \sigma \bigg)^{-1}.
% \end{equs}
% The desired result now follows.
\end{proof}

In the course of proving Theorem \ref{thm:weingarten-recovery} for general values of $N$, we will also need the following technical lemma.

\begin{lemma}\label{lemma:eigenvalue-lower-bound-sum-all-Jucys-Murphy}
Let $N, n \geq 1$. All eigenvalues of
\begin{equs}
\frac{1}{N} \rho(J_n + \cdots + J_1) \in \End((\C^N)^{\otimes n})
\end{equs}
are at least $-\frac{n}{2} + \frac{1}{2}$. More precisely, if $n = mN + r$ with $0 \leq r \leq N-1$, then all eigenvalues are at least
\begin{equs}
-\frac{n}{2} + \frac{m^2}{2} + \frac{r}{2} - \frac{1}{2} \frac{r(r-1)}{N} + \frac{mr}{N}.
\end{equs}
% Finally, when $n = N$, the matrix $\frac{1}{N} \rho(J_N + \cdots + J_1)$ has eigenvalue $-\frac{N}{2} + \frac{1}{2}$ with multiplicity 1, and all other eigenvalues are at least $-\frac{N}{2} + \frac{3}{2}$. 
\end{lemma}
\begin{proof}
As noted in the proof of Lemma \ref{lemma:rho-plus-Jucy-Murphy-eigenvalues}, by Schur-Weyl duality (see e.g. \cite[Theorem 2.1]{collins2006integration}), we have that in the decomposition of $\rho$ into irreps, only those irreps corresponding to Young diagrams $\lambda$ with at most $N$ rows (i.e. $\ell(\lambda) \leq N$) appear. Thus letting $\rho^\lambda$ be the irrep corresponding to $\lambda$, it suffices to show the claim with $\rho$ replaced by $\rho^\lambda$, for any Young diagram $\lambda$ with at most $N$ rows. 

Towards this end, let $\lambda$ be a Young diagram, for example as in Figure \ref{figure:young-diagram}. \revision{Given a box $(i, j) \in \lambda$, its content $c(i, j)$ is the integer $j - i$.}
\begin{figure}[h]
    \centering
    
  \minipage{0.48\textwidth}
    \includegraphics[width=\linewidth]{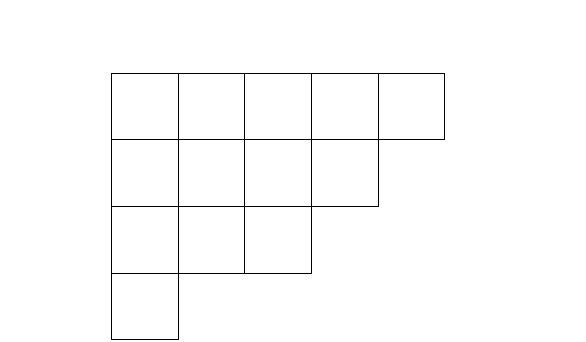}
    \caption{$5 + 4 + 3 + 1 = 13$}\label{figure:young-diagram}
   \endminipage\hfill
  \minipage{0.45\textwidth}
    \includegraphics[page=3,width=\linewidth]{figures/young_diagram.pdf}
    \caption{Young diagram with content labels}\label{figure:young-diagram-content-labeling}
    \endminipage
\end{figure}
As discussed in the proof of Lemma \ref{lemma:rho-plus-Jucy-Murphy-eigenvalues}, for $k \in [n]$ the eigenvalues of $\rho^\lambda(J_k)$ are given by the content of the $k$th box when we range over standard Young tableaux with shape $\lambda$. Even more, \cite{jucys1974symmetric, murphy1981new} show that the $(\rho^\lambda(J_k), k \in [n])$ have a joint eigenbasis indexed by standard Young tableaux with shape $\lambda$, where the eigenvalues corresponding to a given standard Young tableaux are the contents of the boxes of the Young diagram. This discussion shows that on each eigenbasis element, $\rho^\lambda(J_n + \cdots + J_1)$ acts in the same manner, that is as a whole $\rho^\lambda(J_1 + \cdots + J_n)$ acts as a multiple $c_\lambda$ of the identity, where $c_\lambda$ is the sum of contents of all the boxes in $\lambda$.

To envision the computation of $c_\lambda$, we label each box of $\lambda$ with its content, i.e. the number $j - i$, where $(i, j)$ is the row-column coordinate of the box. For the Young diagram in Figure \ref{figure:young-diagram}, we have the labeling in Figure \ref{figure:young-diagram-content-labeling}.
% \begin{figure}[h]
%     \centering
%     \includegraphics[page=3]{figures/young diagram.pdf}
%     \caption{Young diagram with content labels}\label{figure:young-diagram-content-labeling}
% \end{figure}
The constant $c_\lambda$ is then the sum of all box labels. For example, for the Young diagram in Figures \ref{figure:young-diagram} and \ref{figure:young-diagram-content-labeling}, $c_\lambda = 6$. 

Now, fix $n, N \geq 1$. To prove the lemma, we need to understand how negative the content sum $c_\lambda$ may be for a Young diagram with $n$ boxes and at most $N$ rows. Clearly, to minimize $c_\lambda$, we want a Young diagram with as many columns of size $N$ as possible. Thus, if $n = mN + r$ with $0 \leq r \leq N-1$, then the Young diagram in Figure \ref{figure:young-diagram-minimizer} minimizes $c_\lambda$.

\begin{figure}[ht!]
    \centering
    \includegraphics[page=2]{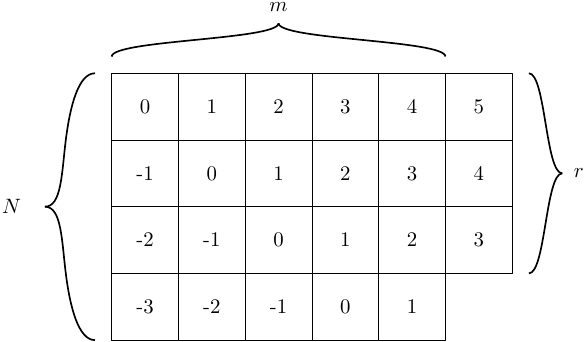}
    \caption{When $n = mN + r$, the ``worst case" Young diagram in terms of smallest content sum.}\label{figure:young-diagram-minimizer}
\end{figure}

The content sum $c_\lambda$ of such a diagram (by first summing the contents along each column) is 
\begin{equs}
c_\lambda &= -\binom{N}{2} + \bigg(-\binom{N}{2} + N\bigg) + \cdots + \bigg(-\binom{N}{2} + (m-1) N\bigg) + \bigg(-\binom{r}{2} + mr\bigg) \\
&= - m \binom{N}{2} + \binom{m}{2} N - \binom{r}{2} + mr.
\end{equs}
% Our goal is to show that $\frac{1}{N} c_\lambda \geq -\frac{1}{2} n + \frac{1}{2}$.
From this, we may obtain
\begin{equs}
\frac{1}{N} c_\lambda &= - \frac{1}{2} m (N - 1) + \frac{1}{2} m(m-1) - \frac{1}{2} \frac{r(r-1)}{N} + \frac{mr}{N} \\
&= -\frac{1}{2}(mN + r) + \frac{1}{2} m^2 + \frac{1}{2} r - \frac{1}{2} \frac{r(r-1)}{N} + \frac{mr}{N}.
\end{equs}
Since $n = mN + r$, this proves the second claim. Moreover, we see that if $m \geq 1$, then the above is at least $-\frac{1}{2}n + \frac{1}{2} m^2 \geq -\frac{1}{2} n + \frac{1}{2}$, as desired. Now, suppose that $m = 0$, so that $n = r$. Then the above is equal to
\begin{equs}
-\frac{1}{2} n + \frac{1}{2} r \bigg(1 + \frac{1}{N} - \frac{r}{N}\bigg).
\end{equs}
One may check that under the restriction $1 \leq r \leq N-1$, the above is minimized at $r = 1$ with a value of $-\frac{1}{2} n + \frac{1}{2}$, as desired.
% Finally, for the last claim, note that $n = N$ corresponds to $(m, r) = (1, 0)$. By our previous formula, this means that the lowest eigenvalue of $\frac{1}{N}\rho(J_N + \cdots + J_1)$ is $-\frac{N}{2}  + \frac{1}{2}$.
\end{proof}

Recall the Jucys-Murphy elements acting on bottom vertices defined in Definition \ref{def:bottom-jucys-murphy}.

\begin{cor}\label{cor:eigenvalue-lower-bound-sum-Jucys-Murphy}
Let $N, n \geq 1$. Let $\rho_+ : \mc{B}_{n, n+1} \ra \End((\C^N)^{\otimes 2n+1})$. All eigenvalues of
\begin{equs}
\frac{1}{N} \rho_+(J_n + \cdots + J_1 + J_{n+1}' + \cdots + J_1')
\end{equs}
are strictly greater than $-n$.
\end{cor}
\begin{proof}
Observe that $\rho_+(J_n + \cdots + J_1)$ acts as the identity on the last $n+1$ coordinates, and on the first $n$ coordinates, $\rho_+(J_n + \cdots + J_1)$ acts as $\rho(J_n + \cdots + J_1)$. In other words, $\rho_+(J_n + \cdots + J_1) = \rho(J_n + \cdots + J_1) \otimes I_{n+1}$, where $I_{n+1} \in \End( (\C^N)^{\otimes (n+1)})$ is the identity. Similarly, we have that $\rho_+(J_{n+1}' + \cdots + J_1') = I_n \otimes \rho(J_{n+1} + \cdots + J_1)$, where $I_n \in \End((\C^N)^{\otimes n})$ is the identity. In general, given two matrices $M_1, M_2$, the eigenvalues of $M_1 \otimes M_2$ are the products of the eigenvalues of $M_1$ and eigenvalues of $M_2$. Combining this fact with Lemma \ref{lemma:eigenvalue-lower-bound-sum-all-Jucys-Murphy}, we obtain that all eigenvalues of $\frac{1}{N}\rho(J_n + \cdots + J_1) \otimes I_{n+1}$ are at least $-\frac{1}{2}n + \frac{1}{2}$, and all eigenvalues of $I_n \otimes \frac{1}{N}\rho(J_{n+1}' + \cdots + J_1')$ are at least $-\frac{1}{2}(n+1) + \frac{1}{2}$. The desired result now follows.
\end{proof}

% We may assume after cyclic reordering that the first edge $(x_1, x_2)$ is the black strand which corresponds to $M_{c(1)}^{\varep(1)}$. In general, given an edge $e = (x_i, x_{i+1})$ corresponding to a black strand, let $i(e) \in [n]$ be the index of the letter that this strand corresponds to -- e.g., $i(x_1, x_2) = 1, i(x_3, x_4) = 2$, $i(x_5, x_6) = 3$, etc.
% We then form the sum
% \begin{equs}
% \Tr(M(\Gamma)) = \prod_{e \in B}  \prod_{e \in R} \delta^{i_j i_{j+1}}
% \end{equs}
% \begin{equs}
% \Tr(M(\Gamma)) = (M_{c(1)}^{\varep(1)})_{i_1 i_2} \delta^{i_2 i_3} (M_{c(2)}^{\varep(2)})_{i_3 i_4} \delta^{i_4 i_5} \cdots (M^{\varep(n)}_{c(n)})_{i_{V-1} i_V} \delta^{i_V i_1}.
% \end{equs}

\section{Surface-sum representation of Wilson loop expectations}\label{section:wilson-loop-expectation-epe}

In this section, we show how to apply the Unitary Weingarten calculus (Theorem \ref{thm:weingarten}) to express Wilson loop expectations as sums over edge-plaquette embeddings (which were introduced in Section \ref{section:edge-plaquette-embedding}). We first prove a more abstract result about expectations of traces of words of Haar distributed Unitary matrices (Proposition \ref{prop:word-expectation-as-sum-over-bipartite-maps}) which has no reference to a lattice, and then apply this result to Wilson loop expectations to obtain Theorem \ref{thm:wilson-loop-expectation-sum-over-epe}. 

\begin{definition}\label{def:normalized-weingarten}
Define the normalized Weingarten function $\ovl{\Wg}_N$ by:
\begin{equs}
\ovl{\Wg}_N(\pi) := N^{n + \revision{|\pi|}} \Wg_N(\pi), ~~ \pi \in \symgrp_n.
\end{equs}
Here, $\revision{|\pi|} := n - \#\mrm{cycles}(\pi)$.
\end{definition}

\begin{remark}\label{remark:weingarten-large-N-limit}
We will see later on that the normalized Weingarten function is the more natural quantity to work with, as it leads to nicer statements of our formulas. Another nice thing about $\ovl{\Wg}_N$ is that with this choice of normalization, the limit as $N \toinf$ exists and depends on $\pi$. Indeed, we in fact have (see e.g. \cite[Corollary 2.7]{collins2006integration})
\begin{equs}
\ovl{\Wg}_N(\pi) = \text{M\"{o}b}(\pi) + O(N^{-2}) \text{ as $N \toinf$},
\end{equs}
where if $\pi$ is decomposed into cycles of lengths $C_1, \ldots, C_k$, then
\begin{equs}\label{eq:mobius}
\text{M\"{o}b}(\pi) :=  \prod_{i \in [k]} (-1)^{C_i - 1}\mrm{Cat}(C_i - 1),
\end{equs}
where $\mrm{Cat}(k) := \frac{(2k)!}{k! (k+1)!}$ is the $k$th Catalan number.
\end{remark}

Recall from Theorem \ref{thm:weingarten}
that expectations of traces of words of Haar-distributed Unitary matrices may be expressed in terms of a diagrammatic sum over matchings, with matchings weighted by the Weingarten function. In this section, we will use this to express Wilson loop expectations in lattice gauge theories as weighted sums over edge-plaquette embeddings. The main step is to describe how to obtain a map from a given balanced collection of words $\bm \Gamma$ (which recall specifies the exterior connections) along with interior connections specified by matchings. \revision{Recall from the beginning of Section \ref{section:edge-plaquette-embedding} that a map is a collection of polygons with specified gluings. Thus in the following, we will describe how, given $\bm \Gamma$ and a choice of matchings, we may obtain a collection of faces along with a complete pairing of the set of edges (i.e. the specified gluings).}

Recall the strand diagram $\mrm{SD}(\bm \Gamma)$ introduced in Section \ref{section:strand-diagrams-and-weingarten-caculus}. As in Figure \ref{figure:strand_diagram_example}, the connected components of the strand diagram are in 1-1 correspondence with the words of $\bm \Gamma$. These connected components may be thought of as polygonal faces. A connected component of size $2k$ can be thought of as a $k$-gon\footnote{Due to the way the strand diagram is defined, the connected components are always of even size, because they are found by alternately traversing a black strand and red exterior connection, until arriving back at the starting point.}, by shrinking away the red exterior connections specified by $\bm \Gamma$, as mentioned in Figure \ref{figure:strand_diagram_example}. In summary, the choice of $\bm \Gamma$ specifies a collection of faces, which can be seen in the strand diagram. \revision{However, as we will see, these are not the only faces which will appear in our map. Rather, these will give the plaquette faces, while the blue faces are still to be obtained.}
% (so a connected component of size $k$ can be thought of as a $k$-gon).

% We begin with a collection of faces corresponding to the words of $\bm \Gamma$. Each word $\Gamma_i$ gives a face whose degree (i.e.~number of boundary edges) is the length of $\Gamma_i$. The boundary edges of each such face are naturally labeled by letters in $\{\lambda_1, \ldots, \lambda_L\}$.
% Recalling the strand diagram $\mrm{SD}(\bm \Gamma)$ introduced in Section \ref{section:strand-diagrams-and-weingarten-caculus}, these faces can be thought of as the connected components of $\mrm{SD}(\bm \Gamma)$. 
% These faces can be obtained from adding the exterior connections specified by $\bm \Gamma$ to the strand diagram, as in Figure \ref{figure:strand_diagram_example} (think of the red exterior connections as being shrunk down to a single vertex).

\begin{figure}[h]%
    \centering
    \includegraphics[width=7cm]{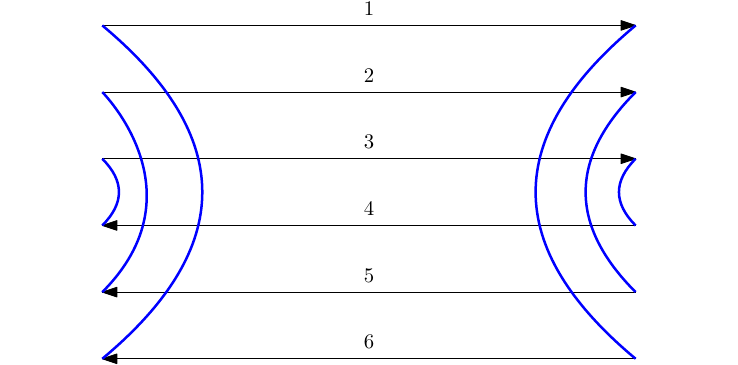}%
    \qquad
    \includegraphics[page=2, width=7cm]{figures/bipartite-map-from-matchings.pdf}%
    \caption{Left: a pair of left and right bijections. Right: the faces obtained by ``shrinking away" the blue matching edges.}\label{figure:bipartite-map-from-matchings}
\end{figure}

Next, consider a pair of left and right bijections of the vertices of a number of strands as displayed in the left of Figure \ref{figure:bipartite-map-from-matchings}. (Recall that the Weingarten calculus involves summing over such pairs.) Think of this as the portion of the diagram corresponding to some letter $\lambda$ in $\{\lambda_1, \ldots, \lambda_L\}$. One can imagine that the two endpoints of each blue line are identified (this corresponds to ``shrinking" each blue line away). In this case, one is then left with a collection of faces as in the right of Figure \ref{figure:bipartite-map-from-matchings}.

\begin{figure}[h]%
    \centering
    {{\includegraphics[page=3, width=7cm]{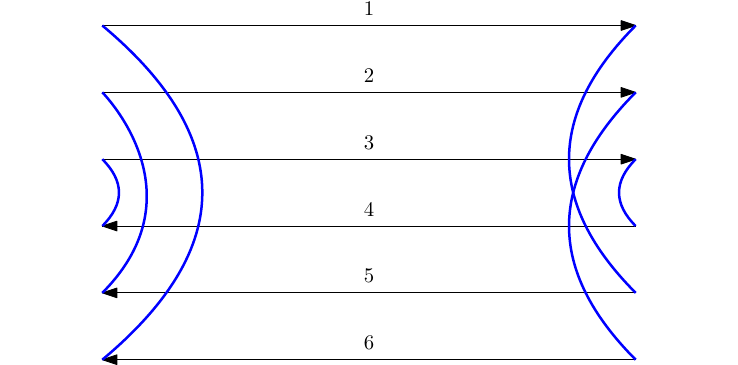} }}%
    \qquad
    {{\includegraphics[page=4, width=7cm]{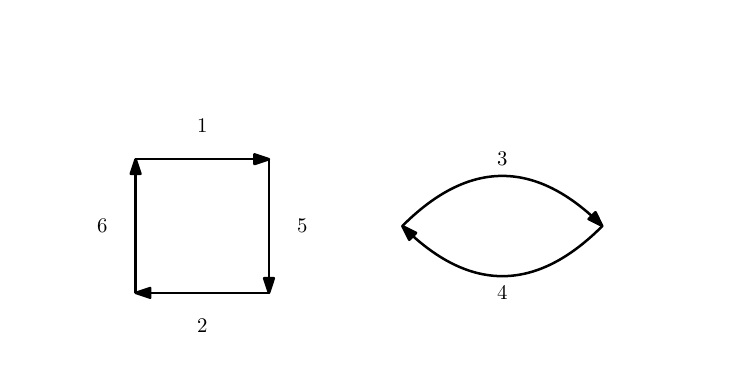} }}%
    \caption{Left: a pair of left and right bijections. Right: the faces obtained by ``shrinking away" the blue matching edges.}\label{figure:bipartite-map-from-matchings-4}
\end{figure}

In Figure \ref{figure:bipartite-map-from-matchings-4}, we give another example of a pair of bijections along with the corresponding collection of faces.

By specifying a pair of left and right bijections $(\sigma_\ell, \tau_\ell)$ for each letter $\lambda_\ell$, we can obtain another collection of faces, in the manner described in Figures \ref{figure:bipartite-map-from-matchings} and \ref{figure:bipartite-map-from-matchings-4}. We thus naturally have two collections of faces: the set of faces which correspond to words in $\bm \Gamma$, and the set of faces obtained from the pairs of left and right bijections.

\begin{convention}\label{conv:faces}
We refer to the faces which correspond to words in $\bm \Gamma$ as ``plaquette-faces", or ``yellow faces". We refer to the faces which are obtained from the pairs of left and right bijections as ``edge-faces", or ``blue faces".
\end{convention}

% \minjae{Bipartite map sounds like the map itself (not the dual) is bipartite. Maybe dual bipartite map (DBM) or face-bipartite map (FBM)?}

Observe that every edge is incident to exactly two faces -- one blue face and one yellow face. This naturally induces a gluing of the faces, and so we obtain a map whose dual is bipartite from the data $(\bm \Gamma, ([\sigma_\ell, \tau_\ell], \ell \in [L]))$. 

We next observe that the number of vertices of this map is precisely the number of components of the strand diagram.
% , the latter of which appears in the formula for the expectation of products of traces of words.
% This relation is captured in the following lemma.

\begin{lemma}\label{lemma:vertices-equals-components}
Let $\bm \Gamma = (\Gamma_1, \ldots, \Gamma_k)$ be a balanced collection of words on $\{\lambda_1, \ldots, \lambda_L\}$. Suppose that for each $\ell \in [k]$, the number of occurrences of $\lambda_\ell$ is $n(\lambda_\ell)$. For $\ell \in [k]$, let $\sigma_\ell, \tau_\ell : [n(\lambda_\ell)] \ra (n(\lambda_\ell) : 2n(\lambda_\ell)]$ be a pair of bijections. Then $\numcomp(\bm \Gamma,\ ([\sigma_\ell, \tau_\ell],\ \ell \in [L]))$ is equal to the number of vertices in the corresponding map.
\end{lemma}
\begin{proof}
To compute the number of vertices in the map, we can proceed as follows. Recalling that the map arises from combining an interior connection with an exterior connection of the strand diagrams, we may begin by giving each vertex in the strand diagram a unique label. Thus for the portion of the strand diagram corresponding to $\lambda_\ell$, we give a total of $4n_\ell$ labels, since there are $2n_\ell$ strands. Each connection (be it interior or exterior) results in the identification of two labels. In terms of the map, labels which have been identified are in fact the same vertex. Therefore the number of vertices in the map corresponds to the number of different equivalence classes of labels, after performing all label identifications indicated by the connections. The equivalence class of a given label may be obtained by starting at the label, and alternately following the exterior and interior connections, until we arrive back at the initial label. Recalling Example \ref{ex:computing-connected-components}, observe that this is precisely the same method for computing the number of connected components of a given strand diagram with interior and exterior connections. Thus the connected components of the strand diagram are in bijection (moreover, there is a canonical identification) with the vertices of the map.
\end{proof}

% \minjae{Shouldn't we define `blue' and `yellow' faces? When explaining the boundary, I guess we only delete blue faces? E.g. define the faces whose boundary only consists of strands as ``yellow/plaquette faces'' and the other faces as ``blue/edge faces''.}
% \sky{I added Convention 3.3. In explaining the boundary in remark 3.10, we only delete certain yellow faces, i.e. those which correspond to loops.}

\revision{Before we make the following definition, we emphasize here that the outcome of the preceding discussion is that given the data $(\bm \Gamma, ([\sigma_\ell, \tau_\ell], \ell \in [L]))$, we have described how to obtain a map whose dual is bipartite. (And moreover, the number of vertices of this map is precisely the number of components in the diagram corresponding to $(\bm \Gamma, ([\sigma_\ell, \tau_\ell], \ell \in [L]))$.) Since Theorem \ref{thm:weingarten} involves a sum over $(\bm \Gamma, ([\sigma_\ell, \tau_\ell], \ell \in [L]))$, we will thus be able to interpret the expectations of products of traces of words in independent Haar Unitaries as a sum over maps. We next define the set of such maps.}

\begin{definition}
Let $\bm \Gamma = (\Gamma_1, \ldots, \Gamma_k)$ be a balanced collection of words on $\{\lambda_1, \ldots, \lambda_L\}$. Define $\bpm(\bm \Gamma)$ (short for ``dual bipartite map") to be the set of all possible maps which can be obtained from adding pairs of left and right bijections to the strand diagram of $\bm \Gamma$. For a given map $\mc{M} \in \bpm(\bm\Gamma)$, and $\ell \in [L]$, let $\mu_\ell(\mc{M})$ be the partition of $n(\lambda_{\ell})$ (the total number of occurrences of $\lambda_\ell$) given by $1/2$ times the degrees of the blue faces which are glued in to the strand diagram of $\lambda_\ell$.
\end{definition}

\begin{remark}
All maps in $\bpm(\bm \Gamma)$ are orientable. The faces coming from words of $\bm \Gamma$ are endowed with a natural orientation. It follows from construction that the edge-faces can then be endowed with orientations that are consistent with the existing orientations. This will be in contrast to the Orthogonal and Symplectic cases in Section \ref{section:orthogonal-and-symplectic}, where the constructed maps are not necessarily orientable.
% The faces corresponding to a word in $\bm \Gamma$ are endowed with a natural orientation (given by traversing the word). The faces coming from interior matchings can then always be endowed with a consistent orientation. For instance, in Figures \ref{figure:bipartite-map-from-matchings} and \ref{figure:bipartite-map-from-matchings-4}, the orientations of these faces should be the reverse of what is drawn. 
\end{remark}

\begin{remark}
Observe that for any $\mc{M} \in \bpm(\bm \Gamma)$,  we have that 
\begin{equs}\label{eq:edge-and-face-identities}
E(\mc{M}) &= 2\sum_{i \in [L]} n(\lambda_i), ~~ F(\mc{M}) = k + \sum_{i \in [L]} \ell(\mu_i(\mc{M})).
\end{equs}
Here, $\ell(\mu_i(\mc{M}))$ is the number of parts of the partition $\mu_i(\mc{M})$. The first identity says that the number of edges is equal to the total number of strands in the strand diagram, and the second identity says that the total number of faces is equal to the number of words plus the total number of cycles of the interior connections of the strand diagram.
\end{remark}

\begin{prop}\label{prop:word-expectation-as-sum-over-bipartite-maps}
Let $\bm \Gamma = (\Gamma_1, \ldots, \Gamma_k)$ be a balanced collection of words on $\{\lambda_1, \ldots, \lambda_L\}$. We have that
% \begin{equs}
% \E[ U(\Gamma)] := \E[ U(\Gamma_1) \cdots U(\Gamma_k) ] = \sum_{\mc{M} \in \bpm(\Gamma_1, \ldots, \Gamma_k)} \bigg(\prod_{\ell \in [L]} \Wg_N(\mu_\ell(\mc{M})) \bigg)N^{V(\mc{M})}.
% \end{equs}
\begin{equs}
\E[ \Tr(U(\bm \Gamma))] = \sum_{\mc{M} \in \bpm(\bm \Gamma)} \bigg(\prod_{\ell \in [L]} \ovl{\Wg}_N(\mu_\ell(\mc{M})) \bigg)N^{\chi(\mc{M}) - k}.
\end{equs}
\end{prop}
\begin{proof}
By Theorem \ref{thm:weingarten}, the definition of $\bpm(\bm \Gamma)$, and Lemma \ref{lemma:vertices-equals-components}, we have that
\begin{equs}
\E[ \Tr(U(\bm \Gamma))] = \sum_{\mc{M} \in \bpm(\bm \Gamma)} \bigg(\prod_{\ell \in [L]} \Wg_N(\mu_\ell(\mc{M})) \bigg)N^{V(\mc{M})}.
\end{equs}
Applying the identities \eqref{eq:edge-and-face-identities}, we further obtain 
\begin{equs}
\E[\Tr( U(\bm \Gamma))] = \sum_{\mc{M} \in \bpm(\bm \Gamma)} \bigg(\prod_{i \in [L]} N^{n(\lambda_i) + \revision{|\mu_i(\mc{M})|}} \Wg_N(\mu_i(\mc{M})) \bigg) N^{V(\mc{M}) - E(\mc{M}) + F(\mc{M}) -k}.
% \end{equs}
% where $\|\mu_i(\mf{m})\| := n_i - \ell(\mu_i(\mf{m}))$. 
\end{equs}
The desired result now follows. \qedhere
% \begin{equs}
% N^k \E[\Gamma_1 \cdots \Gamma_k] = \sum_{\mf{m} \in \bpm(\Gamma_1, \ldots, \Gamma_k)} \bigg(\prod_{i \in [L]} N^{n_i + \|\mu_i(\mf{m})\|} \Wg_N(\mu_i(\mf{m})) \bigg) N^{\chi(\mf{m})},
% \end{equs}
% where $\|\mu_i(\mf{m})\| := n_i - \ell(\mu_i(\mf{m}))$. 
\end{proof}

Next, suppose the letters $\{\lambda_1, \ldots, \lambda_L\}$ are edges of the lattice $\Lambda$. Suppose the collection of words $\bm \Gamma = \bm \Gamma(s, K)$ is made of a string $s$ and a plaquette count $K : \mc{P} \ra \N$. Here, a string $s = (\ell_1, \ldots , \ell_n)$ is a collection of loops. Each loop is a word. The plaquette count $K$ specifies that there should be $K(p)$ copies of the boundary of $p$ in $\bm \Gamma$. We say that $(s, K)$ is balanced if the corresponding $\bm \Gamma$ is balanced.

\begin{figure}[h]%
    \centering
    {{\includegraphics[width=10cm]{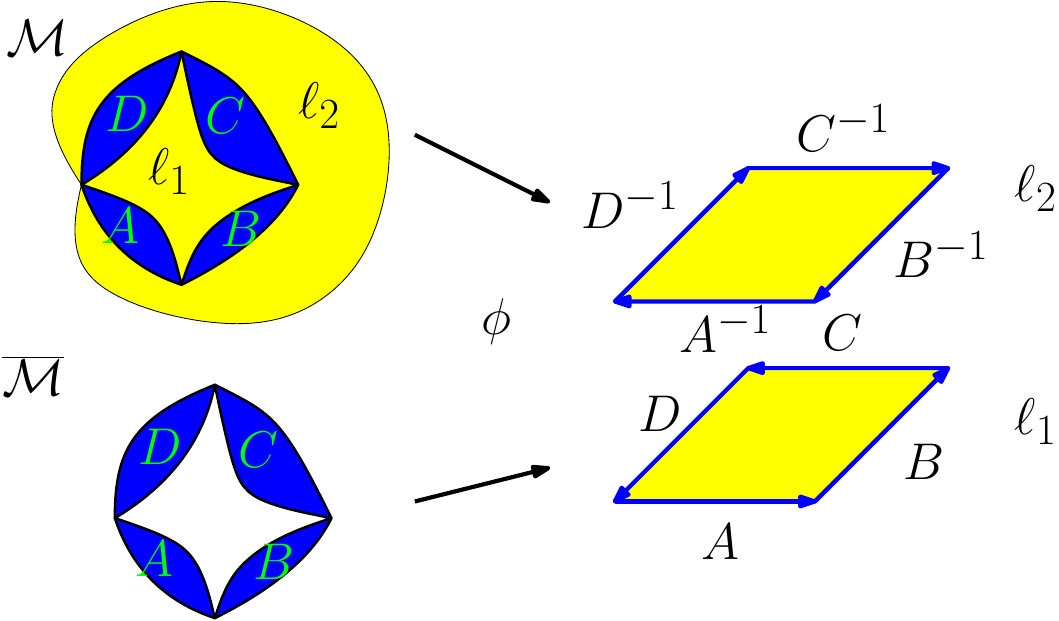} }}
    \caption{An example of an edge-plaquette embedding with boundary when $\ell_1=ABCD$ and $\ell_2=D^{-1}C^{-1}B^{-1}A^{-1}$. The top left sphere is a map $\mc M$ whose dual is bipartite and the bottom left map $\ovl {\mc M}$ is obtained by removing two yellow faces corresponding to $\ell_1$ and $\ell_2$. The Euler characteristic changes from 2 to 0 after removal of the two faces. The boundary of $\ovl{\mc{M}}$ maps onto the edges of $\ell_1$ and $\ell_2$, and thus we can interpret $\ovl{\mc{M}}$ as having boundary given by the union of $\ell_1$ and $\ell_2$.}\label{fig:boundary}
\end{figure}

Recall the definition of an edge-plaquette embedding (Section \ref{section:edge-plaquette-embedding}). For balanced $(s, K)$, a map $\mc{M} \in \bpm(\bm \Gamma(s, K))$ induces an edge-plaquette embedding $(\ovl{\mc{M}}, \phi)$, where the function $\phi$ is determined by the requirement that it maps each edge of $\mc{M}$ (which is canonically labeled by a letter in $\{\lambda_1, \ldots, \lambda_L\}$, which we are assuming to be edges of the lattice $\Lambda$) to the corresponding edge of $\Lambda$. Here, $\ovl{\mc{M}}$ is obtained from $\mc{M}$ by deleting the faces of $\mc{M}$ which correspond to the loops $\ell_1, \ldots \ell_n$ (i.e. the faces whose boundaries are mapped by $\phi$ to the loops). This way, all faces of $\ovl{\mc{M}}$ are mapped to plaquettes of $\Lambda$, as required in the definition of edge-plaquette embeddings. Note that $\ovl{\mc{M}}$ has $n$ boundary components, which are mapped by $\phi$ to the loops $\ell_1, \ldots, \ell_n$. In this way, we interpret $(\ovl{\mc{M}}, \phi)$ as having ``boundary" $\ell_1, \ldots, \ell_n$.  See Figure~\ref{fig:boundary} for an illustration.

We will later use that
\begin{equs}\label{eq:euler-relation-boundary-epe}
\chi(\ovl{\mc{M}}) = \chi(M) - n,
\end{equs}
which follows since $\ovl{\mc{M}}$ is obtained from $\mc{M}$ by deleting $n$ faces.

We now apply these considerations to lattice Yang-Mills. Let $s$ be a string. Recall equation \eqref{eqn::latticeymexpanded}, which we reproduce here:
\begin{equs}\label{eq:wilson-loop-exponential-taylor-expansion}
\langle W_s \rangle_{\Lambda, \beta} = Z_{\Lambda, \beta}^{-1} \sum_{K : \mc{P} \ra \N} \frac{(N\beta)^K}{K!} \int W_s(Q) \prod_{p \in \mc{P}} \Tr(Q_p)^{K(p)} \prod_{e \in E_\Lambda} dQ_e.
\end{equs}
For each fixed $K : \mc{P} \ra \N$, we may apply Proposition \ref{prop:word-expectation-as-sum-over-bipartite-maps} to obtain an expression for the integral above in terms of a sum over edge-plaquette embeddings. We first set some notation.

% As in Section \ref{section:intro}, to compute the Wilson loop expectation
% \begin{equs}
% \langle W_s \rangle_{\Lambda, \beta} := Z_{\Lambda, \beta}^{-1} \int W_s(Q) \exp\bigg(2N \beta \sum_{p \in \mc{P}^+} \mrm{Re} \Tr(Q_p)\bigg) \prod_{e \in E_\Lambda} dQ_e,
% \end{equs}
% we may expand for each $p \in \mc{P}^+$:
% \begin{equs}
% \exp\big(2N \beta \mrm{Re}(\Tr(Q_p))\big) &= \sum_{k=0}^\infty \frac{(N\beta)^k}{k!} \big(\Tr(Q_p) + \Tr(Q_p^{-1})\big)^k \\
% &= \sum_{k_+, k_- = 0}^\infty \frac{(N\beta)^{k_+ + k_-}}{(k_+)! (k_-)!} \Tr(Q_p)^{k_+} \Tr(Q_p^{-1})^{k_-},
% \end{equs}
% and then insert this expansion to obtain:
% \begin{equs}\label{eq:wilson-loop-exponential-taylor-expansion}
% \langle W_s \rangle_{\Lambda, \beta} = Z_{\Lambda, \beta}^{-1} \sum_{K_+, K_- : \mc{P} \ra \N}  \frac{(N\beta)^{|K_+| + |K_-|}}{(K_+)! (K_-)!} \int W_s(Q) \prod_{p \in \mc{P}} \Tr(Q_p)^{K_+(p)} \Tr(Q_p^{-1})^{K_-(p)} \prod_{e \in E_\Lambda} dQ_e,
% \end{equs}
% where recall our notation $(K_+)! = \prod_{p \in \mc{P}^+} (K_+(p))!$, and we set $|K_+| := \sum_{p \in \mc{P}} K_+(p)$, and likewise for $|K_-|$. For each fixed $K_+, K_- : \mc{P} \ra \N$, we may apply Proposition \ref{prop:word-expectation-as-sum-over-bipartite-maps} to obtain an expression for the integral above in terms of a sum over edge-plaquette embeddings. We first set some notation.

\begin{definition}
Let $s = (\ell_1, \ldots, \ell_n)$ be a string, and let $K : \mc{P} \ra \N$ be a plaquette count. Define the set $\epe(s, K)$ of edge-plaquette embeddings associated to $s, K$ as follows. If $(s, K)$ is unbalanced, then $\epe(s, K) := \varnothing$. If $(s, K)$ is balanced, we define $\epe(s, K)$ to be the set of edge-plaquette embeddings $(\ovl{\mc{M}}, \phi)$ obtained from maps $\mc{M} \in \bpm(\bm \Gamma(s, K))$. In words, $\epe(s, K)$ is the set of edge-plaquette embeddings with plaquette counts specified by $K$ and which have ``boundary" $s$.
Next, define
\begin{equs}\label{eq:epe-partition}
\epe(s) := \bigsqcup_{K : \mc{P} \ra \N} \epe(s, K).
\end{equs}
We use the suggestive notation $\ptl (\ovl{\mc{M}}, \phi) = s$ to mean $(\ovl{\mc{M}}, \phi) \in \epe(s)$.

For $(\ovl{\mc{M}}, \phi) \in \epe(s)$, and $e \in E_\Lambda$, let $\mu_e(\phi)$ be the partition of $|\phi^{-1}(e)|/2$ induced by $1/2$ times the degrees of the faces of $\phi^{-1}(e)$. Define 
\begin{equs}
\area(\ovl{\mc{M}}, \phi) &:= \sum_{p \in \mc{P}} |\phi^{-1}(p)|, \\
(\phi^{-1})! &:= \prod_{p \in \mc{P}} |\phi^{-1}(p)|! .
\end{equs}
Note that if $(\ovl{\mc{M}}, \phi) \in \epe(s, K)$, then $\area(\ovl{\mc{M}}, \phi) = \sum_p K(p)$ and $(\phi^{-1})! = K!$. 
% Note also that $\area(\mc{M}, \phi)$ is equal to the total number of plaquette-faces of $(\mc{M}, \phi)$.
\end{definition}

% \begin{remark}[Boundaries of edge-plaquette embeddings]\label{rem:boundary-epe}
% One can also think of $\epe(s)$ as the set of edge-plaquette embeddings which have ``boundary" $s$, by deleting the plaquette-faces (see Convention~\ref{conv:faces}) of $\mc{M}$ which correspond to $\ell_1, \ldots, \ell_n$ (i.e. the plaquette-faces whose boundary is mapped by $\phi$ to one of edges of $\ell_1, \ldots, \ell_n$). Denote by $(\ovl{\mc{M}}, \phi)$ the map obtained in this way from a given $(\mc{M}, \phi)$. Then $\ovl{\mc{M}}$ is a map with $n$ boundary components (whose duals consist of neighboring edge-faces), and the boundary components are mapped by the embedding $\phi$ to $\ell_1, \ldots, \ell_n$. This is the sense in which $(\ovl{\mc{M}}, \phi)$ has ``boundary" given by $s$.
% Also, note that $\chi(\ovl{\mc{M}}) = \chi(\mc{M}) - n$. Thus if one wants to sum over maps $(\ovl{\mc{M}}, \phi)$ with boundary $s$, then the term $\chi(\mc{M}) - 2n$ which appears in what follows should be replaced by $\chi(\ovl{\mc{M}}) - n$.
% \end{remark}

% \begin{remark}
% Observe that for $(\mc{M}, \phi) \in \epe(s, K_+, K_-)$, we have that 
% \begin{equs}\label{eq:number-edges-and-faces}
% E(\mc{M}) &= 2\sum_{e \in E_\Lambda} |\phi^{-1}(e)|,\\
% F(\mc{M}) &= n + |K_+| + |K_-| + \sum_{e \in \Lambda} \ell(\mu_e(\phi)).
% \end{equs}
% \end{remark}

\revision{The following theorem is the main result of Section \ref{section:wilson-loop-expectation-epe}. It is the precise version of Theorem \ref{thm:informal-wilson-loop-expectation-sum-over-epe}.}

\begin{theorem}\label{thm:wilson-loop-expectation-sum-over-epe}
\revision{Let $\Lambda$ be a finite lattice.} Let $s = (\ell_1, \ldots, \ell_n)$ be a string. We have that
\begin{equs}
\langle W_s \rangle_{\Lambda, \beta} = Z_{\Lambda, \beta}^{-1} \sum_{\ptl (\mc{M}, \phi) = s} \frac{\beta^{\area(\mc{M}, \phi)}}{(\phi^{-1})!} \bigg(\prod_{e \in E_\Lambda} \ovl{\Wg}_N(\mu_e(\phi))\bigg)   N^{\chi(\mc{M}) - n}.
\end{equs}
\revision{Here, the sum converges absolutely in the following sense. Recall that the sum is over the set $\epe(s)$, which may be partitioned as in \eqref{eq:epe-partition}. Then
\[ \sum_{K : \mc{P} \ra \N} \sum_{\substack{(\mc{M}, \phi) \in \epe(s, K) }} \Bigg|\frac{\beta^{\area(\mc{M}, \phi)}}{(\phi^{-1})!} \bigg(\prod_{e \in E_\Lambda} \ovl{\Wg}_N(\mu_e(\phi))\bigg)   N^{\chi(\mc{M}) - n}\Bigg| < \infty.  \]}
\end{theorem}
\begin{proof}
\revision{
For the absolute convergence, combining equations \eqref{eq:euler-relation-boundary-epe} and \eqref{eq:wilson-loop-exponential-taylor-expansion}, and using that $|\Tr(Q)| \leq N$ for unitary matrices $Q \in \unitary(N)$, we have that (recall that our Wilson loops are defined using the normalized trace)
\begin{equs}
\sum_{K : \mc{P} \ra \N} \bigg| \frac{(N\beta)^K}{K!} \int W_s(Q) \prod_{p \in \mc{P}} \Tr(Q_p)^{K(p)} \prod_{e \in E_\Lambda} dQ_e \bigg| &\leq \sum_{K : \mc{P} \ra \N} \frac{(N\beta)^K}{K!} \prod_{p \in \mc{P}} N^{K(p)} \\ 
&= \sum_{K : \mc{P} \ra \N} \frac{(N^2\beta)^K}{K!} = \prod_{p \in \mc{P}} e^{N^2 \beta} < \infty.
\end{equs}
Next, by Proposition \ref{prop:word-expectation-as-sum-over-bipartite-maps}, we have that for each fixed $K : \mc{P} \ra \N$, 
% Combining equations \eqref{eq:euler-relation-boundary-epe} and \eqref{eq:wilson-loop-exponential-taylor-expansion} with Proposition \ref{prop:word-expectation-as-sum-over-bipartite-maps}, we have that (recall that our Wilson loops are defined using the normalized trace)
\begin{equs}
\frac{(N\beta)^K}{K!} \int W_s(Q) \prod_{p \in \mc{P}} \Tr(Q_p)^{K(p)} =  \frac{(N\beta)^{K}}{K!} \sum_{(\mc{M}, \phi) \in \epe(s, K)} \bigg(\prod_{e \in E_\Lambda} \ovl{\Wg}_N(\mu_e(\phi))\bigg) N^{\chi(\mc{M}) - n} N^{-K}.
\end{equs}
% \begin{equs}
% \langle W_s \rangle_{\Lambda, \beta} = Z_{\Lambda, \beta}^{-1} \sum_{K : \mc{P} \ra \N} \frac{(N\beta)^{K}}{K!} \sum_{(\mc{M}, \phi) \in \epe(s, K)} \bigg(\prod_{e \in E_\Lambda} \ovl{\Wg}_N(\mu_e(\phi))\bigg) N^{\chi(\mc{M}) - n} N^{-K}.
% \end{equs}
Recalling that $\area(\mc{M}, \phi) = \sum_p K(p)$ and $(\phi^{-1})! = K!$ for $(\mc{M}, \phi) \in \epe(s, K)$, the desired result now follows upon summing in $K$.
% Using the identities \eqref{eq:number-edges-and-faces} and the definition of the normalized Weingarten function $\ovl{\Wg}_N$ (Definition \ref{def:normalized-weingarten}), we may obtain 
% \begin{equs}
% N^{-n} \langle W_s \rangle = 
% Z_{\Lambda, \beta}^{-1} \sum_{(\mc{M}, \phi) \in \epe(s)} \frac{\beta^{\area(\mc{M}, \phi)}}{(\phi^{-1})!} \bigg(\prod_{e \in \Lambda} \ovl{\Wg}_N(\mu_e(\phi)) \bigg) N^{\chi(\mc{M}) - 2n}. \qedhere
% \end{equs}
}
\end{proof}

\revision{In the first few versions of this article, we included a heuristic discussion of the large-$N$ limit of the lattice Yang--Mills theories in terms of surface sums. This work has since been carried out in \cite{BCSK2024}, which established that the large-$N$ limit of Wilson string expectations can be expressed in terms of a certain sum over planar maps, which is very similar to the surface sums which were suggested in previous versions of this work. There was however one surprising aspect of the surface sums, which was that an additional ``non-separability" condition had to be imposed to obtain the result. See \cite[Appendix A]{BCSK2024} for more discussion on this topic.}

\section{Brownian motion and Poisson process exploration}\label{section:poisson-exploration}

In this section, we prove Theorem \ref{thm:weingarten-recovery}. First, in Section \ref{section:poisson-exploration}, we define and analyze a particular exploration process that is central to our proof. We then give the proof of the theorem in the case where $N$ is large, where it is easier to focus on the main ideas. In Section \ref{section:poisson-exploration-general-N-proof}, we extend the argument to the case of general $N$.
% In Section \ref{section:rep-theory}, we cover the representation theory facts that are needed in the course of the proof.

\subsection{Strand-by-strand exploration}\label{section:strand-by-strand}

We begin towards the proof of Theorem \ref{thm:weingarten-recovery}. The main difficulty is that the weights $w_T(\pi)$ appearing in Lemma \ref{lemma:expectations-words-unitary-bm}, when expressed as a series in $T$, do not converge absolutely when $T \toinf$. In fact, the series is of the schematic form $\sum_k \frac{(-T)^k}{k!} c_k$, for some coefficients $c_k$. Clearly, to show convergence as $T \toinf$, we need to take advantage of delicate cancellations which occur, rather than any sort of absolute summability. Uncovering these cancellations is the main technical part of the argument. This will be achieved via a certain exploration of the Poisson point process that we introduced in Section \ref{section:poisson-process-intro} which will provide an alternate form for the weights $w_T(\pi)$ which makes taking the $T \toinf$ limit easy.

% \begin{notation}
% We will often refer to the opposite-direction swaps introduced in Section \ref{section:poisson-process-intro} as ``turnarounds". 
% \end{notation}

In the following, recall the Poisson process and strand diagram material introduced in Section \ref{section:poisson-process-intro}. Because the Poisson processes corresponding to different letters are independent, it will suffice to just analyze the portion of the strand diagram corresponding to a single letter $\lambda$.

\begin{notation}
Recall the notation from Section \ref{section:poisson-process-intro} that $n_+(\lambda)$ (resp. $n_-(\lambda)$) is the number of times $\lambda$ (resp. $\lambda^{-1}$) appears in the collection of words $\bm \Gamma$. In this section, we will almost always have that $n_+(\lambda) = n_-(\lambda)$, and so we simplify notation and just write $n(\lambda) = n_+(\lambda) = n_-(\lambda)$. Moreover, with $\lambda$ fixed, we will just write $n$ instead of $n(\lambda)$. We will also write $\mc{D}_T = \mc{D}_T(\lambda), \Sigma = \Sigma(\lambda), \Sigma(T) = \Sigma_T(\lambda)$. In this section, we change to writing $\Sigma(T)$ instead of making $T$ a subscript because this will visually simplify many of the formulas that we write. (On the other hand, we found it more natural to write $\Sigma_T(\lambda)$ in Section \ref{section:poisson-process-intro}.)
\end{notation}
% \begin{notation}
% In this section, it will be notationally convenient for us to assume that the strand diagram corresponding to $\lambda$ has $2n$ total strands, with $n$ right-directed and left-directed strands each. This corresponds to the case that $\lambda$ and $\lambda^{-1}$ each appear a total of $n$ times in the given collection of words $\bm \Gamma$. This is in contrast to the notation of Definition \ref{defn-strand-diagram}, where $n_\ell$ is the total number of occurrences of a given letter $\lambda_\ell$ and its inverse $\lambda_\ell^{-1}$. To make consistent with this previous notation, we could perhaps introduce $n_+ = n_- = n/2$ and work with $n_+$. However, the parameter $n$ often appears in subscripts or superscripts, and adding a subscript ``$+$" to $n$ will result in iterated subscripts, which will complicate many expressions. Therefore, we decide just to use $n$ to denote the number of right-directed (and left-directed) strands.
% \end{notation}

Recall that by Lemma \ref{lemma:expectations-words-unitary-bm}, expectations of words of Unitary Brownian motion may be expressed in terms of $w_T(\pi)$ (which is defined in \eqref{eq:w-T-def}), for matchings $\pi \in \mc{M}(2n)$. 
% Now that we have introduced the Brauer algebra $\mc{B}_{n}$ in Definition \ref{def:brauer-algebra}, one may in fact view the matchings $\pi$ as elements of $\mc{B}_{2n}$ (even more, as elements of the walled Brauer algebra $\mc{B}_{n, n}$, as we will see). 
We next describe an alternative way to think about $w_T(\pi)$ in terms of the Brauer algebra (introduced in Section \ref{section:rep-theory}) which will be useful.

% Even more, we may restrict to the walled Brauer algebra $\mc{B}_{n,n} \sse \mc{B}_{2n}$, for reasons we next describe.

Let $P \sse \mc{D}_T$ be a finite point set. Recalling that $\mc{D}_T$ is defined as a finite disjoint union of $\binom{2n}{2}$ copies of the interval $[0, T]$, we may view $P$ as a finite subset $\{x_1, \ldots, x_{|P|}\}$ of $[0, T]$, where $0 \leq x_1 \leq \cdots \leq x_{|P|} \leq T$, along with a collection of labels $\mf{l}_1, \ldots, \mf{l}_{|P|}$ where each $\mf{l}_k = \{i_k, j_k\}$, where $i_k < j_k \in [2n]$. The label $\mbf{l}_k$ can be thought of as keeping track of which copy of $[0, T]$ the point $x_k$ belongs to. In the following, we assume that the $x_1, \ldots, x_{|P|}$ are distinct. Note that this is a.s. true for the random point set $\Sigma_T$.

In the following, recall the walled Brauer algebra elements $(i ~ j)$ and $\langle i ~ j \rangle$ introduced in Definition \ref{def:brauer-algebra-generators}.
% We first describe an alternative way to write the Brauer algebra element
% \begin{equs}
% (-1)^{\numswaps(P)} N^{-|P|} N^{\numintloops(P)} \pi(P) \in \mc{B}_{2n}
% \end{equs}
% which essentially appears in the definition of $w_T(\pi)$.

\begin{definition}\label{def:points-to-brauer-fn}
Given a finite collection of points $P \sse \mc{D}_T$, we define $\pointstobrauerfn_N(P) \in \mc{B}_{n, n}$ as follows. We use the representation of $P$ in terms of $0 \leq x_1 < \cdots < x_{|P|} \leq T$ and $\mf{l}_1, \ldots, \mf{l}_{|P|}$. For each $k \in [|P|]$, define $b_k \in \mc{B}_{n, n}$ as follows. If $\mf{l}_k = \{i_k, j_k\}$ is such that $i_k, j_k \in [n]$ or $i_k, j_k \in (n:2n]$, then define $b_k := (-1/N)(i_k ~ j_k)$. Otherwise, if $i_k \in [n]$, $j_k \in (n:2n]$, define $b_k := (1/N)\langle i_k ~ j_k \rangle$. We then define
\begin{equs}
\pointstobrauerfn_N(P) := b_1 \cdots b_{|P|} \in \mc{B}_{n, n}.
\end{equs}
We often omit the dependence of $\pointstobrauerfn_N$ on $N$ and just write $\pointstobrauerfn(P)$.
\end{definition}

The key observation relating $w_T(\pi)$ to the Brauer algebra is the following.

\begin{lemma}\label{lemma:brauer-algebra-F-relation}
For any finite set of points $P \sse \mc{D}_T$, we have that 
\begin{equs}
\pointstobrauerfn(P) = (-1)^{\numswaps(P)} N^{-|P|} N^{\numintloops(P)} \pi(P).
\end{equs}
\end{lemma}
\begin{proof}
In the definition of $\pointstobrauerfn(P)$, each point which corresponds to a swap carries a $-1/N$ factor, while each point which corresponds to a turnaround carries a $1/N$ factor. This explains the factor $(-1)^{\numswaps(P)} N^{-|P|}$. The term $N^{\numintloops(P)} \pi(P)$ comes about because this is precisely how multiplication in the Brauer algebra is defined (recall in the discussion after Definition \ref{def:brauer-algebra} we set $\zeta = N$ so that each loop formed during multiplication contributes a factor $N$).
\end{proof}

With this definition and lemma, we may express
% \begin{equs}
% w_T(\pi) = e^{\binom{2n}{2}T - n T} \E[\pointstobrauerfn(\Sigma(T)) \ind(\pi(\Sigma(T)) = \pi)],
% \end{equs}
% Or, as elements of the walled Brauer algebra $\mc{B}_{n, n}$, we have the equality
\begin{equs}\label{eq:w-T-pi-brauer}
\sum_{\pi \in \mc{M}(n, n)} w_T(\pi) \pi = e^{\binom{2n}{2}T - n T} \E[\pointstobrauerfn(\Sigma(T))].
\end{equs}
In the following, we will mainly focus on understanding the large-$T$ limit of the right\revision{-}hand side above. 
% \[ \lim_{T \toinf} e^{\binom{2n}{2}T - n T} \E[\pointstobrauerfn(\Sigma(T))].\]
Clearly, once we understand this, we will also know 
\[ \lim_{T \toinf} w_T(\pi) \text{ for any $\pi$ }. \]

\begin{remark}
Observe that $\pointstobrauerfn(\Sigma(T))$ can be thought of as a ``random walk" on the Brauer algebra, since due to the Poisson distribution of $\Sigma(T)$, we are multiplying a number of ``increments" of the form $(i ~ j)$ or $\langle i ~ j \rangle$, and each increment is equally likely. Of course, here we ignored the $-1/N$ and $1/N$ factors, which makes the random walk interpretation not exactly correct.
\end{remark}
Our starting point is the next lemma, which observes some cancellation that happens once a turnaround appears. In words, it says that once a blue turnaround appears, the only points which can thereafter appear that touch either of the matched strands must be the turnaround between the same two strands. Ultimately, this cancellation is possible because swaps come with a negative sign and turnarounds come with a positive sign.

\begin{lemma}\label{lemma:cancellation}
Let $i_0 \in [n]$, $j_0 \in (n:2n]$. Let $\Sigma_{-\{i_0, j_0\}}(T)$ be the Poisson point process obtained by deleting all points of $\Sigma(T)$ touching the $i_0$th or $j_0$th strands. Then
\begin{equs}
\langle i_0 ~ j_0 \rangle \E[\pointstobrauerfn(\Sigma(T))] &= e^{-4(n-1) T} \langle i_0 ~ j_0 \rangle \E[\pointstobrauerfn(\Sigma_{-\{i_0, j_0\}}(T))], \\
\E[\pointstobrauerfn(\Sigma(T))] \langle i_0 ~ j_0 \rangle &= e^{-4(n-1) T} \E[\pointstobrauerfn(\Sigma_{-\{i_0, j_0\}}(T))] \langle i_0 ~ j_0 \rangle .
\end{equs}
\begin{remark}
\revision{We remark that if we applied the representation $\rho_+$ to the identities, then these identities would follow directly from the relation \eqref{eq:unitary-rho-plus-identity} between $\rho_+$ and the action of unitary matrices, combined with \eqref{eq:w-T-pi-brauer} which relates $\E [F(\Sigma(T))]$ to expectations of Unitary Brownian motion. However, because $\rho_+$ in general is not injective (see e.g. the discussion in Remark \ref{remark:rho-motivation}), we still must prove this lemma via a direct combinatorial argument.}
\end{remark}
\end{lemma}
\begin{proof}
We only show the first identity as the second follows similarly. Let $A_T$ be the event that the process $\Sigma(T)$ contains no points touching the $i_0$th or $j_0$th strand, besides those which give the turnaround $\langle i_0 ~ j_0 \rangle$. Since each strand is involved in $2n-1$ total Poisson processes, the number of Poisson processes that involve the $i$th or $j$th strand is $2(2n-1) - 1 = 4n-3$. On the event $A_T$, all but one of these processes must have zero points, and thus $\p(A_T) = e^{-4(n-1)T}$. Let $\Sigma_{\langle i_0 ~ j_0 \rangle}(T)$ be the process obtained by keeping only those points which give the turnaround $\langle i_0 ~ j_0 \rangle$. On $A_T$, we may split
\begin{equs}
\Sigma(T) = \Sigma_{-\{i_0, j_0\}}(T) \sqcup \Sigma_{\{i_0, j_0\}}(T),
\end{equs}
and moreover 
\begin{equs}
\pointstobrauerfn(\Sigma(T)) = \pointstobrauerfn(\Sigma_{\{i_0, j_0\}}(T)) \pointstobrauerfn(\Sigma_{-\{i_0, j_0\}}(T)).
\end{equs}
We thus have that
\begin{equs}
\langle i_0 ~ j_0 \rangle \E[\pointstobrauerfn(\Sigma(T)) \ind_{A_T}] &= \p(A_T) \langle i_0 ~ j_0 \rangle \E[\pointstobrauerfn(\Sigma_{\{i_0, j_0\}}(T))] \E[\pointstobrauerfn(\Sigma_{-\{i_0, j_0\}}(T))] \\
&= e^{-4(n-1) T} \langle i_0 ~ j_0 \rangle \E[\pointstobrauerfn(\Sigma_{\{i_0, j_0\}}(T))] \E[\pointstobrauerfn(\Sigma_{-\{i_0, j_0\}}(T))] .
\end{equs}
By explicit calculation, using that $\langle i_0 ~ j_0 \rangle^k = N^{k-1} \langle i_0 ~ j_0 \rangle$, we have that
\begin{equs}
\langle i_0 ~ j_0 \rangle \E[\pointstobrauerfn(\Sigma_{\{i_0, j_0\}}(T))] =  e^{-T} \langle i_0 ~ j_0 \rangle  \sum_{k=0}^\infty \frac{T^k}{k!} \bigg(\frac{\langle i_0 ~ j_0 \rangle}{N}\bigg)^k = \langle i_0 ~ j_0 \rangle e^{-T} \sum_{k=0}^\infty \frac{T^k}{k!} = \langle i_0 ~ j_0 \rangle.
\end{equs}
To finish, it suffices to show that
\begin{equs}
\langle i_0 ~ j_0 \rangle \E[\pointstobrauerfn(\Sigma(T))] = \langle i_0 ~ j_0 \rangle \E[\pointstobrauerfn(\Sigma(T)) \ind_{A_T}],
\end{equs}
or in other words,
\begin{equs}
\langle i_0 ~ j_0 \rangle \E[\pointstobrauerfn(\Sigma(T)) \ind_{A_T^c}] = 0.
\end{equs}
We show that for each $k \geq 1$, we have that
\begin{equs}
\langle i_0 ~ j_0 \rangle \E[\pointstobrauerfn(\Sigma(T)) \ind_{A_T^c} ~|~ |\Sigma(T)| = k] = 0.
\end{equs}
(If $k = 0$ then $A_T^c$ cannot occur.) Let $\Omega_k$ be the set of length-$k$ sequences of elements of the set
\begin{equs}
\{(i ~ j) : 1 \leq i < j \leq n\} \cup \{(i ~ j) : n +1 \leq i < j \leq 2n \} \cup \{\langle i ~ j \rangle : i \in [n], j \in (n : 2n]\},
\end{equs}
such that there exists some element not equal to $\langle i_0 ~ j_0 \rangle$ that involves either $i_0$ or $j_0$. For each $(x_1, \ldots, x_k) \in \Omega_k$, let $n_T(x_1, \ldots, x_k)$ be the number of swaps (i.e. elements of the form $(i ~ j)$) in the sequence. Observe that
\begin{equs}
\langle i_0 ~ j_0 \rangle \E[F(\Sigma(T)) \ind_{A_T^c} ~|~ |\Sigma(T)| = k] = \langle i_0 ~ j_0 \rangle \frac{1}{\binom{2n}{2}^k N^k} \sum_{(x_1, \ldots, x_k) \in \Omega_k} (-1)^{n_T(x_1, \ldots, x_k)} x_1 \cdots x_k.
\end{equs}
We now define a bijection $h : \Omega_k \ra \Omega_k$ such that if $h(x_1, \ldots, x_k) = (y_1, \ldots, y_k)$, then 
\begin{equs}
\langle i_0 ~ j_0\rangle (-1)^{n_T(y_1, \ldots, y_k)} y_1 \cdots y_k = - \langle i_0 ~ j_0 \rangle (-1)^{n_T(x_1, \ldots, x_k)} x_1 \cdots x_k.
\end{equs}
Note that this immediately implies that
\begin{equs}
\langle i_0 ~ j_0 \rangle \sum_{(x_1, \ldots, x_k) \in \Omega_k} (-1)^{n_T(x_1, \ldots, x_k)} x_1 \cdots x_k = - \langle i_0 ~ j_0 \rangle \sum_{(x_1, \ldots, x_k) \in \Omega_k} (-1)^{n_T(x_1, \ldots, x_k)} x_1 \cdots x_k, 
\end{equs}
which implies that the above is zero, which would give the desired result. To define $h$, given a sequence $(x_1, \ldots, x_k)$, let $1 \leq r \leq k$ be index of the first element $x_r$ which causes the sequence $(x_1, \ldots, x_n)$ to be in $\Omega_k$. Then either $x_r$ is a swap of the form $(i_0 ~ k)$ or $(k ~ j_0)$, or $x_r$ is a turnaround of the form $\langle i_0 ~ k \rangle$ or $\langle k ~ j_0 \rangle$. If $x_r$ is a swap, we set 
\begin{equs}
h(x_1, \ldots, x_k) := \begin{cases} 
(x_1, \ldots, x_{r-1}, \langle k ~ j_0 \rangle, x_{r+1}, \ldots, x_k) & x_r = (i_0 ~ k) \\
(x_1, \ldots, x_{r-1}, \langle i_0 ~ k \rangle, x_{r+1}, \ldots, x_k) & x_r = (k ~ j_0). 
\end{cases}
\end{equs}
and if $x_r$ is a turnaround, we set
\begin{equs}
h(x_1, \ldots, x_k) := \begin{cases}
(x_1, \ldots, x_{r-1}, (k ~ j_0), x_{r+1}, \ldots, x_k) & x_r = \langle i_0 ~ k \rangle \\
(x_1, \ldots, x_{r-1}, (i_0 ~ k), x_{r+1}, \ldots, x_k) & x_r = \langle k ~ j_0 \rangle.
\end{cases}
\end{equs}
In words, if $x_r$ is a swap involving $i_0$ (resp. $j_0$), then $h$ switches $x_r$ to a turnaround involving $j_0$ (resp. $i_0$). Similarly, if $x_r$ is a turnaround involving $i_0$ (resp. $j_0$), then $h$ switches $x_r$ to a swap involving $j_0$ (resp. $i_0$). Note that $h$ is an involution, and thus a bijection. Also, we clearly have by construction that
\begin{equs}
(-1)^{n_T(x_1, \ldots, x_k)} = - (-1)^{n_T(h(x_1, \ldots, x_k))}.
\end{equs}
Thus to finish, it suffices to show that with $h(x_1, \ldots, x_k) = (y_1, \ldots, y_k)$, we have that $\langle i_0 ~ j_0 \rangle x_1 \cdots x_k = \langle i_0 ~ j_0 \rangle y_1 \cdots y_k$. By construction of $h$, it just suffices to show that $\langle i_0 ~ j_0 \rangle x_1 \cdots x_r = \langle i_0 ~ j_0 \rangle x_1 \cdots x_{r-1} y_r$. By the assumption on $r$, we have that $x_1 \cdots x_{r-1}$ commutes with $\langle i_0 ~ j_0\rangle$, and so
\begin{equs}
\langle i_0 ~ j_0\rangle x_1 \cdots x_r = x_1 \cdots x_{r-1} \langle i_0 ~ j_0\rangle x_r, ~~ \langle i_0 ~ j_0 \rangle x_1 \cdots x_{r-1}y_r = x_1 \cdots x_{r-1} \langle i_0 ~ j_0\rangle y_r.
\end{equs}
To finish, we claim that $\langle i_0 ~ j_0 \rangle x_r = \langle i_0 j_0 \rangle y_r$, i.e. the switching procedure used to define $h$ does not change the overall matching. This follows by the identities $\langle i ~ k\rangle (i ~ j) = \langle i ~ k \rangle \langle j ~ k \rangle$ and $\langle i ~ k \rangle \langle i ~ j \rangle = \langle i ~ k \rangle (j ~ k)$. For the first identity, observe that the two products of matchings in Figure \ref{figure:brauer-identity} are equal.
\begin{figure}[h]%
    \centering
    \subfloat[\centering $\langle 1 ~ 5 \rangle (5 ~ 4)$]
    {{\includegraphics[page=1, width=5cm]{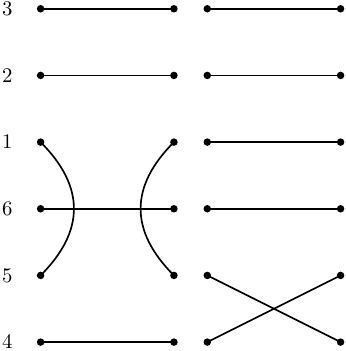} }}%
    \qquad
    \subfloat[\centering $\langle 1 ~ 5\rangle \langle 1 ~ 4 \rangle$]
    {{\includegraphics[page = 2, width=5cm]{figures/brauer-cancellation-identity.pdf} }}%
    \caption{The above two products of matchings are equal.}\label{figure:brauer-identity}
\end{figure}

\noindent The second identity follows similarly. 
\end{proof}

We now finally describe our exploration of the strand diagram corresponding to $\Sigma$. We describe the exploration on the ``infinite-time" process $\Sigma$; the exploration on the ``finite-time" process $\Sigma(T)$ is exactly the same (except since $\Sigma(T)$ contains no points beyond time $T$, nothing happens in the exploration after this time).

\begin{figure}[h]%
    \centering
    {{\includegraphics[width=7cm, page = 1]{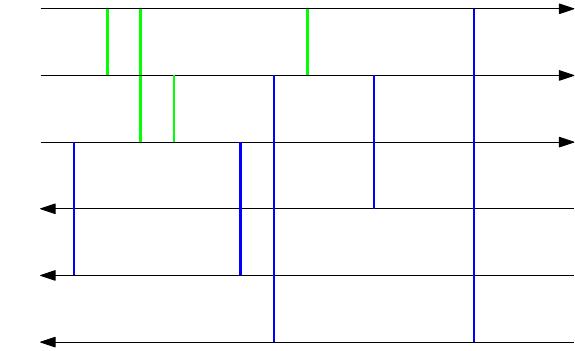} }}%
    \caption{Before the strand-by-strand exploration.}\label{figure:strand-by-strand-example}
\end{figure}

The exploration proceeds strand-by-strand. We first give an informal description with accompanying figures before proceeding to the formal mathematical definition. The main feature of the exploration is that we explore only a single strand at a time, rather than all strands at once. That is, we start at (say) the top strand, and explore left-to-right until we see a swap or a turnaround involving this strand. If we see a swap between the top strand and another strand, then we begin exploring the other strand. If we see a turnaround, then the current exploration era ends, and we begin to explore the next strand. To visualize this exploration, suppose we want to explore the the diagram in Figure \ref{figure:strand-by-strand-example}.

\begin{figure}[h]%
    \centering
    \subfloat[\centering First era (red)]
    {{\includegraphics[width=7cm, page=2]{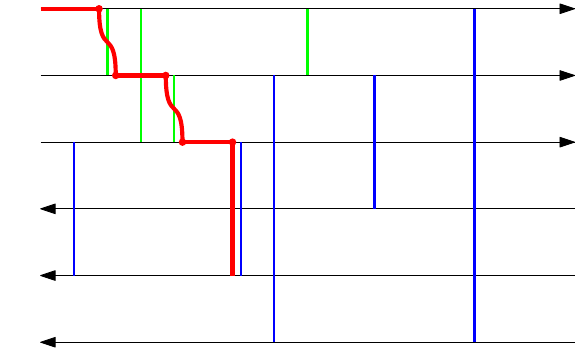} }}%
    \qquad
    \subfloat[\centering Second era (orange)]
    {{\includegraphics[width=7cm, page=3]{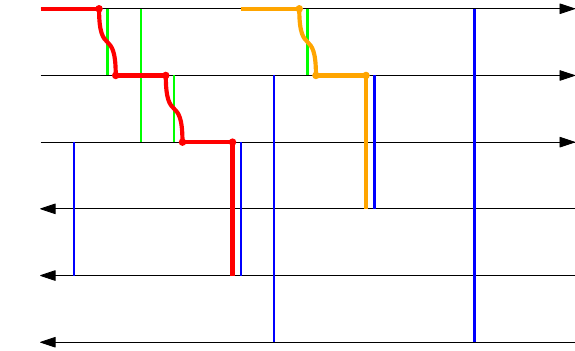} }}%
    \qquad
    \subfloat[\centering Third era (brown)]
    {{\includegraphics[width=7cm, page=4]{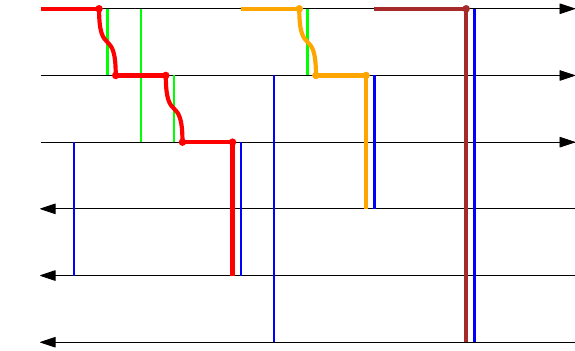} }}%
    \caption{The successive eras of our strand-by-strand exploration. Note that the exploration does not see all of the swaps. One of the key ideas is to show that if we condition on what is seen by the exploration, the contribution from these ``unseen swaps" is in fact zero. This is ultimately due to cancellations that can be uncovered.}\label{figure:strand-by-strand-first-123}
\end{figure}

Our exploration proceeds in three separate eras, drawn as in Figure \ref{figure:strand-by-strand-first-123}.

Note that at the start of the second era, we begin exploring the top strand instead of the second-to-top strand, because of the previous swap between these two strands. Likewise, at the start of the third era, we also begin exploring from the top strand, because this is effectively the bottom strand due to the previously seen swaps. Another thing to note is that in principle, during the first exploration era, it is certainly possible for the point process to have swaps that involve two non-top strands. However, our exploration process does not see these swaps. It turns out that by exploring the random environment in the manner we described, we can in fact assume that in every exploration era, every swap in the point process involves the current exploration strand, so that we don't need to worry about such ``unseen swaps". This property is due to certain cancellations that we may take advantage of, which are very similar in spirit to the cancellations observed in the proof of Lemma \ref{lemma:cancellation}.

\begin{figure}[h]%
    \centering
    {{\includegraphics[width=10cm, page=5]{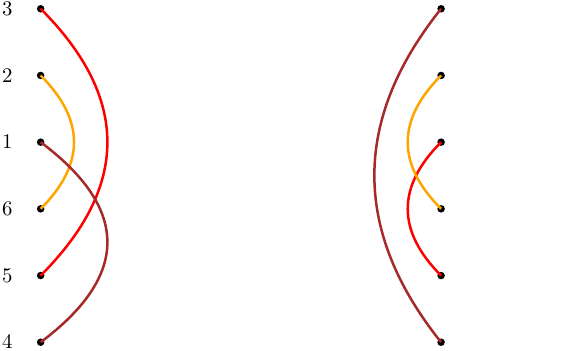} }}%
    \caption{The matching discovered by our strand-by-strand exploration.}\label{figure:strand-by-strand-matching}
\end{figure}

At the end of the last exploration era, we have built up an element of $\mc{B}_{n, n}$ (we have omitted the additional factors of $\pm\frac{1}{N}$ and only drawn the left and right bijections), as displayed in Figure \ref{figure:strand-by-strand-matching}.

Here, the colors are for visual purposes and don't affect the end element of $\mc{B}_{n, n}$: we have colored the matching edges to denote the exploration era in which these pairings were discovered. Now here is why we chose to explore as we did: conditioned on everything we have seen up to the end of the last exploration era, the expectation of $\pointstobrauerfn(\Sigma(T))$ {\it is essentially given\footnote{Technically, this is only true up to some explicit factors $\pm 1/N$ factors, but this is more of a technical detail.} by the matching in Figure \ref{figure:strand-by-strand-matching}}. This property is intimately related to our previous comment that we can assume that there are no unseen swaps, i.e. swaps which do not involve the current exploration strand.
% As we will see, the reason why this is true is due to certain cancellations which we can utilize due to how we explored the strands (these cancellations are very similar to the ones we exploited in the proof of Lemma \ref{lemma:cancellation}). 
This key property of our exploration enables us to give a rather explicit closed-form expression for the overall expectation of $\pointstobrauerfn(\Sigma(T))$. Even more, it is almost trivial to take the $T \toinf$ limit of the closed-form expression, and this allows us to recover the Weingarten calculus.

We now proceed to the precise definition of the exploration. First, for $i \in [2n]$, define
\begin{equs}
\Sigma_i := \bigcup_{j \in [2n] \backslash \{i\}} \Sigma_{\{i, j\}}, 
\end{equs}
i.e. $\Sigma_i$ collects all Poisson processes with which $i$ is involved. In terms of the strands, $\Sigma_i$ collects all swaps and turnarounds touching the $i$th strand. The exploration is described by two processes $(E_t)_{t \geq 0}$, $(\pi_t)_{t \geq 0}$, the first of which takes values in $[n]$, and the second of which takes values in $\symgrp_{n}$ (which we view as the set of bijections of $[n]$). One should think of $E_t$ as tracking the current exploration era, and $\pi_t$ as tracking the current strand of exploration.

We start with $E_0 := 1$, $\pi_0 := \mrm{id}$. We begin exploring $\Sigma_{\pi_0(E_0)} = \Sigma_1$ until we see the first point, which we denote by $U_1$. At time $U_1$, we update $E$ and $\pi$ as follows. Let $j \in [2n]$ be the label of the strand which is matched to the top strand by the first point $U_1$. For $t \in (0, U_1)$, we set $E_t := E_0$, $\pi_t := \pi_0$. Now if $j \in [n]$, then we set $E_{U_1} := E_0$ and $\pi_{U_1} := (\pi_0(1) ~ j) \pi_0$. We then continue exploring $\Sigma_{\pi_{U_1}(E_{U_1})}$ from time $U_1$. Otherwise, if $j \in (n : 2n]$, then we set $E_{U_1} := E_0 + 1$ (i.e. a new exploration era begins) and $\pi_{U_1} := \pi_0$. Additionally, we remove all points of $\Sigma_{\pi_0(E_0)} \cup \Sigma_j$ from $\Sigma$. The exploration then continues on this reduced point process. In terms of the strands, the removal of points corresponds to only looking at those swaps or turnarounds which do not involve the $\pi_0(E_0)$ or $j$th strands. The exploration stops once all exploration eras have ended, i.e. once we have explored all strands up to their first time of turnaround. This is the first time $t$ such that $E_t = n +1$. This is a.s. finite for $\Sigma$ due to properties of Poisson point processes.

For $i \in [n]$, let $T_i := \inf\{t \geq 0 : E_t = i+1\}$, i.e. the time at which the $i$th exploration era ends. Let $\mc{Q}_t$ be the set of points that the exploration has seen up to time $t$. Let $(\mc{F}_t, t \geq 0)$ be the filtration generated by the processes $E, \pi$.

% The first order of business is to show that a reduced exploration process suffices, due to certain cancellations which occur. 

% As we will shortly see, this identity implies that in exploring successive eras, we can treat the strands which have been previously explored as ``out of the game". Towards this end, let us now define the reduced exploration process, which again is described by two processes $(E_t)_{t \geq 0}$, $(\pi_t)_{t \geq 0}$. As before, start with $E_0 = 1$, $\pi_0 = \mrm{id}$. Up to the end of the first exploration era, the processes $E$ and $\pi$ are defined exactly as before. What changes is that once we conclude the first exploration era, we remove the two strands which were matched. In terms of the Poisson process, let $\{E_{\bar{T}_1}, \sigma(E_{\bar{T}_1})\}$ be the label of the point which marks the end of the first exploration era (so $\sigma(E_{\bar{T}_1}) \in \{n+1, \ldots, 2n\}$). We remove from each of the $\mc{P}_i$ all points in $\mc{P}_{E_{\bar{T}_1}}$ and $\mc{P}_{\sigma(E_{\bar{T}_1})}$, and then continue the exploration with the remaining points.  More generally, at the end of the $k$th exploration era, we further remove all points corresponding to the two strands which were matched in this era. One of these strands is labelled by $E_{T_k}$; we suppose the other is labelled by $\sigma(E_{T_k}) \in \{n+1, \ldots, 2n\} \backslash \{\sigma(E_{T_1}), \ldots, \sigma(E_{T_n})\}$.

The following key proposition makes precise the key property of our exploration that we described earlier.

\begin{prop}\label{prop:strand-by-strand-exploration}
We have that
\begin{equs}
e^{\binom{2n}{2} T - nT} \E[\pointstobrauerfn(\Sigma(T))& \ind(T_n \leq T)] = \\
&\E \big[\pointstobrauerfn(\mc{Q}_{T_{n}}) \ind(T_n \leq T) e^{2(n-1) T_1} e^{2(n-2)(T_2 - T_1)} \cdots e^{2(n-n) (T_n - T_{n-1})}\big].
\end{equs}
\end{prop}
\begin{proof}
Fix $N$. We proceed by induction on $n$. When $n = 1$, the result is true for all $T \geq 0$, because then $A_T$ always occurs, and furthermore when $T_1 \leq T$, we have that $\pointstobrauerfn(\Sigma(T)) = \pointstobrauerfn(\mc{Q}_{T_1})$. Now suppose that for some general $n \geq 1$, the result is true for all $T \geq 0$. We proceed to show that the case $n+1$ also holds. We start by conditioning on $\mc{F}_{T_1}$. Pictorially, this corresponds to exploring until the end of the first era, see the left of Figure \ref{figure:exploration-up-to-T1}.
\begin{figure}[ht!]
\centering
\includegraphics[width=7cm]{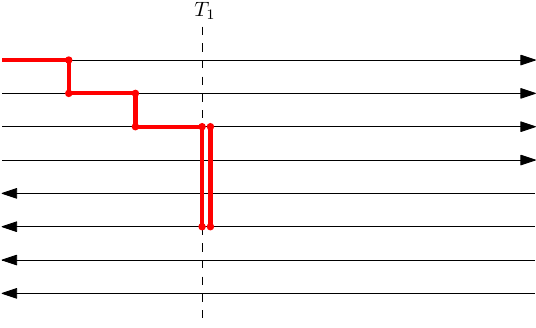}
\qquad
\includegraphics[width=7cm]{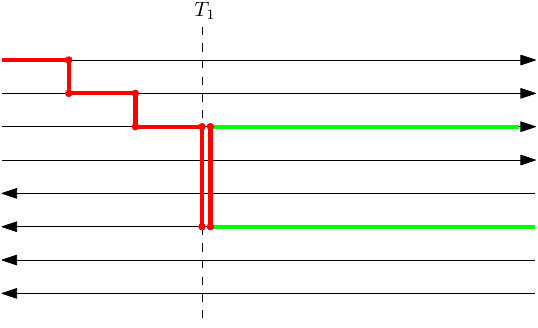}
\caption{Left: We start at the first strand and explore until we see a turnaround. Right: Once we have seen a turnaround, we may treat the two strands involved in the turnaround as ``out of the game".}\label{figure:exploration-up-to-T1}
\end{figure}
One should think of the two parallel vertical red lines as occurring at the same time (namely $T_1$), although for visual purposes we have drawn them to be slightly separated. Next, naturally, we may split the diagram in Figure \ref{figure:exploration-up-to-T1} into two parts: the part to the left of $T_1$, and the part to the right of $T_1$. This corresponds to splitting
\begin{equs}
\Sigma(T) = \Sigma(T_1) \cup (\Sigma(T) \backslash \Sigma(T_1)).
\end{equs}
Since the Poisson processes before $T_1$ and after $T_1$ are conditionally independent, we have that
\begin{equs}
\E[\pointstobrauerfn(\Sigma(T)) & \ind(T_{n+1} \leq T) ~|~ \mc{F}_{T_1}] = \\
&\E[\pointstobrauerfn(\Sigma(T_1)) ~|~ \mc{F}_{T_1}] \E\big[\pointstobrauerfn(\Sigma(T)  \backslash \Sigma(T_1)) \ind(T_{n+1} - T_1 \leq T - T_1)~|~ \mc{F}_{T_1}].
\end{equs}
We first use our inductive assumption to rewrite the second conditional expectation on the right\revision{-}hand side above. By our cancellation lemma (Lemma \ref{lemma:cancellation}), we may assume that there are no swaps or turnarounds which involve either of the two matched strands after $T_1$, as long as we multiply by the explicit exponential factor $e^{-(4(n+1) -1)(T - T_1)} = e^{-4n(T - T_1)}$.
Pictorially, after $T_1$, the two segments which are colored bright green in the right Figure \ref{figure:exploration-up-to-T1} are no longer connected to the other strands in the diagram.
% \begin{figure}
% \begin{center}
% \includegraphics[width=.5\textwidth]{figures/exploration-up-to-T1-with-constraints.pdf}
% \end{center}
% \caption{}\label{figure:exploration-up-to-T1-with-constraints}
% \end{figure}
% Since each individual strand is connected to $2(n+1) - 1 = 2n+1$ other strands, the total number of independent Poisson processes which are forced to be zero is $2(2n+1) - 2 = 4n$. Moreover, any additional turnarounds between the two green strands has no effect, since every additional turnaround creates a new component, which incurs a factor of $N$, which cancels out the factor of $\frac{1}{N}$ that is associated to each turnaround. Thus, we may ignore the two green strands. 
The point now is that after having taken out the two green strands, the expectation of the remainder of the diagram after $T_1$ is exactly given by our inductive assumption. Thus, we have the identity
\begin{equs}
~&e^{\binom{2n}{2} (T - T_1) - n (T - T_1)} \E\big[\pointstobrauerfn(\Sigma(T) \backslash \Sigma(T_1)) \ind(T_{n+1} - T_1 \leq T - T_1)~|~ \mc{F}_{T_1}] = \\
& e^{-4n(T - T_1)} \E\big[\pointstobrauerfn(\mc{Q}_T \backslash \mc{Q}_{T_1}) \ind(T_{n+1} - T_1 \leq T - T_1) e^{2(n-1)(T_2 - T_1)} e^{2(n-2)(T_3 - T_2)} \cdots e^{2(n-n)(T_{n+1} - T_n)} \big].
\end{equs}
Applying this identity, as well as the identity $\binom{2(n+1)}{2} - (n+1) = 4n + \binom{2n}{2} - n$, we obtain
\begin{equs}
~&e^{\binom{2(n+1)}{2} T - (n+1)T} \E[\pointstobrauerfn(\Sigma(T)) \ind(T_{n+1} \leq T)]  = \\
&\E\big[\pointstobrauerfn(\Sigma(T_1)) \pointstobrauerfn(\mc{Q}_{T_{n+1}} \backslash \mc{Q}_{T_1}) \ind(T_{n+1} \leq T) e^{(4n + \binom{2n}{2} - n ) T_1} e^{2(n+1 - 2) (T_2 - T_1)} \cdots e^{2(n+1 - (n+1)) (T_{n+1} - T_n)}\big]. 
\end{equs}
To finish, we now argue that 
\begin{equs}\label{eq:snapping-back-identity}
e^{(4n + \binom{2n}{2} - n) T_1} \E\big[\pointstobrauerfn (\Sigma(T_1)) F(\mc{Q}_{T_{n+1}} \backslash \mc{Q}_{T_1})  ~|~ \mc{F}_{T_{n+1}}] &= e^{2n T_1} F(\mc{Q}_{T_1}) F(\mc{Q}_{T_{n+1}}\backslash \mc{Q}_{T_1}) \\
&= e^{2nT_1} F(\mc{Q}_{T_{n+1}}).
\end{equs}
Note that this would complete the proof of the inductive step. For a picture of what we have in mind when conditioning on $\mc{F}_{T_{n+1}}$, see Figure \ref{figure:exploration-up-to-T1-conditional-on-after}.
\begin{figure}[h]
\centering
\includegraphics[width=7cm]{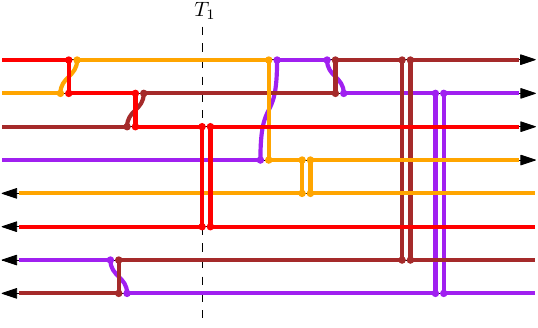}
\qquad
\includegraphics[width=7cm]{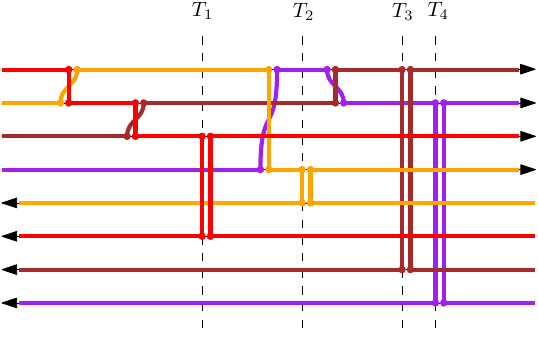}
\caption{Left: we can assume that our exploration process looks like this after applying the inductive assumption. Right: completing the inductive step by arguing that after cancellation, we may assume that there are no other points before $T_1$, besides the previously seen red swaps.}\label{figure:exploration-up-to-T1-conditional-on-after}
\end{figure}

In the left of Figure \ref{figure:exploration-up-to-T1-conditional-on-after}, we treat the portion of the diagram to the right of $T_1$ as fixed, whereas the portions of the strands before $T_1$ which are black have not been fully explored. The identity \eqref{eq:snapping-back-identity} says that after averaging over this randomness, we may simply assume that there are no additional swaps or turnarounds in $[0, T_1]$, so that the expectation is given by the right of Figure \ref{figure:exploration-up-to-T1-conditional-on-after} (which corresponds to the right\revision{-}hand side of the identity).

The identity \eqref{eq:snapping-back-identity} follows by cancellations similar to those exploited in the proof of the cancellation lemma (Lemma \ref{lemma:cancellation}). 
Indeed, observe that the two diagrams in Figure \ref{figure:exploration-up-to-T1-with-cancellation} equal, in the sense that the final matching is the same (the red and orange strands are unchanged, so one only needs to track the brown and purple strands).
\begin{figure}[h]%
    \centering
    % \subfloat[\centering Original]
    {{\includegraphics[width=7cm]{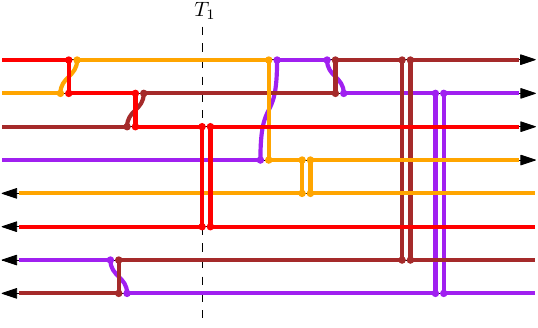} }}%
    \qquad
    % \subfloat[\centering]
    {{\includegraphics[width=7cm]{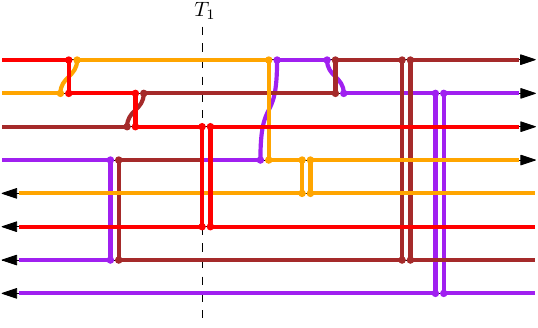} }}%
    \caption{The above two diagrams are equal as elements of $\mc{B}_{n, n}$}\label{figure:exploration-up-to-T1-with-cancellation}
\end{figure}
Note however that the left diagram will have an opposite sign compared to the right diagram, because swaps incur a factor of $-1$ while turnarounds do not. This gives the desired cancellation between swaps and turnaround which do not connect two strands which have been matched by the portion of the diagram after $T_1$. 
Thus the total number of Poisson processes which must have zero points is $\binom{n}{2} + \binom{n+1}{2} + n^2$. Here, $\binom{n}{2}$ counts the possible swaps between two top strands, $\binom{n+1}{2}$ counts the possible swaps between two bottom strands, and $n^2 = n(n+1) - n$ counts the turnarounds which connect a top and bottom strand which are not already connected by the diagram to the right of $T_1$. We now finish by noting the identity
\[
4n + \binom{2n}{2} - n - \binom{n}{2} - \binom{n+1}{2} - n^2 = 2n. \qedhere
\]
\end{proof}

\begin{figure}[h]%
    \centering
    \subfloat[\centering Original]{{\includegraphics[width=7cm, page=1]{figures/exploration-up-to-T1-end-result.pdf} }}%
    \qquad
    \subfloat[\centering Alternative]{{\includegraphics[width=7cm, page=2]{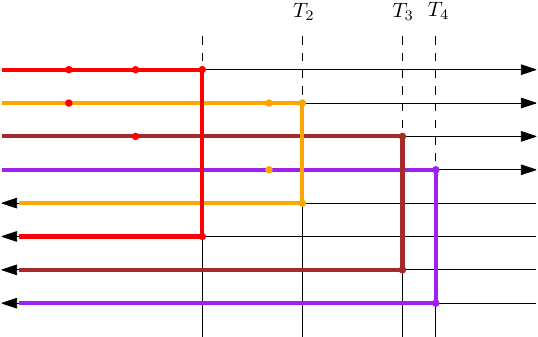} }}%
    % \caption{}%
    % \label{fig:comparison}%
    \caption{If, every time we see a swap, we imagine we ``cut and swap" the two strands which were involved, we go from the left picture the to the right picture. Note that the left bijection is unchanged. The original right bijection can be reconstructed from the left bijection and the swaps. In the right picture, all swaps in the first era involve the top strand, all swaps in the second era involve the second-top strand, etc.}\label{figure:alternative-exploration}
\end{figure}

Next, to extract the Jucys-Murphy elements, it is helpful to think of all the swaps in the $i$th exploration era as involving the $i$th strand. Towards this end, we show that the expectation of $\pointstobrauerfn(\mc{Q}(T_n))$ appearing in Proposition \ref{prop:strand-by-strand-exploration} may be computed by following a slightly different exploration, one in which each exploration era stays on a single strand, and in each era, we keep track of all swaps that touch the strand we are currently exploring. First, we define processes $(\bar{E}_t)_{t \geq 0}$ and $(\bar{\pi}_t)_{t \geq 0}$ as follows. As before, we start with $\bar{E}_0 = 1$ and $\bar{\pi}_0 = \mrm{id}$. We proceed to explore $\mc{P}_{\bar{E}_0}$ (in contrast to before, where we explored $\mc{P}_{\pi_0(E_0)}$). When we see a swap of the form $\{1, j\}$, $j \in [n]$, we update $\bar{\pi} \mapsto \bar{\pi} (1 ~ j) $. When we see a turnaround $\langle 1 ~ j \rangle$, $j \in (n:2n]$, the first exploration era ends, we update $\bar{E}$ to be 2, and we remove from $\mc{P}$ all points in $\mc{P}_{E_0}$. We then continue until the end of the $n$th exploration era. See Figure \ref{figure:alternative-exploration} for how one may visually compare this alternative exploration with our original exploration.

Formally, we may define a bijection on sets of points $P \mapsto P'$, which preserves the Poisson measure, and moreover if we follow our original exploration process on the set $P$, then that amounts to following the alternative exploration on the set $P'$. Under this bijection, the left bijection found by the original exploration is equal to the left bijection found by the alternative exploration, whereas the right bijections of the two explorations differ in a precise way, which is exactly encoded in the process $\bar{\pi}$. As an example, observe that in the previous picture, just before time $T_1$, we have that $\bar{\pi}_t = (4 ~ 3) (4 ~ 2)  = (2 ~ 3 ~ 4)$. Observe that $\bar{\pi}_t(4) = 2$, and on the right\revision{-}hand side of the original exploration, $2$ is matched to $6$. More generally, the rule is as follows. Let $\sigma(\mc{Q})$ be the left bijection found by the alternative exploration process. Then the right bijection $\tau(\mc{Q})$ is given by $\sigma(Q) \bar{\pi}_t$. Finally, because the bijection preserves the Poisson measure, when we apply the two explorations to a Poisson process, then they have the same law. We have thus arrived at the following result.

% Notice here we have not yet drawn in the right\revision{-}hand side of the diagram in the alternative exploration. This is because now the turnaround swap may be different on the left versus on the right, and the matching on the right needs to be determined.

\begin{lemma}\label{lemma:alternative-exploration}
We have that
\begin{equs}
\E \big[&\pointstobrauerfn(\mc{Q}(T_{n})) \ind(T_n \leq T) e^{2(n-1) T_1} e^{2(n-2)(T_2 - T_1)} \cdots e^{2(n-n) (T_n - T_{n-1})}\big] = \\
& \frac{1}{N^n n!}\sum_{\sigma : [n] \ra (n:2n]} \Big[\sigma ~ \sigma \E\big[\bar{\pi}_{T_n}\ind(T_n \leq T)  e^{2(n-1)T_1}  e^{2(n-2)(T_2 - T_1)}\cdots e^{2(n-n)(T_n - T_{n-1})}\big] \Big]
\end{equs}
\end{lemma}

\begin{remark}
The factor of $\frac{1}{N^n}$ arises because each turnaround incurs factor of $\frac{1}{N}$, and there are $n$ total turnarounds on the event $T_n \leq T$. The factor $\frac{1}{n!}$ arises because the first turnaround is equally likely to touch any of the $n$ bottom strands, the second turnaround is equally likely to touch any of the $n-1$ remaining bottom strands, etc.
\end{remark}

\begin{lemma}\label{lemma:expectation-permutation-alternative-exploration}
Let $U_1, \ldots, U_n \stackrel{i.i.d.}{\sim}\mrm{Exp}(1)$. We have that 
\begin{equs}
\E\big[\bar{\pi}_{T_n} \ind(T_n \leq T) e^{2(n-1)T_1}  & e^{2(n-2)(T_2 - T_1)}\cdots e^{2(n-n)(T_n - T_{n-1})}\big] = \\
&n! \E[\exp(- U_n J_n /N) \cdots \exp(-U_1 J_1 / N) \ind(U_1 + \cdots + U_n \leq T)].
\end{equs}
% \begin{equs}
% \E \big[\pointstobrauerfn(\mc{Q}(T_{n}), \labeling_{T_n}) \ind(T_n& \leq T) e^{2(n-1) T_1} e^{2(n-2)(T_2 - T_1)} \cdots e^{2(n-n) (T_n - T_{n-1})}\big] = \\
% & \sum_{\sigma : [n] \ra \{n+1, \ldots, 2n\}} \big[\sigma ~ \E[\exp(- U_n J_n /N) \cdots \exp(-U_1 J_1 / N) \ind(U_1 + \cdots + U_n \leq T)] \sigma \big] .
% \end{equs}
\end{lemma}
\begin{proof}
Note that the duration $T_k - T_{k-1}$ of the $k$th exploration process is an exponential random variable with rate $n-k+1$. We thus have the explicit formula
\begin{equs}
\E\big[\bar{\pi}_{T_n} &\ind(T_n \leq T) e^{2(n-1)T_1} e^{2(n-2)(T_2 - T_1)}\cdots e^{2(n-n)(T_n - T_{n-1})}\big] = \\
&\int_0^T dt_1 \int_{t_1}^T dt_2 \cdots \int_{t_{n-1}}^T dt_n \big(ne^{-nt_1}\big) \big( (n-1) e^{-(n-1)(t_2 - t_1)} \big)\cdots e^{-(t_n - t_{n-1})} ~\times \\
&~~~~~~~~~~~~~~~~~~~~~~~~~~~~~~~~~~ f_n(t_1) f_{n-1}(t_2 - t_1) \cdots f_1(t_n - t_{n-1}) ~\times \\
&~~~~~~~~~~~~~~~~~~~~~~~~~~~~~~~~~~ e^{2(n-1)t_1} e^{2(n-2)(t_2 - t_1)} \cdots e^{2(n-n)(t_n - t_{n-1})}.
\end{equs}
Here, $f_n(t_1)$ is the expected contribution of all swaps in the first exploration era, conditioned on $T_1 = t_1$, $f_{n-1}(t_2 - t_1)$ is the expected contribution of all swaps in the second exploration era, conditioned on $T_2 - T_1 = t_2 - t_1$, etc. Conditioned on $T_1 = t_1$, the number of total swaps is $\mrm{Poi}((n-1) t_1)$, and conditional on the total number of swaps being equal to $k$, the expected contribution is uniformly distributed on all possible sequences of $k$ swaps, i.e. $(-J_n / N (n-1))^k$. We thus have the explicit formula
\begin{equs}
f_n(t_1) &= e^{-(n-1)t_1} \sum_{k=0}^\infty \frac{((n-1) t_1)^k}{k!} \bigg(\frac{-J_n}{N(n-1)}\bigg)^{k} \\
&= e^{-(n-1) t_1} \sum_{k=0}^\infty \frac{(-t_1 J_n / N)^k}{k!} = e^{-(n-1) t_1} e^{-t_1 J_n / N}.
\end{equs}
More generally, we have the formula
\begin{equs}
f_{n-k+1}(t_k - t_{k-1}) = e^{-(n-k) (t_k - t_{k-1})} e^{-(t_k - t_{k-1}) J_{n-k+1}/N}, ~~ k \in [n]. 
\end{equs}
Inserting this into our first display, we obtain
\begin{equs}
\E\big[\bar{\pi}_{T_n} &\ind(T_n \leq T) e^{2(n-1)T_1} e^{2(n-2)(T_2 - T_1)}\cdots e^{2(n-n)(T_n - T_{n-1})}\big] =  \\
&n! \int_0^T dt_1 \int_{t_1}^T dt_2 \cdots \int_{t_{n-1}}^T dt_n ~ e^{-t_1} \cdots e^{-(t_n - t_{n-1})} e^{-t_1 J_n / N} e^{-(t_2 - t_1) J_{n-1}/N} \cdots e^{-(t_n - t_{n-1}) J_1 / N}.
\end{equs}
To finish, observe that the right\revision{-}hand side above is precisely the right\revision{-}hand side of the claimed identity.
\end{proof}

Up to now, we did not need to make any assumption on the size of $N$. We begin to do so here. Later, in Section \ref{section:poisson-exploration-general-N-proof}, we will show how to remove these assumptions, but for now we prefer to work in a simplified setting where the main ideas are more transparent.

\begin{lemma}\label{lemma:limit-permutation-alternative-exploration}
Suppose that $N \geq n$. Then
\begin{equs}
\lim_{T \toinf} \E\big[\bar{\pi}_{T_n} \ind(T_n \leq T) e^{2(n-1)T_1} e^{2(n-2)(T_2 - T_1)}\cdots e^{2(n-n)(T_n - T_{n-1})}\big] = n! N^n (N + J_n)^{-1} \cdots (N + J_1)^{-1}.
\end{equs}
\end{lemma}
\begin{proof}
By Lemma \ref{lemma:expectation-permutation-alternative-exploration}, we may compute
\begin{equs}
\E\big[\bar{\pi}_{T_n} &\ind(T_n \leq T) e^{2(n-1)T_1} e^{2(n-2)(T_2 - T_1)}\cdots e^{2(n-n)(T_n - T_{n-1})}\big]  =\\
&n! \int_0^T du_n \int_0^{T - u_n} du_{n_1} \cdots \int_0^{T - (u_n + \cdots + u_2)} du_1 \big(e^{-u_n} e^{-u_n J_n / N}\big) \cdots \big(e^{-u_1} e^{-u_1 J_1 / N}\big) .
\end{equs}
Since $N \geq n$, we have that $\|J_k / N\| < 1$ for all $k \in [n]$. This implies that the following integral is absolutely convergent (recall Remark \ref{remark:inverse-in-group-algebra}):
\begin{equs}
\int_0^\infty du_n \int_0^\infty du_{n-1} \cdots \int_0^\infty du_1 \big(e^{-u_n} e^{-u_n J_n / N}\big) \cdots \big(e^{-u_1} e^{-u_1 J_1 / N}\big) ,
\end{equs}
and moreover, the limit in question is equal to $n!$ times the above. To finish, simply observe that the above splits into a product of $n$ integrals, where the $k$th integral may be evaluated:
\begin{equs}
\int_0^\infty du_k e^{-u_k} e^{-u_k J_k / N} = \int_0^\infty du_k e^{-u_k(\id + J_k / N)} = \bigg(\id + \frac{J_k}{N}\bigg)^{-1} = N (N + J_k)^{-1}.
\end{equs}
The desired result follows.
\end{proof}

Next, we argue why the contribution to the partition function $e^{\binom{2n}{2} T - nT} \E[\pointstobrauerfn(\Sigma(T))]$ coming from the event $\{T_n > T\}$ vanishes in the $T \toinf$ limit (that is, as $T$ becomes large, we can assume that all exploration eras have finished by time $T$). We first show that when the numbers of top strands and bottom strands are mismatched, the expectation vanishes as $T \toinf$. This will be needed in the proof of Proposition \ref{prop:general-N-zero-contribution-from-unfinished-exploration-eras} later.

\begin{lemma}\label{lemma:mismatched-strands-bounded-in-T}
Suppose that $n \geq 2$ and $N \geq 2n$. Suppose that $\Sigma$ is a Poisson process arising from having $n-1$ top strands and $n$ bottom strands. Then
\begin{equs}
\sup_{T \geq 0} e^{\binom{2n-1}{2} T - (n-1) T} \big\| \E[ \pointstobrauerfn (\Sigma(T))] \big\| < \infty.
\end{equs}
\end{lemma}
% \begin{remark}
% If in the exponential prefactor we had $e^{\binom{2n-1}{2} T - (1/2)(2n-1)T}$ instead of $e^{\binom{2n-1}{2} T - (n-1)T}$, then the result would be immediate from the fact that 
% \begin{equs}
% e^{\binom{2n-1}{2} T - (1/2)(2n-1)T} \E[\pointstobrauerfn(\Sigma(T))] = \E[B_T^{\otimes (n-1)} \otimes \bar{B}_T^{\otimes n}],
% \end{equs}
% combined with the fact that the right\revision{-}hand side above converges (as $T \toinf$) to $\E [U^{\otimes (n-1)} \otimes \bar{U}^{\otimes n}] = 0$. The point of the lemma is that the convergence in fact happens at a slightly faster rate, in that we have convergence to zero even with an additional factor of $e^{T/2}$. This is needed for Proposition \ref{prop:zero-contribution-from-unfinished-exploration-eras}.
% \end{remark}
\begin{proof}
We proceed by induction. First, consider the base case $n = 2$. In this case, by conditioning on the first time of turnaround, we can explicitly compute
\begin{equs}
e^{2T} \E[\pointstobrauerfn(\Sigma(T))] = ~&e^{2T} \int_0^T 2 e^{-2u} X(u) Y Z(u) du + e^{2T} e^{-2T} e^{-T} e^{-TJ_2'/N} ,
\end{equs}
where $X(u)$ is the expected contribution of all swaps up to time $u$, and $Z(u)$ is the expected contribution of all points after time $u$, where both are conditioned on the first turnaround happening at time $u$. Also, $Y = \frac{1}{2} \big(\langle 1 ~ 2 \rangle + \langle 1 ~ 3 \rangle\big)$ is the expected contribution of the turnaround, since each of the two turnarounds is equally likely. Note that the time of first turnaround is exponential of rate $2$, which explains the presence of the $2 e^{-2u}$ term. The second term above corresponds to the case where the first turnaround happens after time $T$.

We have the explicit formulas
\begin{equs}
X(u) = e^{-u} e^{-u J_2' / N}, ~~ Z(u) = e^{-2(T-u)},
\end{equs}
where $J_2'$ is the Jucys-Murphy element which we view as acting on the bottom two strands (recall Definition \ref{def:bottom-jucys-murphy}). This formula follows because the number of swaps up to time $u$ is $\mrm{Poisson}(u)$, and each swap incurs a factor $-J_2' / N$. The fact that $Z(u) = e^{-2(T-u)}$ follows because once a turnaround occurs, we can argue via cancellation as in the proof of Lemma \ref{lemma:cancellation} that the only points which can occur thereafter are turnarounds between the same two strands.
% Conditioned on the first turnaround happening at time $u$, the expected contribution of all the swaps before time $u$ is given by $e^{-u} e^{-u J_2' / N}$, where $J_2'$ is second Jucys-Murphy element, which we view as acting on the bottom two strands (in this case, it is just a single swap). This formula follows because the number of swaps up to time $u$ is $\mrm{Poisson}(u)$, and each swap incurs a factor $-J_2' / N$. Next, each of the two turnarounds is equally likely, which gives the term $\frac{1}{2}\big(\langle 1 ~ 2 \rangle + \langle 1 ~ 3\rangle\big)$. Finally, once a turnaround occurs, we can argue via cancellation as in the proof of Lemma \ref{lemma:cancellation} that the only points which can occur thereafter are turnarounds between the same two strands. This contributes the term $e^{-2(T-u)}$.
% The second term above is when the first turnaround does not happen before time $T$. 
Plugging in the formulas for $X(u), Z(u)$, we may obtain the expression
\begin{equs}
\int_0^T e^{-u (\mrm{id} + J_2'/N)} du \big(\langle 1 ~ 2 \rangle + \langle 1 ~ 3 \rangle\big) + e^{-T (\mrm{id} + J_2'/N)} .
\end{equs}
Since $N \geq 2n$ is sufficiently large, as $T \toinf$ the above stays bounded (in fact, it converges to some explicit expression involving $(\mrm{id} + J_2' / N)^{-1}$, as in the proof of Lemma \ref{lemma:expectation-permutation-alternative-exploration}). This shows the case $n = 2$.

Now suppose the claim is true for some $n$. Suppose also that $N \geq 2(n+1)$. We show that the claim is true for $n+1$. As in the base case, by conditioning on the first time of turnaround, we may express (note that $\binom{2n+1}{2} - n = 2n^2$)
\begin{equs}\label{eq:mismatched-strands-inductive}
e^{\binom{2n+1}{2}T - nT} \E[\pointstobrauerfn(\Sigma(T))&] = e^{2n^2 T} \int_0^T n(n+1) e^{-n(n+1)u} X_n(u) Y_n Z_n(u)  du ~+ \\
&e^{2n^2 T} e^{-n(n+1)T}  e^{-\binom{n}{2} T} e^{-T(J_n + \cdots + J_1)/N}  e^{-\binom{n+1}{2}T} e^{-T(J_{n+1}' + \cdots + J_1')/N},
\end{equs}
where
$X_n(u)$ is the expected contribution of all swaps up to time $u$ and $Z_n(u)$ is the expected contribution of all points after time $u$, where both are conditioned on the first turnaround happening at time $u$. Also, $Y_n =  \frac{1}{n(n+1)} \sum_{i \in [n], j \in (n:2n+1]} \langle i ~ j \rangle$ is the expectation of the first turnaround. Let $J_1', \ldots, J_{n+1}'$ be the Jucys-Murphy elements which act on the bottom $n+1$ strands, as in Definition \ref{def:bottom-jucys-murphy}. Similar to before, we may explicitly compute
\begin{equs}
X_n(u) &= e^{-\binom{n}{2} u} e^{-\binom{n+1}{2}u} e^{-u(J_n + \cdots + J_1)/N} e^{-u(J_{n+1}' + \cdots + J_1')/N}, \\
Z_n(u) &= e^{-2(2n-1)(T-u)}  f_n(T - u),
\end{equs}
where $f_n(T-u)$ is the expected contribution of the points involving the remaining $n-1$ top and $n$ bottom strands after time $u$, conditioned on the first turnaround happening at time $u$. Observe that the $e^{-2(2n-1)(T-u)}$ factor in $Z_n(u)$ arises due to similar cancellations as in the proof of Lemma \ref{lemma:cancellation}, which allows us to restrict to the event that after the first turnaround $\langle i ~ j \rangle$, the only points which can involve either of the two matched strands are exactly the turnarounds of the form $\langle i ~ j \rangle$. This means that a total of $2(2n-1)$ rate-$1$ Poisson processes must have zero points on the interval $[u, T]$.

Plugging in our formulas for $X_n(u), Z_n(u)$, and using the identities $2n^2 - n(n+1) - \binom{n}{2} - \binom{n+1}{2} = -n$, $2n^2 - 2(2n-1) = 2(n-1)^2$, we have that the first term on the right\revision{-}hand side of \eqref{eq:mismatched-strands-inductive} is equal to
% Here, as before $J_n, \ldots, J_1$ are the Jucys-Murphy elements, which as usual act on the top $n$ strands, while $J_{n+1}', \ldots, J_1'$ are Jucys-Murphy elements which act on the bottom $n+1$ strands. Also, $X = \frac{1}{n(n+1)} \sum_{i \in [n], j \in [n:2n+1]} \langle i ~ j \rangle$ is the expectation of the first turnaround, and $f_{n-1}(T-u)$ is the expected contribution of all points after the first turnaround, conditioned on the first turnaround happening at time $u$. The first term in the right\revision{-}hand side above may be expressed as
\begin{equs}
n(n+1) \int_0^T e^{-nu} e^{-u(J_n + \cdots + J_1)/N} e^{-u(J_{n+1}' + \cdots + J_1')/N} Y e^{2(n-1)^2 (T-u)} f_n(T-u) du.
\end{equs}
By the inductive assumption, we have that $\sup_{S \geq 0} e^{2(n-1)^2 S} \|f_n(S)\| < \infty$. Also, since $N \geq 2(n+1)$, we have that $\|(J_n + \cdots + J_1)/N\| < n/2$ and $\|(J_{n+1}' + \cdots + J_1')/N\| < n/2$, which implies
\begin{equs}
\int_0^\infty e^{-nu} \big\| e^{-u(J_n + \cdots + J_1)/N} e^{-u(J_{n+1}' + \cdots + J_1')/N}\big\| du < \infty.
\end{equs}
Combining the two facts, we obtain that the first term on the right\revision{-}hand side of \eqref{eq:mismatched-strands-inductive} is uniformly bounded in $T$. The second term in the right\revision{-}hand side of \eqref{eq:mismatched-strands-inductive} may be expressed
\begin{equs}
e^{-nT} e^{-T(J_n + \cdots + J_1)/N}  e^{-T(J_{n+1}' + \cdots + J_1')/N}.
\end{equs}
By arguing as before, we may show that this stays bounded as $T \toinf$ (in fact, it converges to zero). This completes the proof of the inductive step.
\end{proof}

% \begin{proof}
% This follows because for fixed $T$, we have that \sky{technically, this is only true after applying the representation $\rho_+$ to the LHS}
% \begin{equs}
% e^{\binom{2n-1}{2} T - \frac{1}{2}(2n-1) T} \E[ \pointstobrauerfn (\mc{P}(T), \labeling_T)]  = \E[B_T^{\otimes (n-1)} \otimes \bar{B}_T^{\otimes n}], 
% \end{equs}
% and we have that
% \begin{equs}
% \lim_{T \toinf} \E[B_T^{\otimes (n-1)} \otimes \bar{B}_T^{\otimes n}] = \E[ U^{\otimes (n-1)} \otimes \bar{U}^{\otimes n}] = 0. 
% \end{equs}
% \end{proof}

Combining this lemma with an inductive argument, we can obtain the following.

\begin{prop}\label{prop:zero-contribution-from-unfinished-exploration-eras}
Suppose that $N \geq 2n$. We have that
\begin{equs}
\lim_{T \toinf} e^{\binom{2n}{2} T - nT} \E[\pointstobrauerfn(\Sigma(T)) \ind(T_n > T)] = 0.
\end{equs}
\end{prop}
\begin{proof}
Fix $N$. First, when $n = 1$, we have that
\begin{equs}
\E[\pointstobrauerfn(\Sigma(T)) \ind(T_1 > T)] = e^{-T} \mrm{id},
\end{equs}
where $\id$ here denotes the identity element of $\mc{B}_{n, n}$. The right\revision{-}hand side above clearly goes to zero as $T \toinf$. This shows the base case $n = 1$. Now suppose the result is true for some general $n \geq 1$. Suppose also that $N \geq 2(n+1)$. We proceed to show that the $n+1$ case is true. Towards this end, observe that we may decompose
\begin{equs}
\ind(T_{n+1} > T) = \ind(T_1 > T) + \ind(T_1 \leq T < T_{n+1}).
\end{equs}
% It suffices to show that for all $0 \leq k \leq n$, we have that
% \begin{equs}
% \lim_{T \toinf} e^{\binom{2(n+1)}{2} T - (n+1)T}  \E[\pointstobrauerfn(\Sigma(T)) \ind_{A_T} \ind(T_k \leq T < T_{k+1})] = 0.
% \end{equs}
We split into the two cases indicated above. In the first case, we condition on the exploration at time $T$:
\begin{equs}
e^{\binom{2(n+1)}{2} T - (n+1)T}  \E\big[ \E[\pointstobrauerfn(\Sigma(T))  ~|~ \mc{F}_T]\ind(T_1 > T)] \big].
\end{equs}
To help visualize, imagine we have the situation in Figure \ref{figure:not-all-turnarounds-completed}, where we explore the first strand until time $T$, and we have not yet seen a turnaround.
\begin{figure}[h]
\begin{center}
\includegraphics[width=7cm]{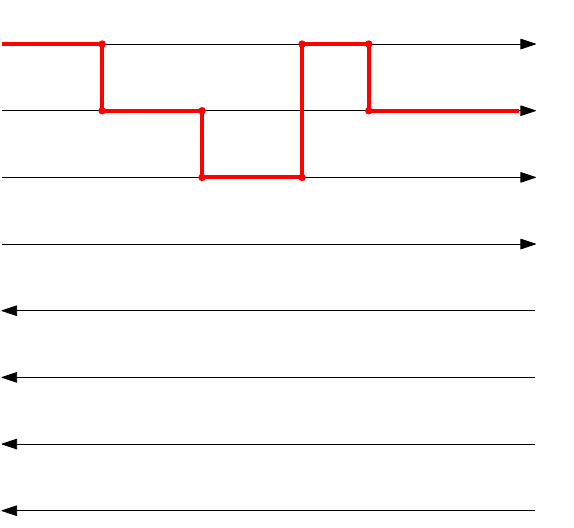}
\end{center}
\caption{We explore the top strand, and do not see a turnaround before time $T$.}\label{figure:not-all-turnarounds-completed}
\end{figure}

Conditioned on this picture, the expectation of the diagram can be computed as follows. First, since we have already explored one strand, the remaining points effectively form a Poisson process corresponding to $n$ top strands and $n+1$ bottom strands. Call this modified process $\bar{\Sigma}(T)$. We visualize this in Figure \ref{figure:not-all-turnarounds-completed-2}, where the top-most strand in the left diagram is dashed, to signify that there are no points touching this strand. 

\begin{figure}[h]
\begin{center}
\includegraphics[page=2, width=7cm]{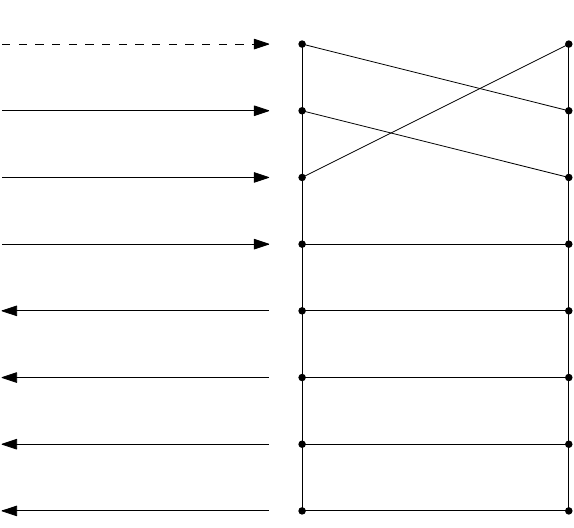}
\end{center}
\caption{Left: $\pointstobrauerfn(\bar{\Sigma}(T))$. Right: $\pointstobrauerfn(\mc{Q}(T))$.}\label{figure:not-all-turnarounds-completed-2}
\end{figure}

Having computed the expectation of the modified diagram in the left of Figure \ref{figure:not-all-turnarounds-completed-2}, to obtain the conditional expectation of $\pointstobrauerfn(\Sigma(T))$ we simply need to multiply by the right diagram in the figure, which captures the effect of all swaps seen by our exploration up to time $T$. This discussion corresponds to the following identity for the conditional expectation:
\begin{equs}
\E[\pointstobrauerfn(\Sigma(T)) ~|~ \mc{F}_T] = \E[ \pointstobrauerfn(\bar{\Sigma}(T))]\pointstobrauerfn(\mc{Q}_T) .
\end{equs}
% where $\bar{\Sigma}$ is a Poisson process corresponding to $n$ top strands and $n+1$ bottom strands.
We may then compute
\begin{equs}
e^{\binom{2(n+1)}{2} T - (n+1)T} & \E\big[ \E[\pointstobrauerfn(\Sigma(T))  ~|~ \mc{F}_T]\ind(T_1 > T)] \big]  \\
&=e^{\binom{2(n+1)}{2} T - (n+1)T} \E[\pointstobrauerfn(\bar{\Sigma}(T))]  \E[\pointstobrauerfn(\mc{Q}_T) \ind(T_1 > T)] \\
&= e^{\binom{2(n+1)}{2} T - (n+1)T} \E[\pointstobrauerfn(\bar{\Sigma}(T))] e^{-(n+1) T} e^{-nT} e^{-T J_{n+1} / N} \\
&= \big( e^{\binom{2n+1}{2} T - nT} \E[\pointstobrauerfn(\bar{\Sigma}(T))]\big)  e^{-T} e^{-T J_{n+1} / N} .
\end{equs}
As $T \toinf$, the right\revision{-}hand side above goes to zero, since by Lemma \ref{lemma:mismatched-strands-bounded-in-T} (and our assumption that $N \geq 2(n+1)$), the term in the parentheses above is $O(1)$, and since $N \geq n+1$, we have that $\|J_{n+1} / N\| < 1$, so that $e^{-T} e^{-T J_{n+1} / N} \ra 0$. This shows the inductive step in the first case.

Next, we consider the case corresponding to $\ind(T_1 \leq T < T_{n+1})$. We condition on the exploration at time $T_1$. Consider the diagram in Figure \ref{figure:not-all-turnarounds-completed-3} which corresponds to $n+1 = 4$.
\begin{figure}[h]
\begin{center}
\includegraphics[page=3, width=7cm]{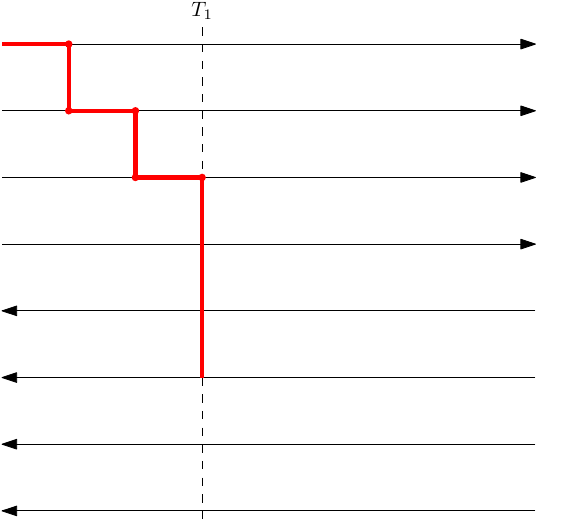}
\end{center}
\caption{We explore the top strand and see a turnaround at time $T_1$.}\label{figure:not-all-turnarounds-completed-3}
\end{figure}

On the event that $T_1 \leq T < T_{n+1}$, the portion of the diagram to the right of $T_1$ can be treated as having $n$ top strands and $n$ bottom strands. By arguing similarly to the previous case, i.e. by splitting our strand diagrams into the portion before $T_1$ and the portion after $T_1$, we may compute the conditional expectation:
\begin{equs}
 \E\big[\pointstobrauerfn(\Sigma(T)) \ind(T_1 \leq T < T_{n+1}) ~|~ \mc{F}_{T_1} = u] = \E\big[\pointstobrauerfn(\bar{\Sigma}(u))] e^{-4n(T - u)} \pointstobrauerfn(\mc{Q}_{T_1}) f_n(T - u),
\end{equs}
where $\bar{\Sigma}$ is a Poisson process corresponding to having $n$ top strands and $n+1$ bottom strands, and $f_n(T-u)$ is the expectation of the remaining $n$ top and $n$ bottom strands after time $u$, on the event that not all $n$ exploration eras end before time is up. Observe that by our inductive assumption, we have that for any $u \geq 0$,
\begin{equs}\label{eq:f-n-inductive-assumption}
\lim_{T \toinf} e^{\binom{2n}{2} (T-u) - n(T-u)} f_n(T-u) = 0.
\end{equs}
Since $T_1$ is an exponential random variable of rate $n+1$, we may compute the expectation:
\begin{equs}
~~ &e^{\binom{2(n+1)}{2} T - (n+1)T}  \E\big[\pointstobrauerfn(\Sigma(T))  \ind(T_1 \leq T < T_{n+1})\big] = \\
& e^{\binom{2(n+1)}{2} T - (n+1)T}\int_0^T du ~(n+1) e^{-(n+1) u} \big(\E\big[\pointstobrauerfn(\bar{\Sigma}(u))]  \big) \big(e^{-nu} e^{-uJ_{n+1}/N}\big)  Y  \big(e^{-4n(T - u)}  f_n(T - u) \big)
\end{equs}
Here, the term $e^{-nu} e^{-u J_{n+1} / N}$ arises from taking the expectation of all swaps in the first exploration era (i.e. $\pointstobrauerfn(\mc{Q}(T_1))$), conditioned on $T_1 = u$, and $Y$ is the expectation of the first turnaround. Since
\begin{equs}
\binom{2(n+1)}{2} - (n+1) &= 2n^2 + 2n \\
\binom{2(n+1)}{2} - (n+1) - 4n &= \binom{2n}{2} - n, 
\end{equs}
we have that the above is further equal to
\begin{equs}
(n+1) \int_0^T du ~ \big( e^{2n^2 u} \E[\pointstobrauerfn(\bar{\Sigma}(u))] \big) \big(e^{-u} e^{-u J_{n+1} / N} \big) Y e^{-4n(T-u)} f_n(T - u).
\end{equs}
Now since $N \geq 2n$, by Lemma \ref{lemma:mismatched-strands-bounded-in-T}, we have that $e^{2n^2 u} \E[\pointstobrauerfn(\bar{\Sigma}(u))]  = O(1)$ and $e^{-u} e^{-u J_{n+1} / N}$ is integrable. Combining this with \eqref{eq:f-n-inductive-assumption} and dominated convergence, we finally obtain
\begin{equs}
\lim_{T \toinf} e^{\binom{2(n+1)}{2} T - (n+1)T}  \E\big[\pointstobrauerfn(\Sigma(T)) \ind(T_1 \leq T < T_{n+1})\big] = 0.
\end{equs}
This completes the proof of the inductive step, and thus the desired result now follows.
\end{proof}

We can now finally take the $T \toinf$ limit. 
% Here is where we impose the condition that $N \geq n$. 
% Later, in Section \ref{section:poisson-exploration-general-N-proof}, we explain a more general argument which does not need this assumption.

\begin{prop}\label{prop:large-T-limit-large-N-case}
Suppose that $N \geq 2n$. Then as $T \toinf$, we have that
\begin{equs}
\lim_{T \toinf} e^{\binom{2n}{2} T - nT} \E[\pointstobrauerfn(\Sigma(T))] = \sum_{\sigma, \tau : [n] \ra (n:2n]} \Wg_N(\sigma^{-1} \tau) [\sigma ~ \tau].
\end{equs}
\end{prop}
\begin{proof}
% First, observe that because $N \geq n$, we have that for $k \in [n]$, $\exp(-U_k J_k / N)$ is integrable with expectation
% \begin{equs}
% \E[\exp(-U_k J_k/N)] = \int_0^\infty e^{-u} e^{-u J_k /N} du = \int_0^\infty e^{-u(\mrm{id} + J_k / N)} du = \bigg(\mrm{id} + \frac{J_k}{N}\bigg)^{-1}.
% \end{equs}
% (Recall that $N + J_k$ is invertible only when $N \geq k$.)
% we may expand
% \begin{equs}
%  \E[\exp(- U_k J_k / N) \ind(U_k \leq T) ] &= \int_0^T  e^{-u} e^{-u J_k / N} du \\
%  &= \sum_{j=0}^\infty \frac{1}{j!}\bigg(\frac{-J_k}{N}\bigg)^j \int_0^T e^{-u} u^j du.
% \end{equs}
% Since $N \geq n$, we have that the series
% $\sum_{j=0}^\infty \big(\frac{-J_k}{N}\big)^j$
% is absolutely convergence, and moreover is equal to $\big(\mrm{id} + J_k / N)^{-1}$. Thus by dominated convergence, we have that 
% \begin{equs}
% \lim_{T \toinf} \E[\exp(- U_k J_k / N) \ind(U_k \leq T) ] = \sum_{j=0}^\infty \frac{1}{j!} \bigg(-\frac{J_k}{N}\bigg)^j j! = \bigg(\mrm{id} + \frac{J_k}{N}\bigg)^{-1}, 
% \end{equs}
% as desired. Next, one may show by similar techniques that when $N \geq n$, we have that
% We can further deduce that (applying \eqref{eq:weingarten-jucy-murphy-large-N} in the last identity)
% \begin{equs} 
% \lim_{T \toinf}  &\E[\exp(- U_n J_n /N) \cdots \exp(-U_1 J_1 / N) \ind(U_1 + \cdots + U_n \leq T)] \\
% &= \E[\exp(-U_n J_n / N)] \cdots \E[\exp(-U_1 J_1 /N)] \\
% &= \bigg(\mrm{id} + \frac{J_n}{N}\bigg)^{-1} \cdots \bigg(\mrm{id} + \frac{J_1}{N}\bigg)^{-1} = N^n \Wg_N.
% \end{equs}
By combining Proposition \ref{prop:strand-by-strand-exploration}, Lemmas \ref{lemma:alternative-exploration} and \ref{lemma:limit-permutation-alternative-exploration}, and Proposition \ref{prop:zero-contribution-from-unfinished-exploration-eras}, we obtain
\begin{equs}
\lim_{T \toinf}  e^{\binom{2n}{2} T - nT} \E[\pointstobrauerfn(\Sigma(T))] &= \sum_{\sigma : [n] \ra (n:2n]} \big[\sigma ~ \sigma  \Wg_N \big] \\
&= \sum_{\sigma : [n] \ra (n:2n]} \sum_{\pi \in \symgrp_n} [\sigma ~ \sigma \pi] \Wg_N(\pi) \\
&= \sum_{\sigma, \tau : [n] \ra (n:2n]} [\sigma ~ \tau] \Wg_N(\sigma^{-1} \tau),
\end{equs}
as desired.
% To finish, recall that $\Wg_N(\sigma^{-1} \tau) = \Wg_N(\sigma \tau^{-1})$, because $\Wg_N$ is a class function.
% Sending $T \toinf$, we have that
% \begin{equs}
% \lim_{T \toinf} \E[\bar{\pi}_{T_n} \ind(T_n \leq T)] &= \lim_{T \toinf} \E[\exp(- U_n J_n /N) \cdots \exp(-U_1 J_1 / N) \ind(U_1 + \cdots + U_n \leq T)]  \\
% &= \E[\exp(-U_n J_n / N)] \cdots \E[\exp(-U_1 J_1 /N)].
% \end{equs}
% For $k \in [n]$, we may compute
% \begin{equs}
% \E[\exp(-U_k J_k/N)] &= \int_0^\infty e^{-u} e^{-u J_k /N} du \\
% &= \int_0^\infty e^{-u(\id + J_k/N)} du \\
% &= \bigg(\id + \frac{J_k}{N}\bigg)^{-1}.
% \end{equs}
% It thus follows that
% \begin{equs}
% \lim_{T \toinf} \E[\bar{\pi}_{T_n} \ind(T_n \leq T)] = \bigg(\id + \frac{J_n}{N}\bigg)^{-1} \cdots \bigg(\id + \frac{J_1}{N}\bigg)^{-1} = \Wg_N. 
% \end{equs}
% From this, we obtain
% \begin{equs}
% \lim_{T \toinf} \sum_{\sigma : [n] \ra [n:2n]} \big[\sigma ~ \sigma \E[\bar{\pi}_{T_n} \ind(T_n \leq T)] \big] &= \sum_{\sigma : [n] \ra [n:2n]} \big[\sigma ~\sigma \Wg_N] \\
% &= \sum_{\sigma : [n] \ra [n:2n]} \sum_{\pi \in \symgrp_n} [\sigma ~ \sigma \pi] \Wg_N(\pi) \\
% &= \sum_{\sigma, \tau : [n] \ra [n:2n]} [\sigma ~ \tau] \Wg_N(\sigma^{-1} \tau).
% \end{equs}
\end{proof}

We can now prove Theorem \ref{thm:weingarten-recovery} in the case $N \geq 2n$.

\begin{proof}[Proof of Theorem \ref{thm:weingarten-recovery} when $N \geq n$]
Recall from \eqref{eq:w-T-pi-brauer} that
\begin{equs}
 e^{\binom{2n}{2} T - nT} \E[\pointstobrauerfn(\Sigma(T))] = \sum_{\pi \in \mc{M}(n, n)} w_T(\pi) \pi.
\end{equs}
By Proposition \ref{prop:large-T-limit-large-N-case}, we obtain
\begin{equs}
\lim_{T \toinf} w_T(\pi) = \sum_{\sigma, \tau : [n] \ra (n:2n]} \ind(\pi = [\sigma ~ \tau]) \Wg_N(\sigma^{-1} \tau).
\end{equs}
Since $w_T(\pi_1, \ldots, \pi_L) = w_T(\pi_1) \cdots w_T(\pi_L)$ by definition, the desired result now follows by Lemma \ref{lemma:expectations-words-unitary-bm}.
\end{proof}

\begin{remark}[Comparison to \cite{Dahlqvist2017}]\label{remark:dahlqvist-comparision}
If one translates Dahlqvist's proof to the language of Poisson point processes, then his strategy amounts to an exploration of the Poisson process which simultaneously explores all strands. This is certainly a natural exploration to try. \cite[Lemma 5.1]{Dahlqvist2017} amounts to the statement that the main contribution comes from the event that all exploration eras end for this ``simultaneous exploration". \cite[Lemma 5.2]{Dahlqvist2017} gives a formula for the limiting contribution on this main event. He then extracts the Weingarten function from this formula by \cite[Lemma 5.3]{Dahlqvist2017}.

We believe that our proof technique via strand-by-strand exploration is intrinsically interesting, because first of all it is rather surprising that such an exploration actually works. Recall that this was Proposition \ref{prop:strand-by-strand-exploration}, whose proof rested on certain cancellations that could be uncovered (Lemma \ref{lemma:cancellation}). Moreover, the strand-by-strand exploration naturally uncovers the Jucys-Murphy elements, thus giving an alternative perspective on the appearance of the Weingarten function. Finally, the strand-by-strand exploration naturally leads to a single-strand recursion that results in a slightly more general version of the Makeenko-Migdal/Master loop/Schwinger-Dyson equation -- see Remarks \ref{remark:schwinger-dyson} and \ref{remark:schwiner-dyson-2} for more discussion.
\end{remark}

\subsection{Extension to general values of \texorpdfstring{$N$}{N}}\label{section:poisson-exploration-general-N-proof}

Recall that in the proof of Proposition \ref{prop:large-T-limit-large-N-case}, we deduced the existence of 
\begin{equs}
\lim_{T \toinf} e^{\binom{2n}{2} T - nT} \E[\pointstobrauerfn(\Sigma(T))] 
\text{ from the existence of }
\lim_{T \toinf} \E[\bar{\pi}_{T_n} \ind(T_n \leq T)] \in \C[\symgrp_n].
\end{equs}
However, when $N \leq n$, the trouble is that the latter limit no longer exists. Thus to prove Theorem \ref{thm:weingarten-recovery} in the case where $N$ is small, we need some alternative argument which does not rely on convergence in the group algebra. Indeed, we will show that although $\lim_{T \toinf} \E[\bar{\pi}_{T_n} \ind(T_n \leq T)]$ does not necessarily exist in $\C[\symgrp_n]$, once we apply the representation $\rho_+$ (Definition \ref{def:rho-plus}), the limit {\it does exist}. Moreover, the limit $\lim_{T \toinf} \rho_+\big(\E[\bar{\pi}_{T_n} \ind(T_n \leq T)]\big)$ already contains enough information in order to compute expectations of traces of words. Once we have built up enough background, the actual proof of Theorem \ref{thm:weingarten-recovery} for general values of $N$ will be a small variation of the proof for large $N$, as the major technical steps were already covered in Section \ref{section:poisson-exploration} (and any additional background covered in Section \ref{section:rep-theory}).

\begin{figure}[h]
\centering
    \includegraphics{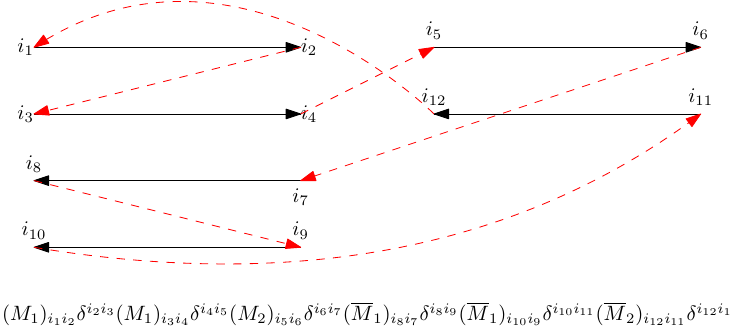}
    \caption{Visualization of the calculation of $\Tr(M_1^2 M_2 M_1^{-2} M_2^{-1})$.}\label{figure:calculating trace}
\end{figure}

Towards this end, it will be useful to recall why expectations of traces of words may be reduced to weighted sums over the Brauer algebra, i.e. why Lemma \ref{lemma:expectations-words-unitary-bm} is true.
Let $\Gamma$ be a word on letters $\{\lambda_1, \ldots, \lambda_L\}$. We may assume $\Gamma = \lambda_{c(1)}^{\varep(1)} \cdots \lambda_{c(n)}^{\varep(n)}$, where $\varep : [n] \ra \{\pm1\}$ and $c : [n] \ra [L]$. Let $M = (M_1, \ldots, M_L)$ be a given collection of $N \times N$ Unitary matrices. 
The computation of $\Tr(M(\Gamma)) = \Tr\big(M_{c(1)}^{\varep(1)} \cdots M_{c(n)}^{\varep(n)}\big)$ may be visualized in terms of the strand diagram as in Figure \ref{figure:calculating trace}, where we consider the concrete case $\Gamma = \lambda_1^2 \lambda_2 \lambda_1^{-2} \lambda_2^{-1}$. 

In Figure \ref{figure:calculating trace}, we can imagine we are traversing the strand diagram. Every black strand contributes a matrix element, and every dashed red strand enforces an identification of indices. In the end we sum over all indices which appear. Of course, we could have written the trace more succinctly as
\begin{equs}
\Tr(M_1^2 M_2 M_1^{-2} M_2^{-1}) = (M_1)_{i_1 i_2} (M_1)_{i_2 i_3} (M_2)_{i_3 i_4} (\ovl{M}_1)_{i_5 i_4} (\ovl{M}_1)_{i_6 i_5} (\ovl{M}_2)_{i_1 i_6},
\end{equs}
but we prefer to keep the $\delta$ functions because they correspond to the dashed red lines. We now want to give an expression as above for general words and strand diagrams. Given the strand diagram of a word $\Gamma$, note that the diagram has a single component, with a unique ordering of its vertices $x_1, \ldots, x_V$ up to cyclic equivalence. This ordering is such that the edges alternate between black strands and dashed red lines. Let $B(\Gamma)$ be the set of black strands, and $R(\Gamma)$ be the set of dashed red lines. Further split $B(\Gamma) = B_+(\Gamma) \cup B_-(\Gamma)$, where $B_+(\Gamma), B_-(\Gamma)$ are the set of positive (i.e. right) and negative (i.e. left)-oriented black strands. In the previous example,
\begin{equs}
B_+(\Gamma) &= \{(x_1, x_2), (x_3, x_4), (x_5, x_6)\},~~ B_-(\Gamma) = \{(x_7, x_8), (x_9, x_{10}), (x_{11}, x_{12})\}, \\
R(\Gamma) &= \{(x_2, x_3), (x_4, x_5), (x_6, x_7), (x_8, x_9), (x_{10}, x_{11})\}.
\end{equs}
Given a collection of indices $i = (i_1, \ldots, i_V) \in [N]^V$, and an edge $e = (x_j, x_{j+1})$, let $i_e = (i_j, i_{j+1})$, $i_{-e} = (i_{j+1}, i_j)$. Let $r(e) \in [L]$ be the index of the letter that $e$ corresponds to. Then the general formula for $\Tr(M(\Gamma))$ in terms of the strand diagram is:
\begin{equs}
\Tr(M(\Gamma)) = \prod_{e \in B_+(\Gamma)} (M_{r(e)})_{i_e} \prod_{e \in B_-(\Gamma)} (\ovl{M}_{r(e)})_{i_{-e}} \prod_{e \in R(\Gamma)} \delta^{i_e},
\end{equs}
where we implicitly sum over $i = (i_1, \ldots, i_V) \in [N]^V$.

More generally, we may extend the definitions $B_+(\bm \Gamma)$, $B_-(\bm \Gamma)$, $R(\bm \Gamma)$ to a collection of words $\bm \Gamma = (\Gamma_1, \ldots, \Gamma_k)$, so that
\begin{equs}
\Tr(M(\Gamma_1)) \cdots \Tr(M(\Gamma_k))  = \prod_{e \in B_+(\bm \Gamma)} (M_{r(e)})_{i_e} \prod_{e \in B_-(\bm \Gamma)} (\ovl{M}_{r(e)})_{i_{-e}} \prod_{e \in R(\bm \Gamma)} \delta^{i_e}.
\end{equs}
Now the point is as follows. If $M_1, \ldots, M_L$ are independent $\UN$-valued Brownian motions, then upon taking expectations of the above, we may obtain that $\E[\Tr(M(\Gamma))]$ is equal to a weighted sum of diagrams as follows. 

First, for $\ell \in [L]$, let $B_+(\bm \Gamma, \lambda_\ell)$, $B_-(\bm \Gamma, \lambda_\ell)$ be the sets of right-directed and left-directed edges corresponding to the letter $\lambda_\ell$. Since the $M_1, \ldots, M_L$ are independent, we have that
\begin{equs}
\E[\Tr(M(\Gamma))] = \prod_{\ell \in [L]} \E\bigg[\prod_{e \in B_+(\bm \Gamma, \lambda_\ell)} (M_{r(e)})_{i_e} \prod_{e \in B_-(\bm \Gamma, \lambda_\ell)} (\ovl{M_{r(e)}})_{i_{-e}} \bigg] \prod_{e \in R(\bm \Gamma)} \delta^{i_e}.
\end{equs}

We recall the following lemma from~\cite[Appendix A]{park2023wilson} (see also~\eqref{eq:w-T-def}) which gives a formula for each of the expectations appearing in the right\revision{-}hand side above. Recall the definition of $w_T(\pi)$ from Definition \eqref{eq:w-T-def} and of $\mc{M}(n, n)$ from Definition \ref{def:walled-brauer-algebra}.

\begin{prop}\label{prop:unitary-brownian-motion-expectation}
Let $i_1, \ldots, i_n$, $i_1', \ldots, i_n'$, $j_1, \ldots, j_n$, $j_1', \ldots, j_n' \in [N]$. We have that
\begin{equs}
\E \big[(B_T)_{i_1j_1} \cdots (B_T)_{i_n j_n} (\ovl{B}_T)_{i_1' j_1'} \cdots (\ovl{B}_T)_{i_n' j_n'} \big] = \sum_{\pi \in \mc{M}(n, n)} w_T(\pi) \ind(\text{indices match with $\pi$}).
\end{equs}
\end{prop}

% Letting $B_T$ be Unitary Brownian motion at time $T$, we have that any product of matrix entries
% \begin{equs}
% \E \big[(B_T)_{i_1j_1} \cdots (B_T)_{i_n j_n} (\ovl{B}_T)_{i_1' j_1'} \cdots (\ovl{B}_T)_{i_n' j_n'} \big]
% \end{equs} 

% is given by the partition function of the Poisson process, restricted to the event that the final diagram is valid for the indices $i_1, \ldots, i_n, i_1', \ldots, i_n'$, $j_1, \ldots, j_n, j_1', \ldots, j_n'$ (in that all indices matched by the diagram are equal). Letting $w_T(\pi)$ denote the weight of a diagram $\pi \in \mc{B}_{n, n}$, we have that
% \begin{equs}
% \E\big[(B_T)_{i_1j_1} \cdots (B_T)_{i_n j_n} (\ovl{B}_T)_{i_1' j_1'} \cdots (\ovl{B}_T)_{i_n' j_n'} \big] = \sum_{\pi \in \mc{B}_{n, n}} w_T(\pi) \ind(\text{indices match with $\pi$}). 
% \end{equs}

\begin{figure}[h]
\centering
    \includegraphics[page=2]{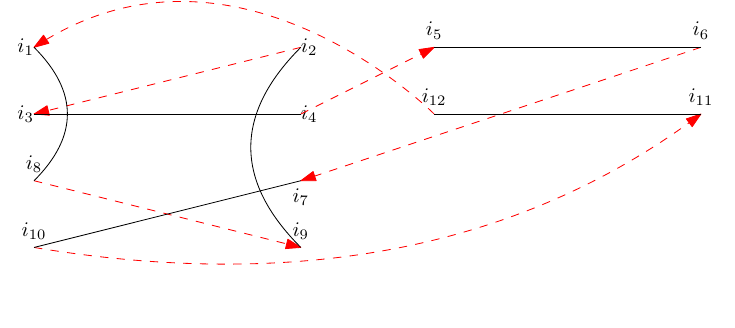}
    \vspace{-10mm}
    \caption{Visualization of the calculation of $\Tr(M_1^2 M_2 M_1^{-2} M_2^{-1})$, where now $M_1 = \rho_+(\pi_1)$, $M_2 = \rho_+(\pi_2)$, for the matchings $\pi_1 \in \mc{M}(4)$, $\pi_2 \in \mc{M}(2)$ displayed in the figure.}\label{figure:calculating trace-2}
\end{figure}

Using this proposition, we may write
\begin{equs}
\E[\Tr(B_T(\bm \Gamma))] = \sum_{\bm \pi = (\pi_1, \ldots, \pi_L)} w_T(\pi_1) \cdots w_T(\pi_L) \prod_{\ell \in L} \prod_{\{a, b\} \in \pi_\ell} \delta_{i_a i_b} \prod_{e \in R(\bm \Gamma)} \delta^{i_e}.
\end{equs}
Now, observe that 
\begin{equs}
\prod_{\ell \in L}  \prod_{\{a, b\} \in \pi_\ell} \delta_{i_a i_b} \prod_{e \in R(\Gamma)} \delta^{i_e} = N^{\numcomp(\bm \Gamma, \bm \pi)},
\end{equs}
where recall $\numcomp(\Gamma, \pi)$ is the number of components obtained by deleting all black strands but including all interior connections specified by the matchings $\pi_1, \ldots, \pi_L$. For instance, in our previous example, suppose our matchings were as in Figure \ref{figure:calculating trace-2}.
Since each edge in Figure \ref{figure:calculating trace-2} (be it red or black) imposes a constraint on the indices, the total number of free summation indices is exactly equal to the number of connected components in the above diagram. Each free summation index may take one of $N$ values, whence the term $N^{\numcomp(\Gamma, \pi)}$. Lemma \ref{lemma:expectations-words-unitary-bm} follows directly from these considerations.

Now recall from Definition \ref{def:rho-plus} that the matrix elements of the representation $\rho_+(\pi)$ are exactly given by
\begin{equs}
(\rho_+(\pi))_{i \sqcup i', j \sqcup j'} = \ind(\text{indices match with $\pi$}).
\end{equs}
Here, $i \sqcup i'$ denotes the length-$2n$ vector of indices given by concatenation: $(i_1, \ldots, i_n, i_1', \ldots ,i_n')$, and similarly for $j \sqcup j'$. Combining this with the previous discussion, we have the following result. First, for some notation, let $\mbf{i}_\ell, \mbf{j}_\ell$ respectively collect all left and right indices which appear in the strand diagram corresponding to $\lambda_\ell$. For the example in Figure \ref{figure:calculating trace-2}, we have that $\mbf{i}_1 = (i_1, i_3, i_8, i_7)$, $\mbf{j}_1 = (i_2, i_4, i_7, i_9)$, $\mbf{i}_2 = (i_5, i_{12})$, $\mbf{j}_2  = (i_6, i_{11})$.

\begin{lemma}\label{lemma:rho-plus-sum-num-comp}
Let $\bm \Gamma$ be a collection of words on letters $\{\lambda_1, \ldots, \lambda_L\}$.
Let $\bm \pi = (\pi_\ell, \ell \in [L])$, where $\pi_\ell \in \mc{M}(n_{\pm}(\lambda_\ell))$ for each $\ell \in [L]$. Then
\begin{equs}
\prod_{\ell \in [L]} \rho_+(\pi_\ell)_{\mbf{i}_\ell \mbf{j}_\ell} \prod_{e \in R(\bm \Gamma)} \delta^{i_e} = N^{\numcomp(\bm \Gamma, \bm \pi)}.
\end{equs}
\end{lemma}

% We thus have that
% \begin{equs}
% \E\big[(B_T)_{i_1j_1} \cdots (B_T)_{i_n j_n} (\ovl{B}_T)_{i_1' j_1'} \cdots (\ovl{B}_T)_{i_n' j_n'} \big] = \sum_{\pi \in \mc{B}_{n, n}} w_T(\pi) \rho_+(\pi)_{ii', jj'}, 
% \end{equs}
% which is another way of saying that
% \begin{equs}
% \E\big[ B_T^{\otimes n} \otimes \ovl{B}_T^{\otimes n}] = \sum_{\pi \in \mc{B}_{n, n}} w_T(\pi) \rho_+(\pi).
% \end{equs}
Using this lemma and the previous discussion, we could have written $\E[\Tr(M(\Gamma))]$ in terms of $\rho_+(\pi_1), \ldots, \rho_+(\pi_L)$, as follows.

\begin{lemma}\label{lemma:rho-plus-expectation-words-unitary-brownian-motion}
Let $\bm \Gamma = (\Gamma_1, \ldots, \Gamma_n)$ be a collection of words with letters $\{\lambda_1, \ldots, \lambda_L\}$. We have that
\begin{equs}
\E[\Tr(B_T(\bm \Gamma))] = \bigg(\sum_{\pi_1 \in \mc{B}_{n_1, n_1}} w_T(\pi_1) \rho_+(\pi_1) \bigg)_{\mbf{i}_1 \mbf{j}_1} \cdots \bigg(\sum_{\pi_L \in \mc{B}_{n_L, n_L}} w_T(\pi_L) \rho_+(\pi_L) \bigg)_{\mbf{i}_L \mbf{j}_L} \prod_{e \in R(\bm \Gamma)} \delta^{i_e},
\end{equs}
\end{lemma}

As mentioned at the beginning of this subsection, we have rewritten expectations of traces of words of Unitary Brownian motion in terms of some function (namely, $\rho_+$) of weighted sums over the Brauer algebra. The point now is that $\lim_{T \toinf} \rho_+\big(\E[\bar{\pi}_{T_n} \ind(T_n \leq T)]\big)$ exists for all $N$. Once we show this, the rest of the proof of Theorem \ref{thm:weingarten-recovery} in the case of general $N$ is exactly the same.

\begin{lemma}[Analog of Lemma \ref{lemma:limit-permutation-alternative-exploration}]\label{lemma:general-N-limit-permutation-alternative-exploration}
We have that
\begin{equs}
\lim_{T \toinf} \rho_+\big(\E\big[\bar{\pi}_{T_n} e^{2(n-1)T_1}  \ind(T_n \leq T) e^{2(n-2)(T_2 - T_1)}\cdots e^{2(n-n)(T_n - T_{n-1})}\big]\big) = n! N^n \rho_+(\Wg_N).
\end{equs}
\end{lemma}
\begin{proof}
By Lemma \ref{lemma:expectation-permutation-alternative-exploration}, we may compute
\begin{equs}
\rho_+\big(&\E\big[\bar{\pi}_{T_n} e^{2(n-1)T_1}  \ind(T_n \leq T) e^{2(n-2)(T_2 - T_1)}\cdots e^{2(n-n)(T_n - T_{n-1})}\big]\big) = \\
&n! \int_0^T du_n \int_0^{T-u_1} du_{n-1}\cdots \int_0^{T-(u_1 + \cdots + u_{n-1})} du_1  e^{-u_n \rho_+(\mrm{id} + J_n / N)} \cdots e^{-u_1 \rho_+ (\mrm{id} + J_1 /N)}  .
\end{equs}
By Lemma \ref{lemma:rho-plus-Jucy-Murphy-eigenvalues}, for all $k \in [n]$, all eigenvalues of $\rho_+(J_k)$ are at least $-N+1$, and thus all eigenvalues of $\rho_+(\mrm{id} + J_k / N)$ are at least $1/N$, and in particular all eigenvalues are strictly positive. Thus as we send $T \toinf$ the above converges to (applying Lemma \ref{lemma:rho-plus-weingarten-is-inverse-rho-plus-cycles} in the final identity)
\begin{equs}
n! \int_0^\infty e^{-u_n \rho_+(\mrm{id} + J_n / N)} du_n \cdots \int_0^\infty  e^{-u_1 \rho_+ (\mrm{id} + J_1 /N)} du_1 &= n!\rho_+(\mrm{id} + J_n / N)^{-1} \cdots \rho_+(\mrm{id} + J_1 / N)^{-1} \\
&= n! N^n \rho_+(N + J_n)^{-1} \cdots \rho_+(N + J_1)^{-1} \\
% &= N^n \rho_+\bigg( \bigg(\mrm{id} + \frac{J_n}{N}\bigg) \cdots \bigg(\mrm{id} + \frac{J_1}{N}\bigg)\bigg)^{-1} \\
% &= N^n \rho_+\bigg(\sum_{\sigma \in \symgrp_n} N^{\#\mrm{cycles}(\sigma)} \sigma \bigg)^{-1} \\
&= n! N^n \rho_+(\Wg_N),
\end{equs}
as desired.
\end{proof}

We also have the following analogs of Lemma \ref{lemma:mismatched-strands-bounded-in-T} and Proposition \ref{prop:zero-contribution-from-unfinished-exploration-eras}.

\begin{lemma}[Analog of Lemma \ref{lemma:mismatched-strands-bounded-in-T}]\label{lemma:general-N-mismatched-strands-bounded-in-T}
Suppose that $\Sigma$ is a Poisson process arising from having $n-1$ top strands and $n$ bottom strands. Then
\begin{equs}
\sup_{T \geq 0} e^{\binom{2n-1}{2} T - (n-1) T} \big\| \rho_+\big( \E[ \pointstobrauerfn (\Sigma(T))] \big)\big\| < \infty.
\end{equs}
\end{lemma}
\begin{proof}
In the proof of Lemma \ref{lemma:mismatched-strands-bounded-in-T}, the condition on $N$ was needed to show that
\begin{equs}
\int_0^\infty e^{-nu} \big\| e^{-u(J_n + \cdots + J_1)/N} e^{-u(J_{n+1}' + \cdots + J_1')/N}\big\| du &< \infty, \\
\sup_{T \geq 0} e^{-nT} \big\| e^{-T(J_n + \cdots + J_1)/N}  e^{-T(J_{n+1}' + \cdots + J_1')/N}\big\| &< \infty.
\end{equs}
When we apply $\rho_+$, we instead need to show that
\begin{equs}
\int_0^\infty e^{-nu} \big\| e^{-u \rho_+(J_n + \cdots + J_1)/N} e^{-u \rho_+(J_{n+1}' + \cdots + J_1')/N}\big\| du &< \infty, \\
\sup_{T \geq 0} e^{-nT} \big\| e^{-T\rho_+(J_n + \cdots + J_1)/N}  e^{-T\rho_+(J_{n+1}' + \cdots + J_1')/N}\big\| &< \infty.
\end{equs}
% Both of the claims follow from the that that the eigenvalues of 
% \begin{equs}
% \exp(- \rho_+(J_k) / N) \exp(-\rho_+(J_k')/N) 
% \end{equs}
% have absolute value strictly less than $e$ for all $2 \leq k \leq n$, as well as the fact that the eigenvalues of $\exp(-\rho_+(J_{n+1}') / N)$ have absolute value strictly less than $e$, and the fact that $J_1 = J_1' = 0$ (by definition).
These claims both follow from Corollary \ref{cor:eigenvalue-lower-bound-sum-Jucys-Murphy}, which gives that the eigenvalues of $\frac{1}{N}\rho_+(J_n + \cdots + J_1) + \frac{1}{N} \rho_+(J_{n+1}' + \cdots + J_1')$ are all strictly greater than $-n$.
\end{proof}

\begin{prop}[Analog of Proposition \ref{prop:zero-contribution-from-unfinished-exploration-eras}]\label{prop:general-N-zero-contribution-from-unfinished-exploration-eras}
We have that
\begin{equs}
\lim_{T \toinf} e^{\binom{2n}{2} T - nT} \rho_+\big( \E[\pointstobrauerfn(\Sigma(T)) \ind(T_n > T)]\big) = 0.
\end{equs}
\end{prop}
\begin{proof}
The points in the proof of Proposition \ref{prop:zero-contribution-from-unfinished-exploration-eras} where we needed $N$ to be large were in the application of Lemma \ref{lemma:mismatched-strands-bounded-in-T} and in arguing that $e^{-u} e^{-u J_{n+1} / N}$ is integrable. For the present proposition, we may apply Lemma \ref{lemma:general-N-mismatched-strands-bounded-in-T} which does not require $N$ to be large. The fact that $e^{-u} e^{-u \rho_+(J_{n+1}) / N}$ is integrable follows from Lemma \ref{lemma:rho-plus-Jucy-Murphy-eigenvalues}, as noted in the proof of Lemma \ref{lemma:general-N-limit-permutation-alternative-exploration}.
\end{proof}

% Note by linearity that
% \begin{equs}
% \sum_{\pi \in \mc{B}_{n, n}} w_T(\pi) \rho_+(\pi) = \rho_+\bigg(\sum_{\pi \in \mc{B}_{n, n}} w_T(\pi) \pi\bigg),
% \end{equs}
% where $\sum_{\pi \in \mc{B}_{n, n}} w_T(\pi) \pi$ is the partition function of the Poisson process strand diagram on $[0, T]$. We may thus rewrite 
% \begin{equs}
% \E[\Tr(M(\Gamma))] = \rho_+\bigg(\sum_{\pi \in \mc{B}_{n_1, n_1}} w_T(\pi) \pi\bigg) \cdots \rho_+\bigg(\sum_{\pi \in \mc{B}_{n_L, n_L}} w_T(\pi) \pi \bigg) \prod_{e \in R(\Gamma)} \delta^{i_e}.
% \end{equs}
% This apparently more complicated way of writing the expectation will come in handy later on, when in certain cases as $T \toinf$ the sum $\sum_{\pi \in \mc{B}_{n, n}} w_T(\pi) \pi$ does not converge, whereas $\rho_+\big(\sum_{\pi \in \mc{B}_{n, n}} w_T(\pi) \pi\big)$ {\it does converge}. This is the key result that will allow us to handle the case where $N$ is small.

\begin{prop}[Analog of Proposition \ref{prop:large-T-limit-large-N-case}]\label{prop:large-T-limit-general-N}
We have that
\begin{equs}
\lim_{T \toinf} e^{\binom{2n}{2} T - nT} \E[\rho_+\big(\pointstobrauerfn(\Sigma(T))\big)] = \sum_{\sigma, \tau : [n] \ra (n:2n]} \rho_+([\sigma ~ \tau]) \Wg_N(\sigma^{-1} \tau).
\end{equs}
\end{prop}
\begin{proof}
We argue exactly as in the proof of Proposition \ref{prop:large-T-limit-large-N-case}, except we replace the applications of Lemma \ref{lemma:limit-permutation-alternative-exploration} and Proposition \ref{prop:zero-contribution-from-unfinished-exploration-eras} with Lemma \ref{lemma:general-N-limit-permutation-alternative-exploration} and Proposition \ref{prop:general-N-zero-contribution-from-unfinished-exploration-eras}.
% First, by arguing as in the proof of Proposition \ref{prop:large-T-limit-large-N-case}, it suffices to show that 
% \begin{equs}
% \lim_{T \toinf} \rho_+\big( \E[\exp(- U_n J_n /N) \cdots \exp(-U_1 J_1 / N) \ind(U_1 + \cdots + U_n \leq T)]\big) = N^n \rho_+ ( \Wg_N).
% \end{equs}
% (This is the only part of the previous proof which does not go through if $N \leq n$.) Towards this end, using that $U_n, \ldots, U_1 \stackrel{i.i.d.}{\sim} \mrm{Exp}(1)$, we may write
% \begin{equs}
% \rho_+&\big( \E[\exp(- U_n J_n /N) \cdots \exp(-U_1 J_1 / N) \ind(U_1 + \cdots + U_n \leq T)]\big) = \\
% &\int_0^T du_n \int_0^{T-u_1} du_{n-1}\cdots \int_0^{T-(u_1 + \cdots + u_{n-1})} du_1  e^{-u_n \rho_+(\mrm{id} + J_n / N)} \cdots e^{-u_1 \rho_+ (\mrm{id} + J_1 /N)} .
% \end{equs}
\end{proof}

\begin{proof}[Proof of Theorem \ref{thm:weingarten-recovery}]
Combining Lemma \ref{lemma:rho-plus-expectation-words-unitary-brownian-motion} and Proposition \ref{prop:large-T-limit-general-N}, we have that
\begin{equs}
\lim_{T \toinf} \E\big[\Tr(B_T(\bm \Gamma))\big] = \prod_{\ell \in L} \bigg(\sum_{\sigma_\ell, \tau_\ell : [n] \ra (n:2n]} \rho_+([\sigma_\ell ~ \tau_\ell]) \Wg_N(\sigma_\ell^{-1} \tau_\ell)\bigg)_{\mbf{i}_\ell \mbf{j}_\ell} \prod_{e \in R(\Gamma)}\delta^{i_e}.
\end{equs}
By Lemma \ref{lemma:rho-plus-sum-num-comp}, the right\revision{-}hand side above may be written
\begin{equs}
\sum_{\pi = ([\sigma_\ell ~ \tau_\ell],\ell \in [L]) } \bigg(\prod_{\ell \in [L]} \Wg_N(\sigma_\ell^{-1} \tau_\ell) \bigg) N^{\numcomp(\Gamma, \pi)},
\end{equs}
as desired.
\end{proof}

\revision{To finish off this section, we use the techniques developed herein to prove a recursion for the Weingarten function for any value of $N$, which will be used in Section \ref{section:master-loop} to prove the Makeenko-Migdal/master loop/Schwinger-Dyson equation. A version of this recursion for large enough $N$ has appeared in \cite{Collins2017}. As we will see in the proof, our recursion directly follows by conditioning on the first step of our exploration process.}

\begin{prop}[Weingarten recursion]\label{prop:weingarten-recursion}
\revision{Let $N, n \geq 1$. We have that
\begin{equs}\label{eq:weingarten-recursion}
\sum_{\sigma, \tau : [n] \ra (n:2n]} \Wg_{N, \revision{n}}(\sigma \tau^{-1}) \rho_+\big([\sigma ~ \tau]\big) = &-\frac{1}{N} \sum_{j=1}^{n-1} \sum_{\sigma, \tau : [n] \ra (n:2n]} \Wg_{N, \revision{n}}(\sigma^{-1} \tau) \rho_+\big((n ~ j) [\sigma ~ \tau] \big)  \\
&+ \frac{1}{N} \sum_{j=n+1}^{2n} \sum_{\substack{\sigma, \tau : [n] \ra (n:2n] \\ \sigma(n) = \tau(n) = j} } \Wg_{N, \revision{n-1}}((\sigma^{-1} \tau)^{\revision{\downarrow}}) \rho_+\big([\sigma ~ \tau]\big).
\end{equs}}\revision{Here, to be clear, we write $\Wg_{N, n}$ (resp. $\Wg_{N, n-1}$) to denote that the Weingarten function is a function on $\symgrp_n$ (resp. $\symgrp_{n-1})$. Also, the notation $(\sigma^{-1} \tau)^{\revision{\downarrow}}$ denotes the element of $\symgrp_{n-1}$ obtained by restricting $\sigma^{-1} \tau$ to $[n-1]$, which is possible because $(\sigma^{-1} \tau)(n) = n$ by the assumption that $\sigma(n) = \tau(n)$.}
\end{prop}
\begin{proof}
\revision{By combining Propositions \ref{prop:strand-by-strand-exploration}, \ref{prop:general-N-zero-contribution-from-unfinished-exploration-eras}, and Proposition \ref{prop:large-T-limit-general-N}, we have that 
\begin{equs}\label{eq:exploration-weingarten-limit}
\sum_{\sigma, \tau : [n] \ra (n:2n]}& \Wg_N(\sigma \tau^{-1}) \rho_+\big([\sigma ~ \tau]\big) = \\
&\lim_{T \toinf}\E \big[\rho_+\big(\pointstobrauerfn(\mc{Q}(T_{n})) \big)\ind(T_n \leq T) e^{2(n-1) T_1} e^{2(n-2)(T_2 - T_1)} \cdots e^{2(n-n) (T_n - T_{n-1})}\big] .
\end{equs}
We will derive a recursion for the left hand side above by looking at the first point seen by our exploration process $\mc{Q}$. For brevity, let
\begin{equs}
f_n(T) := \E \big[\pointstobrauerfn(\mc{Q}(T_{n})) \ind(T_n \leq T) e^{2(n-1) T_1} e^{2(n-2)(T_2 - T_1)} \cdots e^{2(n-n) (T_n - T_{n-1})}\big].
\end{equs}
Let $U_1$ be the time of the first swap \revision{ or turnaroud} seen by $\mc{Q}$. Note that $U_1$ is an exponential random variable with rate $2n-1$ (since there are $n-1$ possible swaps and $n$ possible turnarounds). By conditioning on this time, we may obtain a recursion like
\begin{equs}
f_n(T) = -&\frac{1}{N} \sum_{j=1}^{n-1} \int_0^T e^{-(2n-1) u} (n ~ j) e^{2(n-1) u} f_n(T-u) du + \\
&\frac{1}{N} \sum_{j=n+1}^{2n} \int_0^T e^{-(2n-1)u} \langle n ~ j \rangle e^{2(n-1) u} f_{n-1}(j, T-u) du.
\end{equs}
Note the factor $e^{2(n-1)u}$ comes from the $e^{2(n-1) T_1}$ term. The first sum corresponds to the case that we first see a swap, and the second sum corresponds to the case that we first see a turnaround. Here, $f_{n-1}(j, T- u)$ denotes the corresponding expectation where we take out the top and bottom strand which are matched by the turnaround $\langle n ~ j \rangle$ and continue the exploration on the remaining strands. By applying $\rho_+$ on both sides and sending $T \toinf$, the desired result now follows as a consequence of \eqref{eq:exploration-weingarten-limit}.}
\end{proof}

\section{Makeenko-Migdal/Master loop/Schwinger-Dyson equations}\label{section:master-loop}

In this section, we \revision{utilize the Weingarten recursion (Proposition \ref{prop:weingarten-recursion}) proven using the Poisson point process formulation described in Section \ref{section:poisson-process-intro} and analyzed in Section \ref{section:poisson-exploration} to derive a recursion relation} (Proposition \ref{prop:word-recursion}) on expectations of products of traces of words in independent Haar-distributed Unitary matrices. We then apply this recursion to deduce the Makeenko-Migdal/Master loop/Schwinger-Dyson equations (Theorem \ref{thm:master-loop}) for Wilson loop expectations.

First, we describe the terms which will appear in our recursion. Let $\bm \Gamma = (\Gamma_1, \ldots, \Gamma_k)$ be a collection of words on $\{\lambda_1, \ldots, \lambda_L\}$. We will often refer to the $(i, j)$ location of $\bm \Gamma$, which is meant to be the $j$th letter of $\Gamma_i$. 

% \begin{remark}
% To be clear, in order to talk about specific locations of $\bm \Gamma$, we do not think of the words $\Gamma_i$ as equivalence classes made of cyclically reordered versions of $\Gamma_i$. For example, the word $\lambda_1 \lambda_2 \lambda_3$ is distinguished from $\lambda_2 \lambda_3 \lambda_1$.
% \end{remark}
% \sky{The location is not very well-defined since our words are up to cyclic equivalence, but I want to be able to talk about a specific location of a word.}

\begin{definition}[Splittings and mergers]\label{def:splittings-and-mergers}
Let $\bm \Gamma = (\Gamma_1, \ldots, \Gamma_k)$ be a collection of words on $\{\lambda_1, \ldots, \lambda_L\}$. Let $(i, j)$ be a location of $\bm \Gamma$. Define the set of positive and negative splittings $\splitting_+((i, j), \bm \Gamma)$ and $\splitting_-((i, j), \bm \Gamma)$, as well as the set of positive and negative mergers $\merge_+^U((i, j), \Gamma)$ and $\merge_-^U((i, j), \Gamma)$, as follows.

The set of positive splittings $\splitting_+((i, j), \bm \Gamma)$ is the set of collections of words $\bm \Gamma'$ obtained by splitting $\Gamma_i$ into two words as follows. Let $(i, k)$, $k \neq j$ be another location of $\Gamma_i$ which has the same letter as at location $(i, j)$. Suppose $\Gamma_i$ is of the form $A \lambda B \lambda C$, where $\lambda$ is the letter at locations $(i, j)$ and $(i, k)$. We may split $\Gamma_i$ into $\Gamma_{i, 1} = A\lambda C$ and $\Gamma_{i, 2} = B \lambda$. The set $\splitting_+((i, j), \bm \Gamma)$ is the set of all collections of words that may be obtained this way.

Similarly, the set of negative splittings $\splitting_-((i, j), \bm \Gamma)$ is the set of collections of words $\bm \Gamma'$ obtained by splitting $\Gamma_i$ into two words as follows. Let $(i, k)$, $k \neq j$ be a location of $\Gamma_i$ which has inverse of the letter at location $(i, j)$. We may write $\Gamma_i = A\lambda B \lambda^{-1}C$ or $\Gamma_i = A\lambda^{-1} B \lambda C$. In either case, we split $\Gamma_i$ into $\Gamma_{i, 1} = AC$ and $\Gamma_{i, 2} = B$. The set $\splitting_-((i, j), \bm \Gamma)$ is the set of all collections of words that may be obtained this way.

The set of positive mergers $\merge_+^U((i, j), \bm \Gamma)$ is the set of collections of words $\bm \Gamma'$ obtained by merging $\Gamma_i$ with some $\Gamma_\ell$, $\ell \neq i$, as follows. Let $(\ell, m)$ be a location which has the same letter as at location $(i, j)$. Suppose $\Gamma_i = A \lambda B$ and $\Gamma_\ell = C \lambda D$. Then $\Gamma_i, \Gamma_\ell$ are replaced by their positive merger $A \lambda D C \lambda B$.  The set $\merge_+^U((i, j), \bm \Gamma)$ is the set of all collections of words that may be obtained this way.

Similarly, the set of negative mergers $\merge_-^U((i, j), \bm \Gamma)$ is the set of collections of collections of words $\bm \Gamma'$ obtained by merging $\Gamma_i$ with some $\Gamma_\ell$, $\ell \neq i$, as follows. Let $(\ell, m)$ be a location which has the inverse of the letter at location $(i, j)$. Suppose $\Gamma_i = A \lambda B$ and $\Gamma_\ell = C \lambda^{-1} D$. Then $\Gamma_i, \Gamma_\ell$ are replaced by their negative merger $A D C B$. The set $\merge_-^U((i, j), \bm \Gamma)$ is the set of all collections of words that may be obtained this way. 
\end{definition}

\begin{remark}
The notation $\mbb{M}^U$ is from \cite{shen2022new}. Later on, when we consider other groups in Section \ref{section:orthogonal-and-symplectic}, there will be a larger set of mergers to consider, which we will denote by $\mbb{M}$.
\end{remark}

In the following, let $\tr(U(\bm \Gamma)) = \prod_{i \in [k]} \tr(U(\Gamma_i))$, where $U(\Gamma_i)$ is obtained by substituting into $\Gamma_i$ an independent Haar-distributed Unitary matrix for each letter $\{\lambda_1, \ldots, \lambda_L\}$. Note that in contrast to previous results, we are using the normalized trace here, which we find to be more natural for stating the recursion.

\begin{prop}[Single-location word recursion]\label{prop:word-recursion}
Let $\bm \Gamma = (\Gamma_1, \ldots, \Gamma_k)$ be a collection of words on $\{\lambda_1, \ldots, \lambda_L\}$. For any location $(i, j)$ of $\Gamma$, we have that 
\begin{equs}
\E[\tr(U(\bm \Gamma))] = &- \sum_{\bm\Gamma' \in \splitting_+((i, j), \bm \Gamma)} \E[\tr(U(\bm \Gamma'))] + \sum_{\bm \Gamma' \in \splitting_-((i, j), \bm \Gamma)} \E[\tr(U(\bm \Gamma'))] \\
&- \frac{1}{N^2} \sum_{\bm \Gamma' \in \merge_+^U((i, j), \bm \Gamma)} \E[\tr(U(\bm \Gamma'))] + \frac{1}{N^2} \sum_{\bm \Gamma' \in \merge_-^U((i, j), \bm \Gamma)} \E[\tr(U(\bm  \Gamma'))].
\end{equs}
\end{prop}
\begin{proof}
We can assume that $\bm \Gamma$ is balanced, otherwise both sides are zero. Moreover, without loss of generality, take $(i, j) = (1, 1)$, so that we look at the first letter of $\Gamma_1$. Let $\lambda \in \{\lambda_1, \ldots, \lambda_L\}$ be this letter.
% Recall from Theorem \ref{thm:weingarten} that $\E[\Tr(U(\bm \Gamma))]$ is equal to a sum over pairs of matchings of strand diagrams, weighted by the Weingarten function applied to each pair, as well as $N$ raised to the number of components of the resulting strand diagram. For each strand diagram corresponding to a letter $\lambda' \neq \lambda$, fix a pair of matchings $\sigma_{\lambda'}, \tau_{\lambda'}$. We apply our strand-by-strand Poisson process exploration from Section \ref{section:poisson-exploration} to the strand diagram corresponding to $\lambda$, but stop at the first time we see any point in the first exploration era. This will result in the claimed recursion.
Let $n$ be the number of times that $\lambda$ appears in $\bm \Gamma$, so that the portion of the strand diagram corresponding to $\lambda$ has $n$ right-directed strands and $n$ left-directed strands. We suppose that the top strand (which we imagine as labeled by $n$) corresponds to the first letter of $\Gamma_1$.
We will show that the claimed result follows directly from the Weingarten recursion (Proposition \ref{prop:weingarten-recursion}). By the Weingarten calculus (Theorem \ref{thm:weingarten}) and Lemma \ref{lemma:rho-plus-sum-num-comp}, we have that
\begin{equs}
\E[\Tr(U(\Gamma))] =\sum_{(\sigma_\ell, \tau_\ell),  \ell \in [L]} \prod_{\ell \in [L]} \Wg_N(\sigma_\ell^{-1} \tau_\ell) \prod_{\ell \in [L]} \rho_+\big([\sigma_\ell, \tau_\ell]\big)_{\mbf{i}_\ell, \mbf{j}_\ell} \prod_{e \in R(\bm \Gamma)} \delta^{i_e}.
\end{equs}
Inserting the recursion \eqref{eq:weingarten-recursion} into the sum over $(\sigma_1, \tau_1)$, and then using Lemma \ref{lemma:rho-plus-sum-num-comp} to replace the sum over matrix entries by terms of the form $N^{\numcomp}$, we obtain
\begin{equs}\label{eq:recursion-proof-intermediate}
\E[\Tr(U(\Gamma))] = - \frac{1}{N} \sum_{j=1}^{n-1} I(j) + \frac{1}{N} \sum_{j=n+1}^{2n} I(j),
\end{equs}
where $I(j)$ is defined as follows. For $j \in [n-1]$, we have that
\begin{equs}
I(j) := \sum_{(\sigma_\ell, \tau_\ell),  \ell \in [L]} \prod_{\ell \in [L]} \Wg_{N, \revision{n}}(\sigma_\ell^{-1} \tau_\ell) N^{\numcomp(\bm \Gamma, \bm \pi')},
\end{equs}
where $\bm \pi' = ((n ~ j) [\sigma_1 ~ \tau_1], [\sigma_2 ~ \tau_2], \ldots, [\sigma_L ~ \tau_L])$ (here we hide the dependence of $\bm \pi'$ on $j$)
and for $j \in (n:2n]$, we have that
\begin{equs}
I(j) := \sum_{\substack{(\sigma_\ell, \tau_\ell),  \ell \in [L] \\ \sigma_1(n) = \tau_1(n) = j} } \prod_{\ell \in [L]} \Wg_{N, \revision{n-1}}((\sigma_\ell^{-1} \tau_\ell)^{\revision{\downarrow}}) N^{\numcomp(\bm \Gamma, \bm \pi)}, 
\end{equs}
where $\bm \pi = ([\sigma_\ell ~ \tau_\ell], \ell \in [L])$.
We next show that $I(j)$ for $j \in [n-1]$ gives the positive splittings and mergings, while $I(j)$ for $j \in (n:2n]$ gives the negative splittings and mergings. 

First, suppose that $j \in [n-1]$. In this case, the only difference between $I(j)$ and the sum appearing in the Weingarten calculus (Theorem \ref{thm:weingarten}) is that the collection of matchings $\bm \pi'$ is slightly different. In particular, the first matching is $(n ~ j) [\sigma_1 ~ \tau_1]$ and not $[\sigma_1 ~ \tau_1]$. Diagramatically, the effect of the term $(n ~ j)$ is to swap the incoming red exterior connections at strands $n$ and $j$. See Figure \ref{figure:strand_diagram_recursion} for a visualization.
\begin{figure}[h]%
    \centering
    {{\includegraphics[page=1, width=7cm]{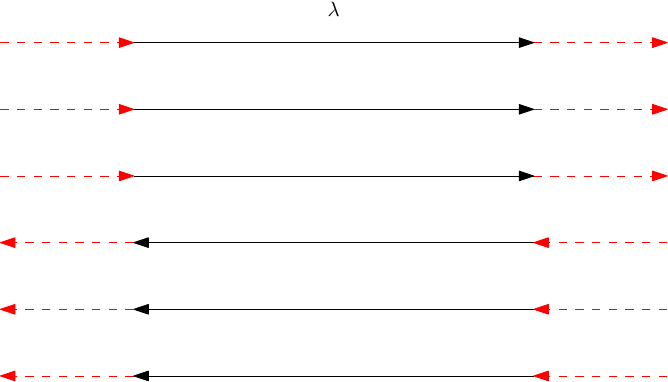} }}%
    \qquad
    {{\includegraphics[page=2, width=7cm]{figures/strand_diagram_recursion.pdf} }}%
    \caption{Left: the portion of the strand diagram of $\bm \Gamma$ corresponding to $\lambda$. Right: the effect of the swap $(n ~ j)$ can be visualized as swapping the two incoming red exterior connections at strands $n$ and $j$.}\label{figure:strand_diagram_recursion}
\end{figure}

Mathematically, Figure \ref{figure:strand_diagram_recursion} corresponds to the following identity for $\numcomp$: there is some modified collection of words $\bm \Gamma'$ such that
\begin{equs}\label{eq:numcomp-identity-recursion-proof}
\numcomp(\bm \Gamma, \bm \pi') = \numcomp(\bm \Gamma', \bm \pi), ~~ \text{ where $\bm \pi = ([\sigma_\ell ~ \tau_\ell], \ell \in [L])$.}
\end{equs}
The modified collection of words $\bm \Gamma'$ is defined as follows. There are two cases to consider. First, if the swap $(n ~ j)$ involves two strands which are in the same word, then since we are assuming that $\lambda$ is the first letter of $\Gamma_1$, this means that the first word may be written $\Gamma_1 = \lambda \Gamma_{1, 1} \lambda \Gamma_{1, 2}$ for some words $\Gamma_{1, 1}, \Gamma_{1, 2}$. We then define the modified collection $\bm \Gamma' = (\lambda \Gamma_{1, 1}, \lambda \Gamma_{1, 2}, \Gamma_2, \ldots, \Gamma_k)$, which ensures that the identity \eqref{eq:numcomp-identity-recursion-proof} is satisfied. Observe that in this case, $\bm \Gamma' \in \splitting_+((1, 1), \bm \Gamma)$.

The other case is if $(n ~ j)$ involves two strands which are originally in different words. For simplicity of notation, suppose that $\Gamma_2$ is the word that contains strand $j$. Then we may write $\Gamma_1 = \lambda \Gamma_1'$, $\Gamma_2 = \Gamma_{2, 1} \lambda \Gamma_{2, 2}$ for some words $\Gamma_1', \Gamma_{2, 1}, \Gamma_{2, 2}$. In this case, we define $\bm \Gamma' = (\lambda \Gamma_1' \lambda \Gamma_{2, 2} \Gamma_{2, 1}, \Gamma_3, \ldots, \Gamma_k)$, which ensures that the identity \eqref{eq:numcomp-identity-recursion-proof} is satisfied. Observe that in this case, $\bm \Gamma' \in \merge_+^U((1, 1), \bm \Gamma)$. To summarize, what we have shown so far is that
\begin{equs}
-\frac{1}{N} \sum_{j=1}^{n-1} I(j) = -\frac{1}{N} \sum_{\bm\Gamma' \in \splitting_+((1, 1), \bm \Gamma)} \E[\Tr(U(\bm \Gamma'))] - \frac{1}{N}\sum_{\bm \Gamma' \in \merge_+^U((1, 1), \bm \Gamma)} \E[\Tr(U(\bm \Gamma'))].
\end{equs}
Since we used the un-normalized trace in \eqref{eq:recursion-proof-intermediate}, we should multiply both sides above by $N^{-k}$, which results in the identity
\begin{equs}\label{eq:positive-operation-intermediate}
-\frac{1}{N} \sum_{j=1}^{n-1} N^{-k} I(j) = -\sum_{\bm\Gamma' \in \splitting_+((1, 1), \bm \Gamma)} \E[\tr(U(\bm \Gamma'))] - \frac{1}{N^2}\sum_{\bm \Gamma' \in \merge_+^U((1, 1), \bm \Gamma)} \E[\tr(U(\bm \Gamma'))],
\end{equs}
(note that splitting results in one additional loop, while merger results in one less loop, which is why the factor of $\frac{1}{N}$ is absorbed in the splitting case, while there is an additional factor of $N$ in the merger case). Observe that the right\revision{-}hand side above is precisely two of the terms which appear in the proposition statement.

We now move on to the case $j \in (n:2n]$. In this case, the only difference between $I(j)$ and the sum appearing in the Weingarten calculus (Theorem \ref{thm:weingarten}) is the restriction that $\sigma_1(n) = \tau_1(n) = j$. Diagramatically, this restriction can be interpreted as deleting the the strands $n$ and $j$, and adding in a dashed red edge which joins the left vertices of strands $n$ and $j$, and another dashed red edge which joins the right vertices of strands $n$ and $j$. See Figure \ref{figure:strand_diagram_recursion_3} for a visualization.

\begin{figure}[h]%
    \centering
    \includegraphics[page=3, width=7cm]{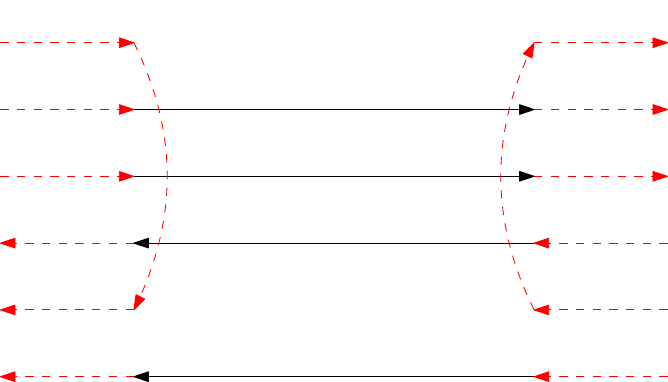} \caption{The restriction $\sigma_1(n) = \tau_1(n) = j$ may be interpreted as deleting the strands $n$ and $j$, and adding in new dashed red edges as indicated.}\label{figure:strand_diagram_recursion_3}
\end{figure}

Mathematically, Figure \ref{figure:strand_diagram_recursion_3} corresponds to the following identity. For any $\bm \pi = ([\sigma_\ell ~ \tau_\ell], \ell \in [L])$ such that $\sigma_1(n) = \tau_1(n) = j$, there exists a modified collection of words $\bm \Gamma'$ such that
\begin{equs}\label{eq:numcomp-identity-recursion-proof-turnaround}
\numcomp(\bm \Gamma, \bm \pi) = \numcomp(\bm \Gamma', \bm \pi'),
\end{equs}
where here $\bm \pi' = ([\sigma_1' ~ \tau_1'], [\sigma_2 ~ \tau_2], \ldots, [\sigma_L ~ \tau_L])$, and $\sigma_1', \tau_1'$ are the bijections $[n-1] : (n-1: 2(n-1)]$ obtained from $\sigma, \tau$ by deleting the vertices of strands $n$ and $j$ (and then perhaps relabeling the remaining strands). 

The modified collection of words $\bm \Gamma'$ is defined as follows. As before, there are two cases to consider. First, suppose that strands $n$ and $j$ are in the same word, which is necessarily $\Gamma_1$ by assumption. This means we may write $\Gamma_1 = \lambda \Gamma_{1, 1} \lambda^{-1} \Gamma_{1, 2}$. We define the modified collection of words by $\bm \Gamma' = (\Gamma_{1, 2}, \Gamma_{1, 1}, \Gamma_2, \ldots, \Gamma_k)$, which ensures that \eqref{eq:numcomp-identity-recursion-proof-turnaround} is satisfied. Note that in this case, $\bm \Gamma' \in \splitting_-((1, 1), \bm \Gamma)$.

The other case is if strands $n$ and $j$ are in different words. For notational simplicity, suppose that strand $j$ is in $\Gamma_2$, so that we may write $\Gamma_1 = \lambda \Gamma_1'$, $\Gamma_2 = \Gamma_{2, 1} \lambda^{-1} \Gamma_{2, 2}$. We define the modified collection of words by $\bm \Gamma' = (\Gamma_1' \Gamma_{2, 2} \Gamma_{2, 1}, \Gamma_3, \ldots, \Gamma_k)$, which ensures that \eqref{eq:numcomp-identity-recursion-proof-turnaround} is satisfied. Note that in this case, $\bm \Gamma' \in \merge_-^U((1, 1), \bm \Gamma)$. In summary, we have shown that
\begin{equs}
\frac{1}{N} \sum_{j=n+1}^{2n} I(j) = \frac{1}{N} \sum_{\bm \Gamma' \in \splitting_-((1, 1), \bm \Gamma)} \E[\Tr(U(\bm \Gamma'))] + \frac{1}{N} \sum_{\bm \Gamma' \in \merge_-^U((1, 1), \bm \Gamma)} \E[\Tr(U(\bm \Gamma'))].
\end{equs}
The desired result now follows by converting the above to using the normalized trace, and then combining this with \eqref{eq:positive-operation-intermediate} and \eqref{eq:recursion-proof-intermediate}.
\end{proof}

Next, we apply the loop recursion Proposition \ref{prop:word-recursion} to obtain a recursion for Wilson loop expectations. 
% In contrast to the notation of Section \ref{section:intro}, we denote collections of loops by $s$ instead of $\mc{L}$, and we refer to $s$ as a string. 
Recall the notation that $W_s(Q) = \prod_{k \in [n]} \tr(Q_{\ell_k})$.

\begin{notation}
Given a string $s = (\ell_1, \ldots, \ell_n)$, let $\phi(s) := \langle W_s\rangle_{\Lambda, \beta}$, where $\langle \cdot \rangle_{\Lambda, \beta}$ denotes expectation with respect to the lattice Yang-Mills measure defined in \eqref{eqn::latticeym}. We omit the dependence of $\phi$ on $\Lambda, \beta, N$.
\end{notation}

Note that Definition \ref{def:splittings-and-mergers} specializes to the case of loops on a lattice: given a string $s$, we have the sets of positive/negative splittings/mergers $\splitting_{\pm}((i, j), s)$ and $\merge_{\pm}^U((i, j), s)$.

\begin{remark}\label{remark:schwinger-dyson}
We remark that our definition of the set of splittings and mergers is slightly different than what appears in \cite{Chatterjee2019a, shen2022new}. In our definition, we consider all possible splittings/mergers that involve the specific location $(i, j)$, whereas in the earlier works, the authors consider any splitting/merger that involves any two locations of the string which correspond to the same lattice edge.
\end{remark}

We need to define another type of string operation which appears for lattice Yang-Mills.

\begin{definition}[Deformations]
Let $s = (\ell_1, \ldots, \ell_n)$ be a string. Let $(i, j)$ be a location in $s$. We define the sets of positive and negative deformations $\deform_+((i, j), s)$ and $\deform_-((i, j), s)$ as follows. 

The set of positive deformations $\deform_+((i, j), s)$ is the set of all possible strings which can be obtained by a positive merger between $s$ at location $(i, j)$ and some oriented plaquette $p \in \mc{P}$. The set of negative deformations $\deform_-((i, j), s)$ is the set of possible strings which can be obtained by a negative merger between $s$ at location $(i, j)$ and some oriented plaquette $p \in \mc{P}$.

Let $e$ be the oriented edge of $\Lambda$ that is at location $(i, j)$ in $s$. Let $p \in \mc{P}$. In order for their to exist a positive merger between $s$ and $p$, note that $p$ must contain $e$. In this case, we denote by $s \oplus_{(i, j)} p$ to be the positive merger of $s$ and $p$ at location $(i, j)$. Similarly, in order for their to exist a negative merger between $s$ and $p$, note that $p$ must contain $-e$. In this case, we denote by $s \ominus_{(i, j)} p$ to be the negative merger of $s$ and $p$ at location $(i, j)$.
\end{definition}

Let $p \succ e$ denote that the plaquette $p$ contains the edge $e$. Note then that (here $e$ is the edge at location $(i, j)$ of $s$)
\begin{equs}\label{eq:defomration-sets}
\deform_+((i, j), s) &= \{s \oplus_{(i, j)} p : p \in \mc{P}, p \succ e\} \\
\deform_-((i, j), s) &= \{s \ominus_{(i, j)} p : p \in \mc{P}, p \succ -e\}.
\end{equs}

% \begin{definition}[Splittings and mergers]
% Let $s = (\ell_1, \ldots, \ell_n)$ be a string. Let $k \in [n]$ and let $i \in [|\ell_k|]$. Define $\splitting_+((i, j), s), \splitting_-((i, j), s)$, the sets of positive and negative splittings of $s$ at location $(i, j)$, as follows. 

% Next, we define the $\merge_+((i, j), s)$, $\merge_-((i, j), s)$, the sets of positive and negative mergers of $s$ at location $(i, j)$, as follows. 
% \end{definition}

\begin{theorem}[Single-location Makeenko-Migdal/Master loop/Schwinger-Dyson equation]\label{thm:master-loop}
Let $s = (\ell_1, \ldots, \ell_n)$ be a string. Let $(i, j)$ be a location in $s$. We have that
\begin{equs}
\phi(s) = &-\sum_{s' \in \splitting_+((i, j), s)} \phi(s') + \sum_{s' \in \splitting_-((i, j), s)} \phi(s') - \frac{1}{N^2} \sum_{s' \in \merge_+^U((i, j), s)} \phi(s') + \frac{1}{N^2} \sum_{s' \in \merge_-^U((i, j), s)} \phi(s') \\
&- 
\beta \sum_{s' \in \deform_+((i, j), s)} \phi(s') + \beta \sum_{s' \in \deform_-((i, j), s)} \phi(s').
\end{equs}
\end{theorem}
\begin{remark}\label{remark:schwiner-dyson-2}
We re-emphasize here that the above recursion is slightly more general than previous literature~\cite{Chatterjee2019a, chatterjee2016, shen2022new, AN2023}, because we defined the string operations appearing on the right\revision{-}hand side of the equation in a slightly more restrictive manner -- recall Remark \ref{remark:schwinger-dyson}. In particular, the right\revision{-}hand side of our formula depends on $j$ while the `unsymmetrized' version stated in~\cite[Theorem 8.1]{Chatterjee2019a} does not. The Makeenko-Migdal/Master loop/Schwinger-Dyson equation of the previous works may be recovered from our equation by summing over all locations of $s$. 

Also, recall Remark \ref{remark:beta-vs-two-beta} that our scaling is so that $\beta$ in our paper corresponds to $2\beta$ in previous papers. This explains why $\beta$ appears in the above recursion, while $\beta/2$ appears in \cite[Equation (1.7)]{shen2022new}. 

\revision{For the sake of comparison, we state here the master loop equation for $\unitary(N)$ derived in \cite{shen2022new, AN2023}, in our notation. Fix an edge $e$ of the lattice, and let $\mathsf{L} = \{(i_1, j_1), \ldots, (i_m, j_m)\}$ be the set of locations of $s$ which contain the edge $e$. Let $\splitting_+(e, s) = \bigsqcup_{(i, j) \in \mathsf{L}} \splitting_+((i, j), s)$, and similarly for the other loop operations. Then
\begin{equs}
|\mathsf{L}|\phi(s) = &-\sum_{s' \in \splitting_+(e, s)} \phi(s') + \sum_{s' \in \splitting_-(e, s)} \phi(s') - \frac{1}{N^2} \sum_{s' \in \merge_+^U(e, s)} \phi(s') + \frac{1}{N^2} \sum_{s' \in \merge_-^U(e, s)} \phi(s') \\
&- 
\beta \sum_{s' \in \deform_+(e, s)} \phi(s') + \beta \sum_{s' \in \deform_-(e, s)} \phi(s').
\end{equs}
It is clear that this follows directly from Theorem \ref{thm:master-loop} upon summing over all locations $(i, j)$ which contain $e$.}
\end{remark}
\begin{proof}
Recall from equation \eqref{eqn::latticeymexpanded} that
\begin{equs}
\phi(s) = Z_{\Lambda, \beta}^{-1} \sum_{K : \mc{P} \ra \N} \frac{(N\beta)^{K}}{K!} \int W_s(Q) \prod_{p \in \mc{P}} \Tr(Q_p)^{K(p)} \prod_{e \in E_\Lambda} dQ_e.
\end{equs}
For brevity, let
\begin{equs}
I(s, K) := \int W_s(Q) \prod_{p \in \mc{P}} \Tr(Q_p)^{K(p)} \prod_{e \in E_\Lambda} dQ_e.
\end{equs}
Fix $K : \mc{P} \ra N$. It may help to keep in mind that $K(p)$ counts the number of copies of $p$ that are present. Before we apply Proposition \ref{prop:word-recursion}, let us set some notation. Let $e$ be the oriented edge of $\Lambda$ that is traversed at location $(i, j)$ in the string $s$. Recall that $p \succ e$ means that $p$ contains $e$, and $p \succ -e$ means that $p$ contains $e$ with the opposite orientation. Recall also that if $p \succ e$ or $p \succ -e$, let $s \oplus_{(i, j)} p$ and $s \ominus_{(i, j)} p$ be the positive and negative deformations of $s$ by $p$ at location $(i, j)$. For $p \in \mc{P}$, let $\delta_p : \mc{P} \ra \N$ be the delta function at $p$. Now applying the word recursion Proposition \ref{prop:word-recursion}, we have that
\begin{equs}
I(s, K) = &-\sum_{s' \in \splitting_+((i, j), s)} I(s', K) +  \sum_{s' \in \splitting_-((i, j), s)} I(s', K) \\
&- \frac{1}{N^2} \sum_{s' \in \merge_+^U((i, j), s)} I(s', K) + \frac{1}{N^2} \sum_{s' \in \merge_-^U((i, j), s)} I(s', K) \\
&-\frac{1}{N} \sum_{\substack{p \in \mc{P} \\ p \succ e}}  K(p) I(s \oplus_{(i, j)} p, K - \delta_p) + \frac{1}{N} \sum_{\substack{p \in \mc{P} \\ p \succ -e}} K(p) I(s \ominus_{(i, j)} p, K - \delta_p).
\end{equs}
(Here, the factor of $K(p)$ arising in the last two terms arises because there are $K(p)$ copies of the plaquette $p$ which can possibly be used to deform $s$.) From this, we obtain 
% (note that in the splitting terms, the $1/N$ factor gets absorbed due to the fact that $s'$ has one more loop than $s$, while in the merging terms, there is an extra $1/N$ factor because $s'$ has one less loop than $s$)
\begin{equs}
\phi(s) = &-\sum_{s' \in \splitting_+((i, j), s)} \phi(s) + \sum_{s' \in \splitting_-((i, j), s)} \phi(s) - \frac{1}{N^2} \sum_{s' \in \merge_+^U((i, j), s)} \phi(s') + \frac{1}{N^2} \sum_{s' \in \merge_-^U((i, j), s)} \phi(s') \\
&+ D_1 + D_2 ,
\end{equs}
where
\begin{equs}
D_1 &:=  -Z_{\Lambda, \beta}^{-1} \frac{1}{N} \sum_{\substack{p \in \mc{P} \\ p \succ e}} \sum_{\substack{K : \mc{P} \ra \N \\ K(p) \geq 1}} \frac{(N\beta)^K}{K!} K(p) I(s \oplus_{(i, j)} p, K - \delta_p), \\
D_2 &:= Z_{\Lambda, \beta}^{-1} \frac{1}{N} \sum_{\substack{p \in \mc{P} \\ p \succ -e}} \sum_{\substack{K : \mc{P} \ra \N \\ K(p) \geq 1}} \frac{(N\beta)^K}{K!} K(p) I(s \ominus_{(i, j)} p, K - \delta_p)
\end{equs}
Observe that we may write (by changing variables $K \mapsto K - \delta_p$ and then recalling \eqref{eq:defomration-sets})
\begin{equs}
D_1 &= -Z_{\Lambda, \beta}^{-1} \frac{1}{N} (N\beta) \sum_{\substack{p \in \mc{P} \\ p \succ e}} \sum_{\substack{K : \mc{P} \ra \N}} \frac{(N\beta)^K}{K!} I(s \oplus_{(i, j)} p, K) \\
&= -\beta \sum_{\substack{p \in \mc{P} \\ p \succ e}} \phi(s \oplus_{(i, j)} p) = -\beta \sum_{s' \in \deform_+((i, j), s)} \phi(s'),
\end{equs}
and similarly
\begin{equs}
D_2 = \beta \sum_{\substack{p \in \mc{P} \\ p \succ -e}} \phi(s \ominus_{(i, j)} p) = \beta \sum_{s' \in \deform_-((i, j), s)} \phi(s') .
\end{equs}
% \begin{equs}
% D_2 = Z_{\Lambda, \beta}^{-1} \frac{1}{N} \frac{N\beta}{2}  \sum_{\substack{p \in \mc{P} \\ p \succ -e}} \sum_{k_1, k_2 : \mc{P} \ra \N} \bigg(\frac{N\beta}{2} \bigg)^{|k_1| + |k_2|} \bigg(\prod_{p' \in \mc{P}} \frac{1}{k_1(p')! k_2(p')!}\bigg)\big(I(s \oplus p, k_1, k_2) - I(s \ominus p, k_1, k_2)\big).
% \end{equs} 
% From this, we see that (recall \eqref{eq:defomration-sets})
% \begin{equs}
% D_1 + D_2 = -\beta\sum_{s' \in \deform_+((i, j), s)} \phi(s') + \beta\sum_{s' \in \deform_-((i, j), s)} \phi(s').
% \end{equs}
The desired result now follows.
\end{proof}

\section{Other groups} \label{sec::othergroups}
In this section, we adapt our results to the cases $G = \orthogonal(N), \SphN, \SUN, \SON$. In Section \ref{section:orthogonal-and-symplectic}, we address the cases $G = \orthogonal(N), \SphN$, and in Section \ref{section:SUN-and-SON}, we address the cases $G = \SUN, \SON$. Define the matrix $J$ by
\begin{equs}\label{eq:J}
J := \begin{pmatrix} 0 & I_{N/2} \\ -I_{N/2} & 0 \end{pmatrix}.
\end{equs}
We quickly recall the definitions of the various groups.
\begin{equs}
\orthogonal(N) &:= \{O \in \mrm{GL}(N, R) : O^T O = I_N \} \\
\SphN &:= \{S \in \UN : S^T J S = J \} \\
\SUN &:= \{U \in \UN : \mrm{det}(U) = 1\} \\
\SON &:= \{O \in \SON : \mrm{det}(O) = 1\}.
\end{equs}

\begin{notation}
Let $\bm \Gamma = (\Gamma_1, \ldots, \Gamma_M)$ be a collection of words on $\{\lambda_1, \ldots, \lambda_L\}$. Given a compact Lie group $G$, we will denote $\Tr(G(\bm \Gamma)) = \Tr(G(\Gamma_1)) \cdots \Tr(G(\Gamma_M))$, where $G(\Gamma_i)$ is obtained by substituting an independent Haar-distributed element of $G$ for each of the letters $\{\lambda_1, \ldots, \lambda_L\}$.
\end{notation}

\subsection{Orthogonal and Symplectic}\label{section:orthogonal-and-symplectic}

In this section, we adapt our previous results to $G = \orthogonal(N), \SphN$. These two cases are at times very similar, and thus we choose to place them in the same section. However, they are also at times very different, which prevents us from handling the two cases completely simultaneously -- there are certain parts which require special attention in the $\orthogonal(N)$ case, and certain parts in the $\SphN$ case.

\begin{notation}
Let $n \geq 1$ be even. In this section, we will denote matchings on $[n]$ (i.e. partitions of $[n]$ into two-element sets) by $\pi, \pi', \pi''$, etc., and often write $\pi : [n] \ra [n]$ (since a matching can be understood as an involution of $[n]$ without any fixed points). In our existing notation, $\pi \in \mc{M}(n/2)$.
\end{notation}

\begin{remark}
In Section \ref{section:orthogonal-and-symplectic}, matchings $\pi$ should be visualized differently than before. Before, a matching $\pi \in \mc{M}(2n)$ was visualized as giving pairings between $2n$ vertices, where we drew $n$ vertices each on the left and right sides. Now, we should only think of a matching $\pi : [n] \ra [n]$ as giving a pairing of $n$ vertices, which are either all drawn on the left, or all on the right. Ultimately, this is due to the fact that matchings $\pi : [n] \ra [n]$ in this section play the role of bijections $\sigma, \tau : [n] \ra (n:2n]$ in previous sections (see e.g. Section \ref{section:strand-diagrams-and-weingarten-caculus}). In the Unitary case, there was a natural division between the top $n$ strands and bottom $n$ strands, and thus also between the top $n$ vertices and bottom $n$ vertices on each side. The bijections $\sigma, \tau$ encoded pairings between the left vertices and between the right vertices which had the special feature that top vertices were always paired to bottom vertices (and vice versa). In the Orthogonal and Symplectic cases, we will see that there is no longer a division between top and bottom, and so the matching $\pi$ encodes an arbitrary pairing of the left vertices (or the right vertices).
\end{remark}

\subsubsection{Orthogonal surface sums}

First, we discuss the surface sums that arise in the $\orthogonal(N)$ case. 
% Because these two cases are very similar, we handle them simultaneously.
We begin by introducing the needed setup in order to state the analog of Theorem \ref{thm:weingarten} (the Unitary Weingarten calculus) for $\orthogonal(N)$.

\begin{definition}[Unoriented-balanced collection of words]
Let $\bm \Gamma = (\Gamma_1, \ldots, \Gamma_M)$ be a collection of words on letters $\{\lambda_1, \ldots, \lambda_L\}$. For $\ell \in [L]$, let $n_\ell$ be the total number of times $\lambda_\ell$ or $\lambda_\ell^{-1}$ occurs in $\bm \Gamma$. We say that $\bm \Gamma$ is unoriented-balanced if $n_\ell$ is even for each $\ell \in [L]$.
\end{definition}

\begin{remark}
By $O \mapsto -O$ distributional symmetry of Haar-distributed $\orthogonal(N)$ matrices, if $\bm \Gamma$ is not unoriented-balanced then $\E[\Tr(O(\bm \Gamma))] = 0$. Thus when computing $\E[\Tr(O(\bm \Gamma))]$, we may assume $\bm \Gamma$ is unoriented-balanced.
\end{remark}

\begin{definition}
Let $n \geq 1$ be even. Let $\pi, \pi' : [n] \ra [n]$ be matchings. Visually, we will think of $\pi, \pi'$ as giving left and right matchings, as in the Figure \ref{figure:O_N_left_right_matching_example}. This defines an element of the Brauer algebra $\mc{B}_n$, which we denote by $[\pi ~ \pi']$.
\begin{figure}[ht!]
    \centering
\includegraphics[width=0.3\textwidth]{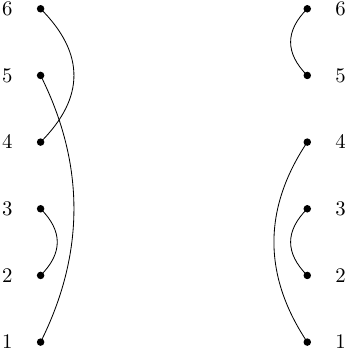}
    \caption{Here, $n = 6$. The left matching is $\pi  = \{\{4, 6\}, \{1, 5\}, \{2, 3\}\}$. The right matching is $\pi' = \{\{5, 6\}, \{1, 4\}, \{2,3\}\}$. The left and right matchings together define an element $[\pi ~ \pi']$ of $\mc{B}_n$.}\label{figure:O_N_left_right_matching_example}
\end{figure}

Let $\cycles(\pi, \pi')$ be the number of connected components in the graph one obtains by adding in the strands connecting the left and right vertices - see Figure \ref{figure:O_N_left_right_matching_example_2} for an example.
\begin{figure}[ht!]
    \centering
\includegraphics[page=2, width=0.3\textwidth]{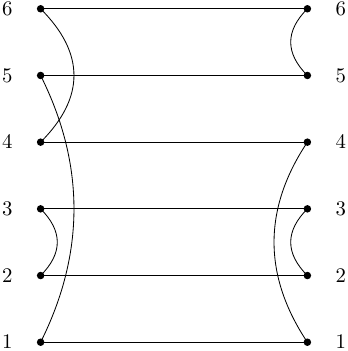}
    \caption{For the left and right matchings $\pi, \pi'$ from Figure \ref{figure:O_N_left_right_matching_example}, there are two connected components in the graph obtained by adding the strands, and thus $\cycles(\pi, \pi') = 2$.}\label{figure:O_N_left_right_matching_example_2}
\end{figure}
\end{definition}

\begin{definition}
Let $n \geq 1$ be even. Given left and right matchings $\pi, \pi' : [n] \ra [n]$, the face profile $\ell(\pi, \pi')$ is the partition of $n$ induced by the cycles of $\pi \pi'$.
\end{definition}

We note that all parts of $\ell(\pi, \pi')$ are even, and thus $\frac{1}{2} \ell(\pi, \pi')$ is a partition of $\frac{n}{2}$. For the matchings $\pi, \pi'$ in Figures \ref{figure:O_N_left_right_matching_example} and \ref{figure:O_N_left_right_matching_example_2}, the face profile $\ell(\pi, \pi') = \{4, 2\}$. Note also that $\cycles(\pi, \pi')$ is exactly the number of parts of $\ell(\pi, \pi')$.

\begin{definition}[Orthogonal Weingarten function]\label{def:orthogonal-weingarten}
Let $\zeta \in \C$. Let $n \geq 1$ be even. We define the Orthogonal Weingarten function $\Wg^{\mrm{O}}_{\zeta, n}$ as follows. The input is a pair of matchings $\pi, \pi' : [n] \ra [n]$, and the output is a number $\Wg^{\mrm{O}}_{\zeta, n}(\pi, \pi') \in \C$. First, define the Gram matrix
\begin{equs}
\mbf{G}^{\mrm{O}}_{\zeta, n}(\pi, \pi') := \zeta^{\cycles(\pi, \pi')}, ~~ \pi, \pi': [n] \ra [n] \text{ matchings}.
\end{equs}
We define $\Wg_{\zeta, n}^{\mrm{O}}$ to be the pseudo-inverse of $\mbf{G} = \mbf{G}^{\mrm{O}}_{\zeta, n}(\pi, \pi') $, that is the symmetric matrix $W$ which satisfies
\begin{equs}
W \mbf{G} W = W \text{ and } \mbf{G} W \mbf{G} = \mbf{G}.
\end{equs}
We typically omit the $n$ variable and write $\Wg_\zeta^\mrm{O}$. The normalized Orthogonal Weingarten function is defined to be
\begin{equs}
\ovl{\Wg}^{\mrm{O}}_\zeta(\pi, \pi') = \zeta^{n - \cycles(\pi, \pi')} \Wg^{\mrm{O}}_\zeta(\pi, \pi').
\end{equs}
\end{definition}

\begin{remark}
From \cite[Theorem 3.13]{collins2006integration}, the normalized Orthogonal Weingarten function has the following large-$N$ asymptotics:
\begin{equs}
\lim_{N \toinf} \ovl{\Wg}^{\mrm{O}}_N(\pi, \pi') = \prod_{a \in \frac{1}{2}\ell(\pi, \pi')} (-1)^{a - 1} \mrm{Cat}(a-1),
\end{equs}
where $\mrm{Cat}(k)$ is the $k$th Catalan number as in \eqref{eq:mobius}, and the product is over all parts in the face profile of $\frac{1}{2} \ell(\pi, \pi')$ (which recall is a partition of $\frac{n}{2}$).

In fact, the proof of the cited theorem extends without change to a general complex parameter $\zeta \toinf$, and thus we have that
\begin{equs}
\lim_{\zeta \toinf} \ovl{\Wg}^{\mrm{O}}_\zeta(\pi, \pi') = \prod_{a \in \frac{1}{2}\ell(\pi, \pi')} (-1)^{a - 1} \mrm{Cat}(a- 1).
\end{equs}
\end{remark}

We state the following lemma which says that the Orthogonal Weingarten function is a function of the face profile of $(\pi, \pi')$. It essentially follows from \cite{matsumoto2013weingarten}, although we haven't found a precise statement in the literature. Thus for the reader's convenience, we give more detail as to why the lemma is true in Appendix \ref{appendix:orthogonal-weingarten-jucys-murphy-relation}.

\begin{lemma}\label{lemma:orthogonal-weingarten-function-face-profile}
The Orthogonal Weingarten function $\Wg_\zeta^{\mrm{O}}(\pi, \pi')$ is a function of the face profile $\ell(\pi, \pi')$ of $\pi, \pi'$.
\end{lemma}

\begin{remark}
We defined the Orthogonal Weingarten function in a slightly different manner than the Unitary Weingarten function (equation \eqref{eq:weingarten-character-sum}). For an expression of $\Wg^{\mrm{O}}_\zeta$ in terms of characters, see \cite[Theorem 3.9]{collins2006integration} or \cite[Proposition 5]{zinn2009jucys}. The interpretation of the Weingarten function which is most relevant for us is as a weight assigned to pairs of left and right matchings, and the most direct definition of $\Wg_\zeta^{\mrm{O}}$ from this point of view is as the pseudo-inverse of the Gram matrix.

Also, note that we defined the Orthogonal Weingarten function for a general complex parameter $\zeta \in \C$. This did not require any extra considerations. For Orthogonal Haar integration, this level of generality is not needed and we could have restricted to $\zeta = N$ a positive integer. However, it turns out that the Symplectic Weingarten function is related to the Orthogonal Weingarten function with $\zeta = -N$ a negative integer -- see Lemma \ref{lemma:os-weingarten-relation}. Moreover, it will be more convenient to work with $\Wg^{\mrm{O}}_{-N}$ rather than the Symplectic Weingarten function, due to a certain sign issue. See Remark \ref{remark:orthogonal-weingarten-more-convenient} for more discussion.
\end{remark}

% Next, we define the symplectic Weingarten function. First, we need the following preliminary definitions.

\begin{definition}\label{def:pi-0}
Let $n \geq 1$ be even. Let $\pi_0 : [n] \ra [n]$ be the matching given by $\{\{n, n-1\}, \{n-2, n-3\}, \ldots, \{2, 1\}\}$. One may visualize $\pi_0$ as in Figure \ref{figure:pi-0}. We omit the dependence of $\pi_0$ on $n$.
\begin{figure}[ht!]
    \centering
\includegraphics[width=.15\textwidth]{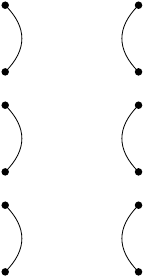}
    \caption{$[\pi_0 ~ \pi_0]$ when $n = 6$}
    \label{figure:pi-0}
\end{figure}
\end{definition}

% \begin{definition}
% For each matching $\pi : [n] \ra [n]$, we fix a permutation $\sigma_\pi$ such that $\sigma_\pi [\pi ~ \pi] \sigma_\pi^{-1} = [\pi_0 ~ \pi_0]$.
% \end{definition}

\begin{definition}\label{def:sigma-pi}
Let $n$ be even. For each matching $\pi : [n] \ra [n]$, we define a permutation $\sigma_\pi \in \symgrp_n$ such that $\sigma_\pi [\pi ~ \pi] \sigma_\pi^{-1} = [\pi_0 ~ \pi_0]$ as follows. We may write $\pi = \{\{\pi(1), \pi(2)\}, \ldots, \{\pi(n-1), \pi(n)\}\}$, where $1 = \pi(1) < \pi(3) < \cdots < \pi(n-1)$, and $\pi(2j-1) < \pi(2j)$ for $j \in [n/2]$. We then define $\sigma_\pi(j) := \pi(j)$.
\end{definition}

\begin{figure}
    \centering
    \includegraphics[width=0.7\textwidth]{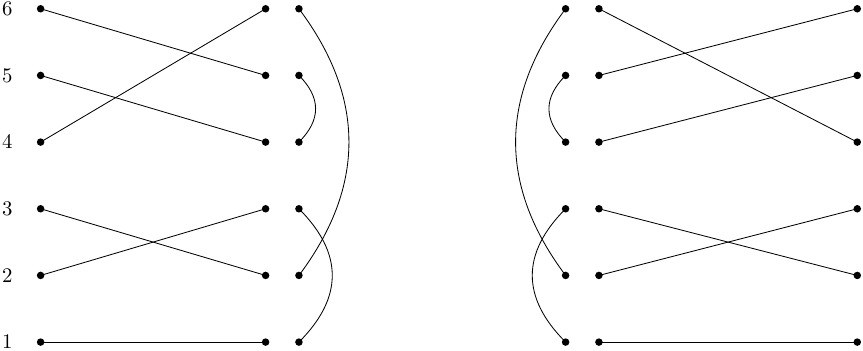}
    \caption{In the middle, we have $[\pi ~ \pi]$ where $\pi = \{\{1, 3\}, \{2, 6\}, \{4, 5\}\}$. Thus $\sigma_\pi = (1 ~ 3 ~ 2 ~ 6 ~ 4 ~ 5)$ (written in one-line notation). We see that upon conjugating $[\pi ~ \pi]$ by $\sigma_\pi$, we get the ``standard form" $[\pi_0 ~ \pi_0]$.}
    \label{figure:os-sigma-pi-example}
\end{figure}

See Figure \ref{figure:os-sigma-pi-example} for an example of $\sigma_\pi$. Visually, $\sigma_\pi$ can be thought of as a permutation of the vertices which takes $[\pi ~ \pi]$ to the ``standard form" $[\pi_0 ~ \pi_0]$. In general, there may be many such permutations; the definition of $\sigma_\pi$ makes a particular choice for each $\pi$. This particular way of choosing the permutation does not matter so much for $\orthogonal(N)$, however for $\SphN$ it is important that $\sigma_\pi$ be defined as it is, due to the fact that $\mrm{sgn}(\sigma_\pi)$ appears in the definition of the Symplectic Gram matrix (see Definition \ref{def:symplectic-weingarten}), and thus also the Symplectic Weingarten function. (Different permutations which take $[\pi ~ \pi]$ to the standard form $[\pi_0 ~ \pi_0]$ may have opposite signs.)

% See Figure \ref{figure:general_pi_to_pi_0} for an example.
% \begin{figure}[ht!]
%     \centering
% \includegraphics[width=0.5\textwidth]{figures/general_pi_to_pi_0.pdf}
%     \caption{The middle element is $[\pi ~ \pi]$ where $\pi = \{\{4, 2\}, \{3, 1\}\}$. By swapping the second and third vertices on both sides, we may change $\pi$ to $\pi_0$. This corresponds to conjugating by the $(2 ~ 3)$, as displayed in the figure.}\label{figure:general_pi_to_pi_0}
% \end{figure}

\begin{figure}[ht!]
    \centering
\includegraphics[width=0.5\textwidth]{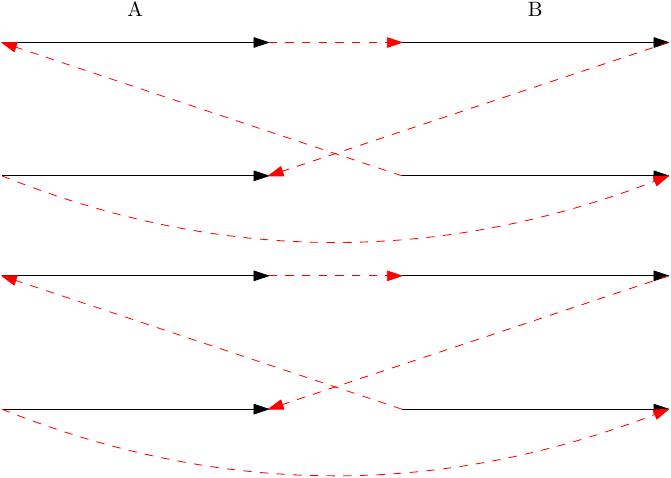}
    \caption{The red exterior connections arising from the collection of words $\bm \Gamma = (A B A^{-1} B^{-1}, A B A^{-1} B^{-1})$.}
    % Compare with the right of Figure \ref{figure:strand diagrams components} in the Unitary case with the exact same $\bm \Gamma$.}
    \label{figure:orthogonal_strand_diagram_from_word}
\end{figure}

Let $\bm \Gamma = (\Gamma_1, \ldots, \Gamma_M)$ be an unoriented-balanced collection of words on $\{\lambda_1, \ldots, \lambda_L\}$. Recall that in the Unitary case, the choice of $\bm \Gamma$ specifies a choice of red exterior connections in our strand diagram. In the orthogonal case, the situation is similar, except now we specify that all strands point in the same direction (right). By doing so, the dashed red strands that we add may not have a consistent orientation with the black strands. This is a reflection of the fact that in the Orthogonal case, the surfaces we obtain may be unorientable. We explain through an example how to obtain the red exterior connections from $\bm \Gamma$ -- see Figure \ref{figure:orthogonal_strand_diagram_from_word}.

\begin{figure}[ht!]
    \centering
\includegraphics[page=2, width=0.5\textwidth]{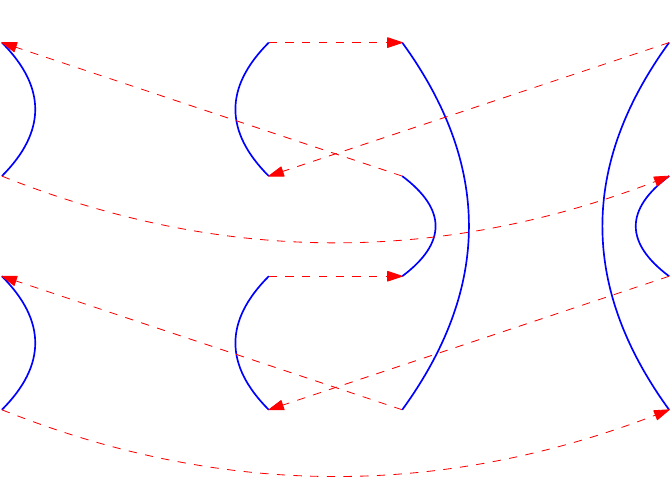}
    \caption{Let $\bm \Gamma$ be the same as in Figure \ref{figure:orthogonal_strand_diagram_from_word}. For some particular choice of $\bm \pi$, we may end up with the blue interior connections as displayed. In this case, $\numcomp(\bm \Gamma, \bm \pi) = 2$.}\label{figure:orthogonal_strand_diagram_from_word_2}
\end{figure}

For each $\ell \in [L]$, let $\pi_\ell, \pi_\ell' : [n_\ell] \ra [n_\ell]$ be matchings. Similar to the Unitary case, we may form the diagram obtained by $\bm \Gamma$ and $\bm \pi = ((\pi_\ell, \pi_\ell'), \ell \in [L])$ by starting with the red exterior connections specified by $\bm \pi$, and then adding in the blue interior connections specified by $\bm \pi$. Let $\numcomp(\bm \Gamma, \pi)$ be the number of components of this diagram. See Figure \ref{figure:orthogonal_strand_diagram_from_word_2} for an example.

We now state the Orthogonal Weingarten calculus.

\begin{prop}[Orthogonal Weingarten calculus]\label{prop:orthogonal-weingarten-calculus}
Let $G = \orthogonal(N)$. Let $\bm \Gamma = (\Gamma_1, \ldots, \Gamma_k)$ be an unoriented-balanced collection of words on $\{\lambda_1, \ldots, \lambda_L\}$. Then
\begin{equs}
\E[\Tr(G(\bm \Gamma))] = \sum_{\bm \pi = ([\pi_\ell, \pi_\ell'], \ell \in [L])} \bigg(\prod_{\ell \in L} \Wg_N^{{\mrm{O}}}(\pi_\ell, \pi_\ell')\bigg) N^{\numcomp(\bm \Gamma, \bm \pi) }.
\end{equs}
Here, the sum in the right\revision{-}hand side is over $\pi$ which is a collection of pairs of matchings $\pi_\ell, \pi_\ell' : [n_\ell] \ra [n_\ell]$, $\ell \in [L]$.
\end{prop}

We proceed towards applying Proposition \ref{prop:orthogonal-weingarten-calculus} to give expressions for Wilson loop expectations of $\orthogonal(N)$ lattice gauge theories. First, we need some setup. Exactly as in the Unitary case, given $(\bm \Gamma, \bm \pi)$, we may obtain a map whose dual is bipartite as follows. We start with one yellow face for each word in $\bm \Gamma$. For each letter $\lambda_\ell$, the left and right matchings $\pi_\ell, \pi_\ell'$ giving the interior connections in the portion of the diagram corresponding to $\lambda_\ell$ then specify an additional collection of blue faces which are glued to the yellow faces which contain the letter $\lambda_\ell$ or its inverse.

\begin{definition}
Let $\bm \Gamma = (\Gamma_1, \ldots, \Gamma_k)$ be an unoriented-balanced collection of words on $\{\lambda_1, \ldots, \lambda_L\}$. Define $\bpm_{\mrm{OS}}(\bm \Gamma)$ to be the set of all possible maps which can be obtained from adding interior left and right matchings to the strand diagram corresponding to $\bm \Gamma$. For a given map $\mc{M} \in \bpm_{\mrm{OS}}(\bm\Gamma)$, and $\ell \in [L]$, let $\mu_\ell(\mc{M})$ be the partition of $n_{\ell}$ (the total number of occurrences of $\lambda_\ell$ and $\lambda_\ell^{-1}$) given by the degrees of the blue faces which are glued in to the strand diagram of $\lambda_\ell$ (this is the same as the face profile of the left and right matchings $\bm \pi = (\pi_\ell, \pi_\ell', \ell \in [L])$ used to construct $\mc{M}$).
\end{definition}

Here, the subscript ``$\mrm{OS}$" is short for Orthogonal and Symplectic, since $\bpm_{\mrm{OS}}$ is the set of maps that one obtains in these cases.

\begin{remark}
Unlike in the Unitary case, the maps in $\bpm_{\mrm{OS}}$ may be unorientable. 
\end{remark}

The next result is the analog of Proposition \ref{prop:word-expectation-as-sum-over-bipartite-maps}.

\begin{prop}
Let $G = \orthogonal(N)$. Let $\bm \Gamma = (\Gamma_1, \ldots, \Gamma_k)$ be an unoriented-balanced collection of words on $\{\lambda_1, \ldots, \lambda_L\}$. We have that
\begin{equs}
\E[\Tr(G(\bm \Gamma))] = \sum_{\mc{M} \in \bpm_{\mrm{OS}}(\bm \Gamma)} \bigg(\prod_{\ell \in [L]} \ovl{\Wg}_N^{\mrm{O}}(\mu_\ell(\mc{M}))\bigg) N^{\chi(M) - k}.
\end{equs}
\end{prop}

As in the Unitary case, when the letters $\{\lambda_1, \ldots, \lambda_L\}$ are edges of the lattice $\Lambda$, and the collections of words $\bm \Gamma = \bm \Gamma(s, K)$ is obtained from a string $s$ and a plaquette count $K$, then any map $\mc{M} \in \bpm_{\mrm{OS}}(\bm \Gamma)$ naturally gives an edge-plaquette embedding $(\ovl{\mc{M}}, \phi)$, where $\phi$ is determined by the requirement that it maps edges of $\mc{M}$ to the corresponding edges of the lattice. As in the Unitary case, the map $\ovl{\mc{M}}$ is obtained from $\mc{M}$ by deleting all faces whose boundary is mapped by $\phi$ to a loop in $s$, and so $(\ovl{\mc{M}}, \phi)$ may be interpreted as having ``boundary" $s$.

\begin{definition}
Let $s = (\ell_1, \ldots, \ell_n)$ be a string, and let $K : \mc{P} \ra \N$. Define the set $\epe_{\mrm{OS}}(s, K)$ of edge-plaquette embeddings associated to $(s, K)$ to as follows. If $\bm \Gamma(s, K)$ is not unoriented-balanced, then $\epe(s, K) := \varnothing$. If $\bm \Gamma(s, K)$ is unoriented-balanced, define $\epe_{\mrm{OS}}(s, K)$ to be the set of edge-plaquette embedding $(\ovl{\mc{M}}, \phi)$ obtained from maps $\mc{M} \in \bpm_{\mrm{OS}}(\bm \Gamma(s, K))$.

Next, define
\begin{equs}\label{eq:epe-OS-partition}
\epe_{\mrm{OS}}(s) := \bigsqcup_{K : \mc{P} \ra \N} \epe_{\mrm{OS}}(s, K).
\end{equs}
For $(\mc{M}, \phi) \in \epe_{\mrm{OS}}(s)$, and $e \in E_\Lambda$, let $\mu_e(\phi)$ be the partition of $|\phi^{-1}(e)|/2$ induced by $1/2$ times the degrees of the faces of $\phi^{-1}(e)$. Define 
\begin{equs}
\area(\mc{M}, \phi) &:= \sum_{p \in \mc{P}} |\phi^{-1}(p)|, \\
(\phi^{-1})! &:= \prod_{p \in \mc{P}} |\phi^{-1}(p)|! .
\end{equs}
Note that if $(\mc{M}, \phi) \in \epe_{\mrm{OS}}(s, K)$, then $\area(\mc{M}, \phi) = \sum_p K(p)$ and $(\phi^{-1})! = K!$. 
\end{definition}

We now arrive at the following theorem, which is the analog of Theorem \ref{thm:wilson-loop-expectation-sum-over-epe}. Since the proof is very similar to the proof of the corollary, it is omitted.

\begin{theorem}
\revision{Let $\Lambda$ be a finite lattice.} Let $s = (\ell_1, \ldots, \ell_n)$ be a string. For $\orthogonal(N)$ lattice gauge theory, we have that
\begin{equs}
\langle W_s \rangle_{\Lambda, \beta} = Z_{\Lambda, \beta}^{-1} \sum_{(\mc{M}, \phi) \in \epe_{\mrm{OS}}(s)} \frac{\beta^{\area(\mc{M}, \phi)}}{(\phi^{-1})!} \bigg(\prod_{e \in E_\Lambda} \ovl{\Wg}_N^{\mrm{O}}(\mu_e(\phi))\bigg)   N^{\chi(\mc{M}) - n}.
\end{equs}
\revision{The sum converges absolutely in the following sense. Recalling that  $\epe_{\mrm{OS}}(s) = \bigsqcup_{K : \mc{P} \ra \N} \epe_{\mrm{OS}}(s, K)$, we have that
\begin{equs}
\sum_{K : \mc{P} \ra \N} \Bigg|\sum_{(\mc{M}, \phi) \in \epe_{\mrm{OS}}(s, K)} \frac{\beta^{\area(\mc{M}, \phi)}}{(\phi^{-1})!} \bigg(\prod_{e \in E_\Lambda} \ovl{\Wg}_N^{\mrm{O}}(\mu_e(\phi))\bigg)   N^{\chi(\mc{M}) - n}\Bigg| < \infty.
\end{equs}}
\end{theorem}

\subsubsection{Symplectic surface sums}

Next, we discuss the surface sums in the Symplectic case. This case is more complicated than before due to a certain sign issue. We start by working towards the definition of the Symplectic Weingarten function.

\begin{definition}[Symplectic Weingarten function]\label{def:symplectic-weingarten}
Define the Symplectic Weingarten function $\Wg_{N, n}^{\mrm{Sp}}$ as follows. First, define the Gram matrix
\begin{equs}
\mbf{G}_{N, n}^{\mrm{Sp}}(\pi, \pi') := (-1)^{n/2} \mrm{sgn}(\sigma_\pi) \mrm{sgn}(\sigma_{\pi'}) (-N)^{\cycles(\pi, \pi')}, ~~ \pi, \pi' : [n] \ra [n] \text{ matchings}.
\end{equs}
We define $\Wg_{N, n}^{\mrm{Sp}}$ to be the pseudo-inverse of $\mbf{G}_{N, n}^{\mrm{Sp}}$. We typically omit the dependence on $n$ and write $\Wg_N^{\mrm{Sp}}$.
\end{definition}

\begin{remark}
This definition of the Symplectic Weingarten function is not so easy to find in the literature. For instance, the first paper on the topic \cite{collins2006integration} does not give an explicit formula for the Symplectic Weingarten function, nor does the recent survey \cite{collins2022weingarten}. The paper \cite{Magee2019a} which applies the Symplectic Weingarten calculus only posits the existence of some function which can be used to compute Symplectic matrix integrals (see \cite[Theorem 3.1]{Magee2019a}). The paper \cite{matsumoto2013weingarten} defines the Symplectic Weingarten function as a certain element $W$ of the group algebra $\C[\symgrp_n]$ (Matsumoto denotes this element by $\mrm{Wg}^{\mrm{Sp}}$). The relation between Matsumoto's definition and our definition via pseudo-inverses is precisely stated in \cite[Lemma 2.5]{matsumoto2013weingarten}, which says that the Weingarten weight $\Wg_N^{\mrm{Sp}}(\pi, \pi')$ assigned to a pair of matchings $\pi, \pi'$ is precisely $W(\sigma_\pi^{-1} \sigma_{\pi'})$. We prefer to give the pseudo-inverse definition in the present paper, because it is the most easy to state and understand. This way, the reader who only wishes to be able to understand the weights that appear in our surface sums can do so without having to spend too much time on background material.
\end{remark}

By comparing the definitions of the Orthogonal (Definition \ref{def:orthogonal-weingarten}) and Symplectic Weingarten functions, the next lemma follows immediately. (Here we also use the uniqueness of the pseudo-inverse of a matrix.)

\begin{lemma}\label{lemma:os-weingarten-relation}
We have that
\begin{equs}
\Wg_{N}^{\mrm{Sp}}(\pi, \pi') = (-1)^{n/2} \mrm{sgn}(\sigma_\pi) \mrm{sgn}(\sigma_{\pi'}) \Wg_{-N}^{\mrm{O}}(\pi, \pi'), ~~ \pi, \pi' : [n] \ra [n].
\end{equs}
\end{lemma}

\begin{remark}
This relation between the Orthogonal and Symplectic Weingarten functions has previously been observed, see for instance the end of \cite[Section 2.3.2]{matsumoto2013weingarten}. When $N \geq n$, this identity is also stated as \cite[Lemma 3.2]{Magee2019a}. We note that by defining Weingarten functions as pseudo-inverses of the appropriate Gram matrices, it is trivial to see that the relation holds for general $N$ (indeed, even general $\zeta \in \C$).
\end{remark}

\begin{remark}\label{remark:symplectic-weingarten-not-a-function-of-face-profile}
Recall that $\Wg^{\mrm{O}}_{-N}(\pi, \pi')$ is a function of the face profile $\ell(\pi, \pi')$. Lemma \ref{lemma:os-weingarten-relation} shows that $\Wg_N^{\mrm{Sp}}(\pi, \pi')$ is {\it not} a function of the face profile $\ell(\pi, \pi')$, because $\mrm{sgn}(\sigma_\pi) \mrm{sgn}(\sigma_{\pi'})$ is not determined by $\ell(\pi, \pi')$. For a simple example, see Figure \ref{figure:symplectic_weingarten_sign_issue_example}. 
\begin{figure}
    \centering
    \includegraphics[width=0.4\textwidth]{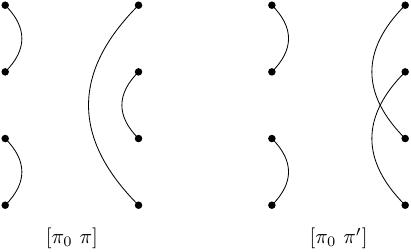}
    \caption{Left: $[\pi_0 ~ \pi]$ with $\pi = \{\{1, 4\}, \{2, 3\}\}$. Thus $\sigma_\pi = (1 ~ 4 ~ 2 ~ 3)$, and so $\mrm{sgn}(\sigma_\pi) = 1$. Right: $[\pi_0 ~ \pi']$ with $\pi' = \{\{1, 3\}, \{2, 4\}\}$. Thus $\sigma_{\pi'} = (1 ~ 3 ~ 2 ~ 4)$, and so $\mrm{sgn}(\sigma_{\pi'}) = -1$. Consequently, we see that $\Wg_N^{\mrm{Sp}}(\pi_0, \pi) = (-1)^2 \mrm{sgn}(\sigma_{\pi_0}) \mrm{sgn}(\sigma_{\pi'}) \Wg_{-N}^{\mrm{O}}(\pi_0, \pi') = \Wg_{-N}^{\mrm{O}}(\pi_0, \pi')$, whereas $\Wg_N^{\mrm{Sp}}(\pi_0, \pi') = - \Wg_{-N}^{\mrm{O}}(\pi_0, \pi')$. Note that the face profiles $\ell(\pi_0, \pi) = \ell(\pi_0, \pi') = \{4\}$, and so $\Wg^{\mrm{O}}_{-N}(\pi_0, \pi) = \Wg^{\mrm{O}}_{-N}(\pi_0, \pi')$. This shows that $\Wg^{\mrm{Sp}}_N(\pi_0, \pi) = -\Wg^{\mrm{Sp}}_N(\pi_0, \pi')$, even though $(\pi_0, \pi)$ has the same face profile as $(\pi_0, \pi')$.}
    \label{figure:symplectic_weingarten_sign_issue_example}
\end{figure}
Thus to obtain weighted sums over surfaces in the Symplectic case, we will use Lemma \ref{lemma:os-weingarten-relation} to replace $\Wg_N^{\mrm{Sp}}$ by $\Wg_{-N}^{\mrm{O}}$, which will allow us to express our weights purely in terms of the surfaces. We note that this was also done in \cite{Magee2019a} -- see Theorem 1.2 and Appendix A of the paper.
\end{remark}

Recall the matrix $J = \begin{pmatrix} 0 & I_{N/2} \\ -I_{N/2} & 0 \end{pmatrix}$ from the definition of $\SphN$.

\begin{definition}\label{def:Delta-pi-prime}
For indices $i_1, i_2 \in [N]$, define $\langle i_1, i_2 \rangle_J := J_{i_1 i_2}$. For $n \geq 1$ even and a permutation $\sigma \in \symgrp_n$, define
\begin{equs}
\Delta_\sigma'(\mbf{i}) := \prod_{k=1}^{n/2} \langle i_{\sigma(2k-1)}, i_{\sigma(2k)} \rangle_J = \langle i_{\sigma(1)}, i_{\sigma(2)} \rangle_J \cdots \langle i_{\sigma(n-1)}, i_{\sigma(n)} \rangle_J, ~~ \mbf{i} = (i_1, \ldots, i_n) \in [N]^n.
\end{equs}
For a matching $\pi : [n] \ra [n]$, we abuse notation and write $\Delta_\pi'$ for $\Delta_{\sigma_\pi}'$.
\end{definition}

We next state the matrix-entry version of the Symplectic Weingarten calculus. This is essentially \cite[Theorem 2.4]{matsumoto2013weingarten} (see also \cite[Lemma 2.5]{matsumoto2013weingarten}).
% while the second form is what our Poisson exploration process most directly proves. Additionally, the second form is better for obtaining weighted surface sums, since it involves the Orthogonal Weingarten function which is a function of the face profile -- recall Remark \ref{remark:symplectic-weingarten-not-a-function-of-face-profile}.

\begin{prop}[Symplectic Weingarten calculus]\label{prop:symplectic-weingarten-matrix-entry}
Let $G = \SphN$. Let $n \geq 1$ be even. For any $\mbf{i} = (i_1, \ldots, i_n), \mbf{j} = ( j_1, \ldots, j_n) \in [N]^{n}$, we have that
\begin{equs}
\E[G_{i_1 j_1} \cdots G_{i_n j_n}] = \sum_{\pi, \pi' : [n] \ra [n]} \Delta_{\pi}'(\mbf{i}) \Delta_{{\pi'}}'(\mbf{j}) \Wg_N^{\mrm{Sp}}(\pi, \pi').
\end{equs}
% Equivalently, we have that
% \begin{equs}
% \E[S^{\otimes n}] = \sum_{\pi, \pi' : [n] \ra [n]} \Wg_{-N}^{\mrm{O}}(\pi, \pi') \rho_-([\pi ~ \pi']).
% \end{equs}
\end{prop}

By applying Proposition \ref{prop:symplectic-weingarten-matrix-entry} and Lemma \ref{lemma:os-weingarten-relation}, one can obtain the following word-expectation version of the Symplectic Weingarten calculus. We remark that going from the matrix-entry version to the word-expectation version of Weingarten caclulus is not as simple as in the Unitary or Orthogonal cases (where one may use the argument described in Section \ref{section:poisson-exploration-general-N-proof}), and one has to carefully handle signs. The proof is omitted -- see \cite[Appendix A]{Magee2019a} for the relevant details.

\begin{prop}\label{prop:symplectic-weingarten-word-expectation}
Let $G = \SphN$. Let $\bm \Gamma = (\Gamma_1, \ldots, \Gamma_k)$ be an unoriented-balanced collection of words on $\{\lambda_1, \ldots, \lambda_L\}$. We have that
\begin{equs}
\E[\Tr(G(\bm \Gamma))] = (-1)^k \sum_{\bm \pi = ([\pi_\ell, \pi_\ell'], \ell \in [L])} \prod_{\ell \in [L]} \Wg^{\mrm{O}}_{-N}(\pi_\ell, \pi_\ell') (-N)^{\numcomp(\bm \Gamma, \bm \pi)}.
\end{equs}
Consequently, we have that
\begin{equs}
\E[\Tr(G(\bm \Gamma))] = (-1)^k \sum_{\mc{M} \in \bpm_{OS}(\bm \Gamma)} \bigg(\prod_{\ell \in [L]} \ovl{\Wg}_{-N}^{\mrm{O}}(\mu_\ell(\mc{M}))\bigg) (-N)^{\chi(M) - k}.
\end{equs}
\end{prop}

\begin{remark}\label{remark:orthogonal-weingarten-more-convenient}
To obtain the second claim in Proposition \ref{prop:symplectic-weingarten-word-expectation}, it was crucial that we used the Orthogonal Weingarten function rather than the Symplectic Weingarten function, since the former is a function of the face profile but the latter is not (recall Remark \ref{remark:symplectic-weingarten-not-a-function-of-face-profile}). In other words, $\Wg_N^{\mrm{Sp}}$ is not a function of $\mu_\ell(\mc{M})$, so we could not have replaced $\Wg_N^{\mrm{Sp}}(\pi_\ell, \pi_{\ell'})$ with $\Wg_N^{\mrm{Sp}}(\mu_\ell(\mc{M}))$.
\end{remark}

We can now give a representation of Wilson loop expectations in $\SphN$ lattice gauge theories as Weingarten-weighted surface sums. 

\begin{theorem}
\revision{Let $\Lambda$ be a finite lattice.} Let $s = (\ell_1, \ldots, \ell_n)$ be a string. For $\SphN$ lattice gauge theory, we have that
\begin{equs}
\langle W_s \rangle_{\Lambda, \beta} = (-1)^n Z_{\Lambda, \beta}^{-1} \sum_{(\mc{M}, \phi) \in \epe_{\mrm{OS}}(s)} \frac{\beta^{\mrm{area}(\mc{M}, \phi)}}{(\phi^{-1})!} \bigg(\prod_{e \in E_\Lambda} \ovl{\Wg}_{-N}^{\mrm{O}}(\mu_e(\phi))\bigg) (-N)^{\chi(M) - n}.
\end{equs}
\revision{The sum converges absolutely in the following sense. Recalling that  $\epe_{\mrm{OS}}(s) = \bigsqcup_{K : \mc{P} \ra \N} \epe_{\mrm{OS}}(s, K)$, we have that
\begin{equs}
\sum_{K : \mc{P} \ra \N} \Bigg|\sum_{(\mc{M}, \phi) \in \epe_{\mrm{OS}}(s, K)} \frac{\beta^{\mrm{area}(\mc{M}, \phi)}}{(\phi^{-1})!} \bigg(\prod_{e \in E_\Lambda} \ovl{\Wg}_{-N}^{\mrm{O}}(\mu_e(\phi))\bigg) (-N)^{\chi(M) - n}\Bigg| < \infty.
\end{equs}}
\end{theorem}

% \begin{lemma}
% Let $\bm \Gamma = (\Gamma_1, \ldots, \Gamma_k)$ be an unoriented-balanced collection of words on $\{\lambda_1, \ldots, \lambda_L\}$. We have that
% \begin{equs}
% \E[\Tr(S(\bm \Gamma))] &= \sum_{\mc{M} \in \bpm_{\mrm{OS}}(\bm \Gamma)} \bigg(\prod_{\ell \in [L]} \ovl{\Wg}_N^S(\mu_\ell(\mc{M}))\bigg)  (-1)^{\chi(M)} N^{\chi(M) - k} \\
% &=  (-1)^k \sum_{\mc{M} \in \bpm_{\mrm{OS}}(\bm \Gamma)} \bigg(\prod_{\ell \in [L]} \ovl{\Wg}_{-N}^O(\mu_\ell(\mc{M}))\bigg)  (-N)^{\chi(M) - k} 
% \end{equs}
% \end{lemma}

\subsubsection{Exploration process}

In this subsection, we detail how to obtain the Orthogonal and Symplectic Weingarten calculus from taking limits of Brownian motion, much as we did in Section \ref{section:poisson-exploration} for the Unitary case. The key is to use a variant of the exploration process we defined in Section \ref{section:strand-by-strand}. This will again allow us to extract the Jucys-Murphy elements, which as before will relate to the Weingarten function. First, we work towards describing the analog of Proposition \ref{prop:unitary-brownian-motion-expectation} for $\orthogonal(N), \SphN$. 

\begin{remark}\label{remark:SON-ON-BM}
Note that $\SON$ is connected while $\orthogonal(N)$ is not. Thus an $\orthogonal(N)$ Brownian motion started at the identity (or more generally, any element of $\SON$) is exactly the same as an $\SON$ Brownian motion. Thus to reprove the Orthogonal Weingarten calculus, we will need to take the initial value $O_0$ of the $\orthogonal(N)$ Brownian motion to lie in the two connected components of $\orthogonal(N)$ with equal probability. This amounts to multiplying an $\SON$ Brownian motion started at the identity by $O_0$.
\end{remark}

\begin{notation}
In the following, to discuss the cases $G = \orthogonal(N), \SphN$ simultaneously, we set the notation $\varep = 1$ when $G = \orthogonal(N)$ and $\varep = -1$ when $G = \SphN$. We found this useful notation from \cite{Dahlqvist2017}.
\end{notation}

Next, we define a representation $\rho_- : \mc{B}_n(-N) \ra \End((\C^N)^{\otimes n})$ which is needed to relate expectations of Symplectic Brownian motion to weighted sums over Brauer algebra elements. Recall the definition of $\langle \cdot, \cdot \rangle_J$ and $\Delta_\pi'$ from Definition \ref{def:Delta-pi-prime}.

\begin{definition}
Define the representation $\rho_{-} : \mc{B}_n(-N) \ra \End((\C^N)^{\otimes n})$ as follows. It suffices to define $\rho_{-}$ on the generating set $\{(i ~ j), \langle i ~ j \rangle, 1 \leq i  < j < n\}$ of $\mc{B}_n$. We let $\rho_-((i ~ j)) := -\rho_+((i ~ j))$. We let $\rho_-(\langle i ~ j\rangle)$ be the matrix whose $(\mbf{k}, \mbf{l})$ (with $\mbf{k} = (k_1, \ldots, k_n), \mbf{l} = (l_1, \ldots, l_n) \in [N]^n$) matrix entry is
\begin{equs}
(\rho_-(\langle i ~ j\rangle))_{\mbf{k} \mbf{l}} := - \langle k_i, k_j \rangle_J \langle l_i, l_j \rangle_J \prod_{r \neq i, j} \delta_{k_r l_r} .
\end{equs}
\end{definition}

\begin{remark}
The minus sign in the definition of $\rho_-$ is crucial, since $\rho_-$ is supposed to be a representation of $\mc{B}_n(-N)$, which implies that $\rho_-(\langle i ~ j \rangle)^2 = (-N) \rho_-(\langle i ~ j \rangle)$ (since $\langle i ~ j \rangle^2 = (-N) \langle i ~ j \rangle$). The minus sign ensures that this is the case.
\end{remark}

\begin{remark}\label{remark:rho-minus-matrix-entries}
For a matching $\pi : [n] \ra [n]$, observe that
\begin{equs}
~[\pi ~ \pi] = \langle \pi(1) ~ \pi(2) \rangle \cdots \langle \pi(2n-1) ~ \pi(2n) \rangle.
\end{equs}
This implies that
\begin{equs}
\rho_-([\pi ~ \pi])_{\mbf{i}, \mbf{j}} = (-1)^{n/2} \Delta_{\pi}'(\mbf{i}) \Delta_{\pi}'(\mbf{j}).
\end{equs}
More generally, we have that
\begin{equs}\label{eq:rho-minus-matrix-entries}
\rho_-([\pi ~ \pi'])_{\mbf{i} \mbf{j}} = (-1)^{n/2} \mrm{sgn}(\sigma_\pi) \mrm{sgn}(\sigma_{\pi'}) \Delta_\pi'(\mbf{i}) \Delta_{\pi'}'(\mbf{j}).
\end{equs}
\end{remark}

\begin{notation}
We will take $\rho_1$ to mean $\rho_+$ and $\rho_{-1}$ to mean $\rho_-$. This way we can write $\rho_\varep$.
\end{notation}

Let $n \geq 1$. Consider a strand diagram with $n$ total strands, all of them right-directed. We define a Poisson point process $\Sigma_{\mrm{OS}}$ as follows. Define
\begin{equs}
\mc{D}_{\mrm{OS}} := \bigsqcup_{\substack{i, j \in [n] \\ i \neq j}} [0, \infty), ~~ \mc{D}_{\mrm{OS}}(T) := \bigsqcup_{\substack{i, j \in [n] \\ i \neq j}} [0, T].
\end{equs}
We let $\Sigma_{\mrm{OS}}$ be a rate-$1$ Poisson process on $\mc{D}_{\mrm{OS}}$. We also define $\Sigma_{\mrm{OS}}(T) := \Sigma_{\mrm{OS}} \cap \mc{D}_{\mrm{OS}}(T)$, which is a rate-$1$ Poisson process on $\mc{D}_{\mrm{OS}}(T)$. 

In terms of the strands, $\mc{D}_{\mrm{OS}}(T)$ is visualized as follows. Between any pair of strands, there are two independent rate-1 Poisson processes (which is why we only have the restriction $i \neq j$ and not $i < j$ as in the Unitary case): one which gives the swaps between the two strands, and one which gives the turnarounds between the two strands. The Poisson processes corresponding to different pairs of strands are also independent. Let $\Sigma_{\mrm{OS}}(T)$
be process obtained by keeping only those points of $\Sigma_{\mrm{OS}}$ which occur before time $T$.

% See Figure \ref{figure:orthogonal_poisson_process_realization} for an example realization of $\Sigma_{\mrm{OS}}(T)$.

\begin{figure}[ht!]
    \centering
\includegraphics[width=0.5\textwidth]{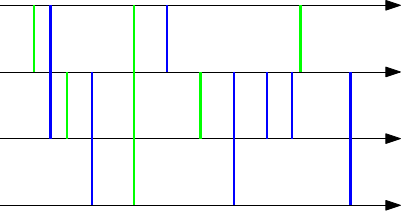}
    \caption{An example realization of $\Sigma_{\mrm{OS}}(T)$ for some finite $T$. The green lines represent swaps and blue lines represent turnarounds.}\label{figure:orthogonal_poisson_process_realization}
\end{figure}

\begin{remark}[Comparison with the Unitary case]
In the Orthogonal and Symplectic cases, all strands point the same direction, and there may be turnarounds between same-direction strands. Whereas in the Unitary case, there were only swaps between same-direction strands, and turnarounds between opposite-direction strands. 
\end{remark}

Let $P \sse \mc{D}_{\mrm{OS}}(T)$ be a finite point set. As in the Unitary case, we may represent $P$ as a finite set of points $\{x_1, \ldots, x_{|P|}\} \sse [0, T]$, $0 \leq x_1 \leq \cdots \leq x_{|P|} \leq T$, along with labels $\mf{l}_1, \ldots, \mf{l}_{|P|}$, where each $\mf{l}_k = (i_k, j_k)$, with $i_k, j_k \in [n]$, $i_k \neq j_k$. Hereafter, we assume that the $x_k$ are all distinct, which is a.s. the case for $\Sigma_{\mrm{OS}}(T)$.

\begin{definition}
Let $P \sse \mc{D}_{\mrm{OS}}(T)$ be a finite subset. Define $\pointstobrauerfn_\varep(P) \in \mc{B}_n(\varep N)$ as follows. We first define $b_k \in \mc{B}_n$ for $k \in [|P|]$. Using the representation just described, let $\mf{l}_k = (i_k, j_k)$. We let
\begin{equs}
b_k := \begin{cases} -\frac{\varep}{N}(i_k ~ j_k) & i_k < j_k \\
\frac{\varep}{N} \langle i_k ~ j_k \rangle & i_k > j_k . \end{cases}
\end{equs}
We then define $\pointstobrauerfn_\varep(P) := b_1 \cdots b_{|P|} \in \mc{B}_n(\varep N)$.
\end{definition}

\begin{remark}
In this definition, the parameter of the Brauer algebra is taken to be $\varep N$, which means that in the Symplectic case, closed loops arising during multiplication contribute factors of $-N$, as opposed to $N$.
\end{remark}

The following proposition is the analog of Proposition \ref{prop:unitary-brownian-motion-expectation}. It says that for $G = \SON, \SphN$, expectations of $G$-valued Brownian motion may be expressed in terms of the Poisson point process $\Sigma_{\mrm{OS}}$. See~\cite[Appendix A]{park2023wilson} for the proof.

\begin{prop}\label{prop:orthogonal-symplectic-brownian-motion-expectation}
Let $n \geq 1$. Let $G = \SON, \SphN$. Let $B_T$ be a $G$-valued Brownian motion at time $T$. We have that
\begin{equs}
\E[B_T^{\otimes n}] = e^{2\binom{n}{2} T - \frac{n}{2}(1 - \frac{\varep}{N}) T} \rho_\varep\big(\E[\pointstobrauerfn_\varep(\Sigma_{\mrm{OS}}(T))]\big).
\end{equs}
Equivalently, for $\mbf{i} = (i_1, \ldots, i_n), \mbf{j} = (j_1, \ldots, j_n) \in [N]^n$, we have that
\begin{equs}
\E[(B_T)_{i_1 j_1} \cdots (B_T)_{i_n j_n}] = e^{2\binom{n}{2} T - \frac{n}{2}(1 - \frac{\varep}{N})T} \rho_\varep\big(\E[\pointstobrauerfn_\varep(\Sigma_{\mrm{OS}}(T))]\big)_{\mbf{i} \mbf{j}}.
\end{equs}
\end{prop}

% \begin{remark}
% A more succinct statement of Proposition \ref{prop:orthogonal-symplectic-brownian-motion-expectation} would be to say that
% \begin{equs}
% \E[B_T^{\otimes n}] = e^{2\binom{n}{2} T - \frac{n}{2}(1 - \frac{\varep}{N})} \rho_\varep\big(\E[\pointstobrauerfn_\varep(\Sigma_{\mrm{OS}}(T))]\big).
% \end{equs}
% \end{remark}

Next, we state the following lemma which relates the Orthogonal Weingarten fuction to Jucys-Murphy elements. This may essentially be found in \cite{matsumoto2013weingarten}, however it may not be so clear why this is the case without actually reading the paper. Thus, for the reader's convenience, we provide some discussion in Appendix \ref{appendix:orthogonal-weingarten-jucys-murphy-relation} of why the following lemma follows from the results of \cite{matsumoto2013weingarten}.

\begin{lemma}[Relation of Weingarten function to Jucys-Murphy elements]\label{lemma:os-weingarten-jucys-murphy-relation}
Let $\pi : [n] \ra [n]$ be a matching. We have that
\begin{equs}
\rho_\varep([\pi ~ \pi_0]) \rho_\varep((\varep N + J_{n-1}) (\varep N + J_{n-3}) \cdots (\varep N + J_1))^{-1} = \sum_{\pi' : [n] \ra [n]} \Wg^{\mrm{O}}_{\varep N}(\pi_0, \pi') \rho_\varep([\pi ~ \pi']) .
\end{equs}
\end{lemma}

% \begin{definition}
% For $\mbf{i} = (i_1, \ldots, i_n) \in [N]^n$ and $\sigma \in \symgrp_n$, define
% \begin{equs}
% \Delta_\sigma'(\mbf{i}) := \langle i_{\sigma(1)}, i_{\sigma(2)} \rangle_J \cdots \langle i_{\sigma(n-1)}, i_{\sigma(n)} \rangle_J.
% \end{equs}
% For a matching $\pi : [n] \ra [n]$, we abuse notation and write $\Delta_\pi'(\mbf{i})$ for $\Delta_{\sigma_\pi}'(\mbf{i})$, where $\sigma_\pi$ is as in Definition \ref{def:sigma-pi}.
% \end{definition}

As in the Unitary case, the main theorem that we will prove using our exploration process is the recovery of the Weingarten calculus, stated as follows.

\begin{theorem}[Weingarten recovery]\label{thm:os-weingarten-recovery}
Let $n \geq 1$ be even. Let $O_0 \in \orthogonal(N)$ be a random matrix which has equal probability of being in the two connected components of $\orthogonal(N)$, or equivalently $\E[\mrm{det}(O_0)] = 0$. For $G = \orthogonal(N)$, we have that
\begin{equs}
\lim_{T \toinf} \E[B_T^{\otimes n}] \E[O_0^{\otimes n}]  = \sum_{\pi, \pi' : [n] \ra [n]} \Wg^{\mrm{O}}_N(\pi, \pi') \rho_+([\pi ~ \pi']).
\end{equs}
For $G = \SphN$, we have that
\begin{equs}\label{eq:symplectic-weingarten-recovery}
\lim_{T \toinf} \E[B_T^{\otimes n}] = \sum_{\pi, \pi' : [n] \ra [n]} \Wg^{\mrm{O}}_{-N}(\pi, \pi') \rho_-([\pi ~ \pi']).
\end{equs}
\end{theorem}

\begin{remark}
To see why equation \eqref{eq:symplectic-weingarten-recovery} is equivalent to the matrix-entry version of the Symplectic Weingarten calculus (Proposition \ref{prop:symplectic-weingarten-matrix-entry}), recall by Remark \ref{remark:rho-minus-matrix-entries} that
\begin{equs}
\big(\rho_-([\pi ~ \pi'])\big)_{\mbf{i} \mbf{j}} = (-1)^{n/2} \mrm{sgn}(\sigma_\pi) \mrm{sgn}(\sigma_{\pi'}) \Delta_\pi'(\mbf{i}) \Delta_{\pi'}'(\mbf{j}).
\end{equs}
Then by the relation between the Orthogonal and Symplectic Weingarten functions (Lemma \ref{lemma:os-weingarten-relation}), it follows that
\begin{equs}
\Wg^{\mrm{O}}_{-N}(\pi, \pi') \big(\rho_-([\pi ~ \pi'])\big)_{\mbf{i} \mbf{j}} &= \big((-1)^{n/2} \mrm{sgn}(\sigma_\pi) \mrm{sgn}(\sigma_{\pi'}) \Wg^{\mrm{O}}_{-N}(\pi, \pi')\big) \Delta_\pi'(\mbf{i}) \Delta_{\pi'}'(\mbf{j}) \\
&= \Wg_N^{\mrm{Sp}}(\pi, \pi) \Delta_\pi'(\mbf{i}) \Delta_{\pi'}'(\mbf{j}).
\end{equs}
\end{remark}

% Next, we define the analog of the event $A_T$ from Definition \ref{def:A-T}.

% \begin{definition}
% Let $T \geq 0$. Let $A_T^{\mrm{OS}}$ be the event that once a blue turnaround appears, the only points which thereafter appear that touch either of the matched strands must be the blue turnaround or green swap between the same two strands.
% \end{definition}

% The following lemma is the analog of Lemma \ref{lemma:cancellation}. The proof is very similar, and thus omitted.

% \begin{lemma}
% We have that
% \begin{equs}
% \E[\pointstobrauerfn(\Sigma_{\mrm{OS}}(T))] = \E[\pointstobrauerfn(\Sigma_{\mrm{OS}}(T)) \ind_{A_T^{\mrm{OS}}}].
% \end{equs}
% \end{lemma}

Next, we define the analog $\mc{Q}^{\mrm{OS}}_T$ of the strand-by-strand exploration process $\mc{Q}_T$ from Section \ref{section:strand-by-strand}. As before, the exploration is encoded by two processes $(E_t)_{t \geq 0}$, $(\pi_t)_{t \geq 0}$. Here, $E$ takes values in $[n/2]$ and tracks the current exploration era, and $\pi$ takes values in $\symgrp_n$. In words, the exploration starts at the top strand, and follows swaps until the first turnaround. At this time, the exploration proceeds to the next-highest strand which hasn't been matched. The exploration continues until all strands have been matched (i.e. until the end of the $n/2$th exploration era).

\begin{notation}
For notational brevity in what follows, define
\begin{equs}
h_n(t_1, \ldots, t_{n/2}) := e^{(2(n-1) - (1 - \frac{\varep}{N})) t_1} e^{(2(n-3) - (1 - \frac{\varep}{N})) (t_2 - t_1)} \cdots e^{(2(1) - (1 - \frac{\varep}{N})) (t_{n/2} - t_{n/2 - 1})}
\end{equs}
\end{notation}

The following is the analog of Proposition \ref{prop:strand-by-strand-exploration} which recall encoded the key cancellation property of our strand-by-strand exploration process.

\begin{prop}\label{prop:os-strand-by-strand}
Let $n \geq 1$ be even. Let $G = \SON, \SphN$. We have that
\begin{equs}
e^{2\binom{n}{2} T - \frac{n}{2}(1 - \frac{\varep}{N})T} &\E[\pointstobrauerfn_\varep(\Sigma_{\mrm{OS}}(T)) \ind(T_{n/2} \leq T)] = \E\big[F_\varep(\mc{Q}_{T_{n/2}}^{\mrm{OS}}) \ind(T_{n/2} \leq T) h_n(T_1, \ldots, T_{n/2}) \big].
\end{equs}
\end{prop}
\begin{proof}
The proof is very similar to the proof of Proposition \ref{prop:strand-by-strand-exploration}, in that we proceed by induction, except now the combinatorics is slightly different. Throughout, we write $\Sigma$ and $\mc{Q}$ instead of $\Sigma_{\mrm{OS}}$ and $\mc{Q}^{\mrm{OS}}$ for brevity.

The base case $n = 2$ may be handled by direct calculation, which we omit. Suppose that the proposition is true for $n -2 \geq 2$ and any $T \geq 0$. As before, we condition on the time $T_1$, which is the time of first turnaround, which results in two strands being matched. After this time, we may assume that any swaps or turnarounds involving either of the matched strands must involve precisely the two matched strands (by essentially the same argument as in the proof of the cancellation lemma, Lemma \ref{lemma:cancellation}). Each strand is involved in $2(n-1)$ independent Poisson processes, and thus the number of independent Poisson processes which must have zero points on the interval $[T_1, T]$ is $2(2(n-1)) - 4 = 4n - 8$. The Poisson process which gives the turnarounds between the two matched strands contributes a factor 1, and the Poisson process which gives the swaps between the two matched strands contributes a factor $e^{-(T - T_1)} e^{-\varep(T - T_1)/N}$ (when we condition on $T_1$). We thus obtain
\begin{equs}
\E[&\pointstobrauerfn_\varep(\Sigma(T))   \ind(T_n \leq T) ~|~ \mc{F}_{T_1}] = \\
&\E[\pointstobrauerfn_\varep(\Sigma(T_1)) ~|~ \mc{F}_{T_1}] \E[\pointstobrauerfn_\varep(\Sigma(T) / \Sigma(T_1)) \ind(T_n \leq T) ~|~ \mc{F}_{T_1}] e^{-(4n-8) (T-T_1)} e^{-(T-T_1)} e^{-\varep(T-T_1)/N} .
\end{equs}
Now observe that
\begin{equs}
2\binom{n}{2} - \frac{n}{2}\bigg(1 - \frac{\varep}{N}\bigg) - (4n-7) - \frac{\varep}{N} &= n^2 - 5n + 6 - \frac{n}{2}\bigg(1 - \frac{\varep}{N}\bigg) + 1 - \frac{\varep}{N} \\
&= 2\binom{n-2}{2} - \frac{n-2}{2} \bigg(1 - \frac{\varep}{N}\bigg).
\end{equs}
From this, we obtain
\begin{equs}
&e^{2\binom{n}{2}T - \frac{n}{2}(1 - \frac{1}{N})T}\E[\pointstobrauerfn_\varep(\Sigma(T))  \ind(T_{n/2} \leq T) ~|~ \mc{F}_{T_1}] = \bigg(e^{2\binom{n}{2} T_1 - \frac {n}{2} (1 - \frac{1}{N})T_1} \E[ \pointstobrauerfn_\varep(\Sigma(T_1)) ~|~ \mc{F}_{T_1}] \bigg) ~\times \\
&\bigg(e^{2\binom{n-2}{2}(T - T_1) - \frac{n-2}{2}(1 - \frac{\varep}{N}) (T - T_1)} \E[\pointstobrauerfn_\varep(\Sigma(T) \backslash \Sigma(T_1)) \ind(T_{n/2} \leq T) ~|~ \mc{F}_{T_1}] \bigg).
\end{equs}
At this point, we recognize that the second factor is exactly given by the inductive assumption:
\begin{equs}
e^{2\binom{n-2}{2}(T - T_1) - \frac{n-2}{2}(1 - \frac{\varep}{N}) (T - T_1)} &\E[\pointstobrauerfn_\varep(\Sigma(T) \backslash \Sigma(T_1)) \ind_{A_T} \ind(T_{n/2} \leq T) ~|~ \mc{F}_{T_1}] = \\
&\E[ \pointstobrauerfn_\varep(\mc{Q}_{T_{n/2}} \backslash \mc{Q}_{T_1}) \ind(T_{n/2} \leq T) h_{n-2}(T_2, \ldots, T_{n/2}) ~|~ \mc{F}_{T_1}]
\end{equs}
Thus, to finish the induction, we just need to show that
\begin{equs}
e^{2\binom{n}{2} T_1 - \frac {n}{2} (1 - \frac{\varep}{N})T_1} \E[ \pointstobrauerfn_\varep(\Sigma(T_1)) \pointstobrauerfn(\mc{Q}_{T_{n/2}} \backslash \mc{Q}_{T_1}) ~|~ \mc{F}_{T_{n/2}}] = e^{(2(n-1) - (1 - \frac{\varep}{N}))T_1} F_\varep(\mc{Q}_{T_1})\pointstobrauerfn_\varep(\mc{Q}_{T_{n/2}} \backslash \mc{Q}_{T_1}) .
\end{equs}
Again, this follows by accounting for the contributions before time $T_1$ of all the swaps and turnarounds not involving the top strand. There are a total of $2\binom{n-1}{2}$ such processes. Out of these, there are $\frac{n-2}{2}$ processes which contribute 1 (the turnarounds between two strands which are matched on the right), there are $\frac{n-2}{2}$ processes which contribute $e^{-T_1} e^{-\varep T_1/N}$ (the swaps between two strands which are matched on the right), and every other process must have zero points on $[0, T_1]$, and thus contributes $e^{-T_1}$. In total, we get
\begin{equs}
\E[ \pointstobrauerfn_\varep(\Sigma(T_1)) \pointstobrauerfn_\varep(\mc{Q}_{T_{n/2}} \backslash \mc{Q}_{T_1}) ~|~ \mc{F}_{T_{n/2}}]  = e^{-(2\binom{n-1}{2} - (n-2))T_1} e^{-\frac{n-2}{2} T_1} e^{-\varep\frac{n-2}{2} T_1/N} F_\varep(\mc{Q}_{T_1}) \pointstobrauerfn_\varep(\mc{Q}_{T_{n/2}} \backslash \mc{Q}_{T_1}) .
\end{equs}
To finish, note that (using that $2\binom{n}{2} - 2\binom{n-1}{2} = 2(n-1)$)
\[
2\binom{n}{2} - \frac{n}{2}\bigg(1 - \frac{\varep}{N}\bigg) - \bigg(2\binom{n-1}{2} - (n-2)\bigg) - \frac{n-2}{2} - \frac{n-2}{2} \frac{\varep}{N} = 2(n-1) - \bigg(1 - \frac{\varep}{N}\bigg). \qedhere
\]
\end{proof}

Next, we give an explicit expression for the right\revision{-}hand side of Proposition \ref{prop:os-strand-by-strand}. First, let $E_{\pi_0}$ be the event that in the exploration process, $n$ gets matched to $n-1$, $n-2$ gets matched to $n-3$, $\ldots$, $2$ gets matched to $1$. In other words, $E_{\pi_0}$ is the event that the left matching discovered by our exploration is $\pi_0$ (which was defined in Definition \ref{def:pi-0}, see also Figure \ref{figure:pi-0}).

% \begin{notation}
% Let $\pi_0 : [n] \ra [n]$ be the matching given by $\{\{n, n-1\}, \{n-2, n-3\}, \ldots, \{2, 1\}\}$. One may visualize $\pi_0$ as in Figure \ref{figure:pi-0}.
% \begin{figure}[ht!]
%     \centering
% \includegraphics[width=.1\textwidth]{figures/pi-0.pdf}
%     \caption{$\pi_0$ when $n = 6$}
%     \label{figure:pi-0}
% \end{figure}
% \end{notation}

% Note that the event $E$ is the event that the left matching discovered by our exploration is $\pi_0$.

For notational brevity, we make the following definition.

\begin{notation}
Define
\begin{equs}
I(T, n) := \int_0^T du_1 \int_0^{T - u_1} du_2 \cdots \int_0^{T - (u_1 + \cdots + u_{n/2 - 1})} du_{n/2} & \big(e^{-u_1} e^{-\varep u_1 J_{n-1} / N}\big) \times \big(e^{-u_2} e^{-\varep u_2 J_{n-3} / N} \big) \\
& ~\times \cdots \times \big( e^{-u_{n/2}} e^{-\varep u_{n/2} J_1 / N}\big).
\end{equs}
\end{notation}

We will need the following lemma.

\begin{lemma}\label{lemma:pi-0-J-n}
We have that
\begin{equs}
~[\pi_0 ~ \pi_0] J_n = [\pi_0 ~ \pi_0] (1 + J_{n-1}).
\end{equs}
\end{lemma}
\begin{proof}
In the case $n = 4$, Figure \ref{figure:pi_0_J_n} contains the proof by explicitly identifying the terms which appear on the left and right\revision{-}hand sides of the claimed identity.
\begin{figure}[ht!]
    \centering
\includegraphics[width=0.7\textwidth]{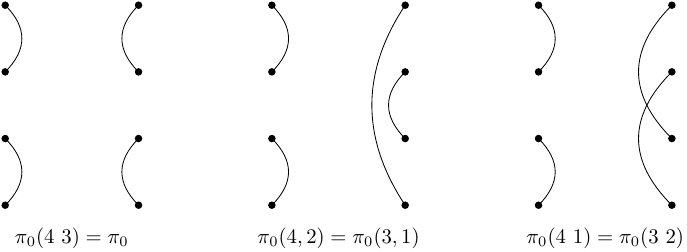}
    \caption{Term-by-term identification of the three terms in each of $\pi_0 J_4$ and $\pi_0(1 + J_3)$.}\label{figure:pi_0_J_n}
\end{figure}
The proof when $n \geq 6$ is essentially same as the case $n = 4$. The case $n =2$ is trivial.
\end{proof}

\begin{lemma}\label{lemma:os-strand-by-strand-expression}
We have that
\begin{equs}
\E\big[\pointstobrauerfn_\varep(\mc{Q}_{T_{n/2}}^{\mrm{OS}})\ind(T_{n/2} \leq T) h_n(T_1, \ldots, T_{n/2})\ind_{E_{\pi_0}} \big] =  (\varep N)^{-n/2} [\pi_0 ~ \pi_0] I (T, n).
\end{equs}
% where
% \begin{equs}
% I = \int_0^T du_1 \int_0^{T - u_1} du_2 \cdots \int_0^{T - (u_1 + \cdots + u_{n/2 - 1})} du_{n/2} & \big(e^{-u_1} e^{-\varep u_1 J_{n-1} / N}\big) \times \big(e^{-u_2} e^{-\varep u_2 J_{n-3} / N} \big) \\
% & ~\times \cdots \times \big( e^{-u_{n/2}} e^{-\varep u_{n/2} J_1 / N}\big).
% \end{equs}
\end{lemma}
\begin{proof}
By considering an alternative exploration as in Section \ref{section:strand-by-strand}, we may explicitly compute
\begin{equs}
\E\big[\pointstobrauerfn_\varep(\mc{Q}_{T_{n/2}}^{\mrm{OS}}) \ind(T_{n/2} \leq T) h_n(T_1, \ldots, T_{n/2}) \ind_{E_{\pi_0}}  \big] = (\varep N)^{-n/2}[\pi_0 ~ \pi_0]  I',
\end{equs}
where
\begin{equs}
I' = \int_0^T& du_1 \int_0^{T-u_1} du_2 \cdots \int_0^{T - (u_1 + \cdots + u_{n/2 - 1})} du_{n/2} \big( e^{-2(n-1) u} e^{-\varep u_1 J_n/N}\big) ~\times \\
&\big( e^{-2(n-3) u_2} e^{-\varep u_2 J_{n-2}/N}\big) \times \cdots \times \big( e^{-2u_{n/2}} e^{-\varep u_{n/2} J_2/N}\big) h_n(u_1, \ldots, u_{n/2}).
\end{equs}
% \begin{equs}
% I' = &\int_0^T du_1 \int_0^{T-u_1} du_2 \cdots \int_0^{T - (u_1 + \cdots + u_{n/2 - 1})} du_{n/2} e^{(2(n-1) - (1-\frac{\varep}{N})) u_1}e^{-(n-1) u_1} e^{-(n-1) u_1} e^{-\varep u_1 J_{n}/N}  \\
% & ~\times e^{(2(n-3) - (1 - \frac{\varep}{N})) u_2} e^{-(n-3) u_2} e^{-(n-3) u_2} e^{-\varep u_2 J_{n-2} / N} \\
% & \times \cdots \times e^{(2(1) - (1 - \frac{\varep}{N})) u_{n/2}} e^{-u_{n/2}} e^{-u_{n/2}} e^{-\varep u_{n/2} J_2 / N}.
% \end{equs}
We explain how $e^{-2(n-1)u} e^{-\varep u_1 J_n/N}$ arises. One factor of $e^{-(n-1) u}$ comes from the density of $T_1$, which is an exponential random variable with rate $n-1$. Conditioned on $T_1 = u_1$, we need to average over all swaps involving the top strand before time $u$, which contributes $e^{-(n-1) u} e^{-u_1\varep J_n/N}$. 

Next, observe that the formula for $I'$ may be simplified to
\begin{equs}
I' = \int_0^T du_1 &\int_0^{T - u_1} du_2 \cdots \int_0^{T - (u_1 + \cdots + u_{n/2-1})} du_{n/2} \big(e^{-(1 - \frac{\varep}{N}) u_1} e^{-\varep u_1 J_n / N}\big) ~\times \\
&\big(e^{-(1 - \frac{\varep}{N}) u_2} e^{-\varep u_2 J_{n-2}/N} \big) \cdots \times \big(e^{-(1 - \frac{\varep}{N}) u_{n/2}} e^{-\varep u_{n/2} J_2/N}\big).
\end{equs}
Next, by Lemma \ref{lemma:pi-0-J-n} we have that $[\pi_0 ~ \pi_0] J_n = [\pi_0 ~ \pi_0] (1 + J_{n-1})$. Since the Jucys-Murphy elements commute, we may then obtain that
\begin{equs}
~ [\pi_0 ~ \pi_0] J_n^k = [\pi_0 ~ \pi_0] (1 + J_{n-1})^k \text{ for all $k \geq 0$.}
\end{equs}
This implies that $[\pi_0 ~ \pi_0] e^{u J_n} = e^u [\pi_0 ~ \pi_0] e^{u J_{n-1}}$ for any $u \in \R$. More generally, by the same argument as in the proof of Lemma \ref{lemma:pi-0-J-n}, we may obtain that $[\pi_0 ~ \pi_0] e^{u J_{2r}} = e^u [\pi_0 ~ \pi_0]  e^{u J_{2r-1}}$ for $1 \leq r \leq n/2$. By applying these identities, we obtain $[\pi_0 ~ \pi_0] I' = [\pi_0 ~ \pi_0] I(T, n)$, and the desired result follows.
\end{proof}

Now suppose that the left matching discovered by our strand-by-strand exploration is some arbitrary $\pi$. Then we may first permute the strands so that $\pi$ becomes $\pi_0$, apply Lemma \ref{lemma:os-strand-by-strand-expression} to compute the expectation of our strand-by-strand-exploration in the case where the left matching is $\pi_0$, and then permute the strands back. This gives a formula for the expectation of our strand-by-strand exploration in the case where the left matching is $\pi$, which then leads to the following lemma. Recall from Definition \ref{def:sigma-pi} that for each matching $\pi : [n] \ra [n]$, we fixed a permutation $\sigma_\pi \in \symgrp_n$ such that $\sigma_\pi [\pi ~ \pi] \sigma_\pi^{-1} = [\pi_0 ~ \pi_0]$.

\begin{lemma}\label{lemma:os-strand-by-strand-expression-general}
We have that
\begin{equs}
\E\big[\pointstobrauerfn_\varep(\mc{Q}_{T_{n/2}}^{\mrm{OS}}) \ind(T_{n/2} \leq T) h_n(T_1, \ldots, T_{n/2})\big] = (\varep N)^{-n/2}\sum_{\pi : [n] \ra [n]} [\pi ~ \pi_0] I(T, n) \sigma_\pi.
\end{equs}
\end{lemma}
\begin{proof}
By summing over all possible left machings $\pi$, we obtain
\begin{equs}
\E\big[\pointstobrauerfn_\varep(\mc{Q}_{T_{n/2}}^{\mrm{OS}}) \ind(T_{n/2} \leq T) h_n(T_1, \ldots, T_{n/2})\big] = (\varep N)^{-n/2} \sum_{\pi : [n] \ra [n]} [\pi ~ \pi] I(T, n, \pi),
\end{equs}
where $I(T, n, \pi)$ is the analog of $I(T, n)$ defined for a general $\pi$. Writing $[\pi ~ \pi] = \sigma_\pi^{-1} [\pi_0 ~ \pi_0] \sigma_\pi$, we may write
\begin{equs}
~[\pi ~ \pi] I(T, n, \pi) = \sigma_\pi^{-1} [\pi_0 ~ \pi_0] (\sigma_\pi I(T, n, \pi) \sigma_\pi^{-1}) \sigma_\pi.
\end{equs}
Since conjugating by $\sigma_\pi$ corresponds to permuting the strands according to $\sigma_\pi$, we have that $\sigma_\pi I(T, n, \pi) \sigma_\pi^{-1} = I(T, n)$. I.e., after permuting the strands according to $\sigma_\pi$, $I(T, n, \pi)$ gets taken to $I(T, n)$. To finish, note that $\sigma_\pi^{-1}[\pi_0 ~ \pi_0] = [\pi ~ \pi_0]$ (since if we only permute the labels on the left, then only the left matching gets changed).
\end{proof}

\begin{lemma}\label{lemma:rho-eps-I-T-n-limit}
We have that
\begin{equs}
\lim_{T \toinf} \rho_\varep(I(T, n)) = \varep^{n/2} N^{n/2} \rho_\varep\big((\varep N + J_{n-1}) (\varep N + J_{n-3}) \cdots(\varep N + J_1)\big)^{-1}. 
\end{equs}
\end{lemma}
\begin{proof}
We have that
\begin{equs}
\lim_{T \toinf} \rho_\varep(I(T, n)) &= \int_0^\infty du_1 e^{-u_1} e^{-\varep u_1 \rho_\varep(J_{n-1})/N} \cdots \int_0^\infty du_{n/2} e^{-u_{n/2}} e^{-\varep u_{n/2} \rho_\varep(J_1)/N} \\
&= \rho_\varep\bigg(\id + \varep \frac{J_{n-1}}{N}\bigg)^{-1} \rho_\varep \bigg(\id + \varep \frac{J_{n-3}}{N}\bigg)^{-1} \cdots \rho_\varep\bigg(\id + \varep \frac{J_1}{N}\bigg)^{-1} \\
&= \varep^{n/2} N^{n/2} \rho_\varep(\varep N + J_{n-1})^{-1} \rho_\varep(\varep N + J_{n-3})^{-1} \cdots \rho_\varep(\varep N + J_1)^{-1}.
\end{equs}
Here, the operators $\rho_\varep(\id + \varep J_{n-2k} / N)$, $1 \leq k <n/2$, have strictly positive eigenvalues (and thus are invertible) by Lemma \ref{lemma:rho-plus-Jucy-Murphy-eigenvalues} (recall that $\rho_-((i ~ j)) = - \rho((i ~ j))$ for any transposition $(i ~ j)$).
\end{proof}

% \begin{lemma}
% We have that
% \begin{equs}
% \rho_\varep((\varep N + J_{n-1}) (\varep N + J_{n-3}) \cdots (\varep N + J_1))^{-1} = \rho_\varep(\Wg_{\varep N}^{\mrm{O}}).
% \end{equs}
% \end{lemma}

\begin{lemma}\label{lemma:weingarten-extraction-step}
For any matching $\pi : [n] \ra [n]$, we have that
\begin{equs}
 \rho_\varep([\pi ~ \pi_0]) \rho_\varep\big((\varep N + J_{n-1}) (\varep N + J_{n-3}) \cdots (\varep N + J_1)\big)^{-1} \rho_\varep(\sigma_\pi) = \sum_{\pi' : [n] \ra [n]} \Wg^{\mrm{O}}_{\varep N}(\pi, \pi') \rho_\varep([\pi ~ \pi']) .
\end{equs}
\end{lemma}
\begin{proof}
Applying Lemma \ref{lemma:os-weingarten-jucys-murphy-relation}, we have that 
\begin{equs}
\rho_\varep([\pi ~ \pi_0]) \rho_\varep\big((\varep N + J_{n-1}) (\varep N + J_{n-3}) \cdots (\varep N + J_1)\big)^{-1} \rho_\varep(\sigma_\pi) = \sum_{\pi' : [n] \ra [n]} \Wg_{\varep N}^{\mrm{O}}(\pi_0, \pi') \rho_\varep([\pi ~ \pi'] \sigma_\pi).
\end{equs}
Since $[\pi ~ \pi] = \sigma_\pi^{-1} [\pi_0 ~ \pi_0] \sigma_\pi$, we have that $[\pi ~ \pi'] \sigma_\pi = \sigma_\pi^{-1} [\pi_0 ~ \pi'] \sigma_\pi$. Changing variables $\pi' = \pi' \sigma_\pi$, we obtain that the above is further equal to
\begin{equs}
\sum_{\pi' : [n] \ra [n]} \Wg_{\varep N}^{\mrm{O}}(\pi_0, \pi' \sigma_\pi^{-1}) \rho_\varep([\pi ~ \pi']).
\end{equs}
To finish, observe that $\Wg_{\varep N}^{\mrm{O}}(\pi_0, \pi' \sigma_\pi^{-1}) = \Wg_{\varep N}^{\mrm{O}}(\sigma_\pi^{-1} \pi_0, \pi') = \Wg_{\varep N}^{\mrm{O}}(\pi, \pi')$. The first identity follows since $\sigma_\pi^{-1} [\pi_0 ~ \pi' \sigma_\pi^{-1}] \sigma_\pi = [\sigma_\pi^{-1} \pi_0 ~ \pi']$, which implies that $(\pi_0, \pi' \sigma_\pi^{-1})$ has the same face profile as $(\sigma_\pi^{-1} \pi_0, \pi')$, and the second identity follows since $\sigma_\pi^{-1} \pi_0 = \pi$. 
\end{proof}

Combining Proposition \ref{prop:os-strand-by-strand} and Lemmas \ref{lemma:os-strand-by-strand-expression-general}, \ref{lemma:rho-eps-I-T-n-limit}, and \ref{lemma:weingarten-extraction-step}, we obtain the following result. 

\begin{prop}\label{prop:os-exploration-all-eras-end-large-time-limit}
We have that
\begin{equs}
\lim_{T \toinf} e^{2\binom{n}{2} T - \frac{n}{2}(1 - \frac{\varep}{N})T} &\rho_\varep\big(\E[\pointstobrauerfn_\varep(\Sigma_{\mrm{OS}}(T))  \ind(T_{n/2} \leq T)]\big) = \sum_{\pi, \pi' : [n] \ra [n]} \Wg^{\mrm{O}}_{\varep N}(\pi, \pi')\rho_\varep([\pi ~ \pi']).
\end{equs}
\end{prop}

To complete the proof of Theorem \ref{thm:os-weingarten-recovery}, we need to show that it suffices to restrict to the event that all exploration eras have finished before time $T$. This turns out to be much harder to show for $\orthogonal(N)$ than for $\SphN$ -- see Remarks \ref{remark:orthogonal-harder-jucys-murphy-estimate} and \ref{remark:orthogonal-harder} for some discussion as to why. We begin by introducing some concepts which are needed to handle the $\orthogonal(N)$ case.
% With this in mind, in the following we will devote most of our attention to the $\orthogonal(N)$ case.
% Unless explicitly mentioned, we will assume that $G = \orthogonal(N)$ in the remainder of this section.

\begin{definition}
Let $G = \orthogonal(N)$. Define $\proj_n := \E[G^{\otimes n}] \in \End((\C^N)^{\otimes n})$.
\end{definition}

By properties of Haar integration, we have that $\proj_n$ is symmetric and $\proj_n^2 = \proj_n$. Thus, $\proj_n$ is the orthogonal projection onto its image, which is precisely the subspace of $\orthogonal(N)$-invariant vectors $\{v \in (\C^N)^{\otimes n} : O^{\otimes n} v = v\}$. 
Observe moreover that with $O_0$ as in Theorem \ref{thm:os-weingarten-recovery}, we have that
\begin{equs}
\lim_{T \toinf} \E[B_T^{\otimes n}] \E[O_0^{\otimes n}] = \E[G^{\otimes n}] = \proj_n.
\end{equs}
We thus have that
\begin{equs}
\lim_{T \toinf} \E[B_T^{\otimes n}] \E[O_0^{\otimes n}] \proj_n = \proj_n^2 = \proj_n.
\end{equs}
This discussion shows that when taking limits of $\SON$ Brownian motion to recover results about $\orthogonal(N)$ Haar integration, we may first project to the space of invariant vectors, and this does not change the limit. This projection is a technical convenience that will make it easier to argue why the contribution from the case where not all eras end by time $T$ goes to zero as $T \toinf$.

With this discussion in mind, we state the following proposition.

\begin{prop}\label{prop:orthogonal-no-contribution-when-not-all-eras-end}
We have that
\begin{equs}
\lim_{T \toinf} \Big\| e^{2\binom{n}{2} T - \frac{n}{2}(1 - \frac{1}{N})T} \rho_\varep\big(\E[\pointstobrauerfn_\varep(\Sigma_{\mrm{OS}}(T))  \ind(T_{n/2} > T)]\big) \proj_n \Big\|_{op} &= 0, ~~ G = \orthogonal(N) \\
\lim_{T \toinf} \Big\| e^{2\binom{n}{2} T - \frac{n}{2}(1 + \frac{1}{N})T} \rho_{\varep}\big(\E[\pointstobrauerfn_\varep(\Sigma_{\mrm{OS}}(T))  \ind(T_{n/2} > T)]\big) \Big\|_{op} &= 0, ~~ G = \SphN.
\end{equs}
\end{prop}

We will prove this proposition by an inductive argument, which rests on the following technical lemmas. The proofs are deferred to Section \ref{section:orthogonal-exploration-technical-proofs}.

\begin{lemma}\label{lemma:orthogonal-jucys-murphy-op-norm-bound}
Let $n$ be even. For any $u \geq 0$, we have that
\begin{equs}
\|e^{-u\rho(J_n)}\|_{op} &\leq e^{(N-1)u}, \\
\| e^{-u \rho(J_n)} \proj_n\|_{op} &\leq e^{(N-2)u}.
\end{equs}
\end{lemma}

\begin{remark}\label{remark:orthogonal-harder-jucys-murphy-estimate}
The first estimate of Lemma \ref{lemma:orthogonal-jucys-murphy-op-norm-bound} immediately follows from Lemma \ref{lemma:rho-plus-Jucy-Murphy-eigenvalues}, which says that all eigenvalues of $\rho_+(J_n) = \rho(J_n)$ are at least $-N+1$. However, this estimate is not good enough for the proof of Proposition \ref{prop:orthogonal-no-contribution-when-not-all-eras-end} when $G = \orthogonal(N)$. The point of the second estimate of Lemma \ref{lemma:orthogonal-jucys-murphy-op-norm-bound} is that if we restrict to the subspace of $\orthogonal(N)$-invariant vectors (which is the effect of adding the $\proj_n$ term), then we can in fact obtain a better estimate for $\|e^{- u \rho(J_n)}\|_{op}$.
\end{remark}

\begin{lemma}\label{lemma:orthogonal-odd-strands-operator-norm-bound}
Let $n$ be even. For any $T \geq 0$, we have that
\begin{equs}
\big\| \big(I \otimes \E[B_T^{\otimes (n-1)}] \big) \proj_n \big\|_{op} &\lesssim_{N, n} T^{\frac{n}{2}-1} e^{-\frac{1}{2}(1 - \frac{1}{N})\revision{T}}, ~~ G = \orthogonal(N), \\
\big\| I \otimes \E[B_T^{\otimes (n-1)}] \big\|_{op} &\lesssim_{N, n} T^{\frac{n}{2}-1} e^{-\frac{1}{2}(1 + \frac{1}{N})\revision{T}}, ~~ G = \SphN.
\end{equs}
\end{lemma}

\begin{remark}\label{remark:orthogonal-harder}
Another reason why $\orthogonal(N)$ is more delicate than $\SphN$ may be seen in the statement of Lemma \ref{lemma:orthogonal-odd-strands-operator-norm-bound}. For $\orthogonal(N)$, we need to add in the additional projection $\proj_n$ in order to obtain the stated estimate. Indeed, in certain cases $\lim_{T \toinf} I \otimes \E[B_T^{\otimes (n-1)}]$ is not even zero -- note that this limit is equal to $I \otimes \E[S^{\otimes (n-1)}]$, where $S$ is a Haar-distributed $\SON$ random matrix. If $N$ is odd and $n-1 \geq N$ is also odd, then $\E[S^{\otimes (n-1)}]$ may be nonzero. The most direct example of this is when $n - 1 = N$, because if $\E[S^{\otimes N}]$ were equal to zero, then this would imply that any matrix entry has expectation zero:
\begin{equs}
\big(\E[S^{\otimes N}]\big)_{\mbf{i} \mbf{j}} = \E[S_{i_1 j_1} \cdots S_{i_n j_n}] = 0, ~~ \mbf{i} = (i_1, \ldots, i_N), \mbf{j} = (j_1, \ldots, j_N)\in [N]^N.
\end{equs}
This would further imply that $\E[\mrm{det}(S)] = 0$. On the other hand, $\mrm{det}(S) = 1$ deterministically. Thus $\E[S^{\otimes N}] \neq 0$. Thus to prove Lemma \ref{lemma:orthogonal-odd-strands-operator-norm-bound}, we will need to argue why we still have convergence to zero at an exponential rate, if we restrict to the subspace of $\orthogonal(N)$-invariant vectors.

On the other hand, since $-I \in \SphN$, we have by parity that $I \otimes \E[S^{\otimes (n-1)}] = 0$ if $S$ is a Haar-distributed $\SphN$ random matrix (and $n$ is even). It then isn't too hard to further prove that the convergence of $I \otimes \E[B_T^{\otimes (n-1)}]$ to zero happens at an exponential rate -- one can argue similar to the Unitary case.
\end{remark}

\begin{lemma}\label{lemma:turnaround-projection}
Let $n$ be even. We have that
\begin{equs}
\rho_+(\langle n ~ n-1 \rangle) \proj_n = I^{\otimes 2} \otimes \proj_{n-2}.
\end{equs}
\end{lemma}

In the following, we will also use without explicit reference the fact that for any $O \in \orthogonal(N)$, $O^{\otimes n}$ commutes with any element of $\rho_+(\mc{B}_n)$. As a consequence, $\proj_n$ also commutes with any element of $\rho_+(\mc{B}_n)$.

\begin{proof}[Proof of Proposition \ref{prop:orthogonal-no-contribution-when-not-all-eras-end}]
First, assume $G = \orthogonal(N)$. We proceed by induction. First, in the base case $n = 2$, we may obtain by explicit calculation
\begin{equs}
e^{2T - (1 - \frac{1}{N})T} \rho_+\big(\E[\pointstobrauerfn_+(\Sigma_{\mrm{OS}}(T)) \ind(T_1 > T)]\big) = e^{-(1 - \frac{1}{N})T} e^{-T \rho(J_2) / N}.
\end{equs}
Now by Lemma \ref{lemma:orthogonal-jucys-murphy-op-norm-bound}, we have that $\|e^{-T \rho(J_2)/N} \proj_2\|_{op} \leq  e^{(1 - \frac{2}{N})T}$. The desired result when $n = 2$ then follows by combining the two estimates.

Next, suppose the result is true for some even $n \geq 2$. Consider the case $n+2$. We first show that there is no contribution when the first era doesn't end, that is
\begin{equs}\label{eq:orthogonal-first-era-must-end}
\lim_{T \toinf} \Big\|e^{2\binom{n+2}{2} T - \frac{n+2}{2}(1 - \frac{1}{N})T} \rho_+\big(\E[\pointstobrauerfn_+(\Sigma_{\mrm{OS}}(T))  \ind(T_1 > T)]\big) \proj_{n+2} \Big\|_{op} = 0.
\end{equs}
Towards this end, consider a realization of $\Sigma_{\mrm{OS}}(T)$ on the event $T_1 > T$, as in the left of Figure \ref{figure:single_strand_first_era_doesnt_end}. By imagining that every time we see a swap involving the current strand of exploration, we ``cut and swap" the current strand and the other strand involved in the swap, we obtain a map on point configurations which preserves the law of $\Sigma_{\mrm{OS}}(T)$. After applying this map (see the right of Figure \ref{figure:single_strand_first_era_doesnt_end}), we obtain another Poisson point process $\tilde{\Sigma}_{\mrm{OS}}(T)$, which has the property that all swaps which involve the first strand of exploration touch the top strand.
\begin{figure}[ht!]
    \centering
\includegraphics[page=1, width=0.7\textwidth]{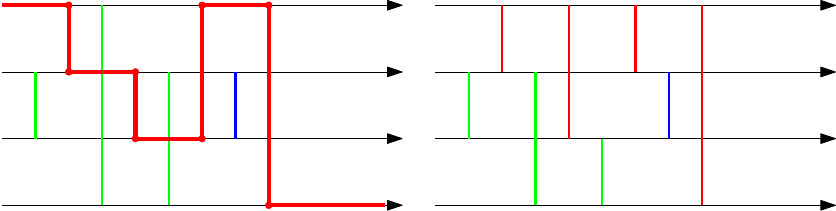}
    \caption{The green lines represent swaps, and the blue lines represent turnarounds. On the event $\{T_1 > T\}$, the exploration of the first strand makes it all the way to the right (see left). We may map the left point process into the right point process, which has the property that during the first exploration era, all swaps which are seen by the exploration involve the top strand.}
    \label{figure:single_strand_first_era_doesnt_end}
\end{figure}

To determine $\pointstobrauerfn_+(\Sigma_{\mrm{OS}}(T))$ from $\tilde{\Sigma}_{\mrm{OS}}(T)$, we split $\tilde{\Sigma}_{\mrm{OS}}(T)$ into two parts: all points not involving the top strand, and all points involving the top strand -- see Figure \ref{figure:single_strand_first_era_doesnt_end-2}. Here, the points involving the top strand must be read in reverse order.
\begin{figure}[ht!]
    \centering
\includegraphics[page=2, width=0.7\textwidth]{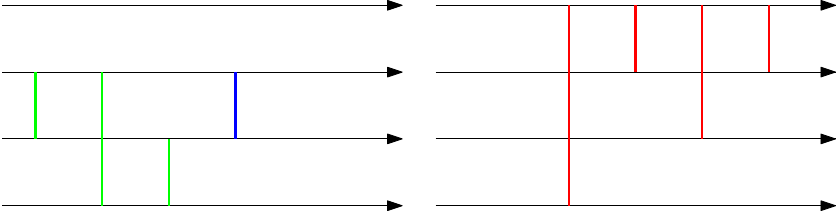}
    \caption{Left: the points of $\tilde{\Sigma}_{\mrm{OS}}(T)$ not involving the top strand. Right: the points of $\tilde{\Sigma}_{\mrm{OS}}(T)$ involving the top strand, arranged in reverse order.}
    \label{figure:single_strand_first_era_doesnt_end-2}
\end{figure}

If we now multiply together the two matchings in Figure \ref{figure:single_strand_first_era_doesnt_end-2}, we obtain the matching in Figure \ref{figure:single_strand_first_era_doesnt_end-3}, which is precisely the same matching one obtains by following all the swaps/turnaround in the original points process $\Sigma_{\mrm{OS}}(T)$ (recall the left of Figure \ref{figure:single_strand_first_era_doesnt_end}).

\begin{figure}[ht!]
    \centering
\includegraphics[page=3, width=\textwidth]{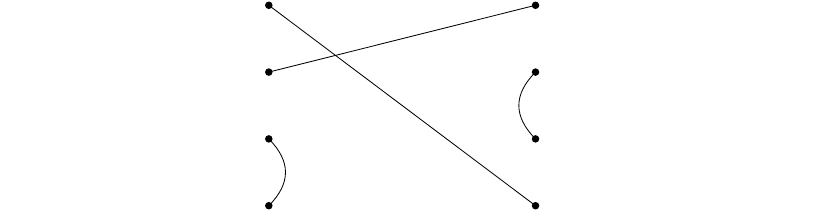}
    \caption{The matching one obtains from following all the swaps/turnarounds in $\Sigma_{\mrm{OS}}(T)$, or equivalently by first following all swaps/turnarounds not involving the top strand in $\tilde{\Sigma}_{\mrm{OS}}(T)$, and then following in reverse order all swaps involving the top strand in $\tilde{\Sigma}_{\mrm{OS}}(T)$.}
    \label{figure:single_strand_first_era_doesnt_end-3}
\end{figure}

Let $\Sigma_{\mrm{OS}}^{\mrm{top}}(T)$ be the process obtained by keeping only those points of $\Sigma_{\mrm{OS}}(T)$ which involve the top strand. Let $\Sigma^{\mrm{rest}}(T)$ be the process made of all other points, i.e. $\Sigma^{\mrm{rest}}(T) = \Sigma_{\mrm{OS}}(T) / \Sigma_{\mrm{OS}}^{\mrm{top}}(T)$.
The preceding discussion shows that 
\begin{equs}
\E[\pointstobrauerfn_+(\Sigma_{\mrm{OS}}(T))  \ind(T_1 > T)] = e^{-(n+1)T} \E[\pointstobrauerfn_+(\Sigma_{\mrm{OS}}^{\mrm{rest}}(T))] \E[\pointstobrauerfn_+(\Sigma_{\mrm{OS}}^{\mrm{top}}(T))],
\end{equs}
% where on the right\revision{-}hand side we use $\Sigma_{\mrm{OS}}$ since it has the same law as $\tilde{\Sigma}_{\mrm{OS}}(T)$. 
By an explicit calculation, we have that
\begin{equs}
e^{(n+1)T - \frac{1}{2}(1 - \frac{1}{N})T} \rho_+\big(\E[\pointstobrauerfn_+(\Sigma_{\mrm{OS}}^{\mrm{top}}(T))]\big) = e^{-\frac{1}{2}(1 - \frac{1}{N})T} e^{-T J_{n+2}/N}.
\end{equs}
We also have that
\begin{equs}
e^{(2 \binom{n+2}{2} - 2(n+1))T - \frac{n+1}{2}(1 - \frac{1}{N})T} \rho_+\big(\E[\pointstobrauerfn_+(\Sigma_{\mrm{OS}}^{\mrm{rest}}(T))]\big) = I \otimes \E[B_T^{\otimes (n+1)}].
\end{equs}
Combining, we thus obtain
\begin{equs}
e^{2\binom{n+2}{2} T - \frac{n+2}{2}(1 - \frac{1}{N})T} \rho_+\big(&\E[\pointstobrauerfn_+(\Sigma_{\mrm{OS}}(T))  \ind(T_1 > T)]\big) \proj_{n+2} \\
&= \big(I \otimes \E[B_T^{\otimes (n+1)}]\big) e^{-\frac{1}{2}(1 - \frac{1}{N})T} e^{-T \rho(J_{n+2})/N} \proj_{n+2}  \\
&= \big( \big(I \otimes \E[B_T^{\otimes (n+1)}]\big) \proj_{n+2} \big) \big(  e^{-\frac{1}{2}(1 - \frac{1}{N})T} e^{-T \rho(J_{n+2})/N} \proj_{n+2}  \big).
\end{equs}
The second identity follows since $\proj_{n+2}^2 = \proj_{n+2}$, and $\proj_{n+2}$ commutes with $\rho(J_{n+2})$. By applying Lemmas \ref{lemma:orthogonal-jucys-murphy-op-norm-bound} and \ref{lemma:orthogonal-odd-strands-operator-norm-bound}, the last term above has operator norm which is bounded by
\begin{equs}
\big\|\big(I \otimes \E[B_T^{\otimes (n+1)}]\big) \proj_{n+2}\big\|_{op} \big\| e^{-\frac{1}{2}(1 - \frac{1}{N})T} e^{-T \rho(J_{n+2})/N} \proj_{n+2}\big\|_{op} &\lesssim T^{n+2} e^{-\frac{1}{2}(1 - \frac{1}{N})T} e^{-\frac{1}{2}(1 - \frac{1}{N})T}  e^{(1 - \frac{2}{N})T} \\
&\lesssim T^{n+2} e^{-\frac{1}{N} T},
\end{equs}
which converges to zero as $T \toinf$. This shows the claim \eqref{eq:orthogonal-first-era-must-end}. 

Thus to finish, it suffices to show that
\begin{equs}
\lim_{T \toinf} \Big\| e^{2\binom{n+2}{2} T - \frac{n+2}{2}(1 - \frac{1}{N})T} \rho_+\big(\E[\pointstobrauerfn_+(\Sigma_{\mrm{OS}}(T))  \ind(T_1 \leq T, T_{(n+2)/2} > T)]\big) \proj_{n+2}\Big\|_{op} = 0. 
\end{equs}
On the event $T_1 \leq T$, we may follow the exploration until the end of the first era. Let $E$ be the event that the first era ends with the turnaround $\langle n+2 ~ n+1 \rangle$. We will focus on this case, as the case of a general turnaround may either be reduced to the case by permuting the strands, or may be similarly argued, just with more notation. By a discussion similar to that outlined in Figures \ref{figure:single_strand_first_era_doesnt_end} - \ref{figure:single_strand_first_era_doesnt_end-3}, we may compute
\begin{equs}
&e^{2\binom{n+2}{2} T - \frac{n+2}{2}(1 - \frac{1}{N})T} \rho_+\big(\E[\pointstobrauerfn_+(\Sigma_{\mrm{OS}}(T))  \ind(T_1 \leq T, T_{(n+2)/2} > T) \ind_E] \big) \proj_{n+2} \\
&=\int_0^T du ~ (I \otimes \E[B_u^{\otimes (n+1)}]) \rho_+(\langle n +2 ~ n+1 \rangle) (I^{\otimes 2} \otimes f_n(T-u)) e^{-u \rho(J_{n+2})/N} \proj_{n+2}, 
\end{equs}
where here $f_n(T-u)$, is the total partition function for a system with $n$ strands, not all exploration eras end by time $T-u$. For brevity, let $I(T)$ denote the term on the right\revision{-}hand side above. Observe that our inductive assumption implies that for any $u \geq 0$, 
\begin{equs}
\lim_{T \toinf} \|f_n(T - u) \proj_n\|_{op} = 0.
\end{equs}
To insert $\proj_n$, note by Lemma \ref{lemma:turnaround-projection} that $\rho_+(\langle n+2 ~ n+1 \rangle) \proj_{n+2} = I^{\otimes 2} \otimes  \proj_n$. Using this and the fact that $\proj_{n+2}^2 = \proj_{n+2}$, we have that
\begin{equs}
I(T) = \int_0^T du \big((I \otimes \E[B_u^{\otimes (n+1)}]) \proj_{n+2} \big) \rho_+(\langle n +2 ~ n+1 \rangle) (I^{\otimes 2} \otimes (f_n(T-u) \proj_n) ) e^{-u \rho(J_{n+2})/N} \proj_{n+2}.
\end{equs}
By the inductive assumption, the operator norm of the integrand above converges pointwise to zero as $T \toinf$. Recall also the previously obtained bound (via Lemmas \ref{lemma:orthogonal-jucys-murphy-op-norm-bound} and \ref{lemma:orthogonal-odd-strands-operator-norm-bound})
\begin{equs}
\big \| (I \otimes \E[B_u^{\otimes (n+1)}]) \proj_{n+2}  \big\|_{op} \big\| e^{-u \rho(J_{n+2})/N} \proj_{n+2}\big\|_{op} \lesssim u^{n+2} e^{-\frac{1}{N} u}.
\end{equs}
We may thus apply dominated convergence to conclude that $\lim_{T \toinf} \|I(T)\|_{op} = 0$. This finishes the proof of the inductive step. Thus the case $G = \orthogonal(N)$ is proven.

The case $G = \SphN$ follows in a similar (and indeed, simpler) fashion. By a similar discussion, in the inductive step we may obtain the following identity when the first exploration era does not end:
\begin{equs}
e^{2\binom{n+2}{2} T - \frac{n+2}{2}(1 + \frac{1}{N})T} \rho_-\big(\E[\pointstobrauerfn_-(\Sigma_{\mrm{OS}}(T)) \ind(T_1 > T)]\big) &= \big(I \otimes \E[B_T^{\otimes(n+1)}]\big) e^{-\frac{1}{2}(1 + \frac{1}{N})T} e^{T \rho_-(J_{n+2})/N} \\
&= \big(I \otimes \E[B_T^{\otimes(n+1)}]\big) e^{-\frac{1}{2}(1 + \frac{1}{N})T} e^{-T \rho(J_{n+2})/N},
\end{equs}
where in the second identity we used that (by definition) $\rho_-((i ~ j)) = -\rho((i ~ j))$ for transpositions $(i ~ j)$. Then applying Lemmas \ref{lemma:orthogonal-jucys-murphy-op-norm-bound} and \ref{lemma:orthogonal-odd-strands-operator-norm-bound}, we may bound
\begin{equs}
\Big\|\big(I \otimes \E[B_T^{\otimes(n+1)}]\big) e^{-\frac{1}{2}(1 + \frac{1}{N})T} e^{-T \rho(J_{n+2})/N}\Big\|_{op} \lesssim T^{\frac{n}{2}} e^{-(1 + \frac{1}{N})T} e^{(1 - \frac{1}{N})T} \ra 0 \text{ as $T \toinf$}.
\end{equs}
Thus as before, we may work on the event $\{T_1 \leq T, T_{(n+2)/2} > T\}$. The contribution from this event may be bounded similar to before. We omit the details.
\end{proof}

Before we combine everything and prove Theorem \ref{thm:os-weingarten-recovery}, we state the following lemma which is needed for the case $G = \orthogonal(N)$, whose proof is deferred to Section \ref{section:orthogonal-exploration-technical-proofs}.

\begin{lemma}\label{lemma:image-in-invariant-subspace}
For every pair of matchings $\pi, \pi' : [n] \ra [n]$, $\rho_+([\pi ~ \pi'])$ maps into the subspace of $\orthogonal(N)$-invariant vectors, i.e. $\mrm{Im}(\rho_+([\pi ~ \pi'])) \sse \mrm{Im}(\proj_n)$.
\end{lemma}

% We may now combine everything and prove Theorem \ref{thm:os-weingarten-recovery}.

\begin{proof}[Proof of Theorem \ref{thm:os-weingarten-recovery}]
First, consider the case $G = \orthogonal(N)$. By combining Propositions \ref{prop:os-exploration-all-eras-end-large-time-limit} and \ref{prop:orthogonal-no-contribution-when-not-all-eras-end}, we obtain
\begin{equs}
\lim_{T \toinf} \E[B_T^{\otimes n}] \E[O_0^{\otimes n}] &= \lim_{T \toinf} \E[B_T^{\otimes n}] \E[O_0^{\otimes n}] \proj_n \\
&= \lim_{T \toinf} e^{2\binom{n}{2} T - \frac{n}{2}(1 - \frac{1}{N})T}\rho_\varep\big(\E[\pointstobrauerfn(\Sigma_{\mrm{OS}}(T))  \ind(T_{n/2} \leq T)]\big) \E[O_0^{\otimes n}] \proj_n \\
&= \sum_{\pi, \pi' : [n] \ra [n]} \Wg^G_N(\pi, \pi')\rho_\varep([\pi ~ \pi'])  \proj_n.
\end{equs}
Let $A = \sum_{\pi, \pi' : [n]} \Wg_N^{\mrm{O}}(\pi, \pi') \rho_\varep([\pi ~ \pi'])$. To conclude that $\lim_{T \toinf} \E[B_T^{\otimes n}] = A$, use that (by Lemma \ref{lemma:image-in-invariant-subspace}) $\mrm{Im}(A) \sse \mrm{Im}(\proj_n)$, and $A \proj_n = \proj_n A$. This implies $A \proj_n = \proj_n A  = A$. Thus the case $G = \orthogonal(N)$ is proven. The case $G = \SphN$ follows similarly (without the extra considerations involving $\proj_n$).
\end{proof}

\subsubsection{Technical proofs}\label{section:orthogonal-exploration-technical-proofs}

In this section, we prove Lemmas \ref{lemma:orthogonal-jucys-murphy-op-norm-bound} and \ref{lemma:orthogonal-odd-strands-operator-norm-bound}. The main difficulty is in proving the estimates that involve the projection $\proj_n$, because as mentioned in Remarks \ref{remark:orthogonal-harder-jucys-murphy-estimate} and \ref{remark:orthogonal-harder}, the addition of the $\proj_n$ term leads to better estimates.
% First, we note that from Lemma \ref{lemma:rho-plus-Jucy-Murphy-eigenvalues}, all eigenvalues of $\rho(J_n)$ are at least $-N+1$, which would imply that $\|e^{-T \rho(J_n)}\|_{op} \leq e^{(N-1)T}$. However, this estimate is not good enough for the proof of Proposition \ref{prop:orthogonal-no-contribution-when-not-all-eras-end}. The point of Lemma \ref{lemma:orthogonal-jucys-murphy-op-norm-bound} is that if we restrict to the subspace of $\orthogonal(N)$-invariant vectors, then we can in fact obtain a better estimate for $\|e^{- u \rho(J_n)}\|_{op}$. In a similar vein, when $n$ is even and $N$ is odd it may be the case that $\|I \otimes \E[B_T^{\otimes (n-1)}]\|_{op}$ does not even converge to zero as $T \toinf$. On the other hand, Lemma \ref{lemma:orthogonal-odd-strands-operator-norm-bound} says that if we restrict to the subspace of $\orthogonal(N)$-invariant vectors, then the operator norm does converge to zero at an exponential rate.
We proceed to introduce the additional representation theory elements that are needed to see why these improved estimates hold. We note that everything we introduce is classical.

We first describe a spanning set for the space of $\orthogonal(N)$-invariant vectors. From classical representation theory (see e.g. \cite[Section 3]{Dahlqvist2017}), when $n$ is even, the space of $\orthogonal(N)$-invariants $\{v \in (\C^N)^{\otimes n} : O^{\otimes n} v = v\}$ is spanned by a family of vectors $\{u_\pi, \pi : [n] \ra [n]\}$ which are indexed by matchings $\pi$. The vector $u_\pi \in (\C^N)^{\otimes n}$ is given by (with implicit summation over repeated indices)
\begin{equs}
u_\pi := \prod_{\{a, b\} \in \pi} \delta^{i_a i_b} e_{i_1} \otimes \cdots \otimes e_{i_n}.
\end{equs}

\begin{remark}
Dahqlvist \cite{Dahlqvist2017} uses Brownian motion to prove this fact that $\{u_\pi, \pi : [n] \ra [n]\}$ is a spanning set for the space of $\orthogonal(N)$-invariants (i.e. the First Fundamental Theorem of invariant theory). Thus one may wonder if we are cheating a bit in using this explicit knowledge of $\orthogonal(N)$-invariants in order to Proposition \ref{prop:orthogonal-no-contribution-when-not-all-eras-end}. We don't think our argument is circular, because our focus is not to re-prove representation theory results using Brownian motion, but rather to show that our particular strand-by-strand exploration process indeed suffices to recover the Weingarten calculus.
Moreover, we find our strand-by-strand exploration intrinsically interesting, for the reasons given in Remark \ref{remark:dahlqvist-comparision}.
% because it naturally leads to a more general form of the Schwinger-Dyson/Master loop equation than has previously appeared in \cite{Chatterjee2019a, chatterjee2016, jafarov2016, shen2022new}.
\end{remark}

Observe that for any $\pi$, there exists a permutation $\sigma \in \symgrp_n$ such that $\rho(\sigma) u_{\pi_0} = u_\pi$. Indeed, recall that we previously fixed $\sigma_\pi$ such that $\sigma_\pi [\pi ~ \pi] \sigma_\pi^{-1} = [\pi_0 ~ \pi_0]$, and that visually, this had the interpretation that $[\pi ~ \pi]$ may be taken to $[\pi_0 ~ \pi_0]$ by permuting the left labels according to $\sigma_\pi$ and the right labels by $\sigma_\pi^{-1}$ -- recall Figure \ref{figure:os-sigma-pi-example}. 
% If we only look at the right\revision{-}hand side of the figure, then this says that $\sigma_\pi \pi = \pi_0$, or $\pi = \sigma_\pi^{-1} \pi_0$. Therefore $\rho(\sigma_\pi^{-1}) u_{\pi_0} = u_\pi$.
From this, we can obtain that
\begin{equs}\label{eq:rho-sigma-pi-u-pi-0}
\rho(\sigma_\pi) u_\pi = u_{\pi_0}, \text{ or } \rho(\sigma_\pi^{-1}) u_{\pi_0} = u_\pi.
\end{equs}

For matchings $\pi, \pi' : [n] \ra [n]$, the matrix elements of $\rho_+([\pi ~ \pi'])$ are given by
\begin{equs}
\rho_+([\pi ~ \pi'])_{\mbf{i} \mbf{j}} = \prod_{\{a, b\} \in \pi} \delta_{i_a i_b} \prod_{\{a, b\} \in \pi'} \delta_{j_a j_b}, ~~ \mbf{i} = (i_k, k \in [n]), \mbf{j} = (j_k, k \in [n]) \in [N]^n.
\end{equs}
The right\revision{-}hand side above is precisely $\langle u_\pi, e_{\mbf{i}} \rangle \langle u_{\pi'}, e_{\mbf{j}} \rangle$.
In other words, we have that $\rho_+([\pi ~ \pi'])$ is the rank-one matrix given by
\begin{equs}
\rho_+([\pi ~ \pi']) = u_\pi u_{\pi'}^T
\end{equs}

\begin{proof}[Proof of Lemma \ref{lemma:image-in-invariant-subspace}]
The preceding discussion shows $\mrm{Im}(\rho_+([\pi ~ \pi'])) \sse \mrm{span}(u_\pi) \sse \mrm{Im}(\proj_n)$.
\end{proof}

\begin{definition}\label{def:H-n}
Let $\mc{H}_n$ be the subgroup of $\symgrp_n$ such that $\sigma [\pi_0 ~ \pi_0] \sigma^{-1} = [\pi_0 ~ \pi_0]$. In words, $\mc{H}_n$ is the subgroup of $\symgrp_n$ which leaves $\pi_0$ fixed upon permutation of the vertices. Let $P_{\mc{H}_n} := \frac{1}{|\mc{H}_n|} \sum_{h \in \mc{H}_n} h \in \C[\symgrp_n]$.
\end{definition}

Next, we recall the following classic results from the representation theory of the symmetric group. We closely follow the discussion from \cite[Section 1.3]{zinn2009jucys}.
There is a family of group algebra elements $e_T$ indexed by standard Young tableau $T$ with $n$ boxes such that 
\begin{equs}
e_T e_{T'} = \delta_{T T'} e_T, ~~\sum_{T : |T| = n} e_T = 1.
\end{equs} 
The $e_T$ are known as Young's orthogonal idempotents. These elements have the additional property that they diagonalize the Jucys-Murphy elements. That is, 
\begin{equs}
J_k e_T = e_T J_k = c(T, k) e_T, ~~ k \in [n],
\end{equs}
where $c(T, k)$ is the content of box $k$ in $T$, i.e. $c(T, k) = j - i$ if box $k$ has coordinates $(i, j)$ in $T$. For a Young diagram $\lambda$, let $\mrm{SYT}(\lambda)$ be the set of all standard Young tableau with shape $\lambda$. Define
\begin{equs}
P_\lambda := \sum_{T \in \mrm{SYT}(\lambda)} e_T \in \C[\symgrp_n].
\end{equs}
From the given properties of $e_T$, $P_\lambda$ acts on $\C[\symgrp_n]$ as the projection onto the subspace $V_\lambda$ corresponding to the irrep $\lambda$. An explicit formula for this projection is given by
\begin{equs}\label{eq:P-lambda-character-formula}
P_\lambda = \frac{\chi_\lambda(\id)}{n!} \sum_{\sigma \in \symgrp_n} \chi_\lambda(\sigma) \sigma.
\end{equs}
Since $\chi_\lambda$ is constant on conjugacy classes, $P_\lambda$ is central, i.e. it commutes with all elements of $\C[\symgrp_n]$.

We note that for any Young diagram $\lambda \vdash n$, the matrix $\rho(P_\lambda) \in \End((\C^N)^{\otimes n})$ is the orthogonal projection onto its image. Similarly, for any Young tableau with $n$ boxes, $\rho(e_T)$ is the orthogonal projection onto its image. Moreover, the subspaces $\mrm{Im}(\rho(P_\lambda))$ and $\mrm{Im}(\rho(P_{\lambda'}))$ are orthogonal for $\lambda \neq \lambda'$. Similarly, the subspaces $\mrm{Im}(\rho(e_T)), \mrm{Im}(\rho(e_{T'}))$ are orthogonal for $T \neq T'$.

\begin{notation}
Given $\lambda$, let $2\lambda$ be the Young tableau obtained by ``doubling", i.e. by multiplying each part in the partition by $2$. 
\end{notation}

The following lemma is the key observation which leads to improved estimates for $e^{-u \rho(J_n)} \proj_n$. 

\begin{lemma}[Proposition 4 of \cite{zinn2009jucys}]\label{lemma:restriction-to-doubled-tableau}
In order for $e_T P_{\mc{H}_n} \neq 0$, $T$ must have shape $2 \lambda$ for some $\lambda \vdash \frac{n}{2}$.
\end{lemma}

\begin{lemma}\label{lemma:orthogonal-invariant-projection-doubled-young-diagram}
For all $\orthogonal(N)$-invariant vectors $v \in (\C^N)^{\otimes n}$, we have that
\begin{equs}
v = \sum_{\lambda  \vdash \frac{n}{2}} \rho(P_{2\lambda}) v.
\end{equs}
\end{lemma}
\begin{proof}
In general, we may write (recall \eqref{eq:rho-sigma-pi-u-pi-0})
\begin{equs}
v = \sum_\pi \alpha_\pi u_\pi = \sum_\pi \alpha_\pi \rho(\sigma_\pi^{-1}) u_{\pi_0} = \rho\bigg(\sum_\pi \alpha_\pi \sigma_\pi^{-1} \bigg) u_{\pi_0}.
\end{equs}
For brevity, let $X := \sum_\pi \alpha_\pi \sigma_\pi^{-1} \in \C[\symgrp_n]$. Now, since $\mc{H}_n$ stabilizes $u_{\pi_0}$, we have that 
\begin{equs}
\rho(X) u_{\pi_0} = \rho(X) \rho(P_{\mc{H}_n}) u_{\pi_0} &= \rho(X) \sum_T \rho(e_T P_{\mc{H}_n}) u_{\pi_0} = \rho(X) \sum_{\lambda \vdash \frac{n}{2}} \sum_{\lambda \in \mrm{SYT}(2\lambda)} \rho (e_T P_{\mc{H}_n}) u_{\pi_0} \\
&= \rho(X) \rho\bigg(\sum_{\lambda \vdash \frac{n}{2}} \sum_{T \in \mrm{SYT}(2\lambda)} e_T \bigg) \rho( P_{\mc{H}_n}) u_{\pi_0} \\
&= \rho(X) \rho\bigg(\sum_{\lambda \vdash \frac{n}{2}} P_{2\lambda} \bigg) u_{\pi_0} = \rho\bigg(\sum_{\lambda \vdash \frac{n}{2}} P_{2\lambda} \bigg)  \rho(X) u_{\pi_0}.
\end{equs}
In the second identity, we used Lemma \ref{lemma:restriction-to-doubled-tableau}, and in the last identity, we used that $P_{2\lambda}$ is central.
\end{proof}

\begin{proof}[Proof of Lemma \ref{lemma:orthogonal-jucys-murphy-op-norm-bound}]
The first estimate follows immediately from the fact that all eigenvalues of $\rho(J_n)$ are at least $-N + 1$ (by Lemma \ref{lemma:rho-plus-Jucy-Murphy-eigenvalues}). We proceed to prove the second estimate. By Lemma \ref{lemma:orthogonal-invariant-projection-doubled-young-diagram}, we have that
\begin{equs}
e^{-u \rho(J_n)} \proj_n = e^{-u \rho(J_n)} \sum_{\lambda \vdash \frac{n}{2}} \rho(P_{2\lambda}) \proj_n.
\end{equs}
It suffices to show that
\begin{equs}
\bigg\|e^{-u \rho(J_n)} \sum_{\lambda \vdash \frac{n}{2}} \rho(P_{2\lambda}) \bigg\|_{op} \leq e^{(N-2)u}.
\end{equs}
Since $(\rho(P_{2\lambda}), \lambda \vdash \frac{n}{2})$ is a family of projections onto orthogonal subspaces, it suffices to show that for each $\lambda \vdash \frac{n}{2}$, we have that
\begin{equs}
\big\|e^{-u \rho(J_n)} \rho(P_{2\lambda})\big\|_{op} \leq  e^{(N-2)u}.
\end{equs}
To see this, first note that in order for $\rho(P_{2\lambda}) \neq 0$, $2\lambda$ must have at most $N$ rows. Thus, we will assume that this is the case. Recalling that $P_{2\lambda} = \sum_{T \in \mrm{SYT}(2\lambda)} e_T$, and $(\rho(e_T), T \in \mrm{SYT}(2\lambda))$ is a family of projections onto orthogonal subspaces, it suffices to show that for each $T \in \mrm{SYT}(2\lambda)$, we have that
\begin{equs}
\big\|e^{-u \rho(J_n)} \rho(e_T)\big\|_{op} \leq  e^{(N-2)u}.
\end{equs}
Since $J_n e_T = c(T, n) e_T$, we have that 
\begin{equs}
e^{- u \rho(J_n)} \rho(e_T) = e^{-u c(T, n)} \rho(e_T),
\end{equs}
and so it is enough to argue that $c(T, n) \geq -(N-2)$. This follows because $T$ has shape $2\lambda$, and $2\lambda$ has at most $N$ rows, which implies that the location of $n$ in $T$ cannot be $(N, 1)$. Any other location in $T$ must have content at least $-(N-2)$.
\end{proof}

\begin{proof}[Proof of Lemma \ref{lemma:turnaround-projection}]
Since $\rho_+(\langle n ~ n - 1 \rangle)$ acts as the identity on the last $n-2$ tensor coordinates, it is enough to assume $n = 2$ and prove $\rho_+(\langle 2 ~ 1 \rangle) O^{\otimes 2} = \rho_+(\langle 2 ~ 1 \rangle)$. Since $O^{\otimes 2}$ commutes with $\rho_+(\langle 2 ~ 1\rangle)$, we have that $\rho_+(\langle 2 ~ 1 \rangle) O^{\otimes 2} = O^{\otimes 2} \rho_+(\langle 2 ~ 1 \rangle)$. Now observe that when $n = 2$, we have that $\langle 2 ~ 1 \rangle = [\pi_0 ~ \pi_0]$. Since $\rho_+([\pi_0 ~ \pi_0])$ maps into the subspace of $\orthogonal(N)$-invariants (by Lemma \ref{lemma:image-in-invariant-subspace}), it follows that $O^{\otimes 2} \rho_+([\pi_0 ~ \pi_0]) = \rho_+([\pi_0 ~ \pi_0])$.
\end{proof}

\begin{definition}
Following the notation of \cite{Dahlqvist2017}, let $\varep_N \in \C[\symgrp_N]$ be given by
\begin{equs}
\varep_N := \frac{1}{N!} \sum_{\sigma \in \symgrp_N} \mrm{sgn}(\sigma) \sigma.
\end{equs}
\end{definition}

\begin{remark}
Observe that $\varep_N$ is precisely $P_{\lambda_{\mrm{min}}}$, where $\lambda_{\mrm{min}} = (1, \ldots, 1)$ is the Young tableau corresponding to the sign representation of $\symgrp_N$.
\end{remark}

\begin{lemma}\label{lemma:kernel-annihilated}
Suppose $N \geq 3$ is odd. We have that $(I \otimes \rho(\varep_N)) \proj_{N+1}= 0$. Also, for any $1 \leq i < j \leq N$, $\varep_N \langle i ~ j \rangle = 0 \in \mc{B}_N$. Here, to be clear $\rho : \mc{B}_N \ra \End((\C^N)^{\otimes N})$.
\end{lemma}
\begin{proof}
It suffices to show that for any matching $\pi : [N+1] \ra [N+1]$, the corresponding invariant vector $u_\pi$ is annihilated by $I \otimes \rho(\varep_N)$, i.e. $(I \otimes \rho(\varep_N)) u_\pi = 0$. To see this, note that for any $\pi$, there is some pair of vertices $\{i, j\}$ matched by $\pi$, with both $i, j \leq N$ (here we use the assumption that $N \geq 3$). Since these vertices are matched, swapping them does not change the matching, and so we have that $(I \otimes \rho((i ~ j))) u_\pi = u_\pi$. On the other hand, we have that $\varep_N (i ~ j) = -\varep_N$, and thus $(I \otimes \rho(\varep_N)) \rho((i ~ j)) = I \otimes \rho(\varep_N (i ~ j)) = -I \otimes \rho(\varep_N)$. We thus have
\begin{equs}
(I \otimes \rho(\varep_N)) u_\pi = (I \otimes \rho(\varep_N)) \rho((i ~ j)) u_\pi = -(I \otimes \rho(\varep_N)) u_\pi, 
\end{equs}
and thus $(I \otimes \rho(\varep_N)) u_\pi = 0$. The second claim follows by the a similar argument, i.e. we start from the observation $(i ~ j) \langle i ~ j \rangle = \langle i ~ j \rangle$.
\end{proof}

\begin{lemma}\label{lemma:laplacian-second-minimal-eigenvalue}
We have that $\frac{1}{N} \rho(J_N + \cdots + J_1) \in \End((\C^N)^{\otimes N})$ has eigenvalue $-\frac{N}{2} + \frac{1}{2}$ with eigenspace $\mrm{Im}(\rho(\varep_N))$. All other eigenvalues are at least $-\frac{N}{2} + \frac{3}{2}$.
\end{lemma}
\begin{proof}
From the discussion in Lemma \ref{lemma:eigenvalue-lower-bound-sum-all-Jucys-Murphy}, recall that for each $\lambda \vdash N$, $\rho(P_\lambda)$ projects onto an eigenspace of $\rho(J_N + \cdots + J_1)$ with eigenvalue given by the content sum $c_\lambda$. The minimal content sum in this case is achieved when $\lambda = \lambda_{\mrm{min}} = (1, \ldots, 1)$, i.e. the Young diagram with $N$ parts of size $1$, or equivalently a single column of height $N$. The content sum in this case is
\begin{equs}
-(N-1) - (N-2) - \cdots - 1 = -\frac{N(N-1)}{2}.
\end{equs}
Thus the minimal eigenvalue of $\frac{1}{N} \rho(J_N + \cdots + J_1)$ is $-\frac{N-1}{2} = -\frac{N}{2} + \frac{1}{2}$. The associated eigenspace is $\mrm{Im}(\rho(P_{\lambda_{\mrm{min}}}))$. The first claim now follows upon recalling that $P_{\lambda_{\mrm{min}}} = \varep_N$. The next smallest eigenvalue is given by moving the box $(N, 1)$ to $(1, 2)$, i.e. by the Young diagram $\lambda = (2, 1, \ldots, 1)$. The content sum in this case is
\begin{equs}
-(N-2) - \cdots - 1 + 1 = -\frac{N(N-1)}{2} + N.
\end{equs}
The second claim now follows.
\end{proof}

\begin{definition}
Define
\begin{equs}
\Delta^n_\varep(n-1) := -\frac{(n-1)}{2}\bigg(1 - \frac{\varep}{N}\bigg) I^{\otimes n} - \frac{1}{N} I \otimes \rho_{N, n-1}(J_{n-1} + \cdots + J_1) \in \End((\C^N)^{\otimes n}).
\end{equs}
Here, we write the subscripts $\rho_{N, n-1}$ to be clear that $\rho_{N, n-1} : \symgrp_{n-1} \ra \End((\C^N)^{\otimes (n-1)})$.
\end{definition}

\begin{lemma}\label{lemma:laplacian-estimate}
For any $n \geq 2$ such that $n-1 \neq N$, we have that
\begin{equs}
\| e^{u \Delta^n_\varep(n-1)}\|_{op} \leq e^{-\frac{1}{2}(1 - \frac{\varep}{N})u}, ~~ u \geq 0.
\end{equs}
If $n-1 = N$, we have that
\begin{equs}
\| e^{u \Delta^n_1(n-1)} \proj_n \|_{op} \leq e^{-u},  ~ \| e^{u \Delta^n_1(n-1)} \rho_+(\langle 2 ~ 1 \rangle)\|_{op}  \leq Ne^{-u}, ~~ \|e^{u \Delta_{-1}^n(n-1)}\|_{op} \leq e^{-u}, ~~ u \geq 0.
\end{equs}
\end{lemma}
\begin{proof}
For brevity, write $\rho$ instead of $\rho_{N, n-1}$. For the first estimate, note that the case $\varep = -1$ readily follows from the case $\varep = 1$ because
\begin{equs}
\Delta_{-1}^n(n-1) = \Delta_1^n(n-1) - \frac{n-1}{N}.
\end{equs}
Thus, we focus on the case $\varep = 1$. Define
\begin{equs}\label{eq:Delta-plus-Delta-minus-relation}
\Delta_1^{n-1}(n-1) := -\frac{(n-1)}{2}\bigg(1 - \frac{1}{N}\bigg) I^{\otimes (n-1)} - \frac{1}{N} \rho(J_{n-1} + \cdots + J_1) \in \End((\C^N)^{\otimes (n-1)})
\end{equs}
Then $\Delta_1^n(n-1) = I \otimes \Delta_1^{n-1}(n-1)$.
Thus, it suffices to just look at $\Delta_1^{n-1}(n-1)$. 
% By the discussion in Lemma \ref{lemma:eigenvalue-lower-bound-sum-all-Jucys-Murphy}, the eigenvalues of $\Delta_1^{n-1}(n-1)$ are given by the content sums $c(\lambda)$ of Young diagrams $\lambda \vdash n-1$ with at most $N$ rows. When $n -1 = mN + r$, the minimal value of $\frac{1}{N} c(\lambda)$ across all such Young diagrams is
By Lemma \ref{lemma:eigenvalue-lower-bound-sum-all-Jucys-Murphy}, the eigenvalues of $\rho(J_{n-1} + \cdots + J_1)$ are lower-bounded by
\begin{equs}
- \frac{1}{2} (n-1) + \frac{1}{2} m^2 + \frac{1}{2} r - \frac{1}{2} \frac{r(r-1)}{N} + \frac{mr}{N}.
\end{equs}
From this, it follows that all eigenvalues of $\Delta^{n-1}_1(n-1)$ are at most
\begin{equs}
-\frac{n-1}{2} \bigg(1 - \frac{1}{N}\bigg) + \frac{1}{2}(n-1) - \frac{1}{2} m^2 - \frac{1}{2} r + \frac{1}{2} \frac{r(r-1)}{N} - \frac{mr}{N}.
\end{equs}
Using that $n-1 = mN + r$, this may be simplified to
\begin{equs}
\frac{1}{2}m(1-m) - \frac{1}{2} r \bigg(1 - \frac{r}{N}\bigg) - \frac{mr}{N}.
\end{equs}
If $(m, r) \neq (1, 0)$, then the above is easily seen to be at most $-\frac{1}{2}(1 - \frac{1}{N})$, which implies
\begin{equs}
\|e^{u \Delta^n_1(n-1)} \|_{op} = \|e^{u (I \otimes \Delta^{n-1}_1(n-1))}\|_{op} = \|e^{u \Delta_1^{n-1}(n-1)}\|_{op} \leq e^{-\frac{1}{2}(1 - \frac{1}{N})u}.
\end{equs}
If $(m, r) = (1, 0)$, then $n-1 = N$. 
% and the minimal content sum corresponds to the Young diagram given by $\lambda_{\mrm{min}} = (1, \ldots, 1)$ (i.e. there is a single column of length $N$). The representation of $\symgrp_N$ corresponding to this diagram is the sign representation, and thus $\chi_{\lambda_{\mrm{min}}}(\sigma) = \mrm{sgn}(\sigma)$. Thus $\varep_N = \frac{d_{\lambda_{\mrm{min}}}}{N!} \sum_{\sigma \in \symgrp_N} \chi_{\lambda_{\mrm{min}}}(\sigma) \sigma = P_{\lambda_{\min}}$, which implies that $\rho(\varep_N)$ is the orthogonal projection onto the subspace of $(\C^N)^{\otimes N}$ on which $\rho(\sigma)$ acts as $\mrm{sgn}(\sigma)$ for each $\sigma \in \symgrp_N$. 
We may split
\begin{equs}
e^{T \Delta^{N}_1(N)} = e^{T \Delta^{N}_1(N)} (1 - \rho(\varep_N)) + \rho(\varep_N).
\end{equs}
By Lemma \ref{lemma:laplacian-second-minimal-eigenvalue}, we have that on $\mrm{Im}(1 - \rho(\varep_N))$, all eigenvalues of $\Delta^N_1(N)$ are at most $-1$, and thus
% $e^{T \Delta^{n-1}_1(n-1)}(1 - \rho(\varep_N))$ 
% The eigenvalue of $\Delta^{n-1}_1(n-1)$ on any eigenspace corresponding to $\lambda \neq \lambda_0$ is at most $-1$. This follows because besides after $\lambda_{\mrm{min}}$, the Young diagram with the next smallest content sum is the one by shifting the box $(N, 1)$ (i.e. the bottom box) to $(1, 2)$. This results in a net gain of $N$, so upon multiplying by $-\frac{1}{N}$, we decrease the eigenvalue from $0$ to $-1$. 
% This implies that 
\begin{equs}
\big\| e^{T \Delta^{N}_1(N)} (1 - \rho(\varep_N)) \big\|_{op} \leq e^{-T}.
\end{equs}
By Lemma \ref{lemma:kernel-annihilated}, we have that for $M = \proj_{N+1}$ or $\rho_+(\langle i ~ j \rangle)$, $(I \otimes \rho(\varep_N)) M = 0$. Combining these two, it follows that
\begin{equs}
\big\|e^{T \Delta^{n}_1(n-1)} M\big\|_{op} \leq e^{-T} \|M\|_{op}.
\end{equs}
We have that $\|\proj_{N+1}\|_{op} \leq 1$ since $\proj_{N+1}$ is an orthogonal projection. Since $\langle i ~ j \rangle^2 = N \langle i ~ j \rangle$, we obtain $\rho_+(\langle i ~ j\rangle)^2 = N \rho_+(\langle i ~ j \rangle)$, which implies $\|\rho_+(\langle i ~ j \rangle) \|_{op} = N$.

Finally, when $n - 1 = N$, we have by \eqref{eq:Delta-plus-Delta-minus-relation} that $\Delta_{-1}^n(n-1) = \Delta_1^n(n-1) - 1$. By the preceding discussion, all eigenvalues of $\Delta_1^n(n-1)$ are at most $0$, and thus by equation \eqref{eq:Delta-plus-Delta-minus-relation} all eigenvalues of $\Delta_{-1}^n(n-1)$ are at most $-1$, and thus the estimate $\|e^{u \Delta^n_{-1}(n-1)}\|_{op} \leq e^{-u}$ immediately follows.
\end{proof}

\begin{proof}[Proof of Lemma \ref{lemma:orthogonal-odd-strands-operator-norm-bound}]
First, consider the case $G = \orthogonal(N)$. We proceed by induction. First, in the base case $n = 2$, we have that $\E[B_u] = e^{-\frac{1}{2}(1 - \frac{1}{N})u}$, and so
\begin{equs}
I \otimes \E[B_u] = e^{-\frac{1}{2} ( 1- \frac{1}{N})u} I^{\otimes 2}.
\end{equs}
The desired estimate in this case immediately follows.

Now, suppose that the result is true for some even $n \geq 2$. Consider the case $n+2$. We have that
\begin{equs}
I \otimes \E[B_u^{\otimes (n+1)}] = e^{2\binom{n+2}{2}T - \frac{n+2}{2}(1 - \frac{1}{N})T} \rho_+\big(\E[\pointstobrauerfn_+(\Sigma_{\mrm{OS}}(T)) \ind_{E_1}]\big),
\end{equs}
where $E_1$ is the event that there are no points touching the top strand. Let $E_2$ be the event that there is some turnaround in $\Sigma_{\mrm{OS}}(T)$. Then on the complement of $E_2$, there are only swaps, and we may compute
\begin{equs}
e^{2\binom{n+2}{2}T - \frac{n+2}{2}(1 - \frac{1}{N})T} \rho_+\big(\E[\pointstobrauerfn_+(\Sigma_{\mrm{OS}}(T)) \ind_{E_1} \ind_{E_2^c}]\big) = e^{T \Delta^{n+2}_1(n+1)}.
\end{equs}
% If $n+1 \neq N$, then $\Delta_1^{n+2}(n+1)$ has spectrum less than or equal to $-\frac{1}{2}(1 - \frac{1}{N})$, and so
% \begin{equs}
% \big\| e^{T \Delta_1^{n+2}(n+1)} \proj_{n+2}\big\|_{op} \leq \big\| e^{T \Delta_1^{n+2}(n+1)}\big\|_{op} \leq e^{-\frac{1}{2}(1 - \frac{1}{N})T}.
% \end{equs}
% If $n+1 = N$, then $\Delta^{N+1}_1(N)$ has a kernel given by $\mrm{Im}(\rho(\varep_N))$. All other eigenvectors have eigenvalue at most $-\frac{1}{2}$. Since $\rho(\varep_N) \proj_{N+1} = 0$, it follows that all eigenvalues of $e^{T \Delta^{N+1}_1(N)} \proj_{N+1}$ are at most $e^{-\frac{T}{2}}$. 
By Lemma \ref{lemma:laplacian-estimate}, we have that
\begin{equs}
\big\| e^{T \Delta^{n+2}_1(n+1)} \proj_{n+2} \big\|_{op} \leq e^{-\frac{1}{2}(1 - \frac{1}{N}) T}.
\end{equs}
Combining, we thus obtain
\begin{equs}
\big\| e^{2\binom{n+2}{2}T - \frac{n+2}{2}(1 - \frac{1}{N})T} \rho_+\big(\E[\pointstobrauerfn_+(\Sigma_{\mrm{OS}}(T)) \ind_{E_1} \ind_{E_2^c}]\big) \proj_{n+2}\big\|_{op} \leq e^{-\frac{1}{2}(1 - \frac{1}{N})T}.
\end{equs}
To finish, it suffices to show a similar estimate with $\ind_{E_2^c}$ replaced by $\ind_{E_2}$. Let $E_2^0 \sse E_2$ be the event that the first turnaround in $\Sigma_{\mrm{OS}}(T)$ is $\langle 2 ~ 1 \rangle$. We will show the estimate with $\ind_{E_2}$ replaced by $\ind_{E_2^0}$. The general estimate will follow by the same argument, just with more notation. On the event $E_2^0$, we may condition on the time of the first turnaround to obtain
\begin{equs}
e^{2\binom{n+2}{2}T - \frac{n+2}{2}(1 - \frac{1}{N})T} \rho_+\big(\E[\pointstobrauerfn_+(&\Sigma_{\mrm{OS}}(T)) \ind_{E_1} \ind_{E_2^0}]\big) = \\
&\int_0^T du ~ e^{u \Delta^{n+2}_1(n+1)} \rho_+(\langle 2 ~ 1 \rangle) \big(I \otimes \E[B_{T-u}^{\otimes (n-1)}]  \otimes I^{\otimes 2}\big).
\end{equs}
Here, the $e^{u \Delta_1^{n+2}(n+1)}$ term arises because given that the first turnaround happens at time $u$, we average over the contribution from all swaps which happen before $u$. The term $I \otimes \E[B_{T-u}^{\otimes (n-1)}] \otimes I^{\otimes 2}$ arises because once we see the turnaround $\langle 2 ~ 1 \rangle$ at time $u$, we can ignore those strands after time $u$, and only look at the top $n-2$ strands on the interval $[u, T]$.  

Now, observe that (by a variant of Lemma \ref{lemma:turnaround-projection})
\begin{equs}
\rho_+(\langle 2 ~ 1 \rangle) \proj_{n+2} = \proj_n \otimes I^{\otimes 2}. 
\end{equs}
From this, we obtain
\begin{equs}
e^{2\binom{n+2}{2}T - \frac{n+2}{2}(1 - \frac{1}{N})T} &\rho_+\big(\E[\pointstobrauerfn_+(\Sigma_{\mrm{OS}}(T)) \ind_{E_1} \ind_{E_2^0}]\big) \proj_{n+2} \\
&=\int_0^T du ~ e^{u \Delta^{n+2}_1(n+1)} \rho_+(\langle 2 ~ 1 \rangle) I \otimes \E[B_{T-u}^{\otimes (n-1)}]  \otimes I^{\otimes 2} \proj_{n+2} \\
&= \int_0^T du e^{u \Delta_1^{n+2}(n+1)} \rho_+(\langle 2 ~ 1 \rangle) \Big(\big((I \otimes \E[B_{T-u}^{\otimes (n - 1)}]) \proj_n\big) \otimes I^{\otimes 2} \Big). 
\end{equs}
By our inductive assumption, we have that
\begin{equs}
\big\|\big((I \otimes \E[B_{T-u}^{\otimes (n - 1)}]) \proj_n\big\|_{op} \lesssim (T-u)^{\frac{n}{2}-1} e^{-\frac{1}{2}(1 - \frac{1}{N})(T-u)}.
\end{equs}
By Lemma \ref{lemma:laplacian-estimate}, we have that
\begin{equs}
\big\| e^{u \Delta^{n+2}_1(n+1)} \rho_+(\langle 2 ~ 1 \rangle)\big\|_{op} \lesssim e^{-\frac{1}{2}(1 - \frac{1}{N})u}.
\end{equs}
Putting our two estimates, together, we obtain
\begin{equs}
\big\| e^{2\binom{n+2}{2}T - \frac{n+2}{2}(1 - \frac{1}{N})T} \rho_+\big(\E[\pointstobrauerfn_+(\Sigma_{\mrm{OS}}(T)) \ind_{E_1} \ind_{E_2^0}]\big) \big\|_{op} &\lesssim \int_0^T du ~(T-u)^{\frac{n}{2}-1} e^{-\frac{1}{2}(1 - \frac{1}{N})T} \\
&\lesssim T^{\frac{n+2}{2} - 1} e^{-\frac{1}{2}(1 - \frac{1}{N})T},
\end{equs}
which proves the inductive step. Thus the case $G = \orthogonal(N)$ is proven.

The case $G = \SphN$ is similar (and indeed, simpler). We sketch the changes. In the first part of the inductive step, we may compute
\begin{equs}
e^{2\binom{n+2}{2}T - \frac{n+2}{2}(1 + \frac{1}{N})T} \rho_-\big(\E[\pointstobrauerfn_-(\Sigma_{\mrm{OS}}(T)) \ind_{E_1} \ind_{E_2^c}]\big) &= e^{-\frac{n+2}{2}(1 + \frac{1}{N})T} e^{\frac{T}{N} \rho_-(J_{n+2} + \cdots + J_1)} \\
&= e^{-\frac{n+2}{2}(1 + \frac{1}{N})T} e^{-\frac{T}{N} \rho(J_{n+2} + \cdots + J_1)} = e^{T \Delta_{-1}^{n+2}(n+1)},
\end{equs}
where we used that (by definition) $\rho_-((i ~ j)) = -\rho((i ~ j))$ for transpositions $(i ~ j)$. By Lemma \ref{lemma:laplacian-estimate}, we have that $\|e^{T \Delta_{-1}^{n+2}(n+1)}\|_{op} \leq e^{-\frac{1}{2}(1 + \frac{1}{N})T}$. The contribution from the case $\ind_{E_1} \ind_{E_2}$ may be handled similar to before. We omit the details.
\end{proof}

\subsubsection{Makeenko-Migdal/Master loop/Schwinger-Dyson equation}

We next discuss the Makeenko-Migdal/Master loop/Schwinger-Dyson equation for $G = \orthogonal(N), \SphN$. First, we introduce additional string operations which appear for these groups.

\begin{definition}[Mergers, Twistings]\label{def:mergers-twistings}
Let $\bm \Gamma = (\Gamma_1, \ldots, \Gamma_k)$ be a collection of words on $\{\lambda_1, \ldots, \lambda_L\}$. Let $(i, j)$ be a location of $\bm \Gamma$. Define the set of positive and negative mergers $\merge_+((i, j), \bm \Gamma)$ and $\merge_-((i, j), \bm \Gamma)$, as well as the set of positive and negative twistings $\mbb{T}_+((i, j), \Gamma)$ and $\mbb{T}_-((i, j), \Gamma)$, as follows. Throughout, denote the letter at location $(i, j)$ by $\lambda$, and suppose $\Gamma_i = A \lambda B$.

The set of positive mergers $\merge_+((i, j), \bm \Gamma)$ is the set of collections of words $\bm \Gamma'$ obtained by merging $\Gamma_i$ with some $\Gamma_\ell$, $\ell \neq i$, in one of two ways. The first way: let $(\ell, m)$ be a location which also has letter $\lambda$. Suppose $\Gamma_\ell = C \lambda D$. Then $\Gamma_i, \Gamma_\ell$ are replaced by $A \lambda D C \lambda B$. The second way: let $(\ell, m)$ be a location which has $\lambda^{-1}$. Suppose $\Gamma_\ell = C \lambda^{-1} D$. Then $\Gamma_i, \Gamma_\ell$ are replaced by $A \lambda C^{-1} D^{-1} \lambda B$.

The set of negative mergers $\merge_-((i, j), \bm \Gamma)$ is the set of collections of words $\bm \Gamma'$ obtained by merging $\Gamma_i$ with some $\Gamma_\ell$, $\ell \neq i$, in one of two ways. The first way: let $(\ell, m)$ be a location which also has letter $\lambda$. Suppose $\Gamma_\ell = C \lambda D$. Then $\Gamma_i, \Gamma_\ell$ are replaced by $A C^{-1} D^{-1} B$. The second way: let $(\ell, m)$ be a location which has $\lambda^{-1}$. Suppose $\Gamma_\ell = C \lambda^{-1} D$. Then $\Gamma_i, \Gamma_\ell$ are replaced by $A D C B$.

The set of positive twistings $\mbb{T}_+((i, j), \bm \Gamma)$ is the set of collections of words $\bm \Gamma'$ obtained by replacing $\Gamma_i$ with another word as follows. If $\lambda^{-1}$ does not appear in $\Gamma_i$, the set $\mbb{T}_-((i, j), \bm \Gamma)$ is empty. Thus, suppose $\lambda^{-1}$ also appears in $\Gamma_i$. Let $(i, k)$ be a location which has $\lambda^{-1}$. If $k > j$ then recalling that $\Gamma_i = A \lambda B$, we may write $B = C \lambda^{-1} D$. We then replace $\Gamma_i = A \lambda C \lambda^{-1} D$ by $A \lambda C^{-1} \lambda^{-1} D$. If $k < j$ then we may write $A = E \lambda^{-1} F$. We then replace $\Gamma_i = E \lambda^{-1} F \lambda B$ by $E \lambda^{-1} F^{-1} \lambda B$.

The set of negative twistings $\mbb{T}_- ((i, j), \bm \Gamma)$ is the set of collections of words $\bm \Gamma'$ obtained by replacing $\Gamma_i$ with another word as follows. If $\lambda$ appears only once in $\Gamma_i$, the set $\mbb{T}_-((i, j), \bm \Gamma)$ is empty. Thus, suppose $\lambda$ appears at least twice in $\Gamma_i$. Denote $(i, k)$ be another location which has $\lambda$. If $k > j$ then recalling that $\Gamma_i = A \lambda B$, we may write $B = C \lambda D$. We then replace $\Gamma_i = A \lambda C \lambda D$ by $A C^{-1} D$. If $k < j$ then we may write $A = E \lambda F$. We then replace $\Gamma_i = E \lambda F \lambda C$ by $E F^{-1} C$.
\end{definition}

\begin{remark}
From the perspective of our Poisson point process on strand diagrams, the reason why the $\orthogonal(N)$ and $\SphN$ cases result in additional loop operations is because there may now be turnarounds between two same-direction strands, and swaps between two opposite-direction strands. (Recall that in the Unitary case, same-direction strands only had swaps and opposite-direction strands only had turnarounds.) 
% Note that as defined, $\merge_{\pm}((i, j), \bm \Gamma)$ contains the set $\merge_{\pm}^U((i, j), \bm \Gamma)$ from Definition \ref{def:splittings-and-mergers}. Again, this is due to the fact that there may now be turnarounds between two same-direction strands, and swaps between two opposite-direction strands.
\end{remark}

\begin{prop}[Single-location $\orthogonal(N), \SphN$ word recursion]
Let $G = \orthogonal(N), \SphN$. Let $\bm \Gamma = (\Gamma_1, \ldots, \Gamma_k)$ be a collection of words on $\{\lambda_1, \ldots, \lambda_L\}$. For any location $(i, j)$ of $\Gamma$, we have that 
\begin{equs}
\varep \bigg(1 - \frac{\varep}{N}\bigg)\E[\tr(G(\bm \Gamma))] = &- \sum_{\bm\Gamma' \in \splitting_+((i, j), \bm \Gamma)} \E[\tr(G(\bm \Gamma'))] + \sum_{\bm \Gamma' \in \splitting_-((i, j), \bm \Gamma)} \E[\tr(G(\bm \Gamma'))] \\
&- \frac{1}{N^2} \sum_{\bm \Gamma' \in \merge_+((i, j), \bm \Gamma)} \E[\tr(G(\bm \Gamma'))] + \frac{1}{N^2} \sum_{\bm \Gamma' \in \merge_-((i, j), \bm \Gamma)} \E[\tr(G(\bm  \Gamma'))] \\
&- \frac{1}{N} \sum_{\bm \Gamma' \in \mbb{T}_+((i, j), \bm \Gamma)} \E[\tr(G(\bm \Gamma'))] + \frac{1}{N} \sum_{\bm \Gamma' \in \mbb{T}_-((i, j), \bm \Gamma)} \E[\tr(G(\bm \Gamma'))]
\end{equs}
\end{prop}
\begin{proof}[Proof (sketch)]
The proof proceeds by stopping our strand-by-strand exploration process at the time of the first point, as in the proof of the $\UN$ word recursion (Proposition \ref{prop:word-recursion}). This gives a recursion for the Orthogonal Weingarten function very much analogous to the key identity \eqref{eq:weingarten-recursion}, which recall directly implied Proposition \ref{prop:word-recursion}. The main ideas are very similar but the details are a bit different -- we sketch out where the differences lie. When we explore the strand-by-strand exploration until the first point, we see either a turnaround or a swap, which may connect same-direction or opposite-direction strands. Moreover, the two strands may be part of the same word or different words. We present the two tables in Figure \ref{figure:ON_SPN_recursion_table} which indicate which of the loop operations each of these cases contributes to.
\begin{figure}[ht!]
    \centering
\includegraphics[width=0.7\textwidth]{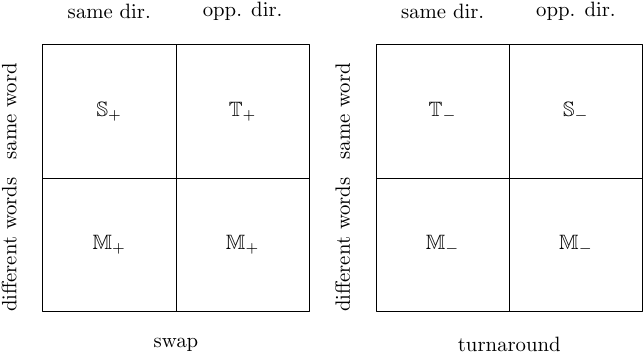}
    \caption{Left: the various cases when the first point is a swap. Right: the various cases when the first point is a turnaround.}\label{figure:ON_SPN_recursion_table}
\end{figure}
\end{proof}

The word recursion then immediately implies the Makeenko-Migdal/Master loop/Schwinger-Dyson equation. The proof is omitted, as it is very similar to the proof of the Unitary Makeenko-Migdal/Master loop/Schwinger-Dyson equation using the Unitary word recursion (see Section \ref{section:master-loop}).

\begin{theorem}[Single-location $\orthogonal(N)$ and $\SphN$ Makeenko-Migdal/Master loop/Schwinger -Dyson equation]\label{thm:os-master-loop}
Let $s = (\ell_1, \ldots, \ell_n)$ be a string. Let $(i, j)$ be a location in $s$. For $G = \orthogonal(N), \SphN$ lattice Yang-Mills theory, we have that
\begin{equs}
\varep \bigg(1 - \frac{\varep}{N}\bigg)\phi(s) = &-\sum_{s' \in \splitting_+((i, j), s)} \phi(s') + \sum_{s' \in \splitting_-((i, j), s)} \phi(s') \\
&- \frac{1}{N^2} \sum_{s' \in \merge_+((i, j), s)} \phi(s') + \frac{1}{N^2} \sum_{s' \in \merge_-((i, j), s)} \phi(s') \\
&- \frac{1}{N} \sum_{s' \in \mbb{T}_+((i, j), s)} \phi(s') + \frac{1}{N} \sum_{s' \in \mbb{T}_-((i, j), s)} \phi(s') \\
&- 
2\beta \sum_{s' \in \deform_+((i, j), s)} \phi(s') + 2\beta \sum_{s' \in \deform_-((i, j), s)} \phi(s').
\end{equs}
\end{theorem}

\begin{remark}
In the Unitary case, we had the factor $\beta$ in front of the deformation terms, whereas in the Orthogonal and Symplectic cases, we have the factor $2\beta$. This difference is ultimately due to the fact that there may be swaps and turnaround between any two strands, no matter their directions. 
\end{remark}

\subsection{Special Unitary and Special Orthogonal}\label{section:SUN-and-SON}

The Weingarten calculus for $\SUN$ and $\SON$ is far less developed than for $\UN, \orthogonal(N)$, and $\SphN$. The only formula we have seen for the $\SUN$ Weingarten function is in the physics literature \cite[Equation (20)]{borisenko2020n}. We have not seen a formula for the $\SON$ Weingarten function. Therefore, in this paper we will not do as much for $\SUN$ and $\SON$ as we did for the previous three groups. In particular, we will not recover the Weingarten calculus via large-time limits of Brownian motion. Instead, we will focus on giving surface-sum representations of Wilson loop expectations and proving the Makeenko-Migdal/Master loop/Schwinger-Dyson equations.

We first show that although we don't have explicit formulas for the $\SUN$ and $\SON$ Weingarten functions, we can still relate (via soft arguments) $\SUN$ and $\SON$ Haar expectations to some elements of the Brauer algebra. These ``Weingarten elements" will then provide the weights that appear in our surface-sum representations.

\begin{definition}
Given a matching $\pi \in \mc{M}(n)$, let $\pi^T$ be the reflection of $\pi$, or i.e. the matching obtained by swapping the vertices on the left with the vertices on the right.
\end{definition}

\begin{prop}\label{prop:sun-weingarten-exists}
For $n, m \geq 0$, there exist elements $\Wg^{\mrm{SU}}_N \in \mc{B}_{n, m}, \Wg^{\mrm{SO}}_N \in \mc{B}_{n}$ such that
\begin{equs}
\E[G^{\otimes n} \otimes \bar{G}^{\otimes m}] &= \rho_+(W^{\mrm{SU}}_N), ~~ G = \SUN, \\
\E[G^{\otimes n}] &= \rho_+(\Wg^{\mrm{SO}}_N), ~~ G = \SON.
\end{equs}
Moreover, these elements are invariant under reflection:
\begin{equs}
\Wg_N^{\bullet}(\pi) = \Wg_N^{\bullet}(\pi^T), ~~ \bullet \in \{\SUN, \SON\},
\end{equs}
as well as invariant under conjugation:
\begin{equs}
\sigma \Wg^{\mrm{SU}}_N \sigma^{-1} &= \Wg^{\mrm{SU}}_N \text{ for all $\sigma \in \symgrp_n \times \symgrp_m \sse \mc{B}_{n, m}$}, \\
\sigma \Wg^{\mrm{SO}}_N \sigma^{-1} &= \Wg^{\mrm{SO}}_N \text{ for all $\sigma \in \symgrp_n$}.
\end{equs}
\end{prop}

There are various ways one can prove this proposition. The representation-theoretic way would be to note that in the $\SUN$ case, $\E[S^{\otimes n} \otimes \bar{S}^{\otimes m}]$ commutes with $U^{\otimes n} \otimes \bar{U}^{\otimes m}$ for any $U \in \SUN$, and then to use the fact that any such operator must be of the form $\rho(W)$ for some $W \in \mc{B}_{n, m}$. The fact that $W$ may be assumed to be invariant under conjugation can be ensured by averaging over all possible conjugations, since this does not change $\E[S^{\otimes n} \otimes \bar{S}^{\otimes m}]$. The $\SON$ case can be handled similarly.

Another way to show the proposition is via $\SUN$ and $\SON$ Brownian motion. We have already introduced how $\SON$ Brownian motion is related to the Brauer algebra (Proposition \ref{prop:orthogonal-symplectic-brownian-motion-expectation}), and we will need to introduce $\SUN$ Brownian motion in order to derive the Makeenko-Migdal/Master loop/Schwinger-Dyson equations for $\SUN$. Thus we will supply a proof of Proposition \ref{prop:sun-weingarten-exists} using Brownian motion. The first step is to introduce the analog of Proposition \ref{prop:orthogonal-symplectic-brownian-motion-expectation} for $\SUN$, i.e. to state how expectations of $\SUN$ Brownian motion are related to Brauer algebra elements. We begin with the necessary setup.

% \begin{definition}[Poisson point process for $\SUN$]
We proceed to define the Poisson point process which relates to $\SUN$ Brownian motion. Let $n, m \geq 0$. Define the spaces
\begin{equs}
\mc{D}_{\mrm{SU}} := \mc{D}_{\mrm{U}} \sqcup \bigsqcup_{i, j \in [n+m]} [0, \infty), ~~ \mc{D}_{\mrm{SU}}(T) := \mc{D}_{\mrm{U}}(T) \sqcup \bigsqcup_{i, j \in [n+m]} [0, T],
\end{equs}
where here $\mc{D}_{\mrm{U}}, \mc{D}_{\mrm{U}}(T)$ are as in the Unitary case:
\begin{equs}
\mc{D}_{\mrm{U}} = \bigsqcup_{\substack{i, j \in [n+m] \\ i < j}} [0, \infty), ~~ \mc{D}_{\mrm{U}}()T = \bigsqcup_{\substack{i, j \in [n+m] \\ i < j}} [0, T].
\end{equs}
Let $\Sigma_{\mrm{SU}}$ be a rate-$1$ Poisson process on $\mc{D}_{\mrm{SU}}$. Let $\Sigma_{\mrm{SU}}(T) := \Sigma_{\mrm{SU}} \cap \mc{D}_{\mrm{SU}}(T)$ and note that $\Sigma_{\mrm{SU}}(T)$ is a rate-$1$ Poisson process on $\mc{D}_{\mrm{SU}}(T)$.

We may naturally split $\Sigma_{\mrm{SU}}(T) = \Sigma_{\mrm{U}}(T) \sqcup \Sigma_{\mrm{SU}}'(T)$, where $\Sigma_{\mrm{U}}(T)$ is a rate-$1$ Poisson process on $\mc{D}_U(T)$. Note that $\Sigma_{\mrm{U}}(T)$ is precisely the Poisson process which arises in the Unitary case. One can think of $\Sigma_{\mrm{SU}}(T)$ as starting with the processes of swaps and turnarounds as in the Unitary case, and then adding in additional independent rate-$1$ Poisson processes between any two pairs of strands, not necessarily distinct\footnote{Whereas the processes for swaps and turnarounds are always between distinct strands.}. The process $\Sigma_{\mrm{SU}}'(T)$ gives these additional points.

\begin{figure}[ht!]
    \centering
\includegraphics[width=0.5\textwidth]{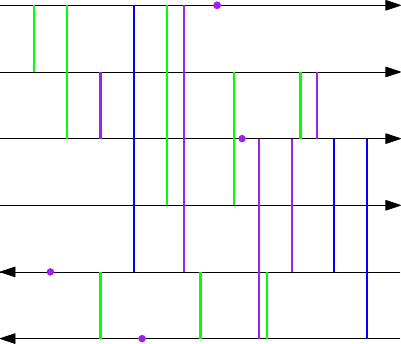}
    \caption{Example realization of $\Sigma_{\mrm{SU}}(T)$. The green lines represent swaps, blue lines represent turnarounds, and purple lines/points represent the points of $\Sigma_{\mrm{SU}}'$. In particular, a purple point on a given strand belongs to the Poisson process that encodes points between that strand and itself.}\label{figure:SUN_point_process_realization}
\end{figure}

To visualize $\Sigma_{\mrm{SU}}(T)$, consider a strand diagram with $n$ right-directed strands and $m$ left-directed strands. The points of $\Sigma_{\mrm{U}}(T)$ give swaps and turnarounds exactly as in the Unitary case. The points of $\Sigma_{\mrm{SU}}'(T)$ are represented by purple lines or points. See Figure \ref{figure:SUN_point_process_realization} for an example realization of $\Sigma_{\mrm{SU}}$.

% Define a Poisson point process $\Sigma_{\mrm{SU}}$ as follows. We imagine we have $n$ right-directed strands and $m$ left-directed strands. Let $\Sigma_{\mrm{U}}$ be the Poisson point process on this strand diagram corresponding to the Unitary case, i.e. between any two same-direction strands, there is a rate-1 Poisson process giving swaps, and between any two opposite-direction strands, there is a rate-1 Poisson process giving turnarounds. We define $\Sigma_{\mrm{SU}} := \Sigma_{\mrm{U}} \sqcup \Sigma_{\mrm{SU}}'$, where $\Sigma_{\mrm{SU}}'$ is an independent Poisson process which is made of independent rate-1 Poisson processes between any two pairs of strands, not necessarily distinct\footnote{Whereas the processes for swaps and turnarounds are always between distinct strands.}. We differentiate between the various types of points by assigning the color green to swaps, blue to turnarounds, and purple to the points of $\Sigma'_{\mrm{SU}}$.
% \end{definition}

\begin{definition}
For finite subsets $P \sse \mc{D}_{\mrm{SU}}(T)$, define $\pointstobrauerfn^{\mrm{SU}}(P) \in \mc{B}_{n, m}$ as follows. First, we decompose $P = P_{\mrm{U}} \sqcup P'$, where $P_U = P \cap \mc{D}_{\mrm{U}}(T)$. We define
\begin{equs}
\pointstobrauerfn^{\mrm{SU}}(P) := \pointstobrauerfn(P_{\mrm{U}}) \pointstobrauerfn'(P'), 
\end{equs}
where $\pointstobrauerfn(P_{\mrm{U}}) \in \mc{B}_{n, m}$ is exactly as in the Unitary case, and $\pointstobrauerfn'(P') \in \R$ is a scalar defined as follows. Let $K$ be the number of points in $P'$ between same-direction strands, and $K'$ be the number of points in $P'$ between opposite-direction strands. Then
\begin{equs}
\pointstobrauerfn'(P') := \bigg(\frac{1}{N^2}\bigg)^K \bigg(-\frac{1}{N^2}\bigg)^{K'}.
\end{equs}
% Define $\pointstobrauerfn^{\mrm{SU}}$ which maps point process realizations $P$ to elements of the Brauer algebra $\mc{B}_{n, m}$ as follows. We may split our point process realization $P = P_{\mrm{U}} \sqcup P'$, where $P_{\mrm{U}}$ collects all swaps and turnarounds, and $P'$ collects all purple points. We define $\pointstobrauerfn^{\mrm{SU}}(P) := \pointstobrauerfn(P_{\mrm{U}}) \pointstobrauerfn'(P') \in \mc{B}_{n, m}$,
% One should think of this as saying that each point in $P'$ between same-direction strands incurs a factor of $N^{-2}$, while each point in $P'$ between opposite-direction strands incurs a factor of $-N^{-2}$.
\end{definition}

\begin{remark}
One should think of this definition as saying that each point in $P'$ between same-direction strands incurs a factor of $N^{-2}$, while each point in $P'$ between opposite-direction strands incurs a factor of $-N^{-2}$.
\end{remark}

We can now state how $\SUN$ Brownian motion is related to  the Brauer algebra $\mc{B}_{n,m}$. See~\cite[Appendix A]{park2023wilson} for the proof.

\begin{prop}\label{prop:SUN-brownian-motion-expectation}
We have that
\begin{equs}
\E[B_T^{\otimes n} \otimes \bar{B}_T^{\otimes m}] = e^{(n+m)^2 T - \frac{n+m}{2}(1 + \frac{1}{N^2}) T} \rho_+\big(\E[\pointstobrauerfn^{\mrm{SU}}(\Sigma_{\mrm{SU}}(T)) ]\big), ~~ T \geq 0.
\end{equs}
\end{prop}

Proposition \ref{prop:SUN-brownian-motion-expectation} immediately implies Proposition \ref{prop:sun-weingarten-exists}, as we next show.

\begin{proof}[Proof of Proposition \ref{prop:sun-weingarten-exists}]
First, consider the $\SUN$ case. Note that the map $\rho_+ : \mc{B}_{n, m} \ra \End((\C^N)^{\otimes (n+m)})$ is a linear map into a finite-dimensional vector space. This implies that its image $\rho_+(\mc{B}_{n, m})$ is a closed subspace of $\End((\C^N)^{\otimes (n+m)})$ (as very subspace of a finite-dimensional vector space is closed).

By Proposition \ref{prop:SUN-brownian-motion-expectation}, we have that $\E[B_T^{\otimes n} \otimes \bar{B}_T^{\otimes m}] \in \rho_+(\mc{B}_{n, m})$ for all $T \geq 0$. Thus by the preceding discussion, we also have that $\E[G^{\otimes n} \otimes \bar{G}^{\otimes m}] = \lim_{T \toinf} \E[B_T^{\otimes n} \otimes \bar{B}_T^{\otimes m}] \in \rho_+(\mc{B}_{n, m})$. Therefore there exists $W \in \mc{B}_{n, m}$ such that $\E[S^{\otimes n} \otimes \bar{S}^{\otimes m}] = \rho_+(W)$. Using $W$, we construct $\Wg^{\mrm{SU}}_N$ which possesses the claimed symmetries.

Let $W^T \in \mc{B}_{n, m}$ be defined by $W^T(\pi) := W(\pi^T)$. We have that
\begin{equs}
\rho_+(W^T) &= \big(\rho_+(W)\big)^T = \big(\rho_+(W)\big)^* = \big(\E[G^{\otimes n} \otimes \bar{G}^{\otimes m}]\big)^* \\
&=\E[(G^*)^{\otimes n} \otimes \ovl{G^*}^{\otimes m}] = \E[(G^{-1})^{\otimes m} \otimes \ovl{G^{-1}}^{\otimes m}]  \\
&= \E[G^{\otimes n} \otimes \bar{G}^{\otimes m}] = \rho_+(W),
\end{equs}
where the first identity follows by the definition of $\rho_+$, the second follows because $\rho_+$ has real-valued matrix entries, the fourth follows by linearity, the fifth follows since $G^* = G^{-1}$ when $G \in \SUN$, and the sixth follows by the inversion-invariance of Haar measure on compact groups.

Next, for any $\sigma \in \symgrp_n \times \symgrp_m$, we have that
\begin{equs}
\rho_+(\sigma W \sigma^{-1}) = \rho_+(\sigma) \E[S^{\otimes n} \otimes \bar{S}^{\otimes m}]\rho_+(\sigma)^{-1} = \E[S^{\otimes n} \otimes \bar{S}^{\otimes m}] .
\end{equs}
We may thus define
\begin{equs}
\Wg^{\mrm{SU}}_N := \frac{1}{n! m!} \sum_{\sigma \in \symgrp_n \times \symgrp_m} \sigma \bigg(\frac{W + W^T}{2}\bigg) \sigma^{-1} \in \mc{B}_{n, m}. 
\end{equs}
Then $\Wg^{\mrm{SU}}_N$ satisfies all the required properties. This shows the $\SUN$ case. The $\SON$ case follows in the exact same manner.
\end{proof}

Suppose we have a collection of words $\bm \Gamma = (\Gamma_1, \ldots, \Gamma_n)$ on letters $\{\lambda_1, \ldots, \lambda_L\}$, along with a collection of Brauer algebra elements $\bm \pi = (\pi_\ell, \ell \in [L])$, where $\pi_\ell \in \mc{B}_{n_\ell + m_\ell}$, where $n_\ell, m_\ell$ are the respective number of times $\lambda_\ell, \lambda_\ell^{-1}$ appears in $\bm \Gamma$. In the $\SUN$ case, we may further assume $\pi_\ell \in \mc{B}_{n_\ell, m_\ell} \sse \mc{B}_{n_\ell + m_\ell}$. Now as noted in previous sections, the choice of $\bm \Gamma$ specifies the exterior connections of the strand diagram, while the choice of $\bm \pi$ specifies the interior connections. Let $\numcomp(\bm \Gamma, \bm \pi)$ be the number of components of the graph one obtains by including both the exterior and exterior connections. This slightly generalizes our previous definition of $\numcomp$ to the case where the $\pi_\ell$ are not of the special form of a combined left and right matching.

Proposition \ref{prop:sun-weingarten-exists} implies the following proposition about $\SUN$ and $\SON$ word expectations. The proof is essentially the same as the discussion in Section \ref{section:poisson-exploration-general-N-proof}, and thus it is omitted.

\begin{prop}\label{prop:SUN-word-expectation}
Let $\bm \Gamma = (\Gamma_1, \ldots, \Gamma_n)$ be a collection of words on letters $\{\lambda_1, \ldots, \lambda_L\}$. Then
\begin{equs}
\E[\Tr(G(\bm \Gamma))] &= \sum_{\bm \pi = (\pi_\ell, \ell \in [L])} \bigg( \prod_{\ell \in L} \Wg^{\mrm{SU}}_N(\pi_\ell) \bigg) N^{\numcomp(\bm \Gamma, \bm \pi)}, ~~ G = \SUN, \\
\E[\Tr(G(\bm \Gamma))] &=  \sum_{\bm \pi = (\pi_\ell, \ell \in [L])} \bigg( \prod_{\ell \in L} \Wg^{\mrm{SO}}_N(\pi_\ell) \bigg) N^{\numcomp(\bm \Gamma, \bm \pi)}, ~~ G = \SON.
\end{equs}
Here, the first sum is over $\pi_\ell \in \mc{B}_{n_\ell, m_\ell}$, $\ell \in [L]$, where $n_\ell, m_\ell$ are the respective number of times that $\lambda_\ell, \lambda_\ell^{-1}$ appear in $\bm \Gamma$, and the second sum is over $\pi_\ell \in \mc{B}_{n_\ell + m_\ell}$, $\ell \in [L]$.
\end{prop}

\begin{remark}
Unlike in the Unitary case, the collection of words $\bm \Gamma$ is not required to be balanced when $G = \SUN$, or unoriented-balanced when $G = \SON$. Ultimately, this is due to the fact that $\E[G^{\otimes n} \otimes \bar{G}^{\otimes m}]$ may be nonzero even if $m \neq n$ when $G = \SUN$, and $\E[G^{\otimes n}]$ may be nonzero for odd $n$ when $G = \SON$. Recall Remark \ref{remark:orthogonal-harder} for an example of the latter. Ultimately, the reason for this is because the elements of $\SUN, \SON$ must have determinant $1$.
\end{remark}

\begin{figure}[ht!]
    \centering
\includegraphics[width=.9\textwidth]{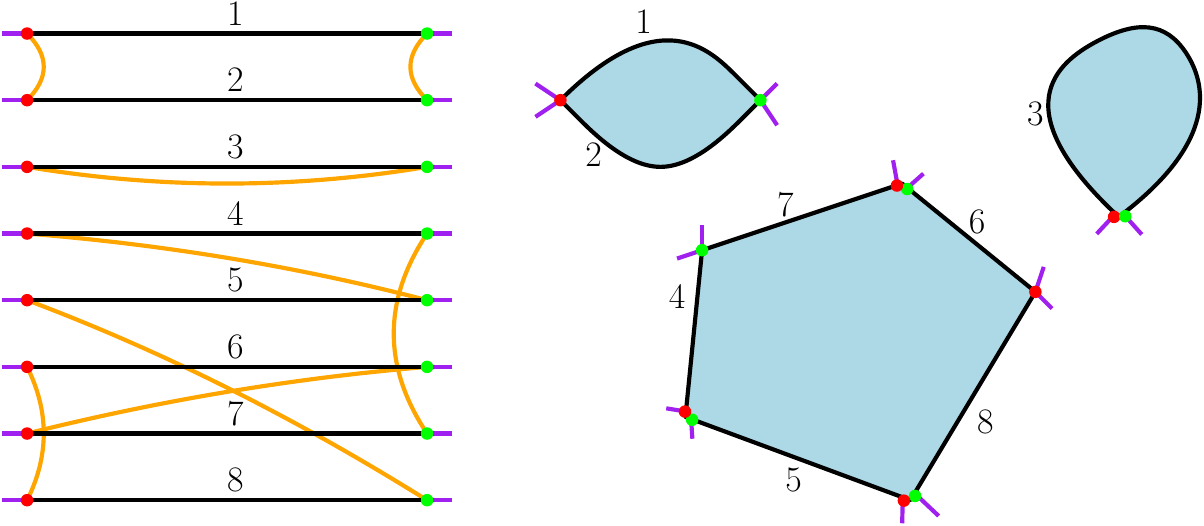}
    \caption{{\bf {\em A priori} setting:} Imagine that at first, the oranges paths at the left are not present and each black edge represents the boundary of a (not shown) yellow face. In this {\em a priori} picture every purple segment on the left is connected to a purple segment on the right by a horizontal black line.  {\bf Constructing blue faces from orange matching:} Let the orange curves indicate an arbitrary matching of the 16 red and green vertices.  There are certain cycles obtained by alternating between orange paths and black paths; if we shrink the orange paths to points, these cycles become the polygons shown on the right, whose interiors are shaded blue. If an orange edge on the left connects a green and red vertex, then the corresponding vertex on the right is colored \revision{both} red and green. {\bf Gluing interpretation:} If we start with the {\em a priori} set up and then {\em glue the blue faces into the diagram} this has the effect of changing the purple-to-purple matching from the black one to the orange one.}
    \label{fig::oddbluefaces}
\end{figure}

Next, we apply Proposition \ref{prop:SUN-word-expectation} to give a surface-sum expression for $\SUN$ and $\SON$ lattice gauge theories. To do this, we need to explain how an arbitrary element of $\mc{B}_{n, m}$ or $\mc{B}_{n+m}$ can be interpreted as giving a collection of blue faces that are glued in to the existing yellow faces. This contrasts with all the previous cases, where we could restrict to those elements of $\mc{B}_{n, n}$ (Unitary) or $\mc{B}_n$ (Orthogonal and Symplectic) which only pair vertices on the same side. We cannot do the same here, because the element $\Wg^{\mrm{SU}}_N$ (resp. $\Wg^{\mrm{SO}}_N$ may in general give nonzero weight to elements of $\mc{B}_{n, m}$ (resp. $\mc{B}_{n+m}$) which are not of this special form. In Figure \ref{fig::oddbluefaces}, we explain how to go from the interior connections specified by $\bm \pi$ to a collection of faces with specified gluings.

For $\pi \in \mc{B}_{n+m}$, we define the ``face profile" of $\pi$ to be the collection of faces that one obtains from $\pi$, as described in Figure \ref{fig::oddbluefaces}, additionally with the coloring of the vertices by red, green, or red and green, as specified in the figure. 

\begin{remark}\label{remark:SUN-SON-face-profile}
The invariance of $\Wg^{\mrm{SU}}_N$ under conjugation implies that it is a function of the face profile of $\pi \in \mc{B}_{n, m} \sse \mc{B}_{n+m}$, and similarly the invariance of $\Wg^{\mrm{SO}}_N$ under conjugation implies that it is a function of the face profile of $\pi \in \mc{B}_{n+m}$. To see this, note that invariance under conjugation is the same as invariance under permutation of the strands (where in the $\SUN$ case, we mean invariance under separate permutations of the top right-directed strands and bottom left-directed strands). If $\pi, \pi'$ have the same face profile, then there exists a permutation of the strands which takes one to the other. 

Put another way, starting only from the blue faces in Figure \ref{fig::oddbluefaces}, we may reconstruct a matching $\pi'$ which will be related the the displayed orange matching by a reflection and permutation of the strands. The vertices of the blue face which are red and green indicate that the corresponding orange edge connects a left vertex to a right vertex, while vertices which are only red or only green indicate that the corresponding orange matching edge connects same-side vertices.
\end{remark}

We now make the following definition which captures the types of surfaces that one obtains from the gluing procedure described in Figure \ref{fig::oddbluefaces}.

\begin{definition}[Flexible edge-plaquette embeddings]
Consider a pair $(\mathcal M, \phi)$ where $\mathcal M$ is a planar (or higher genus) map and $\phi:\mathcal M \to \Lambda$ is a map from the edges of $\mc{M}$ to the edges of $\Lambda$, and from the faces of $\mc{M}$ to the plaquettes of $\Lambda$. We call this pair a {\bf flexible edge-plaquette embedding} if the following hold:
\begin{enumerate}
\item The dual graph of $\mathcal M$ is bipartite. The faces of $\mathcal M$ in one partite class are designated as ``edge-faces'' (shown blue in figures) and those in the other class are called ``plaquette-faces'' (shown yellow in figures).
\item $\phi$ maps each plaquette-face of $\mathcal M$ isometrically {\em onto} a plaquette in $\mathcal P$.
\item $\phi$ maps each edge-face of $\mathcal M$ onto a single edge of $\Lambda$.
\end{enumerate}
\end{definition}

\begin{figure}[ht!]
    \centering
\includegraphics[width=.6\textwidth]{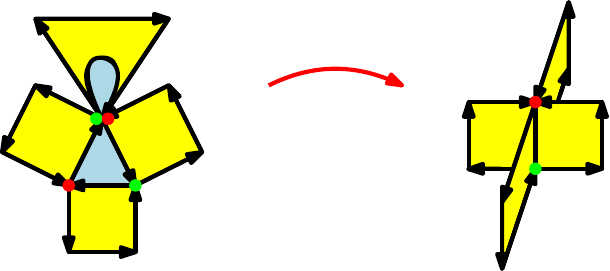}
    \caption{{\bf Flexible edge-plaquette embedding example:} Shown left is part of an oriented planar map. In a flexible edge-plaquette embedding, the embedding function $\phi$ maps directed edges of the map to directed edges of the lattice, but it is not required that $\phi$ extends to a single-valued function on vertices of the map. Here the edges of the blue triangle and blue $1$-gon all map to the red-green edge on the right; the vertex shared by the triangle and $1$-gon is colored both red and green to illustrate that it does not map to a single vertex on the right. Recall that when $\UN$ is replaced by $\orthogonal(N)$ the corresponding surfaces become {\em non-orientable}. When $\UN$ or $\orthogonal(N)$ is replaced by $\SUN$ or $\SON$ the corresponding surfaces become {\em flexible} in the sense illustrated here.}
    \label{fig::flexiblesemifold}
\end{figure}

\begin{figure}[ht!]
    \centering
\includegraphics[width=.6\textwidth]{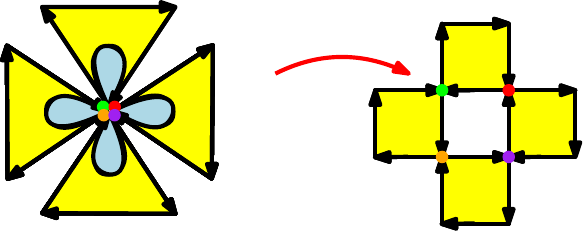}
    \caption{{\bf Flexible edge-plaquette embedding example:} In a flexible edge-plaquette embedding, a single vertex in the map (left) can in principle correspond to several vertices in the lattice (right). But each plaquette (directed edge) on the left has  a uniquely defined image plaquette (directed edge) on the right, and the boundary edges of any single blue face on the left all map to the same undirected blue edge on the right. In some sense, the image of a single vertex on the left is a closed cycle on the right, because as one ``moves around the vertex on the left clockwise'' one passes through a sequence of plaquette corners whose images on the right trace a cycle (possibly with some repeated vertices).}
    \label{fig::flexiblesemifolds}
\end{figure}

\begin{remark}
Comparing the definitions of flexible edge-plaquette embedding and edge-plaquette embedding, the main difference is that in the flexible case, $\phi$ is not necessarily a graph homomorphism. See Figures \ref{fig::flexiblesemifold} and \ref{fig::flexiblesemifolds}
for examples and intuition.
\end{remark}

In anticipation of the eventual application to lattice gauge theory, now suppose that the letters $\{\lambda_1, \ldots, \lambda_L\}$ are edges of the lattice $\Lambda$. Suppose the word $\bm \Gamma = \bm \Gamma(s, K)$ arises from a string $s = (\ell_1, \ldots, \ell_n)$ and a plaquette count $K : \mc{P} \ra \N$. Then the preceding discussion shows that $(\bm \Gamma, \bm \pi)$ is equivalent to a pair $(\mc{M}, \phi)$ which satisfies all the conditions of a flexible edge-plaquette embedding, except there are some faces of $\mc{M}$ whose boundaries are mapped to the loops in $s$ (rather than plaquettes). By deleting such faces, we obtain a flexible edge-plaquette embedding $(\ovl{\mc{M}}, \phi)$.
In summary, $(\bm \Gamma, \bm \pi)$ is equivalent to a flexible edge-plaquette embedding $(\ovl{\mc{M}}, \phi)$.

% Moreover, by removing the faces of $\mc{M}$ whose boundaries are mapped by $\phi$ to the words in $\bm \Gamma$, we obtain a map $(\ovl{\mc{M}}, \phi)$ with ``boundary" given by $\bm \Gamma$.

\begin{definition}
Let $s = (\ell_1, \ldots, \ell_n)$ be a string in $\Lambda$. Let $K : \mc{P} \ra \N$. Define $\fsfm_{\mrm{SU}}(s, K)$ to be the set of flexible edge-plaquette embeddings $(\ovl{\mc{M}}, \phi)$ that one can obtain when there are $K(p)$ copies of the plaquette $p$, for each $p \in \mc{P}$. I.e., $\fsfm_{\mrm{SU}}(s, K)$ is the set of flexible edge-plaquette embeddings that one may obtain by ranging over $\bm \pi = (\pi_e, e \in E_\Lambda^+)$, where $\pi_e \in \mc{B}_{n_e(+), n_e(-)}$, using our correspondence between $(\bm \Gamma(s, K), \bm \pi) \leftrightarrow (\ovl{\mc{M}}, \phi)$. Here, $n_e(+), n_e(-)$ are the respective number of times that $e, e^{-1}$ appear in $\bm \Gamma(s, K)$. Let
\begin{equs}
\fsfm_{\mrm{SU}}(s) := \bigsqcup_{K : \mc{P} \ra \N} \fsfm_{\mrm{SU}}(s, K).
\end{equs}
For $(\ovl{\mc{M}}, \phi) \in \fsfm_\mrm{SU}(s)$ and $e \in E_\Lambda$, let $\mu_e(\phi)$ be the profile of edge-faces of $(\ovl{\mc{M}}, \phi)$ at the edge $e$.

Define $\fsfm_{\mrm{SO}}(s, K)$ and $\fsfm_{\mrm{SO}}(s)$ in the same way, except we only require $\pi_e \in \mc{B}_{n_e(+) + n_e(-)}$ for each $e$. 
\end{definition}

\begin{definition}
Let $n, m \geq 0$. Define the normalized Weingarten function
\begin{equs}
\ovl{\Wg}^{\mrm{SU}}_N(\pi) &:= N^{n+m - \cycles(\pi)} \Wg^{\mrm{SU}}_N(\pi), ~~ \pi ~ \in \mc{B}_{n, m}, \\
\ovl{\Wg}^{\mrm{SO}}_N(\pi) &:= N^{n - \cycles(\pi)} \Wg^{\mrm{SO}}_N(\pi), ~~ \pi ~ \in \mc{B}_{n}
\end{equs}
where $\cycles(\pi)$ is the number of faces in the face profile of $\pi$.
\end{definition}

We can now state the following theorem which expresses Wilson loop expectations in $\SUN$ lattice gauge theory as sums over flexible edge-plaquette embeddings.

\begin{theorem}
\revision{Let $\Lambda$ be a finite lattice.} Let $s = (\ell_1, \ldots, \ell_n)$ be a string. For $\SUN$ lattice gauge theory, we have that
\begin{equs}
\langle W_s \rangle_{\Lambda, \beta} = Z_{\Lambda, \beta}^{-1} \sum_{(\mc{M}, \phi) \in \fsfm_{\mrm{SU}}(s)} \frac{\beta^{\mrm{area}(\mc{M}, \phi)}}{(\phi^{-1})!} \prod_{e \in E_\Lambda} \ovl{\Wg}^{\mrm{SU}}_N(\mu_e(\phi)) N^{\chi(M) - n}. 
\end{equs}
For $\SON$ lattice gauge theory, we have that
\begin{equs}
\langle W_s \rangle_{\Lambda, \beta} = Z_{\Lambda, \beta}^{-1} \sum_{(\mc{M}, \phi) \in \fsfm_{\mrm{SO}}(s)} \frac{\beta^{\mrm{area}(\mc{M}, \phi)}}{(\phi^{-1})!} \prod_{e \in E_\Lambda} \ovl{\Wg}^{\mrm{SO}}_N(\mu_e(\phi)) N^{\chi(M) - n}. 
\end{equs}
\revision{The sums converge absolutely in the following sense. Let $G = \mrm{SU}, \mrm{SO}$. Recalling that $\fsfm_G(s) = \bigsqcup_{K : \mc{P} \ra \N} \fsfm_G(s, K)$, we have that
\[ \sum_{K : \mc{P} \ra \N} \Bigg|\sum_{(\mc{M}, \phi) \in \fsfm_{G}(s, K)} \frac{\beta^{\mrm{area}(\mc{M}, \phi)}}{(\phi^{-1})!} \prod_{e \in E_\Lambda} \ovl{\Wg}^{G}_N(\mu_e(\phi)) N^{\chi(M) - n}\Bigg| < \infty.\]}
\end{theorem}

\subsubsection{Makeenko-Migdal/Master loop/Schwinger-Dyson equation}

To obtain the single-strand Makeenko-Migdal/Master loop/Schwinger-Dyson equation, we will need to modify the previous argument for $\UN, \orthogonal(N), \SphN$, because in those cases we had the strand-by-strand exploration, while in the $\SUN$ and $\SON$ cases we do not. Ultimately, this is due to the fact that when $G = \SUN, \SON$, there is some nonzero contribution from the event that all exploration eras do not end by time $T$ (even when we send $T \toinf$). On the other hand, the delicate cancellation properties that we took advantage of when $G = \UN, \orthogonal(N), \SphN$ were only on the event that all exploration eras have ended. Thus when $G = \SUN, \SON$, we cannot expect that the same strand-by-strand exploration will suffice -- in particular, the key property of the strand-by-strand exploration (Propositions \ref{prop:strand-by-strand-exploration} and \ref{prop:os-strand-by-strand}) no longer holds for $\SUN, \SON$.

We begin our alternate approach by introducing analogs of Jucys-Murphy elements for the Brauer and walled Brauer algebras. See \cite{nazarov1996young, jung2020supersymmetric} for more discussion on these elements.

\begin{definition}
Let $n \geq 1$. Define the Brauer algebra elements $x_1, \ldots, x_n \in \mc{B}_n$ by
\begin{equs}
x_k := J_k - \sum_{i=1}^{k-1} \langle k ~ i \rangle = \sum_{i=1}^{k-1} \big((k ~ i) - \langle k ~ i \rangle\big), ~~ k \in [n].
\end{equs}
These elements are the generalizations of the Jucys-Murphy elements for the Brauer algebra. In particular, they are mutually commuting (see \cite[Corollary 2.2]{nazarov1996young}).

Additionally, let $m \geq 1$. Define the walled Brauer algebra elements $z_1, \ldots, z_{n+m} \in \mc{B}_{n, m}$ by
\begin{equs}
z_k := \begin{cases} 
\sum_{i=1}^{k-1} (k ~ i ) - \sum_{i =n+1}^{n+m} \langle k ~ i \rangle & ~~ k \in [n] \\
\sum_{i=n+1}^{k-1} (k ~ i) & k \in (n:n+m].
\end{cases}
\end{equs}
These elements are the generalizations of the Jucys-Murphy elements for the walled Brauer algebra. In particular, they are mutually commuting (see \cite[Proposition 2.6]{jung2020supersymmetric}).
\end{definition}

\begin{lemma}\label{lemma:su-u-brownian-motion}
Let $Z := z_1 + \cdots + z_{n+m}$ for brevity. For $G = \SUN$, we have that
\begin{equs}
\E[B_T^{\otimes n} \otimes \bar{B}_T^{\otimes m}] &= e^{\frac{(n-m)^2}{2N^2}T} e^{-\frac{n+m}{2}T} e^{-\frac{T}{N}\rho_+(Z)} \\
&= e^{\frac{(n-m)^2}{2N^2}T} e^{-\frac{n+m}{2}T} e^{-\frac{T}{N} \rho_+(z_n)} e^{-\frac{T}{N}\rho_+(Z - z_n)}.
\end{equs}
Let $X := x_1 + \cdots + x_n$ for brevity. For $G = \SON$, we have that
\begin{equs}
\E[B_T^{\otimes n}] &= e^{-\frac{n}{2}(1 - \frac{1}{N})T} e^{-\frac{T}{N}X} \\
&= e^{-\frac{n}{2}(1 - \frac{1}{N})T} e^{-\frac{T}{N} \rho_+(x_n)} e^{-\frac{T}{N}\rho_+(X - x_n)}
\end{equs}
\end{lemma}
\begin{proof}
In both cases, the second identity follows from the first by mutual commutativity of the Jucys-Murphy elements. The first identity for $G = \SUN$ (resp. $\SON$) follows by Proposition \ref{prop:SUN-brownian-motion-expectation} (resp. \ref{prop:orthogonal-symplectic-brownian-motion-expectation}) and an explicit Poisson calculation.
% The second identity follows from the first identity and the mutual commutativity of $z_1, \ldots, z_{n+m}$, so we just prove the first identity. Recall that a $\UN$ Brownian motion $U_T$ may be written as $U_T = e^{\icomplex W_T / N} B_T$, where $B$ is a $\SUN$ Brownian motion, $W$ is a $\R$-valued Brownian motion, and $B, W$ are independent. This implies that
% \begin{equs}
% \E[U_T^{\otimes n} \otimes \bar{U}_T^{\otimes m}] &= \E[e^{\icomplex \frac{(n-m)}{N} W_T}] \E[B_T^{\otimes n} \otimes \bar{B}_T^{\otimes m}] \\
% &= e^{-\frac{(n-m)^2}{2N^2} T} \E[B_T^{\otimes n} \otimes \bar{B}_T^{\otimes m}].
% \end{equs}
% To finish, note that
% \begin{equs}
% \E[U_T^{\otimes n} \otimes \bar{U}_T^{\otimes m}] = e^{-\frac{n+m}{2} T} e^{-\frac{T}{N} \rho_+(Z)}.
% \end{equs} \qedhere
\end{proof}

\begin{remark}
In terms of our Poisson point process, this lemma has the following interpretation: we may first explore all points involving the top strand, and then all points which do not involve the top strand. 
\end{remark}

\begin{definition}
Let $\Sigma_{\mrm{SU}}^{\mrm{top}}(T) \sse \Sigma_{\mrm{SU}}(T)$ be the process defined by keeping only those points which involve the top strand. Define $\Sigma_{\mrm{SU}}^{\mrm{rest}}$ to be the complement of $\Sigma_{\mrm{SU}}^{\mrm{top}}$, i.e. the process defined by keeping only those points which do not involve the top strand. Define $\Sigma_{\mrm{SO}}^{\mrm{top}}, \Sigma_{\mrm{SO}}^{\mrm{rest}}$ in the same manner.
\end{definition}

\begin{remark}
By Poisson thinning, $\Sigma_{\mrm{SU}}^{\mrm{top}}$ and $\Sigma_{\mrm{SU}}^{\mrm{rest}}$ are independent Poisson processes. 
\end{remark}

By Lemma \ref{lemma:su-u-brownian-motion} and explicit calculation, we have the following identity, which states Lemma \ref{lemma:su-u-brownian-motion} in terms of our Poisson point process. The proof is omitted.

\begin{lemma}\label{lemma:su-BM-point-process-top-rest}
For $G = \SUN$, we have that
\begin{equs}
\E[B_T^{\otimes n} \otimes \bar{B}_T^{\otimes m}] = e^{-\frac{n+m}{2} (1 + \frac{1}{N^2}) T} \rho_+\Big( e^{(2(n+m)-1) T} \E[\pointstobrauerfn_{\mrm{SU}}(\Sigma_{\mrm{SU}}^{\mrm{top}}(T))] e^{(n+m-1)^2T} \E[\pointstobrauerfn_{\mrm{SU}}(\Sigma_{\mrm{SU}}^{\mrm{rest}}(T))] 
\Big).
\end{equs}
For $G = \SON$, we have that
\begin{equs}
\E[B_T^{\otimes n}] = e^{-\frac{n}{2}(1 - \frac{1}{N})T} \rho_+\Big(e^{2(n-1)T} \E[\pointstobrauerfn(\Sigma_{\mrm{SO}}^{\mrm{top}}(T))] e^{2\binom{n-1}{2}T} \E[\pointstobrauerfn(\Sigma_{\mrm{OS}}^{\mrm{rest}}(T))] \Big).
\end{equs}
\end{lemma}

\begin{remark}
In the above, when $G = \SUN$ we choose to split the exponential prefactor $(n+m)^2 T = (2(n+m) - 1) T + (n+m-1)^2 T$, since $2(n+m) - 1$ is the number of independent rate 1 Poisson processes contributing to $\Sigma_{\mrm{SU}}^{\mrm{top}}(T)$, and $(n+m-1)^2$ is the number of independent rate 1 Poisson processes contributing to $\Sigma_{\mrm{SU}}^{\mrm{rest}}(T)$. Similar considerations hold when $G = \SON$.
\end{remark}

\begin{notation}
For notational brevity in what follows, define
\begin{equs}
X_{\mrm{SU}}(T) &:= e^{-\frac{1}{2}(1 + \frac{1}{N^2}) T} e^{(2(n+m)-1) T} \E[\pointstobrauerfn(\Sigma_{\mrm{SU}}^{\mrm{top}}(T))] \in \mc{B}_{n, m}, \\
Y_{\mrm{SU}}(T) &:= e^{-\frac{n+m-1}{2}(1 + \frac{1}{N^2})T} e^{(n+m-1)^2T} \E[\pointstobrauerfn(\Sigma_{\mrm{SU}}^{\mrm{rest}}(T))] \in \mc{B}_{n, m}, \\
X_{\mrm{SO}}(T) &:= e^{-\frac{1}{2}(1 - \frac{1}{N})T} e^{2(n-1)T} \E[\pointstobrauerfn(\Sigma_{\mrm{SO}}^{\mrm{top}}(T))] \in \mc{B}_n, \\
Y_{\mrm{SO}}(T) &:= e^{-\frac{n-1}{2}(1 - \frac{1}{N})T} e^{2\binom{n-1}{2}T} \E[\pointstobrauerfn(\Sigma_{\mrm{OS}}^{\mrm{rest}}(T))] \in \mc{B}_n.
\end{equs}
% We omit the dependence of $X, Y$ on $m, n, N$.
\end{notation}

By Lemma \ref{lemma:su-BM-point-process-top-rest}, we have that
\begin{equs}\label{eq:su-BM-X-Y}
\E[B_T^{\otimes n} \otimes \bar{B}_T^{\otimes m}] &= \rho_+(X_{\mrm{SU}}(T) Y_{\mrm{SU}}(T)), ~~ G = \SUN \\
\E[B_T^{\otimes n}] &= \rho_+(X_{\mrm{SO}}(T) Y_{\mrm{SO}}(T)), ~~ G = \SON.
\end{equs}

The starting point to deriving an eventual recursion for $\SUN$ or $\SON$ Haar measure is the following recursion for $X$.

\begin{lemma}\label{lemma:X-recursion}
We have that
\begin{equs}
X_{\mrm{SU}}(T) &= e^{-\frac{1}{2}(1 + \frac{1}{N^2})T} + \bigg(\frac{n-m}{N^2} - \frac{z_n}{N}\bigg) \int_0^T du ~ e^{-\frac{1}{2}(1 + \frac{1}{N^2})u} X_{\mrm{SU}}(T-u), \\
X_{\mrm{SO}}(T) &= e^{-\frac{1}{2}(1 - \frac{1}{N})T} - \frac{x_n}{N} \int_0^T du ~ e^{-\frac{1}{2}(1 - \frac{1}{N}) u} X_{\mrm{SO}}(T-u).
\end{equs}
\end{lemma}
\begin{proof}
This follows by considering the time of the first point in $\Sigma_{\mrm{SU}}^{\mrm{top}}(T)$ (resp. $\Sigma_{\mrm{SO}}^{\mrm{top}}(T)$).
\end{proof}

For $\bullet \in \{\mrm{SU}, \mrm{SO}\}$, if we substitute the identities given by Lemma \ref{lemma:X-recursion} into $X_{\bullet}(T) Y_{\bullet}(T)$, the term $X_{\bullet}(T-u) Y_{\bullet}(T)$ for $U \in [0, T]$ appears. The following lemma interprets this term as an appropriate Brownian motion expectation.

\begin{lemma}\label{lemma:X-T-minus-u-Y-T}
For any $0 \leq u \leq T$, we have that
\begin{equs}
\rho_+(X_{\mrm{SU}}(T-u) Y_{\mrm{SU}}(T)) &=\E[(B_T B_u^{-1}) \otimes B_T^{\otimes (n-1)} \otimes \bar{B}_T^{\otimes m}], \\
\rho_+(X_{\mrm{SO}}(T-u) Y_{\mrm{SO}}(T)) &= \E[(B_T B_u^{-1}) \otimes B_T^{\otimes (n-1)}].
\end{equs}
\end{lemma}
\begin{proof}
We only prove the $\SUN$ case as the $\SON$ is very similar. We may write
\begin{equs}
(B_T B_u^{-1}) \otimes B_T^{\otimes (n-1)} \otimes \bar{B}_T^{\otimes m} = \Big( (B_T B_u^{-1})^{\otimes n} \otimes \ovl{B_T B_u^{-1}}^{\otimes m}\Big) ~\big(I \otimes B_u^{\otimes (n-1)} \otimes \bar{B}_u^{\otimes m}\big) .
\end{equs}
Since $B_u$ and $B_T B_u^{-1}$ are independent, upon taking expectations we obtain
\begin{equs}
\E[(B_T B_u^{-1}) \otimes B_T^{\otimes (n-1)} \otimes \bar{B}_T^{\otimes m}] &= \Big( \E[(B_T B_u^{-1})^{\otimes n} \otimes \ovl{B_T B_u^{-1}}^{\otimes m}]\Big) ~\big(I \otimes \E[ B_u^{\otimes (n-1)} \otimes \bar{B}_u^{\otimes m}]\big)\\
&=\E[B_{T - u}^{\otimes n} \otimes \bar{B}_{T-u}^{\otimes m}] ~\big(I \otimes \E[B_u^{\otimes (n-1)} \otimes \bar{B}_u^{\otimes m}]\big) .
\end{equs}
Writing $Z = z_1 + \cdots + z_{n+m}$ for brevity, we have by Lemma \ref{lemma:su-u-brownian-motion} that
\begin{equs}
\E[B_{T - u}^{\otimes n} \otimes \bar{B}_{T-u}^{\otimes m}] &= e^{\frac{(n-m)^2}{2N^2} (T-u)} e^{-\frac{n+m}{2}(T-u)} e^{-\frac{T-u}{N} \rho_+(Z)}, \\
I \otimes \E[B_u^{\otimes (n-1)} \otimes \bar{B}_u^{\otimes m}] &= e^{\frac{(n - m-1)^2}{2N^2} u}e^{-\frac{n+m-1}{2}u} e^{-\frac{u}{N}\rho_+(Z - z_n)}.
\end{equs}
Combining, we obtain (using the mutual commutativity of $z_1, \ldots, z_{n+m}$)
\begin{equs}
\E[(B_T B_u^{-1}) \otimes B_T^{\otimes (n-1)}] = e^{\frac{2(n-m) - 1}{2N^2} (T-u)} e^{-\frac{1}{2}(T - u)} e^{-\frac{T-u}{N} \rho_+(z_n)} e^{\frac{(n-m-1)^2}{2N^2} T}  e^{-\frac{n+m-1}{2}T} e^{-\frac{T}{N}\rho_+(Z - z_n)}.
\end{equs}
By an explicit calculation, we have that the right\revision{-}hand side above is exactly
\begin{equs}
\rho_+\bigg( e^{-\frac{1}{2}(1 + \frac{1}{N^2}) (T-u)} e^{(2(n+m)-1)(T -u)} \E[\pointstobrauerfn(\Sigma_{\mrm{SU}}^{\mrm{top}}(T - u))] e^{-\frac{n+m-1}{2}(1 + \frac{1}{N^2}) T}e^{(n+m-1)^2T} \E[\pointstobrauerfn(\Sigma_{\mrm{SU}}^{\mrm{rest}}(T))] \bigg),
\end{equs}
which is exactly $\rho_+(X_{\mrm{SU}}(T-u) Y_{\mrm{SU}}(T))$, and thus the desired result follows.
\end{proof}

Next, we show that the Brownian motion expectation which appears in Lemma \ref{lemma:X-T-minus-u-Y-T} has a nice limit as $T \toinf$. We prove this for more general $G$ than needed, as the argument is exactly the same.

\begin{lemma}\label{lemma:B-T-B-u-inverse-limit}
Let $G = \UN, \SUN, \SON, \SphN$. For any $u \geq 0$, we have that
\begin{equs}
\lim_{T \toinf} \E[(B_T B_u^{-1}) \otimes B_T^{\otimes (n-1)} \otimes \bar{B}_T^{\otimes m}] = e^{\frac{c_{\mathfrak{g}}}{2} u}\E[G^{\otimes n} \otimes \bar{G}^{\otimes m}].
\end{equs}
\end{lemma}
\begin{proof}
We may write
\begin{equs}
(B_T B_u^{-1}) \otimes B_T^{\otimes (n-1)} \otimes \bar{B}_T^{\otimes m} = \big(B_T^{\otimes n} \otimes \bar{B}_T^{\otimes m} \big) ~ \big(B_u^{-1} \otimes I^{\otimes (n+m-1)} \big).
\end{equs}
For fixed $u$, the conditional distribution $B_T ~|~ B_u$ converges to normalized Haar measure on $G$ as $T \toinf$. We thus obtain
\begin{equs}
\lim_{T \toinf} \E[(B_T B_u^{-1}) \otimes B_T^{\otimes (n-1)} \otimes \bar{B}_T^{\otimes m}] = \E[G^{\otimes n} \otimes \bar{G}^{\otimes m}] \big(\E[B_u^{-1}] \otimes I^{\otimes (n+m-1)}\big).
\end{equs}
By an explicit calculation, we have that
\begin{equs}
\E[B_u^{-1}] = \E[B_u^*] = \big(\E[B_u]\big)^* = \big(e^{\frac{c_{\mathfrak{g}}}{2} u} I\big)^* =e^{\frac{c_{\mathfrak{g}}}{2} u} I .
\end{equs}
The desired result now follows.
\end{proof}

Now by combining the previous few preliminary results, we obtain the following recursions for expectations with respect to $\SUN$ or $\SON$ Haar measure.

\begin{prop}\label{prop:SUN-SON-haar-measure-recursion}
For $n \geq 1, m \geq 0$, we have that
\begin{equs}
\E[G^{\otimes n} \otimes \bar{G}^{\otimes m}] &= \rho_+\bigg(\frac{n-m}{N^2} - \frac{z_n}{N}\bigg) \E[G^{\otimes n} \otimes \bar{G}^{\otimes m}], ~~ G = \SUN, \\
\bigg(1 - \frac{1}{N}\bigg) \E[G^{\otimes n}] &= - \frac{\rho_+(x_n)}{N} \E[G^{\otimes n}], ~~ G = \SON.
\end{equs}
\end{prop}
\begin{proof}
We prove the $\SUN$ case as the $\SON$ is very similar. For brevity, let $E_T(n, m) = \E[B_T^{\otimes n} \otimes \bar{B}_T^{\otimes m}]$. Combining equation \eqref{eq:su-BM-X-Y} with Lemma \ref{lemma:X-recursion}, we obtain
\begin{equs}
E_T(n, m) = e^{-\frac{1}{2}(1 + \frac{1}{N^2}) T} \rho_+(Y_{\mrm{SU}}(T)) + \bigg(\frac{n-m}{N^2} - \frac{\rho_+(z_n)}{N}\bigg) \int_0^T e^{-\frac{1}{2}(1 + \frac{1}{N^2})u} \rho_+(X_{\mrm{SU}}(T-u) Y_{\mrm{SU}}(T)).
\end{equs}
Note that $\rho_+(Y_{\mrm{SU}}(T)) = I \otimes E_T(n-1, m)$, which is $O(1)$ as $T \toinf$. Thus as $T \toinf$, the first term in the right\revision{-}hand side above is $o(1)$. Combining this with Lemmas \ref{lemma:X-T-minus-u-Y-T} and \ref{lemma:B-T-B-u-inverse-limit}, we obtain upon taking $T \toinf$,
\begin{equs}
\E[G^{\otimes n} \otimes \bar{G}^{\otimes m}] =  \bigg(\frac{n-m}{N^2} - \frac{\rho_+(z_n)}{N}\bigg) \int_0^\infty du e^{-\frac{1}{2}(1 + \frac{1}{N^2})u} e^{\frac{c_{\mathfrak{su}(N)}}{2}u} \E[G^{\otimes n} \otimes \bar{G}^{\otimes m}].
\end{equs}
To finish, we use that $c_{\mathfrak{su}(N)} = - 1 + \frac{1}{N^2}$.
\end{proof}

\begin{definition}
Let $\bm \Gamma = (\Gamma_1, \ldots, \Gamma_k)$ be a collection of words on $\{\lambda_1, \ldots, \lambda_L\}$. For each letter $\lambda_i$, $i \in [L]$, let $t(\lambda_i)$ be equal to the number of occurrences of $\lambda_i$ in $\bm \Gamma$ minus the number of occurrences of $\lambda_i^{-1}$ in $\bm \Gamma$. Let $t(\lambda_i^{-1}) := -t(\lambda_i)$. Given a location $(i, j)$ of $\bm \Gamma$, let $t(i, j) := t(\lambda_k^s)$, where $\lambda_k^s$, $s \in \{\pm 1\}$, is the letter at location $(i, j)$.
\end{definition}

Proposition \ref{prop:SUN-SON-haar-measure-recursion} leads immediately (by similar considerations as in the proof of Proposition \ref{prop:word-recursion}) to the following recursions for expectations of words. The proof is omitted. In the following, recall the various string operations defined in Definitions \ref{def:splittings-and-mergers} and \ref{def:mergers-twistings}.

\begin{prop}[Single-location $\SUN$ and $\SON$ word recursion]\label{prop:su-word-recursion}
Let $\bm \Gamma = (\Gamma_1, \ldots, \Gamma_k)$ be a collection of words on $\{\lambda_1, \ldots, \lambda_L\}$. For any location $(i, j)$ of $\Gamma$, we have that for $G = \SUN$,
\begin{equs}
\bigg(1 - \frac{t(i, j)}{N^2}\bigg)\E[\tr(G(\bm \Gamma))] = &- \sum_{\bm\Gamma' \in \splitting_+((i, j), \bm \Gamma)} \E[\tr(G(\bm \Gamma'))] + \sum_{\bm \Gamma' \in \splitting_-((i, j), \bm \Gamma)} \E[\tr(G(\bm \Gamma'))] \\
&- \frac{1}{N^2} \sum_{\bm \Gamma' \in \merge_+^U((i, j), \bm \Gamma)} \E[\tr(G(\bm \Gamma'))] + \frac{1}{N^2} \sum_{\bm \Gamma' \in \merge_-^U((i, j), \bm \Gamma)} \E[\tr(G(\bm  \Gamma'))].
\end{equs}
For $G = \SON$, we have that
\begin{equs}
\bigg(1 - \frac{1}{N}\bigg)\E[\tr(G(\bm \Gamma))] = &- \sum_{\bm\Gamma' \in \splitting_+((i, j), \bm \Gamma)} \E[\tr(G(\bm \Gamma'))] + \sum_{\bm \Gamma' \in \splitting_-((i, j), \bm \Gamma)} \E[\tr(G(\bm \Gamma'))] \\
&- \frac{1}{N^2} \sum_{\bm \Gamma' \in \merge_+((i, j), \bm \Gamma)} \E[\tr(G(\bm \Gamma'))] + \frac{1}{N^2} \sum_{\bm \Gamma' \in \merge_-((i, j), \bm \Gamma)} \E[\tr(G(\bm  \Gamma'))] \\
&- \frac{1}{N} \sum_{\bm \Gamma' \in \mbb{T}_+((i, j), \bm \Gamma)} \E[\tr(G(\bm \Gamma'))] + \frac{1}{N} \sum_{\bm \Gamma' \in \mbb{T}_-((i, j), \bm \Gamma)} \E[\tr(G(\bm \Gamma'))]
\end{equs}
\end{prop}

% \begin{prop}[$\SON$ word recursion]\label{prop:so-word-recursion}
% Let $G = \SON$. Let $\bm \Gamma = (\Gamma_1, \ldots, \Gamma_k)$ be a collection of words on $\{\lambda_1, \ldots, \lambda_L\}$. For any location $(i, j)$ of $\Gamma$, we have that 
% \begin{equs}
% \bigg(1 - \frac{1}{N}\bigg)\E[\tr(G(\bm \Gamma))] = &- \sum_{\bm\Gamma' \in \splitting_+((i, j), \bm \Gamma)} \E[\tr(G(\bm \Gamma'))] + \sum_{\bm \Gamma' \in \splitting_-((i, j), \bm \Gamma)} \E[\tr(G(\bm \Gamma'))] \\
% &- \frac{1}{N^2} \sum_{\bm \Gamma' \in \merge_+((i, j), \bm \Gamma)} \E[\tr(G(\bm \Gamma'))] + \frac{1}{N^2} \sum_{\bm \Gamma' \in \merge_-((i, j), \bm \Gamma)} \E[\tr(G(\bm  \Gamma'))] \\
% &- \frac{1}{N} \sum_{\bm \Gamma' \in \mbb{T}_+((i, j), \bm \Gamma)} \E[\tr(G(\bm \Gamma'))] + \frac{1}{N} \sum_{\bm \Gamma' \in \mbb{T}_-((i, j), \bm \Gamma)} \E[\tr(G(\bm \Gamma'))]
% \end{equs}
% \end{prop}

By applying Proposition \ref{prop:su-word-recursion} to lattice Yang-Mills theories, we may obtain the single-location $\SUN$ and $\SON$ Makeenko-Migdal/Master loop/Schwinger-Dyson equation. The proof is entirely analogous to the proof of the $\UN$ Makeenko-Migdal/Master loop/Schwinger-Dyson equation (Theorem \ref{thm:master-loop}) using the $\UN$ word recursion (Proposition \ref{prop:word-recursion}), and thus it is omitted. Before we state the theorem, we first define the following new string operation which appears for $\SUN$.

\begin{definition}[Expansion]
Let $s = (\ell_1, \ldots, \ell_n)$ be a string. Let $(i, j)$ be a location in $s$. We define the sets of positive and negative expansions $\mbb{E}_+((i, j), s)$ and $\mbb{E}_-((i, j), s)$ as follows. Denote by $e$ the oriented edge of the lattice at location $(i, j)$ in $s$. 

The set of positive expansions $\mbb{E}_+((i, j), s)$ is the set of all possible strings $s'$ which can be obtained by adding an oriented plaquette $p \in \mc{P}$ which contains $e^{-1}$ to the collection of loops $s$.

The set of negative expansions $\mbb{E}_+((i, j), s)$ is the set of all possible strings $s'$ which can be obtained by adding an oriented plaquette $p \in \mc{P}$ which contains $e$ to the collection of loops $s$.
\end{definition}

\revision{In the following, recall the twistings $\mbb{T}_{\pm}$ which were defined in Definition \ref{def:mergers-twistings}.}

\begin{theorem}[Single-location $\SUN$ and $\SON$ Makeenko-Migdal/Master loop/Schwinger -Dyson equation]\label{thm:SUN-SON-master-loop}
Let $s = (\ell_1, \ldots, \ell_n)$ be a string. Let $(i, j)$ be a location in $s$. For $\SUN$ lattice Yang-Mills theory, we have that
\begin{equs}
\bigg(1 - \frac{t(i, j)}{N^2}\bigg)\phi(s) = &-\sum_{s' \in \splitting_+((i, j), s)} \phi(s') + \sum_{s' \in \splitting_-((i, j), s)} \phi(s') \\
&- \frac{1}{N^2} \sum_{s' \in \merge_+^{\revision{U}}((i, j), s)} \phi(s') + \frac{1}{N^2} \sum_{s' \in \merge_-^{\revision{U}}((i, j), s)} \phi(s') \\
&- 
\beta \sum_{s' \in \deform_+((i, j), s)} \phi(s') + \beta \sum_{s' \in \deform_-((i, j), s)} \phi(s') \\
& -\beta \sum_{s' \in \mbb{E}_+((i, j), s)} \phi(s') + \beta \sum_{s' \in \mbb{E}_-((i, j), s)} \phi(s').
\end{equs}
For $\SON$ lattice Yang-Mills theory, we have that
\begin{equs}
\bigg(1 - \frac{1}{N}\bigg)\phi(s) = &-\sum_{s' \in \splitting_+((i, j), s)} \phi(s') + \sum_{s' \in \splitting_-((i, j), s)} \phi(s') \\
&- \frac{1}{N^2} \sum_{s' \in \merge_+((i, j), s)} \phi(s') + \frac{1}{N^2} \sum_{s' \in \merge_-((i, j), s)} \phi(s') \\
&- \frac{1}{N} \sum_{s' \in \mbb{T}_+((i, j), s)} \phi(s') + \frac{1}{N} \sum_{s' \in \mbb{T}_-((i, j), s)} \phi(s') \\
&- 
2\beta \sum_{s' \in \deform_+((i, j), s)} \phi(s') + 2\beta \sum_{s' \in \deform_-((i, j), s)} \phi(s').
\end{equs}
\end{theorem}

\begin{remark}
Observe that the $\SON$ Makeenko-Migdal/Master loop/Schwinger-Dyson equation is exactly the $\orthogonal(N)$ one. This is natural since $\orthogonal(N)$ Brownian motion is essentially $\SON$ Brownian motion (recall Remark \ref{remark:SON-ON-BM}).
\end{remark}

\section{Open problems} \label{sec::open}
Although lattice gauge theory has been very thoroughly studied in physics, there are many simple ideas about the relationship between random surfaces and Yang-Mills theory that have not been so thoroughly explored on the math side. There is also room for innovation: producing clever variants and toy models whose limits might be easier to describe in terms of continuum random surfaces (including those related to Liouville quantum gravity and conformal field theory). If the ultimate goal is to get a handle on a continuum theory, there is a good deal of flexibility in how one sets up the discrete models that are meant to approximate that theory.  We present a series of open problems along those lines, ranging from very general and open-ended to very technical and specific.

\begin{enumerate}
\item For which lattice models can we establish a version of the ``area law''  using the surface sum point of view? Recall that the area law states that the Wilson loop expectation decays exponentially in the minimal area spanned by the loop, at least for reasonably nice loops; see the definitions and discussion in \cite{chatterjee2016a} about the relationship between the area law and ``quark confinement.'' Many such results are known (from various points of view) for small $\beta$, and these results apply in any dimension $d\geq 2$, see e.g.\ \cite{osterwalderseiler1978} for the proof of the $\SUN$ area law for small $\beta$ and general $N$, and the discussion in \cite{Chatterjee2019a} which explains a string-trajectory-based derivation of such a result in the $N \to \infty$ limit for small $\beta$. For general $\beta$, the known results are dimension-dependent:
\begin{enumerate}
\item When $d=2$ the area law is well-known for general groups for any $\beta$ \cite{levy2003, levy2010a, Levy2011a}.
\item When $d=3$ and $N=1$, the area law holds for all $\beta$, see \cite{gopfertmack1981}. Because $\mrm{U}(1)$ is the center of $\UN$ for general $N$, this implies that the $\UN$ area law holds for all $\beta$ and all $N \geq 1$, see \cite{FROHLICH1979}. It is not known whether the $\SUN$ area law holds for large $\beta$ when $N>1$.
\item When $d=4$, interestingly enough, the $\mrm{U}(1)$ area law holds for small $\beta$ but {\em fails} for large $\beta$, see \cite{frohlich1982massless}. It remains a major open problem to prove the area law for any non-commutative group when $\beta$ is large, $N \geq 2$ and $d=4$.  
\end{enumerate}

\item For which lattice models can we establish exponential decay of correlations for the Wilson loop traces using the surface sum point of view? This is related to the so-called ``mass gap'' problem, see e.g.\ discussion in \cite{chatterjee2016a}.  In the settings above, one is usually able to prove exponential decay of correlations in the same settings where one is able to prove the area law. (See \cite{chatterjee2021probabilistic} for an argument that certain strong forms of exponential decay imply the area law.) In particular, it remains a major open problem to prove exponential decay of correlations for any non-commutative group when \revision{$\beta$ is large and $d=4$}.  

\item In the $\UN$ setting, what can we say about the {\em conditional} law of the surface {\em given} the number and type of blue plaquettes at each edge?  \revision{Note that since the Weingarten weights depend only on the blue faces, if we condition on seeing a given collection of blue faces, then every surface we see will have the exact same Weingarten weight, and so} we no longer need to consider the Weingarten function, and the remaining combinatorics becomes simpler: in fact one obtains precisely the sort of model used to study words in GUE matrices using Wick's formula \cite{zvonkin1997matrix}. In this setting all ways of hooking up yellow to blue along edges are allowed and all contribute with the same sign, but there is still a weighting according to the genus, which leads the surface to concentrate around minimal genus configurations in the large $N$ limit. As a simplified model, we could even imagine that we fix the number of blue faces of each type to be exactly the same at each edge. Can we say anything about the scaling limits in this setting? Is the GUE correspondence at all helpful here?

\item Within a three-dimensional lattice like $\mathbb Z^3$, one way to try to understand the scaling limit of an oriented random surface (which could become space-filling in the fine mesh limit, with genus tending to infinity) is to try to understand the limit of the ``height function'' on the dual lattice that changes by $\pm 1$ (depending on orientation) each time one crosses a layer of the surface. Is there a setting in which such a limit can be obtained? The gradient of such a function is in some sense the normal vector field corresponding to the surface. (It is a flow in which one unit of current is assigned for each face of the surface, in the direction orthogonal to that face; the flow is not divergence free but it is curl-free except along the boundary loops.) Is there a qualitative difference between $N=1$ and general $N$ in the limit? The $N=1$ case has been understood by Fr\"olich and Spencer \cite{frohlich1982massless} and has an interesting $\beta$-dependent phase transition (from area law to perimeter law, as mentioned above) that we would not expect to see for larger $N$.

\item Can we prove anything interesting about the variants in which there are many plaquettes but only three can meet along any given edge? For example $\Lambda$ might be the truncated octahedron tessellation (one example of a tessellation by cells where only three cells ever meet along the same edge, see \cite{sheffieldyadin2014}) and $\mathcal P$ can be the collection of of square and hexagonal faces in the tessellation. If we require that each plaquette appears zero times or once, then the only non-zero terms in the surface expansion involve surface in which either zero or two of the three plaquettes contain each given interior edge (i.e.\ each edge not on the Wilson loop).  In this case the surfaces we obtain are simpler: all of the blue faces are 2-gons and the surfaces are self-avoiding. There is no need to consider the Weingarten function in this simplified setting. We remark that this would be the surface analog of the loop $\mrm{O}(n)$ model, studied for instance in \cite{duminil2017exponential}. (Requiring the number of copies of a given plaquette to be small---here either $0$ or $1$---is somehow related to taking a small $\beta$ in the unrestricted-plaquette-number setting.)

\item Recall that in certain contexts it is enough to consider {\em connected surfaces}, such as when there is a single Wilson loop and $N \to \infty$ (recall the discussion just after Theorem \ref{thm:wilson-loop-expectation-sum-over-epe}). Are there other contexts in which it is sufficient to consider connected surfaces?

\item A surface sum like the one in Theorem~\ref{thm:informal-wilson-loop-expectation-sum-over-epe} includes many terms of both signs. Our intuition is that most of these surfaces somehow ``cancel each other out.'' For example, there may be local changes one can make to a surface that change the sign of the associated Weingarten product but do not change the genus of the surface. Is there a clean way to group together the surfaces in this sum that makes this cancellation more transparent? \revision{In the case $d = 2$ and $N = \infty$, this is the subject of the upcoming work \cite{BCSK2025}, which uses the master loop equation to ``explore" the surfaces described in \cite{BCSK2024} in order to find cancellations. The work in particular establishes that for simple loops, there is essentially only one surface which contributes, which is the natural one that is spanned by the loop (recall we are in 2D). All other possible surfaces cancel each other out. \cite{BCSK2025} is also able to handle more complicated loops which are not necessarily simple.}
% One could begin with the case $d=2$, and aim to show that the surfaces that are not locally flat somehow cancel each other out.

%. In this case, one can show that the Wilson loop expectation of a loop should not change if one adds or removes plaquettes that are not surrounded by the loop---in some sense all surfaces that go outside of 

%\item In a surface sum like the one in Theorem~\ref{thm:informal-wilson-loop-expectation-sum-over-epe}, we can imagine that instead of weighting by $N$ to the Euler characteristic, we weight by $N$ to the number of vertices: equivalently we can imagine we sum over labeled surface where each vertex has a label between $1$ and $N$. The labeled surface is more complex than the unlabeled surface, but it also has Markov properties that the unlabeled surface lacks. Can one make sense of a continuum version of this Markov property in any continuum model?

\item Is there a simpler expression (or at least asymptotic expression in the limit of a large number of plaquettes) for the Weingarten function in the case that $N$ is a small integer? Recall that in this case, the sum over representations in \eqref{eq:weingarten-character-sum} involves only those corresponding to Young tableaux with at most $N$ rows. \revision{We remark that when $N = 1$, the Weingarten function simplifies greatly, and is constant across all permutations.}

\item\label{item:finite-T} What is the most natural way to express the finite-$T$ (i.e. Brownian motion at time $T$, as in Section \ref{section:poisson-process-intro}) analog of the Weingarten function and the corresponding random planar maps? Note that adding a few single-edge loops may have a similar effect to switching to finite $T$. This is because weighting Haar measure on $\UN$ by a power of the real part of the trace {\em biases} the measure toward matrices that are near the identity; Brownian motion on a Lie group stopped at a finite time $T$ is also (compared to Haar measure) biased toward matrices that are near the identity.

\item \label{c1to25q} Are there any natural random surface models emerging in the lattice Yang-Mills framework that lead to planar maps similar to those whose limits (can be conjectured to) correspond to Liouville quantum gravity surfaces with $c \in (1,25)$?  Those surfaces are multi-ended and infinite, see e.g.\ \cite{gwynne2020liouville, ding2021regularity, ang2022brownian, ding2023uniqueness}.

\item There have been many recent results about random planar maps of high genus and/or random hyperbolic planar maps, see, e.g.\ \cite{angel2013local,curien2016planar, budzinski2021local,budzinski2022local, delecroix2022large,janson2023unicellular}. Which of these results can be can be extended to embedded random planar maps of the type that emerge in our analysis? 

\item Can we interpret Wilson loops in terms of Liouville quantum gravity at least in the critical $c=1$ setting where we have a ``ladder graph'' and have gauge fixed so that we have the identity on the left and right sides of the ladder, and each yellow plaquette can be treated as a $2$-gon (since the left and right edges can be shrunk to points)?  Since the yellow plaquettes are all $2$-gons, we can interpret them as edges between blue faces: each blue face comes with a ``height'' (the height of the ladder rung) and its neighboring blue faces have heights that are one unit higher or one unit lower.  Essentially one has a planar map of blue faces decorated by a one-dimensional height function, which one might expect to converge to Liouville quantum gravity with parameter $c=1$ in the $N=\infty$ limit. We note there are some physics connections between the large-$N$ 2D Yang-Mills and $c=1$ matrix model~\cite{minahan1993equivalence} whose double scaling limit is related to the Liouville gravity with matter central charge $c=1$.

\item Two-dimensional lattice gauge theory can also be reduced to a ladder graph (if one gauge fixes along a spiral, one essentially obtains a ladder) as in the setting of the previous question. But is there any sense in which the Liouville quantum gravity surfaces for $c = 2$ (as mentioned above in the \ref{c1to25q}th question) can be recovered in these models?

\item What can we learn from models in which there is some correlation between the noise defining distinct edges, so that the analogs of the blue faces are perhaps not mapped to a single edge?  In an extreme case, one can take different edges to be perfectly correlated, so that one has the same random matrix at different locations. For example, one could assign the same random matrix to all edges that are parallel to each other, as is done in \cite{Eguchi1982}.

\item What is the fine-mesh scaling limit of the random surface we obtain when we fix exactly $b$ yellow plaquettes of each type and take $N=\infty$ (so that the surface is simply connected)?  Does it look like a continuum random tree (a.k.a.\ Brownian tree or branched polymer) conditioned to fill out $\Lambda$ in some even way?

\item We alert the reader that the ``spin-foam'' constructions in \cite{oeckl2001dual,conrady2005geometric} provide another approach for converting non-abelian lattice gauge theory into a statistical physical model. We can then pose a general question: what new properties of lattice Yang-Mills theory and/or its continuum scaling limits can be deduced from the spin-foam perspective?

\item The abelian versions of ``spin foam'' are simpler and were used e.g.\ by Fr\"olich and Spencer \cite{frohlich1982massless} to understand the phase transition structure of $\mrm{U}(1)$ lattice gauge theory.  Can an alternative proof of these results be given using the surface expansion described in this paper?

\item Adding extra single-edge faces in both directions has the effect of changing the underlying measure from Haar measure to another conjugation-invariant measure on $\UN$ (which can be a signed measure if we add associate sign weights to different edge configurations). Can one obtain a natural connection between a signed-measure variant of Yang-Mills theory and the sort of random surfaces that arise in conformal field theory?

\item What can one say about supersymmetric variants of this question?  Can a super-symmetric version of Yang-Mills theory be connected to random planar maps whose scaling limits can be understood in terms of Liouville quantum gravity or some other probabilistic continuum random surface model? What about fermionic variants or variants involving Higgs fields? On the latter point, let us remark that the introduction to \cite{CCHS3D} contains a list of references about the lattice Yang-Mills-Higgs model. A configuration in this context assigns a vector to each lattice vertex (in addition to assigning a matrix to each directed edge). In this context, one also considers open Wilson paths (whose endpoints are lattice vertices) in addition to closed Wilson loops.

\end{enumerate}

\begin{appendix}

\section{Properties of the Orthogonal Weingarten function}\label{appendix:orthogonal-weingarten-jucys-murphy-relation}

In this appendix, we give more detail on why Lemmas \ref{lemma:orthogonal-weingarten-function-face-profile} and \ref{lemma:os-weingarten-jucys-murphy-relation} are true, and in particular why it essentially follows from \cite{matsumoto2013weingarten}. Fix $n \geq 1$ even and $\zeta \in \C$. To help the reader, we indicate how to translate between our notation and the notation of \cite[Section 2.2.2]{matsumoto2013weingarten}. Our $\zeta$ translates to $z$. Our $n$ is the equivalent of $2k$. The subgroup $\mc{H}_n \sse \symgrp_n$ we defined in Definition \ref{def:H-n} is $H_k$ in \cite{matsumoto2013weingarten}. One can show that $|\mc{H}_n| = 2^{n/2} (n/2)!$, which translates to $|H_k| = 2^k k!$.

Matsumoto defines the Orthogonal Weingarten function $\Wg^{\mrm{O}}(\cdot; \zeta)$ as an element of the group algebra $\C[\symgrp_n]$. As part of its definition, this element is $\mc{H}_n$ bi-invariant, i.e. 
\begin{equs}\label{eq:weingarten-H-n-bi-invariant}
\Wg^{\mrm{O}}(h \sigma; \zeta) = \Wg^{\mrm{O}}(\sigma; \zeta) = \Wg^{\mrm{O}}(\sigma h; \zeta) \text{ for all $\sigma \in \C[\symgrp_n]$, $h \in \mc{H}_n$}.
\end{equs} 
The relation between Matsumoto's definition and our definition via pseudo-inverses is as follows:
\begin{equs}\label{eq:weingarten-two-def-relation}
\Wg^{\mrm{O}}_\zeta(\pi, \pi') = \Wg^{\mrm{O}}(\sigma_\pi^{-1} \sigma_{\pi'}; \zeta),
\end{equs}
where $\sigma_\pi$ is the permutation associated to $\pi$ as in Definition \ref{def:sigma-pi}. Here and in the following, we will write $\Wg_\zeta^{\mrm{O}}$ for definition of the Weingarten function as a pseudo-inverse, and $\Wg^{\mrm{O}}(\cdot; \zeta)$ for Matsumoto's definition of the Weingarten function as a group algebra element. Now, one can show that the face profile $\ell(\pi, \pi')$ is precisely the coset-type of $\sigma_\pi^{-1} \sigma_{\pi'}$ (which is defined in \cite[Section 2.2.1]{matsumoto2013weingarten}). As mentioned in in \cite[Section 2.2.1]{matsumoto2013weingarten}, two permutations $\sigma, \sigma'$ have the same coset-type if and only if they are part of the same double $\mc{H}_n$ coset, i.e. $\mc{H}_n \sigma \mc{H}_n = \mc{H}_n \sigma' \mc{H}_n$. By the $\mc{H}_n$ bi-invariance \eqref{eq:weingarten-H-n-bi-invariant}, it follows that $\Wg^{\mrm{O}}(\sigma; \zeta)$ is a function of the coset-type of $\sigma$, and then by \eqref{eq:weingarten-two-def-relation}, it follows that $\Wg^{\mrm{O}}_\zeta(\pi, \pi')$ is a function of the face profile $\ell(\pi, \pi')$ of $\pi, \pi'$. This shows Lemma \ref{lemma:orthogonal-weingarten-function-face-profile}.

Next, we discuss Lemma \ref{lemma:os-weingarten-jucys-murphy-relation}. Recall we defined (Definition \ref{def:H-n}) $P_{\mc{H}_n} = \frac{1}{|\mc{H}_n|} \sum_{h \in \mc{H}_n} h$. This translates to $(2^k k!)^{-1} \mbf{1}_k$. The ``zonal spherical function" $\omega^\lambda$ from the paper is for us $\chi^{2\lambda} P_{\mc{H}_n} \in \C[\symgrp_n]$ (where here $\lambda \vdash \frac{n}{2}$). We have that (by Lemma \ref{lemma:restriction-to-doubled-tableau}, as argued in the proof of Lemma \ref{lemma:orthogonal-invariant-projection-doubled-young-diagram})
\begin{equs}\label{eq:P-H-n-doubled-young-tableau}
P_{\mc{H}_n}  = P_{\mc{H}_n} \sum_{\lambda \vdash \frac{n}{2}} P_{2\lambda},
\end{equs}
where recall that (equation \eqref{eq:P-lambda-character-formula}) $P_{2\lambda} = \frac{\chi_{2\lambda}(\id)}{n!} \sum_{\sigma \in \symgrp_n} \chi_{2\lambda}(\sigma) \sigma \in \C[\symgrp_n]$.

Next, as in \cite{matsumoto2013weingarten}, we define for $\lambda \vdash \frac{n}{2}$ the quantity $D_\lambda(\zeta)$ as
\begin{equs}
D_\lambda(\zeta) := \prod_{(i, j) \in \lambda} (\zeta + 2j - i - 1).
\end{equs}
This quantity relates to Jucys-Murphy elements as follows. Define $X_\varep := (\varep N + J_{n-1})(\varep N + J_{n-3}) \cdots (\varep N + J_1)$. 

\begin{lemma}
For any $\lambda \vdash \frac{n}{2}$, we have that 
\begin{equs}\label{eq:jucys-murphy-action-on-P-2lambda}
P_{\mc{H}_n} P_{2\lambda} X_\varep = D_\lambda(\varep N) P_{\mc{H}_n} P_{2\lambda}.
\end{equs}
\end{lemma}
\begin{proof}
This is proven towards the end of \cite[Section 3]{zinn2009jucys}. For the reader's convenience, we reproduce the argument here. Recalling the discussion of Young's orthogonal idempotents from Section \ref{section:orthogonal-exploration-technical-proofs}, we may expand
\begin{equs}
P_{2\lambda} = \sum_{\lambda \in \mrm{SYT}(2\lambda)} e_T.
\end{equs}
By \cite[Proposition 4]{zinn2009jucys}, $P_{\mc{H}_n} e_T \neq 0$ implies that $T$ is obtained by the ``doubling" procedure described on \cite[Page 7]{zinn2009jucys}. As noted in the paper, by direct calculation, for any such $T$, we have that $e_T X_\varep = D_\lambda(\varep N) e_T$. The desired result now follows by combining these observations.
% By the math display just before \cite[Proposition 5]{zinn2009jucys}, we have that
% \begin{equs}
% P_{\mc{H}_n}  X_\varep = \sum_{\lambda \vdash \frac{n}{2}} D_\lambda(\varep N) P_{\mc{H}_n} P_{2\lambda}.
% \end{equs}
% For the reader's notation, we translate our notation to that of \cite{zinn2009jucys}: our $X_\varep$ is the equivalent of $G$ (see \cite[Equation (3.2)]{zinn2009jucys}), and our $D_\lambda(\varep N)$ is the equivalent of $c_\lambda$. Now finish by equation \eqref{eq:P-H-n-doubled-young-tableau}
\end{proof}

In our notation, Matsumoto defines the Orthogonal Weingarten function as an element $\Wg^{\mrm{O}}(\cdot; \zeta) \in \C[\symgrp_n]$ given by the formula
\begin{equs}\label{eq:weingarten-matsumoto-def}
\Wg^{\mrm{O}}(\cdot; \zeta) := |\mc{H}_n| \sum_{\substack{\lambda \vdash \frac{n}{2} \\ D_\lambda(\zeta) \neq 0}} D_\lambda(\zeta)^{-1} P_{2\lambda} P_{\mc{H}_n}. 
\end{equs}
This element is $\mc{H}_n$ bi-invariant, that is $h \Wg^{\mrm{O}}(\cdot; \zeta) = \Wg^{\mrm{O}}(\cdot; \zeta) = \Wg^{\mrm{O}}(\cdot; \zeta) h$ for all $h \in \mc{H}_n$ (these identities are equivalent to \eqref{eq:weingarten-H-n-bi-invariant}). The second identity follows since $P_{\mc{H}_n} h = P_{\mc{H}_n}$ for any $h \in \mc{H}_n$. The first identity follows since $P_{2\lambda}$ is central, so that $P_{2\lambda} P_{\mc{H}_n} = P_{\mc{H}_n} P_{2\lambda}$, combined with $h P_{\mc{H}_n} = P_{\mc{H}_n}$ for all $h \in \mc{H}_n$. 

Combining the $\mc{H}_n$ bi-invariance with the fact that the collection $(\sigma_\pi, \pi : [n] \ra [n])$ forms a complete set of coset representatives of $\mc{H}_n$ as a subgroup of $\symgrp_n$ (as mentioned in the beginning of \cite[Section 2.2.1]{matsumoto2013weingarten}), we may express
\begin{equs}
\Wg^{\mrm{O}}(\cdot; \zeta) = \sum_{\pi : [n] \ra [n]} \Wg^{\mrm{O}}(\sigma_\pi; \zeta ) \mc{H}_n \sigma_\pi.
\end{equs}
From this, we obtain (using that by definition, $\mc{H}_n$ stabilizes $\pi_0$ for the first identity, and equation \eqref{eq:weingarten-two-def-relation} for the second)
\begin{equs}
~[\pi ~ \pi_0] \Wg^{\mrm{O}}(\cdot; \zeta) = |\mc{H}_n| \sum_{\pi' : [n] \ra [n]} [\pi ~ \pi_0] \sigma_{\pi'} \Wg^{\mrm{O}}(\sigma_{\pi'}; \zeta) = |\mc{H}_n| \sum_{\pi' : [n] \ra [n]} [\pi ~ \pi'] \Wg_\zeta^{\mrm{O}}(\pi_0, \pi').
\end{equs}
On the other hand, inserting equation \eqref{eq:weingarten-matsumoto-def}, we have the formula
\begin{equs}
~[\pi ~ \pi_0] \Wg^{\mrm{O}}(\cdot; \zeta) = |\mc{H}_n| [\pi ~ \pi_0] \sum_{\substack{\lambda \vdash \frac{n}{2} \\ D_\lambda(z) \neq 0}}  D_\lambda(\zeta)^{-1}P_{2\lambda}  P_{\mc{H}_n}  .
\end{equs}
Upon equating the previous two identities (and using that $P_{2\lambda}$ is central), we obtain
\begin{equs}
\sum_{\pi' : [n] \ra [n]}\Wg^{\mrm{O}}_\zeta(\pi_0, \pi')  [\pi ~ \pi']  = [\pi ~ \pi_0] \sum_{\substack{\lambda \vdash \frac{n}{2} \\ D_\lambda(\zeta) \neq 0}}  D_\lambda(\zeta)^{-1}P_{\mc{H}_n} P_{2\lambda}  .
\end{equs}
Setting $\zeta = \varep N$ and applying the representation $\rho_\varep$ to both sides of the identity, we obtain
\begin{equs}\label{eq:jucys-murphy-weingarten-relation-intermediate}
\sum_{\pi' : [n] \ra [n]} \Wg^{\mrm{O}}_{\varep N}(\pi_0, \pi') \rho_\varep([\pi ~ \pi'] ) = \rho_\varep([\pi ~ \pi_0]) \sum_{\substack{\lambda \vdash \frac{n}{2} \\ D_\lambda(\varep N) \neq 0}}  D_\lambda(\varep N)^{-1} \rho_\varep( P_{\mc{H}_n}) \rho_\varep( P_{2\lambda}).
\end{equs}
We now claim that
\begin{equs}\label{eq:os-product-jucys-murphy-inverse}
\sum_{\substack{\lambda \vdash \frac{n}{2} \\ D_\lambda(\varep N) \neq 0}}  D_\lambda(\varep N)^{-1} \rho_\varep( P_{\mc{H}_n}) \rho_\varep( P_{2\lambda}) = \sum_{\substack{\lambda \vdash \frac{n}{2}}}  \rho_\varep( P_{\mc{H}_n}) \rho_\varep( P_{2\lambda}) \rho_\varep(X_\varep)^{-1}.
\end{equs}
Given this claim, we obtain that \eqref{eq:jucys-murphy-weingarten-relation-intermediate} is further equal to 
\begin{equs}
\rho_\varep([\pi ~ \pi_0]) \sum_{\substack{\lambda \vdash \frac{n}{2}}}  \rho_\varep( P_{\mc{H}_n}) \rho_\varep( P_{2\lambda}) \rho_\varep(X_\varep)^{-1} &= \rho_\varep([\pi ~ \pi_0]) \rho_\varep(P_{\mc{H}_n})\rho_\varep(X_\varep)^{-1} \\
&= \rho_\varep([\pi ~ \pi_0]) \rho_\varep(X_\varep)^{-1},
\end{equs}
where we used \eqref{eq:P-H-n-doubled-young-tableau} in the second equality and the fact that $\mc{H}_n$ by definition stabilizes $\pi_0$ in the second. Combining the previous few identities, we see that
\begin{equs}
\sum_{\pi' : [n] \ra [n]}\Wg^{\mrm{O}}_\zeta(\pi_0, \pi')  \rho_\varep([\pi ~ \pi']) = \rho_\varep([\pi ~ \pi_0]) \rho_\varep(X_\varep)^{-1},
\end{equs}
which is precisely Lemma \ref{lemma:os-weingarten-jucys-murphy-relation}.

To see the claim \eqref{eq:os-product-jucys-murphy-inverse}, first note that by \eqref{eq:jucys-murphy-action-on-P-2lambda}, we have that
\begin{equs}
\sum_{\substack{\lambda \vdash \frac{n}{2} \\ D_\lambda(\varep N) \neq 0}}  D_\lambda(\varep N)^{-1} \rho_\varep( P_{\mc{H}_n}) \rho_\varep(  P_{2\lambda}) \rho_\varep(X_\varep) = \sum_{\substack{\lambda \vdash \frac{n}{2} \\ D_\lambda(\varep N) \neq 0}} \rho_\varep( P_{\mc{H}_n}) \rho_\varep( P_{2\lambda}).
\end{equs}
As mentioned in the proof of Lemma \ref{lemma:rho-eps-I-T-n-limit}, $\rho_\varep(X_\varep)$ is always invertible, and thus the above implies
\begin{equs}
\sum_{\substack{\lambda \vdash \frac{n}{2} \\ D_\lambda(\varep N) \neq 0}}  D_\lambda(\varep N)^{-1} \rho_\varep( P_{\mc{H}_n}) \rho_\varep(  P_{2\lambda}) = \sum_{\substack{\lambda \vdash \frac{n}{2} \\ D_\lambda(\varep N) \neq 0}} \rho_\varep( P_{\mc{H}_n}) \rho_\varep( P_{2\lambda}) \rho_\varep(X_\varep)^{-1}.
\end{equs}
To finish, it suffices to show that for any $\lambda \vdash \frac{n}{2}$ such that $D_\lambda(\varep N) = 0$, we have that $\rho_\varep(P_{2\lambda}) = 0$. In the case $\varep = 1$, this follows because (as observed in the proof of Lemma \ref{lemma:eigenvalue-lower-bound-sum-all-Jucys-Murphy}) $\rho(P_{2\lambda}) = 0$ unless $\ell(2\lambda) \leq N$, and one may directly check that $\ell(2\lambda) \leq N$ implies $D_\lambda(N) \neq 0$ (the worst case is the box at location $(i, j) = (\ell(2\lambda), 1)$). 

Next, suppose $\varep = -1$. We claim that $\rho_-(P_{2\lambda}) = \rho(P_{(2\lambda)'})$, where $(2\lambda)'$ is the conjugate partition to $2\lambda$. Given this claim, we obtain that $\rho_-(P_{2\lambda}) = 0$ unless $\ell((2\lambda)') \leq N$. Note that $\ell((2\lambda)') = w(2\lambda) = 2w(\lambda)$, where $w(\lambda)$ is the number of columns of $\lambda$. By direct calculation, $2w(\lambda) \leq N$ implies that that $D_\lambda(- N) \neq 0$ (the worst case is the box at location $(i, j) = (1, w(\lambda))$).

To see why $\rho_-(P_{2\lambda}) = \rho(P_{(2\lambda)'})$, note that $\rho_-(\sigma) = \mrm{sgn}(\sigma) \rho(\sigma)$, and so
\begin{equs}
\rho_-(P_{2\lambda}) = \rho\bigg(\frac{\chi_{2\lambda}(\id)}{n!} \sum_{\sigma \in \symgrp_n} \mrm{sgn}(\sigma) \chi_{2\lambda}(\sigma) \sigma\bigg). 
\end{equs}
Using the classical fact that $\chi_{(2 \lambda)'}(\sigma) = \mrm{sgn}(\sigma) \chi_{2\lambda}(\sigma)$, the above is seen to be equal to $\rho(P_{(2\lambda)'})$, as desired.
% (as observed in the proof of Lemma \ref{lemma:eigenvalue-lower-bound-sum-all-Jucys-Murphy}) $\rho(P_{2\lambda}) = 0$ unless $\ell(2\lambda) \leq N$. From the definition of $\rho_-$, we have that $\rho_-(\sigma) = \mrm{sgn}(\sigma) \rho(\sigma)$ for $\sigma \in \symgrp_n$, and thus

% Moreover, if $\ell(2\lambda) \leq N$, then $D_\lambda(\varep N) \neq 0$. Combining these considerations with \eqref{eq:jucys-murphy-action-on-P-2lambda}, it follows that $\rho_\varep(X_\varep)$ is always invertible and moreover
% \begin{equs}
% \sum_{\substack{\lambda \vdash \frac{n}{2} \\ D_\lambda(N) \neq 0}}  D_\lambda(N)^{-1} \rho_\varep( P_{\mc{H}_n}) \rho_\varep( P_{2\lambda}) = \sum_{\substack{\lambda \vdash \frac{n}{2}}}  \rho_\varep( P_{\mc{H}_n}) \rho_\varep( P_{2\lambda}) \rho_\varep(X)^{-1}.
% \end{equs}

\end{appendix}

\bibliographystyle{hmralphaabbrv}
\bibliography{references}

\end{document}